\pdfoutput=1
\date{January 20, 2023}

\documentclass{article}
\usepackage{arxiv}
\usepackage{subfig}
\usepackage{graphicx}
\usepackage{graphics}
\usepackage{amssymb,amsfonts,amsmath,textcomp}
\usepackage{color}
\usepackage{array}
\usepackage{multirow}
\usepackage{mathptmx}
\usepackage{tikz}
\usepackage{gnuplot-lua-tikz}
\usepackage{algorithm}
\usepackage{algpseudocode}
\usepackage{scalerel}
\usepackage{caption}
\usepackage{soul}
\definecolor{mylightblue}{RGB}{55, 113, 200}
\definecolor{myred}{RGB}{255, 0, 0}
\definecolor{mycyan}{RGB}{0, 255, 255}
\graphicspath{{img/}}
\newcommand{\noter}[1]{
	\colorbox{myred} {#1}
}
\newcommand{\noteg}[1]{
	\colorbox{green} {#1}
}
\newcommand{\noteb}[1]{
	\colorbox{blue} {\color{white}#1}
}
\newcommand{\notey}[1]{
	\colorbox{yellow} {#1}
}
\newcommand{\Noter}[1]{
	{\color{red} #1}
}
\newcommand{\Noteg}[1]{
	{\color{green} #1}
}
\newcommand{\Noteb}[1]{
	{\color{blue} #1}
}
\newcommand{\Notey}[1]{
	{\color{yellow} #1}
}
\newcommand{\tS}[0]{{two-scale }
}
\newcommand{\TS}[0]{{Two-scale }
}
\newcommand{\TSD}[0]{{Two-scale$^{gld}$ }
}
\newcommand{\TSDn}[0]{{Two-scale$^{gld}$}
}
\newcommand{\TSI}[0]{{Two-scale$^{gli}$ }
}
\newcommand{\TSIn}[0]{{Two-scale$^{gli}$}
}
\newcommand{\GL}[0]{%
{GFEM$^{gl}$ }
}
\newcommand{\SP}[0]{{super-patch }
}
\newcommand{\SPn}[0]{{super-patch}
}
\newcommand{\tens}[1]{{\ensuremath{\mathsf{#1}}}}
\newcommand{\stretchint}[1]{
	{\vcenter{\hbox{\stretchto[440]{\displaystyle\int}{#1}}}}
}
\newcommand{\scaleint}[1]{
	{\vcenter{\hbox{\scaleto[3ex]{\displaystyle\int}{#1}}}}
}
\newcommand{\vm}[1]{
	{\ensuremath{\mathbf{#1}}}
}
\newcommand{\algva}[0]{
   \rule{\textwidth}{0.4pt} \vspace{-2em}
}
\DeclareCaptionFormat{bfcn}{\textbf{#1} #3}

\title{A two-scale solver for linear elasticity problems in the context of parallel message passing}
\author{Alexis Salzman\thanks{\textit{E-mail address of the corresponding author :} \texttt{alexis.salzman@ec-nantes.fr}}\\  Nantes universit\'e, Ecole centrale de Nantes,\\ GeM Institute, UMR CNRS 6180,\\ \textit{\small 1 rue de la No\"{e}, 44321, Nantes, France} \And  Nicolas	 Mo\"{e}s\\ Nantes universit\'e,   Ecole centrale de Nantes,  \\ GeM Institute, UMR CNRS 6180\\ \textit{ \small 1 rue de la No\"{e}, 44321, Nantes, France}\\
\\
Institut Universitaire de France (IUF)\\ \textit{ \small 1 rue Descartes, 75231 Paris, France}
 }
\newcommand*\lineheight[1]{}
\begin{document}
	\maketitle
\begin{abstract} 
This paper pushes further the intrinsic capabilities of the \GL global-local approach  introduced initially in \cite{Duarte2007}.
We develop a distributed computing approach using MPI (Message Passing Interface) both for the global and local problems.
Regarding local problems, a specific scheduling strategy is introduced.
Then, to measure correctly the convergence of the iterative process, we introduce a reference solution that
revisits the product of classical and enriched functions.
As a consequence, we are able to propose a purely matrix-based implementation of the global-local problem.
The distributed approach is then compared to other parallel solvers
either direct or iterative with domain decomposition.
The comparison addresses the scalability as well as the elapsed time.
Numerical examples deal with linear elastic problems: a polynomial exact solution problem, a complex micro-structure, and, finally, a pull-out test (with different crack extent).
\end{abstract}

\section{Introduction}
Computational methods are addressing more and more complex  problems.
This paper is concerned with elasticity problems having roughly two-scales.
A global scale at the structural level and a local scale representing a micro-structure, a crack or a reinforcement for example.
Note that the local scale  appears over the whole domain 
in the case of a micro-structure and only in narrow zones for cracks and reinforcements.
A direct finite element approach for these problems is to create a unique mesh and solve.
This approach does not take advantage of the two-scale nature.
This leads to difficulties: the matrix can have a huge size and can be badly conditioned.
Iterative domain decomposition resolution is one way to face huge matrix size (see \cite{Saad03}) but does not solve the conditioning issue unless a specific preconditioner is used.
A direct solver can cope with poor conditioning, but matrix size is still an issue.

In order to take into account the two-scale nature of a problem,  the first existing  global to local method (see for example \cite{NOOR1986}) use the raw  global scale field as a boundary condition for local problems (refined encapsulated meshes), which provide the desired accuracy to the simulation.
But since the raw global field does not take into account the local behavior, it does not provide a good boundary condition to the local problems. 
The \GL method, first introduced in \cite{Duarte2007}, adds the idea of local to global interaction via enrichment functions.
More precisely,  local problems, which can overlap, yield solutions used to build enrichment functions.
These functions are  used through the partition of unity method to compute at global scale an improved solution.
The overlapping aspect of the local problems and their intermixing use in the  enriched global problem permit to partially address the boundary condition error of the raw global problem field.
In \cite{Duarte2008} this problem of boundary condition has been addressed by increasing the size of the local problems (the boundary condition at the local scale is pushed far enough) or by increasing the polynomial interpolation order  of the  global scale problem.
But it is in \cite{OHara2009} that a more practical and robust solution to this boundary problem has been proposed by introducing iterations between scales.
From the raw results, boundary conditions are imposed on the local problems, whose solutions are used as enrichment functions for a new enriched global problem, whose solution is again used as boundary conditions for local problems and finally their solutions reused for a last new enriched global problem.

Later in \cite{Pereira2011a} the idea arises that in evolutionary phenomena such as crack propagation, this \GL scale loop can be intertwined with a propagation loop.
But it was only in \cite{Gupta2012} that the \tS looping strategy was further investigated for a larger number of scale iterations as a remedy for the boundary problem, in a context of a crack propagation step with a single local problem.
This work shows that in a few iterations, the \GL method provides a significant error reduction at both scales.
But, as in many other publication, the way this error is obtained is problematic because it is a relative error based on a known solution for this problem (analytical or numerical fine reference).
Another approach to deal with this boundary problem has been studied in more detail in \cite{Kim2010a} by examining how boundary conditions are imposed at the local scale: Dirichlet, Neumann and Cauchy boundary conditions have been compared, the last one providing a better convergence rate in the treated case (fracture mechanics problems). But in \cite{Plews2015a} it turned out that the Dirichlet boundary condition was the best in the given context.

The question of how the local problem solution is used to construct an enrichment function is another topic widely studied in the literature.
In \cite{Babuska2012} the SGFEM technique proposes to remove  from the fine scale field (obtained here with more conventional analytic functions) its piecewise linear interpolant on the coarse\footnote{In this article, the terms  "fine" or "micro" are sometimes used instead  of the original "local"   designation in the \GL method.
And  "coarse" or "macro" are sometimes used instead of  "global".}-scale field.
This operation ensures a  matrix conditioning number for the enriched system equivalent to that of the non enriched system.
In \cite{Gupta2013,Gupta2015} the SGFEM technique, mentioned to potentially boost the \GL method, was adapted to fracture mechanics with an  analytical enrichment  function improving both matrix conditioning and result quality.

The \GL method has also  been successfully evaluated  in a nonlinear context \cite{Kim2012,Plews2016}, in other mechanical domains \cite{Plews2015,Geelen2020} and even integrated in a commercial software \cite{Li2021}.
But since its inception until today, it is its intrinsic parallel capacities that have been put forward by its creators.
In \cite{Kim2011}, since  the local problems (one per patch) are independent of each other, they can be solved (with a direct solver) in parallel (in multithreaded  paradigm, each thread computes a set of patches).
In \cite{Plews2015a}, the authors successfully reused this parallel implementation and addressed, among many aspects, the issue of memory consumption related to local problems.
Especially when the scale loop is used to cope with the boundary condition issue, storing the factorization of local problems reduces the cost of iterations, at the price of a significant increase in the memory footprint.
In \cite{Li2018} the local problems themselves have been treated with a parallel solver  keeping the idea of computing each of them in parallel.
This dual level of parallelism  was also studied at the same time in \cite{Salzman2019} (in a message passing context).

In \cite{Kim2011,Plews2015a,Li2018} the local problems called "sub-local" domains were constructed from one or a few non-overlapping master local domains,  in order to apply the same discretization for all these "sub-local" domains.
This artifact greatly simplify the integration process of the enriched macro-scale problem since all its enrichment function share the same discretization  in this case, thus  the same integration points.

In this paper, the proposed \tS method retains the idea of a loop between scales to deal with boundary condition problems, but adds the calculation of  a new error criterion  to terminate this loop at a given accuracy. 
This criterion is based on the notion of reference solution which introduces a new way to formulate the enrichment function and to create the global problem.
Otherwise, to avoid any limitation imposed by computer memory, the proposed \tS method  uses the parallel message passing interface (MPI) paradigm  to access distributed memory resources.
All the meshes are distributed over the processors, which induces an original scheduling algorithm, described in  section  \ref{scheduling}, to optimize the load balancing of the processors at the fine scale.
The distributed nature of the mesh at the global scale  also leads to a parallel resolution at this scale.  
Compared to the "sub-local" domain approach, the  proposed \tS method  uses the more systematic approach of a local problem per patch while keeping the notion of a single discretization at the fine scale in order to simplify the parallel implementation, the definition of the enrichment functions and the construction of the problem at the global scale. 
And finally, the chosen boundary conditions imposed on these local problems  are  of Dirichlet type.

The paper is organized as follows.
The next section recalls the classical  global-local formulation.
Then in section \ref{TSMethod}, we introduce a reference solution and reformulate the two-scale approach using a pure matrix format.
Section \ref{theoSeqAlgPerf}  studies the issue of scale difference between the global and local level and its impact on performance.
Section \ref{parallel_paradigm} deals with the parallel global-local approach in which both the global 
and local problems are distributed.
Section \ref{NSOLV} deals with the three numerical experiments and the paper ends with conclusion and possible future works.

\section{\GL ingredients}\label{TSMethodOrig}
\subsection{Continuous mechanical problem}
In this work, we will only consider mechanical elasticity but other fields of application are possible.
 In the  domain of interest $\Omega$,  a  region of $\mathbb{R}^3$, the strong form of the equilibrium  equation of the mechanical problem is defined as follows:
\begin{equation}
		\nabla . \tens{\sigma} +\vm{f}= 0
		\label{equilibrium}
\end{equation}  
with  $\vm{f}$ being a prescribed volume loading and the stress tensor is  $\tens{\sigma}=\tens{C}:\tens{\epsilon}$, $\tens{C}$ being the Hook's tensor and $\tens{\epsilon}$ being the strain tensor.
The Neumann and Dirichlet conditions applied on the  $\partial\Omega$ boundary  of the $\Omega$ domain are:
\begin{equation}
\tens{\sigma}.\vm{n} = \vm{t}~ \text{on}~\partial\Omega^N, \vm{u}=\overline{\vm{u}} ~\text{on} ~\partial\Omega^D, \partial\Omega^N\cap \partial\Omega^D=\oslash
\label{BC}
\end{equation}
where $\vm{t}$ is the prescribed tensile load on the   $\partial\Omega^N$ part of $\partial\Omega$, $\vm{n}$ the outgoing normal vector of $\Omega$, $\overline{\vm{u}}$ the prescribed displacements on the  $\partial\Omega^D$ part of $\partial\Omega$ and $\vm{u}$ the displacement field.
Let $\mathcal{M}^{\mathsf{C}}$  be the continuous space of the problem defined on $\Omega$ and compatible with the Dirichlet boundary conditions ($\mathcal{M}_0^{\mathsf{C}}$ being $\mathcal{M}^{\mathsf{C}}$ with null Dirichlet boundary conditions).
The solution $\vm{u}^C \in \mathcal{M}^{\mathsf{C}}$ of the continuous problem defined by \eqref{equilibrium} and \eqref{BC} in their weak form, is given by solving:
\begin{equation}
	\forall \vm{v}^* \in \mathcal{M}_0^{\mathsf{C}},~~ A\left( \vm{u}^C,\vm{v}^* \right)_{\Omega} = B\left(\vm{v}^*\right)_{\Omega},~~\vm{u}^C=\overline{\vm{u}} ~\text{on} ~\partial\Omega^D
	\label{equilib}
\end{equation}
where  the bi-linear and linear  forms are:
	\begin{equation}
		\begin{array}{l}
				A\left( \vm{u},\vm{v} \right)_{\Omega} =\stretchint{5ex}_{\Omega}
				\tens{\epsilon} \left( \vm{u} \right):\tens{C}:\tens{\epsilon}\left( \vm{v}\right)\mathrm{d}V \\
		 B\left( \vm{v} \right)_{\Omega} =\stretchint{5ex}_{\partial\Omega^N}
		\vm{t}\cdot \vm{v} \mathrm{d}S-\stretchint{5ex}_{\Omega}
		\vm{f}\cdot \vm{v} \mathrm{d}V
	\end{array}
		\text{with}~ \vm{u}\in \mathcal{M}^{\mathsf{C}} ~\text{and} ~\vm{v}\in \mathcal{M}_0^{\mathsf{C}}
		\label{bilinform}
	\end{equation}
and the strain is given by the  kinematic equation 
$\tens{\epsilon}\left( \vm{u} \right)=\frac{1}{2}\left(\nabla\vm{u} +\left(\nabla \vm{u}\right)^t\right)$. 
\subsection{Global problem}\label{glob_prob}
To solve \eqref{equilib},  the \tS approach  discretizes the continuous problem into two steps.
It starts by discretizing $\Omega$ with a mesh composed of  macro (or global or coarse) elements capable of representing the general behavior of the problem. 
On this global mesh,  where $I^g$ denotes its set of  nodes (the superscript g stands for "global scale"), we consider the standard first-order finite element approximation functions $N^i:\vm{x}\in \mathbb{R}^3\rightarrow N^i(\vm{x}) \in \left[0,1\right]\subset\mathbb{R}$, $i\in I^g$, associated to a coarse scale node $\vm{x}^i$, with a support  given by the union of all  finite elements sharing the node $\vm{x}^i$ ( this group of macro-elements is called "patch" hereafter).
These $N^i$ are used to linearly interpolate the degrees of freedom (dof), also called values in this work, which represent the classical values of the discretized displacement field at the  mesh nodes. 

Some nodes of this global mesh are then enriched when their support covers certain local behaviors that must be finely described.
Let $I_{e}^g\subseteq I^g$ be the set of enriched global nodes ($card\left( I_{e}^g\right)\leqslant card\left( I^g\right)$ and the subscript e stands for "enriched").
For each node $p\in I_{e}^g$  a \tS enrichment function $\vm{F}^p(\vm{x})$ (described in the next section) enriches the classical basis  ($N^i$) of the discretization space.
The  kinematic equation used to describe the discrete global field $\vm{U}$  at a point $\vm{x}$ is then:
\begin{equation}
	\vm{U}(\vm{x})=\displaystyle\sum_{i\in I^g}  N^i(\vm{x})~ \vm{U}^{i} + \displaystyle\sum_{p\in I_{e}^g}   N^p(\vm{x})~\vm{E}^{p}*\vm{F}^p(\vm{x})
	\label{kinematic}
\end{equation}
where:
\begin{itemize}
	\item $\vm{U}^{i}$ is the vector of classical values  for the node $\vm{x}^i$ 
	\item $\vm{E}^{i}$ is the vector of enriched values  for the node $\vm{x}^i$ 
	\item $*$ operator is the component  by component  multiplication of two vectors 
\end{itemize} 
Thus \eqref{equilib} becomes to find $\vm{U}^G\in \mathcal{M}^{\mathsf{G}}\subset \mathcal{M}^{\mathsf{C}}$ such that:
\begin{equation}
	\forall \vm{v}^* \in \mathcal{M}_0^{\mathsf{G}},~~ A\left( \vm{U}^G,\vm{v}^* \right)_{\Omega} = B\left(\vm{v}^*\right)_{\Omega},~~\vm{U}^G=\overline{\vm{u}} ~\text{on} ~\partial\Omega^D	
	\label{equilibg}
\end{equation}
with $
\mathcal{M}^{\mathsf{G}}=\left\{\vm{U}(\vm{x}):\vm{U}(\vm{x})=\displaystyle\sum_{i\in I^g}  N^i(\vm{x})~ \vm{U}^{i} + \displaystyle\sum_{p\in I_{e}^g}   N^p(\vm{x})~\vm{E}^{p}*\vm{F}^p(\vm{x}), \vm{U}(\vm{x})=\overline{\vm{u}}(\vm{x}) ~\text{for}~ \vm{x}\in \partial\Omega^D  \right\}$.
The integration of \eqref{equilibg} using appropriate Gauss  quadrature points (see \cite{Kim2011} for example) leads to the resolution of a linear system computed with a direct solver.
This system is constructed  by eliminating  the rows of the elementary matrices related to the Dirichlet values and  moving the coupling term (columns related to the Dirichlet values) to the right side.
\subsection{Local problems and enrichment functions}\label{local_prb}
For a node $p\in I_{e}^g$, the \tS enrichment function $\vm{F}^p(\vm{x})$ is constructed numerically using the solution of a micro (or local or fine) problem.
This local problem  follows the same equations as \eqref{equilibrium} and \eqref{BC}, and is discretized by a micro-mesh corresponding to the refinement (in a nested manner) of the macro elements of the  $p$ patch. 
This is the second discretization step of the \tS method. 
Considering $\omega_p$ the region corresponding to the patch $p$, the solution of the local problem $p$ is obtained by finding $\vm{u}^{Q^p}\in \mathcal{M}^{{\mathsf{Q}}^p}\subset \mathcal{M}^{\mathsf{C}}$ such that:
\begin{equation}
	\forall \vm{v}^* \in \mathcal{M}_0^{{\mathsf{Q}}^p},~~ A\left( \vm{u}^{Q^p},\vm{v}^* \right)_{\omega_p} = B\left(\vm{v}^*\right)_{\omega_p},
	~~	\vm{u}^{Q^p}=\vm{U}^G~ \text{on}~\partial\omega_p
	\label{equilibl}
\end{equation}
where: 
\begin{itemize}
	\item $\vm{U}^G$ is the solution of \eqref{equilibg}
	\item $\mathcal{M}^{{\mathsf{Q}}^p}=\left\{\vm{u}(\vm{x}):\vm{u}(\vm{x}) =\displaystyle\sum_{m\in I_p^l} n^m(\vm{x})~\vm{u}^{m^p},
	~  \vm{u}(\vm{x})=\vm{U}^G(\vm{x})~ \text{for}~\vm{x}\in \partial\omega_p
	\right\}$
	\item $I_p^l$ is the set of nodes associated with the fine mesh discretization of $\omega_p$ (superscript l stands for "local scale")
	\item $n^m:\vm{x}\in \mathbb{R}^3\rightarrow n^m(\vm{x}) \in \left[0,1\right]\subset\mathbb{R}$ is the first order finite element shape function associated to the fine scale node $\vm{x}^m$ 
	\item  $\vm{u}^{m^p}$ is the vector of classical dofs  for the node $\vm{x}^m$ corresponding to the fine scale problem $p$
\item $A\left( \vm{u},\vm{v} \right)_{\omega_p} =\stretchint{5ex}_{\omega_p}
\tens{\epsilon} \left( \vm{u} \right):\tens{C}:\tens{\epsilon}\left( \vm{v}\right)\mathrm{d}V$ with $\vm{u}\in \mathcal{M}^{{\mathsf{Q}}^p}$ and $\vm{v}\in \mathcal{M}_0^{{\mathsf{Q}}^p}$
\item $B\left( \vm{v} \right)_{\omega_p} =\stretchint{5ex}_{\partial\omega_p\cap\partial\Omega^N}
\vm{t}\cdot \vm{v} \mathrm{d}S-\stretchint{5ex}_{\omega_p}
\vm{f}\cdot \vm{v} \mathrm{d}V$ with $\vm{v}\in \mathcal{M}_0^{{\mathsf{Q}}^p}$
\end{itemize}
Regarding the boundary condition of the fine-scale problems, different approaches have been tested (in particular in \cite{Kim2010a}), all of them using the global-scale solution for  Dirichlet, Neumann or  a mix of both.
The idea is to impose the global behavior of the problem on the boundary of the fine-scale problems so that the local behaviors that must be finely described can be revealed by the solutions of the fine-scale problems.
In this work, the global-scale solution is only imposed by the Dirichlet boundary condition as presented in \eqref{equilibl}.
Note that the boundary conditions at the global scale are inherited at the fine scale. In particular the coarse-scale Dirichlet boundary condition  are imposed via $\vm{U}^G$ (since $\vm{U}^G$ is the solution of \eqref{equilibg}, $\vm{U}^G=\overline{\vm{u}} ~\text{on} ~\partial\Omega^D$ and  $\vm{u}^{Q^p}=\vm{U}^G=\overline{\vm{u}} ~\text{on} ~\partial\omega_p\cap \partial\Omega^D$). 

The integration of \eqref{equilibl} for all $p$, using the standard Gauss  quadrature points, lead to the construction of $card\left( I_{e}^g\right)$ linear systems (reduced at assembly by elimination of Dirichlet values).
Their resolution can be done in different ways, but,  as already indicated in the introduction, in \cite{Kim2011} they are solved in parallel using threads.
The distribution of patches to threads uses dynamic scheduling: each thread takes a local problem from a list and processes it (with a direct solver) until the list is empty. 
Kim et al. proved that local problems should be sorted in this list by decreasing cost (calculated  a priori with some metric) in order to maintain a good load balance.
They also showed that the most expensive local problem  limits the number of processes that can be used for a given problem.
Beyond this limit, parallel  scaling efficiency drops.
Section \ref{parallel_paradigm} will explain how we go from a multi-thread patch resolution to a full MPI implementation.

When $\vm{u}^{Q^p}$ are available, the \tS enrichment functions $\vm{F}^p$ can be constructed from them.
A first crud approach is to consider that:
\begin{equation}
\forall p\in I_{e}^g,~\forall \vm{x}\in \omega_p ~\vm{F}^p(\vm{x})=\vm{u}^{Q^p}(\vm{x})~ \text{and}~ \forall \vm{x}\notin \omega_p~\vm{F}^p(\vm{x})=\vm{0}
\label{crudenrich}
\end{equation}
This solution can works if the terms of $N^p(\vm{x})~\vm{F}^p(\vm{x})$  are not to close to $N^p(\vm{x})$.
Otherwise,  the created system is badly conditioned because the discretization basis is almost redundant. 
But in the literature,  many authors use the SGFEM technique \cite{Babuska2012}  which removes  from $\vm{u}^{Q^p}$ its piecewise linear interpolant over the coarse-scale field:
\begin{equation}
	\forall p\in I_{e}^g,~ \forall \vm{x}\in \omega_p ~\vm{F}^p(\vm{x})=\vm{u}^{Q^p}(\vm{x})-\displaystyle\sum_{j\in I_{p}^{g}}  N^j(\vm{x})\vm{u}^{Q^p}(\vm{x}^j)~ \text{and}~ \forall \vm{x}\notin \omega_p~\vm{F}^p(\vm{x})=\vm{0}
	\label{sgfem}
\end{equation}
where $I_{p}^g$ is the set of macro-scale nodes associated with  patch $p$.
Eliminating the projection first removes from the enrichment function the part that can be represented by the coarse-scale discretization and thus improves the matrix conditioning of the coarse-scale problem.
Second, it forces the enrichment function to be zero at coarse-scale  nodes location, which helps solve the blending \footnote{Blending appears in a coarse level element when  all  its classical dofs are not enriched.} element problem and also greatly simplifies\footnote{With a non-zero enrichment function, imposing a Dirichlet condition on an enriched node adds the complexity of setting up an equation linking the classical and enriched dofs.} the application of the global Dirichlet boundary conditions, if any, on those enriched nodes. 
In this work, a rather different approach is proposed, as presented in the section \ref{TS_enrich_strategy}.

\subsection{The scale loop}
The main idea of the \tS solver is to  iteratively solve the two types of problems: 
the global problem \eqref{equilibg} and 
the set of micro  problems \eqref{equilibl} over patches solved independently. 
The global and fine problems are linked together.
The fine problems over the patches deliver the enrichment functions which are used in the global problem.
The global problem delivers the boundary conditions for the patches.

The initial step of the \tS loop is a key element of the method.
In general, a global-scale solution is used as the starting point of the loop (it provides a boundary condition for fine scale problems).
This solution can be the result of a global-scale problem without enrichment or the  result of a previous  computation.
In particular,  when evolutionary phenomena are involved, the \tS solver can effectively take advantage of using the last solution found (from the previous evolutionary step)  as the starting solution of the  \tS loop used for the current evolution step.
Some test cases in the section \ref{NSOLV} illustrate this point.

In any case, all the fine-scale problems are mainly handled during this initialization phase. 
The integration of \eqref{equilibl}  yields a set of linear system which  have  constant matrices  and  partially constant right-hand side vectors during iterations. 
Only the Dirichlet coupling term changes right-hand side vector during the loop as $\vm{U}^G$ varies.
Thus, the  linear systems matrices and constant right-hand side vectors are  created only once during the initialization.
And by choosing a direct solver to solve systems, we can also compute the matrix factorization only once.
If preserved (which has an impact on memory, as shown in \cite{Plews2015a}), these factorizations can be reused for backward/forward resolution of fine-scale problems in the remaining steps of the \tS loop.

\section{The proposed \tS method: Matrix-based resolution with a reference target}\label{TSMethod}
\subsection{Overview}\label{TSoverview}
\begin{figure}[!h]
		\subfloat[][Global-scale problem]{
	\includegraphics[width=0.25\textwidth]{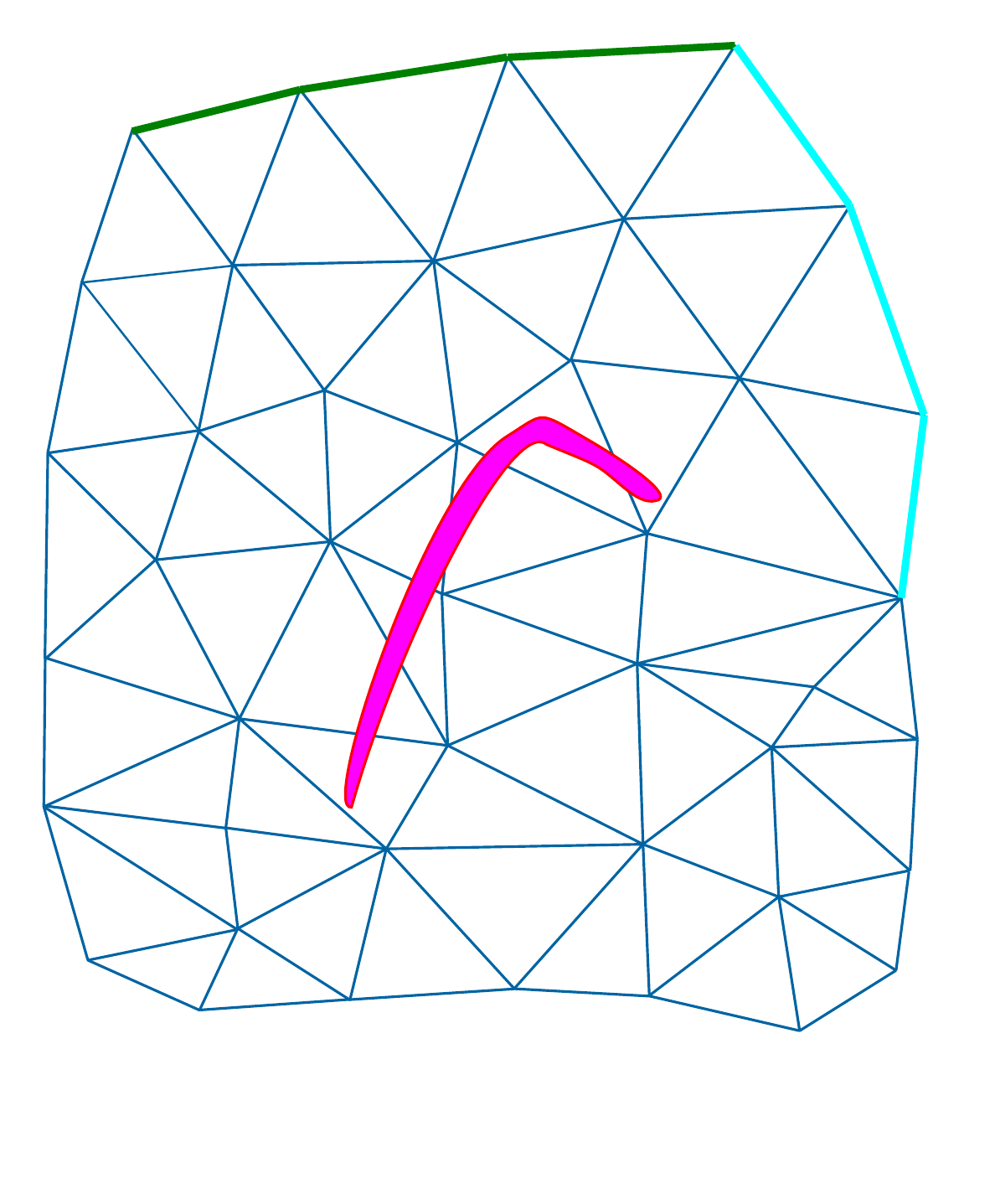}
	\label{ts_coarse_scale}
    }
	\subfloat[][The 25 enriched nodes with 4 of their supports colored]{
		\includegraphics[width=0.25\textwidth]{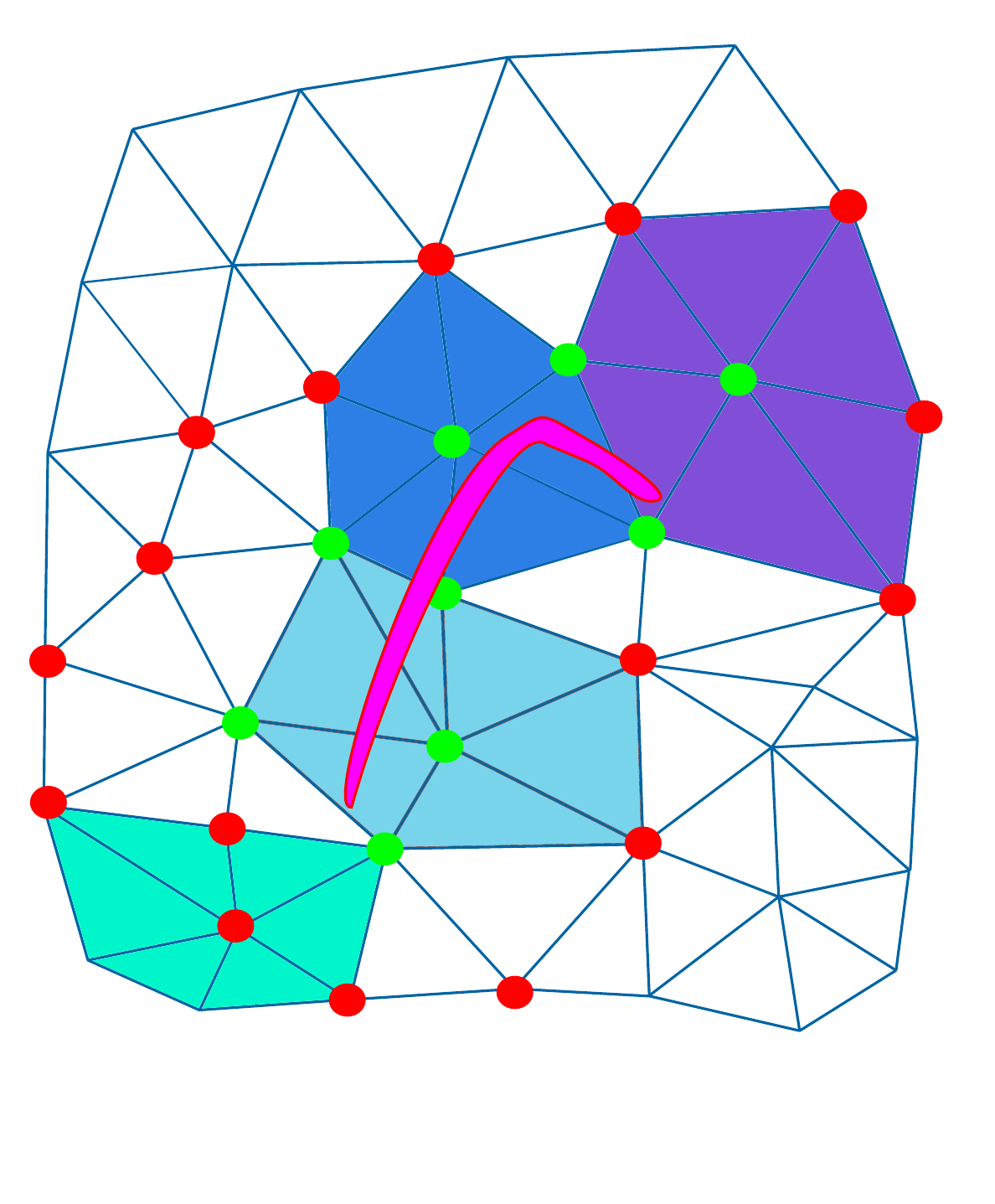}
		\label{ts_enrich}
	}
\subfloat[][\SPn: support union  ]{
	\includegraphics[width=0.25\textwidth]{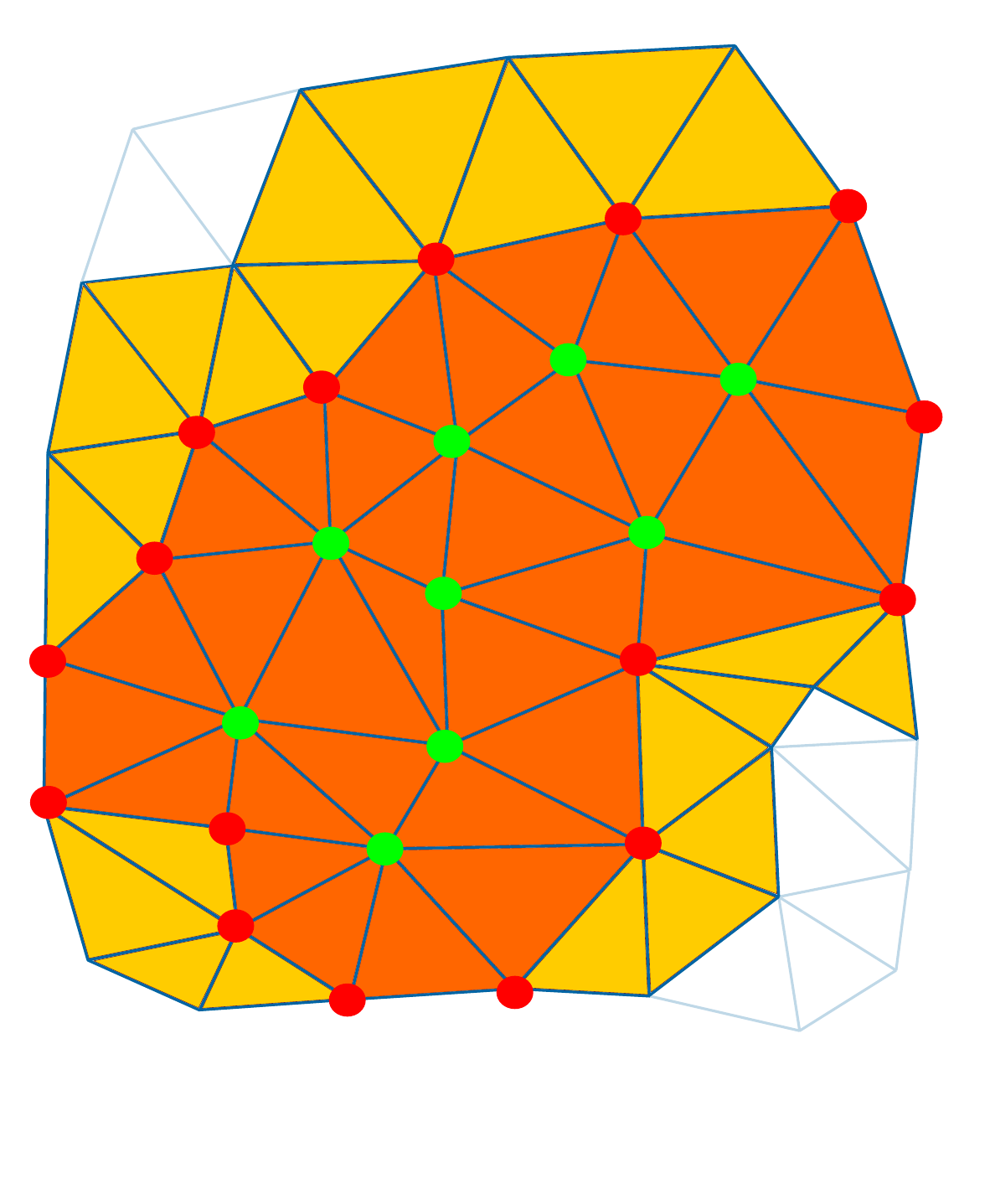}
	\label{ts_sp}
}
\\
	\subfloat[][Fine scale discretisation]{
	\includegraphics[width=0.25\textwidth]{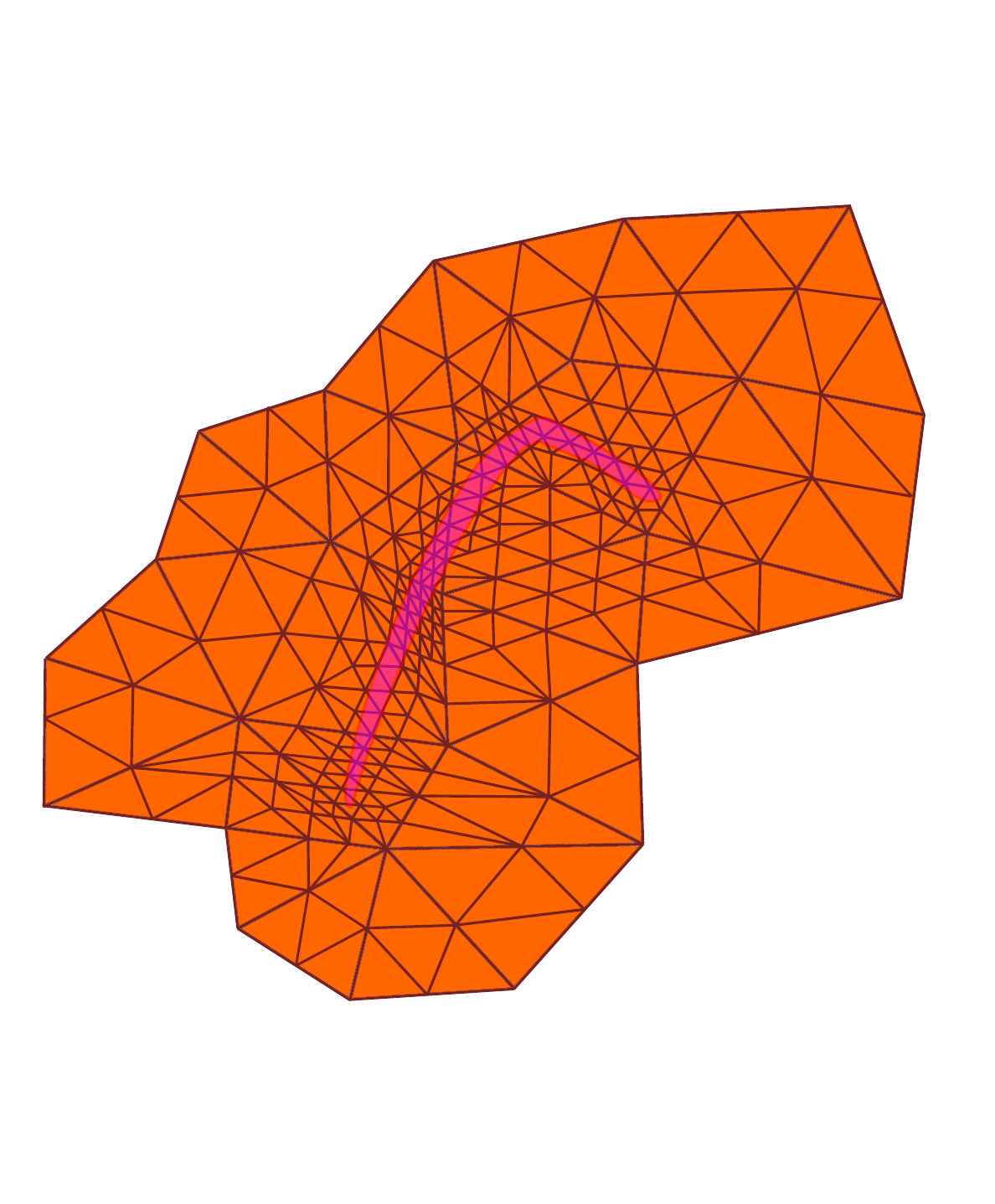}
	\label{ts_fine_scale}
    }
    \subfloat[][Reference field discretization]{
    	\includegraphics[width=0.25\textwidth]{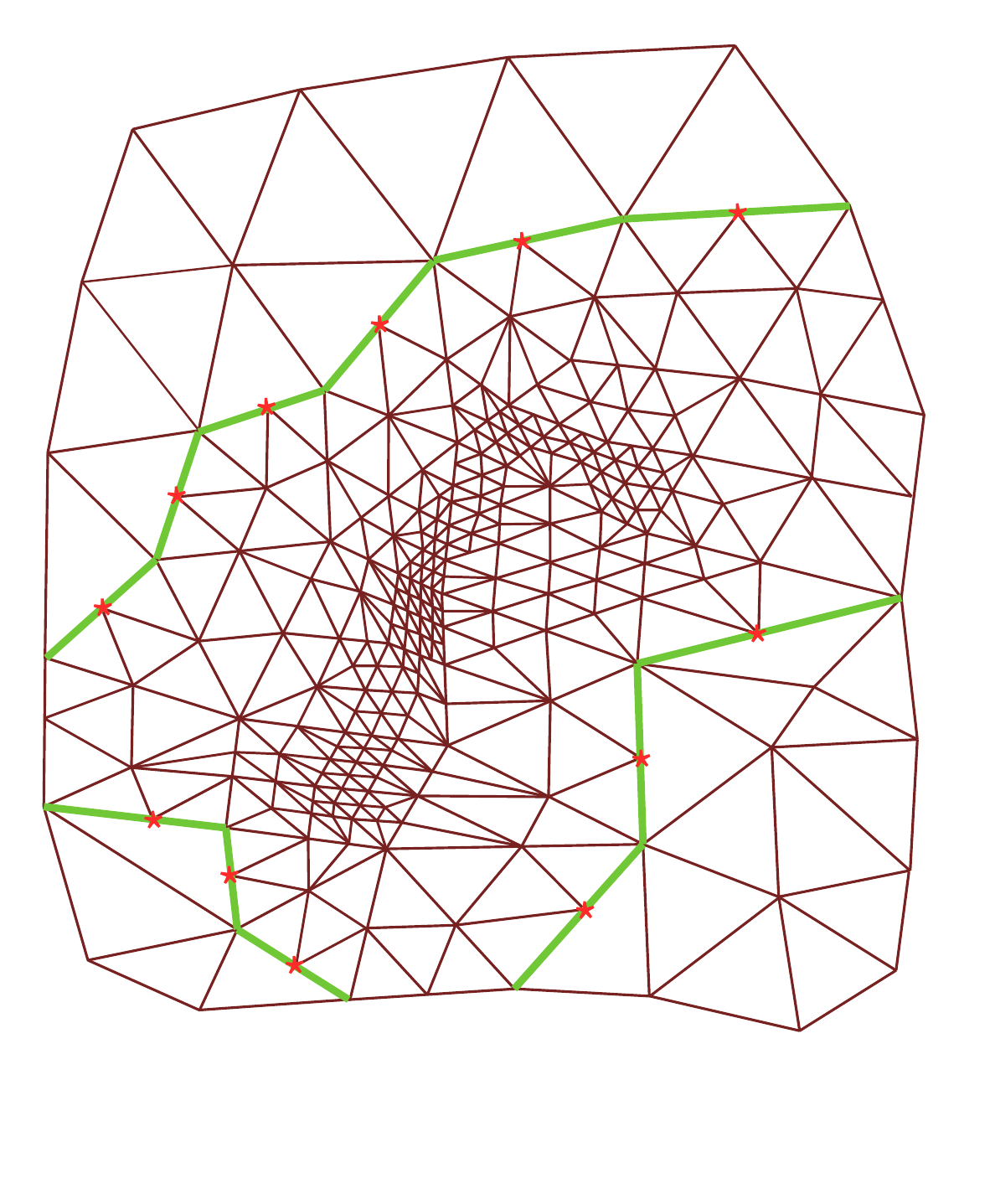}
    	\label{ts_reference_field}
    }
    \subfloat[][Fine scale problems (2 out of 25)]{
	\includegraphics[width=0.25\textwidth]{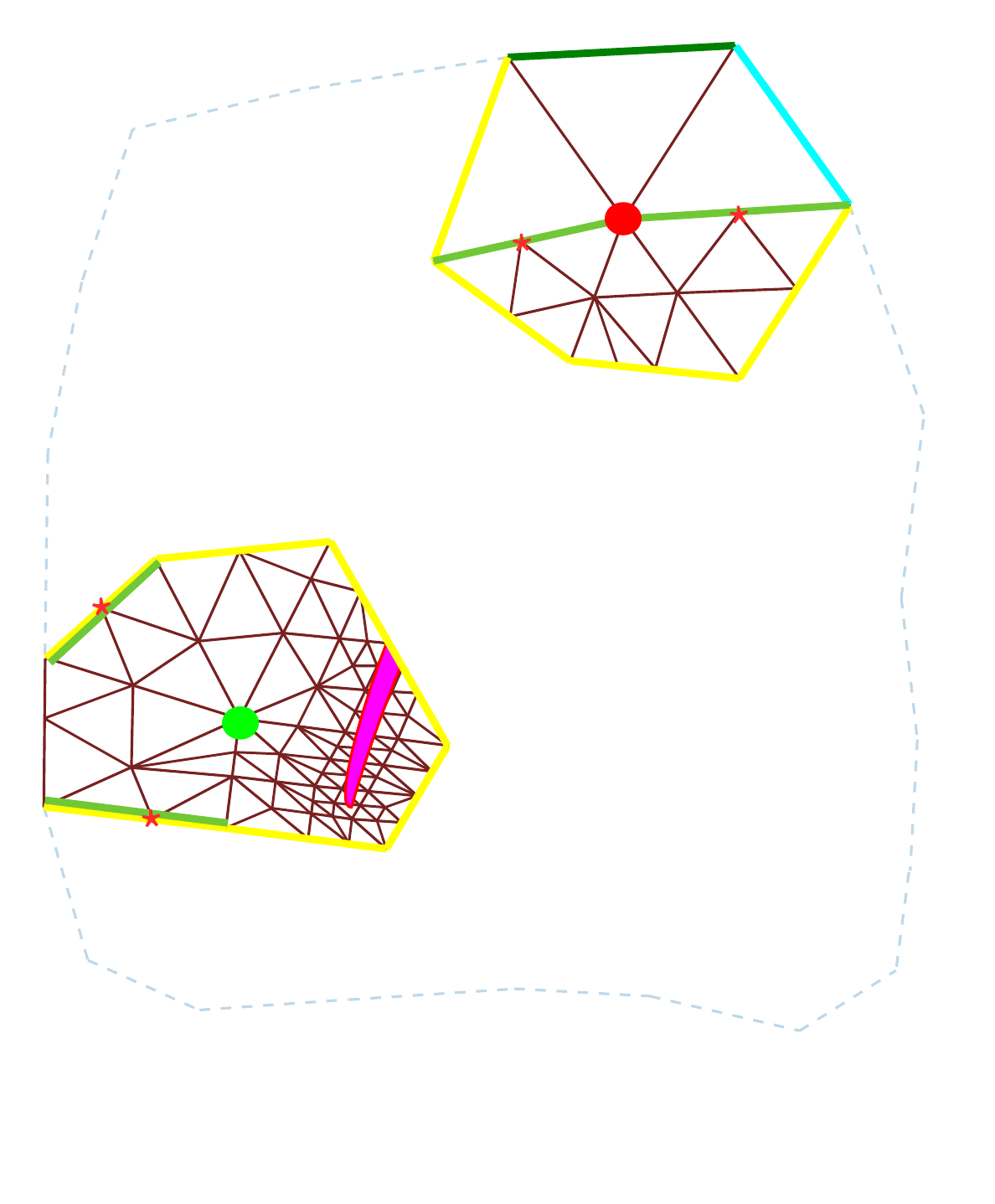}
	\label{ts_fine_prb}
}\\
\subfloat[][Legend]{
\def\svgwidth{0.46\textwidth}
\begingroup%
  \makeatletter%
  \providecommand\color[2][]{%
    \errmessage{(Inkscape) Color is used for the text in Inkscape, but the package 'color.sty' is not loaded}%
    \renewcommand\color[2][]{}%
  }%
  \providecommand\transparent[1]{%
    \errmessage{(Inkscape) Transparency is used (non-zero) for the text in Inkscape, but the package 'transparent.sty' is not loaded}%
    \renewcommand\transparent[1]{}%
  }%
  \providecommand\rotatebox[2]{#2}%
  \ifx\svgwidth\undefined%
    \setlength{\unitlength}{349.16327559bp}%
    \ifx\svgscale\undefined%
      \relax%
    \else%
      \setlength{\unitlength}{\unitlength * \real{\svgscale}}%
    \fi%
  \else%
    \setlength{\unitlength}{\svgwidth}%
  \fi%
  \global\let\svgwidth\undefined%
  \global\let\svgscale\undefined%
  \makeatother%
  \begin{picture}(1,0.32480512)%
    \lineheight{1}%
    \setlength\tabcolsep{0pt}%
    \put(0,0){\includegraphics[width=\unitlength,page=1]{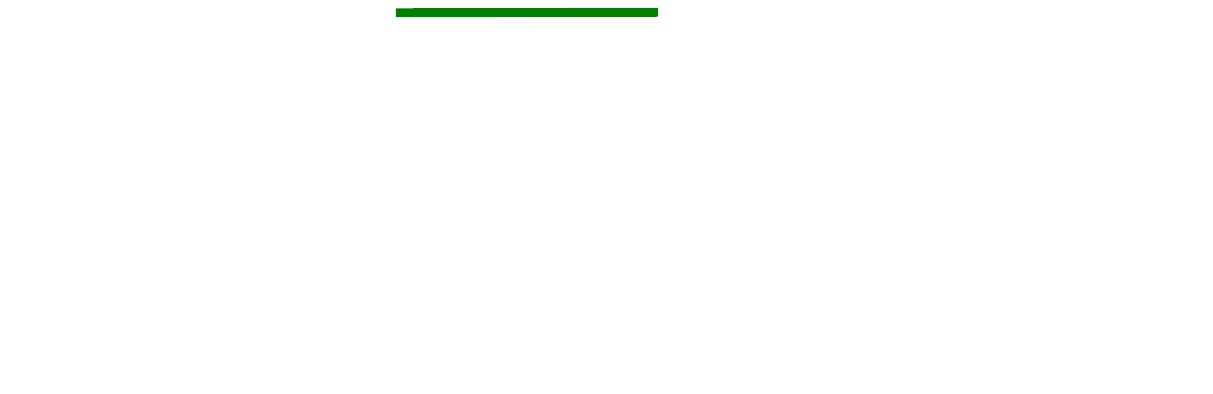}}%
    \put(0.5750374,0.30915628){\color[rgb]{0,0,0}\makebox(0,0)[lt]{\lineheight{0}\smash{\begin{tabular}[t]{l}\scriptsize{Neumann}\end{tabular}}}}%
    \put(0,0){\includegraphics[width=\unitlength,page=2]{Fig01g.pdf}}%
    \put(0.57597232,0.2643017){\color[rgb]{0,0,0}\makebox(0,0)[lt]{\lineheight{0}\smash{\begin{tabular}[t]{l}\scriptsize{Dirichlet}\end{tabular}}}}%
    \put(0,0){\includegraphics[width=\unitlength,page=3]{Fig01g.pdf}}%
    \put(0.57702894,0.21777677){\color[rgb]{0,0,0}\makebox(0,0)[lt]{\lineheight{0}\smash{\begin{tabular}[t]{l}\scriptsize{Dirichlet ($\vm{u}_d^p$)}\end{tabular}}}}%
    \put(0,0){\includegraphics[width=\unitlength,page=4]{Fig01g.pdf}}%
    \put(0.27095745,0.145){\color[rgb]{0,0,0}\makebox(0,0)[lt]{\lineheight{0}\smash{\begin{tabular}[t]{l}\scriptsize{Visual effect to show the location of the local} \end{tabular}}}}%
    \put(0.27095745,0.115){\color[rgb]{0,0,0}\makebox(0,0)[lt]{\lineheight{0}\smash{\begin{tabular}[t]{l}\scriptsize{ behavior to be  captured} \end{tabular}}}}%
    \put(0.03399547,0.30154357){\color[rgb]{0,0,0}\makebox(0,0)[lt]{\lineheight{0}\smash{\begin{tabular}[t]{l}\scriptsize{Boundary conditions: }\end{tabular}}}}%
    \put(0,0){\includegraphics[width=\unitlength,page=5]{Fig01g.pdf}}%
    \put(0.27095745,0.0574839){\color[rgb]{0,0,0}\makebox(0,0)[lt]{\lineheight{0}\smash{\begin{tabular}[t]{l}\scriptsize{Linear relation}\end{tabular}}}}%
    \put(0,0){\includegraphics[width=\unitlength,page=6]{Fig01g.pdf}}%
    \put(0.27479608,0.00599376){\color[rgb]{0,0,0}\makebox(0,0)[lt]{\lineheight{0}\smash{\begin{tabular}[t]{l}\scriptsize{Hanging nodes}\end{tabular}}}}%
  \end{picture}%
\endgroup%
\label{ts_legend}
}
\subfloat{
\def\svgwidth{0.46\textwidth}
\begingroup%
  \makeatletter%
  \providecommand\color[2][]{%
    \errmessage{(Inkscape) Color is used for the text in Inkscape, but the package 'color.sty' is not loaded}%
    \renewcommand\color[2][]{}%
  }%
  \providecommand\transparent[1]{%
    \errmessage{(Inkscape) Transparency is used (non-zero) for the text in Inkscape, but the package 'transparent.sty' is not loaded}%
    \renewcommand\transparent[1]{}%
  }%
  \providecommand\rotatebox[2]{#2}%
  \ifx\svgwidth\undefined%
    \setlength{\unitlength}{399.05079811bp}%
    \ifx\svgscale\undefined%
      \relax%
    \else%
      \setlength{\unitlength}{\unitlength * \real{\svgscale}}%
    \fi%
  \else%
    \setlength{\unitlength}{\svgwidth}%
  \fi%
  \global\let\svgwidth\undefined%
  \global\let\svgscale\undefined%
  \makeatother%
  \begin{picture}(1,0.35919328)%
    \lineheight{1}%
    \setlength\tabcolsep{0pt}%
    \put(0,0){\includegraphics[width=\unitlength,page=1]{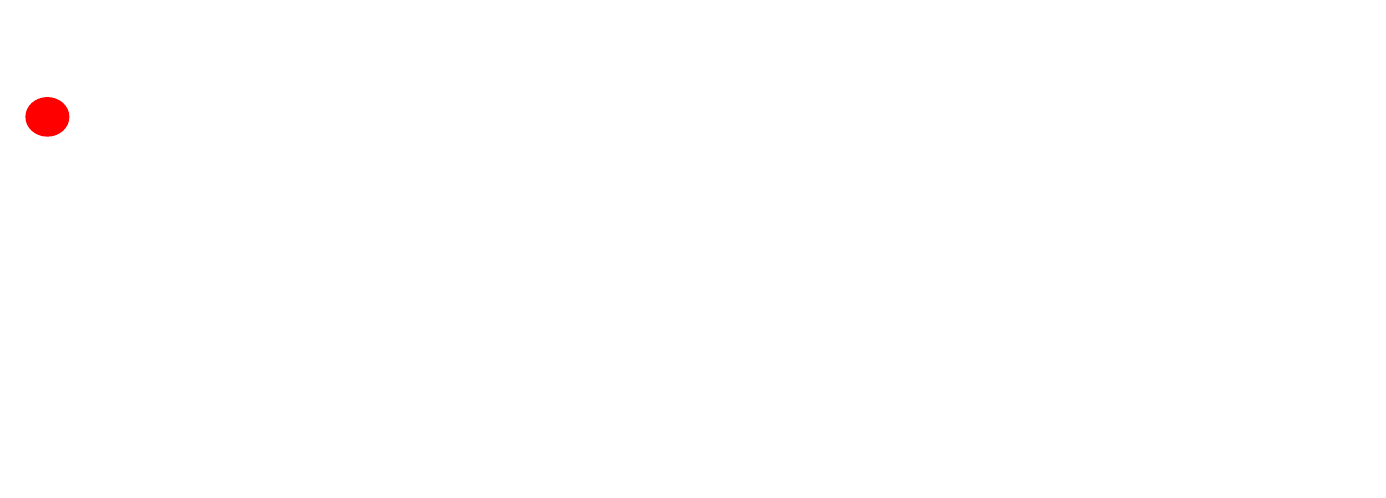}}%
    \put(0.08816039,0.26898148){\color[rgb]{0,0,0}\makebox(0,0)[lt]{\lineheight{0}\smash{\begin{tabular}[t]{l}\scriptsize{Global scale enriched node ($\in I_{e}^g$) with a partially refined patch}\end{tabular}}}}%
    \put(0,0){\includegraphics[width=\unitlength,page=2]{Fig01gb.pdf}}%
    \put(0.08816039,0.33843793){\color[rgb]{0,0,0}\makebox(0,0)[lt]{\lineheight{0}\smash{\begin{tabular}[t]{l}\scriptsize{Global scale enriched node ($\in I_{e}^g$) with a fully refined patch}\end{tabular}}}}%
    \put(0,0){\includegraphics[width=\unitlength,page=3]{Fig01gb.pdf}}%
    \put(0.17760806,0.15915917){\color[rgb]{0,0,0}\makebox(0,0)[lt]{\lineheight{0}\smash{\begin{tabular}[t]{l}\scriptsize{Global-scale mesh}\end{tabular}}}}%
    \put(0,0){\includegraphics[width=\unitlength,page=4]{Fig01gb.pdf}}%
    \put(0.17760806,0.04575297){\color[rgb]{0,0,0}\makebox(0,0)[lt]{\lineheight{0}\smash{\begin{tabular}[t]{l}\scriptsize{Reference mesh}\end{tabular}}}}%
    \put(0,0){\includegraphics[width=\unitlength,page=5]{Fig01gb.pdf}}%
    \put(0.72510056,0.15262199){\color[rgb]{0,0,0}\makebox(0,0)[lt]{\lineheight{0}\smash{\begin{tabular}[t]{l}\scriptsize{\SPn: refined}\end{tabular}}}}%
    \put(0,0){\includegraphics[width=\unitlength,page=6]{Fig01gb.pdf}}%
    \put(0.72510056,0.02882823){\color[rgb]{0,0,0}\makebox(0,0)[lt]{\lineheight{0}\smash{\begin{tabular}[t]{l}\scriptsize{\SPn: not refined }\end{tabular}}}}%
  \end{picture}%
\endgroup%
}
	\caption{\TS method presented in a fictitious 2D problem (for sake of visibility fine-scale mesh is not refined as it should)}
	\label{ts_fig}
\end{figure}

Based on the ingredients of the section \ref{TSMethodOrig}, the figure \ref{ts_fig} illustrates in 2D, but the same principle applies in 3D, the  general concept of the proposed \tS method with a fictitious problem on a $\Omega$ region. 
Figure \ref{ts_coarse_scale} shows the discretization of the global scale problem  with its boundary conditions.
The magenta area is a simple visual effect to show the position of the local behavior to be captured.
This area is just a visual artifact where the mesh needs to be refined  to appropriately simulate a crack, a specific material inclusion, a damaged material, a rivet, ... 
All macro nodes for which a macro element of their support has been refined  are enriched (figure \ref{ts_enrich} where red and green dots form the $I_e^g$ set introduced in the section \ref{glob_prob}).
Their enrichment functions are obtained from the solutions of the fine-scale problems  associated with their patches (figure \ref{ts_fine_prb} where global and local boundary condition  are imposed, $\vm{u}_d^p$ being  introduced in the section \ref{TS_local_prb} and defined precisely  in \ref{anexe_algo}).
The union of all the elements of the enriched patches (yellow + orange in the figure \ref{ts_sp}) is called the \SPn, noted SP thereafter.
Conversely, the set of macro-element that are not part of the SP (shaded in the figure \ref{ts_sp})  will be noted as NSP in what follows.
One of the goals of the \SP is to impose a single fine mesh discretization for all fine-scale problems.
More precisely, the idea is to impose a fine discretization for all patches covering the local behavior to be captured.
Thus their union (in orange on the figure \ref{ts_sp}) is the starting mesh of the adaptation strategy described in  \ref{rsplit_anexe} (see this appendix for details on how the elements are cut/split/refined) which guarantees that all the elements of this area are at least split once and that the local behavior is well refined (figure \ref{ts_fine_scale}).
Using only  orange elements prevents the refinement of any other element, in order to respect the "one hanging node per edge" rule (avoiding the enrichment of more nodes).
The yellow part of the SP is then the set of macro elements not reﬁned and therefore used without modification  at the fine-scale level.
Thus, in the SP, at the boundary between two macro-elements, the micro-meshes are either compatible or leave hanging nodes (visible in figure \ref{ts_reference_field} and presented in section \ref{TS_reffield}). 
This unified discretization imposed in the SP, allows to process  all quantities related to fine-scale grid  (integration at Gauss point, level set, damage, mapping, ... ) only once.
In this work, the union of the SP refined mesh and the NSP macro-elements (which naturally connect  to the yellow fine scale grid elements) is called the reference mesh (figure \ref{ts_reference_field}).
It represented a new vision of what the \GL method intends to discretize at the local and global level alternatively.
The section  \ref{TS_reffield} will detail the problem associated with this mesh.

A rich naming convention is implemented in \ref{anex_sets} to clarify  the different discretizations and the status of the values involved in the different problems solved by the \tS solver.
Thus, all subscribe letters appearing in the matrix or vector in the following algorithms or discussions come from the table \ref{tab_dof_set}.

The general scheme  of this new \tS version is given by the algorithm \ref{TS_algebra}.
\begin{algorithm}[h]
	\footnotesize
	\begin{algorithmic}
		\State Create $\vm{U}_g^0$ macro-scale dofs \Comment{~}
		\State $\left(\vm{A}_{gg}^{ini},\vm{B}_{g}^{ini},\vm{A}_{FF},\vm{B}_F,\left\|\vm{B}_r\right\|,\vm{A}_{qq}^{patches},\vm{BI}_q^{patches},\vm{D}_{qd}^{patches} \right) \gets$ \Call{INIT}{}\Comment{~}
		\State Retrieve macro-scale classical dofs to initialize $\vm{U}_g$ (enriched dof set to zero)
		\State $ \vm{u}_F \gets$ \Call{UPDATE\_MICRO\_DOFS}{$\vm{U}_{g}$}
		\Repeat
		\State $ \vm{u}_q^{patches} \gets$\Call{MICRO-SCALE\_RESOLUTION}{$\vm{u}_F$,  $\vm{A}_{qq}^{patches}$, $\vm{BI}_q^{patches}$, $\vm{D}_{qd}^{patches}$}\Comment{~}
		\State  $\left( \vm{A}_{gg},\vm{B}_{g}\right)\gets$\Call{UPDATE\_MACRO\_PRB}{$\vm{A}_{gg}^{ini}$, $\vm{B}_{g}^{ini}$, $\vm{u}_q^{patches}$}
		\State $\vm{U}_{g}\gets \vm{A}_{gg}^{-1}\cdot \vm{B}_{g}$  {\large$\dagger$} \Comment{~}
		\State $ \vm{u}_F \gets$\Call{UPDATE\_MICRO\_DOFS}{$\vm{U}_{g}$}
		\State $resi\gets$\Call{COMPUTE\_RESIDUAL}{$\vm{u}_F$, $\vm{A}_{FF}$, $\vm{B}_F$, $\left\|\vm{B}_r\right\|$}\Comment{~}
		\Until{ $resi< \epsilon$ }
	\end{algorithmic}
	\caption{\TS algorithm: general procedure.   Procedures INIT, UPDATE\_MACRO\_PRB, UPDATE\_MICRO\_DOFS, MICRO-SCALE\_RESOLUTION and COMPUTE\_RESIDUAL are  respectively depicted in algorithm \ref{TS_algebra_init}, \ref{TS_algebra_macro_update} , \ref{TS_algebra_micro_update}, \ref{TS_algebra_patch} and \ref{TS_algebra_residual} of \ref{anexe_algo}.
	The $\epsilon$ value is the user-defined  target accuracy that the criterion given in the equation \eqref{residual_error} must meet to end the scaling loop. See \ref{convention_anexe} for notation convention : subscripted letters,superscript  $patches$,$\triangleright$, ...}
	\label{TS_algebra}
\end{algorithm}
After an initialization step,
the alternate resolution between scales starts with the MICRO-SCALE\_RESOLUTION procedure (algorithm \ref{TS_algebra_patch}) which computes the solutions of fine-scale problems.
These solutions are used by the  UPDATE\_MACRO\_PRB procedure (algorithm \ref{TS_algebra_macro_update})  to update the linear system at the global-scale ($\vm{A}_{gg}$ and $\vm{B}_{g}$).
Then, a conventional  solver ($\dagger$ in the algorithm \ref{TS_algebra} described in section \ref{gsolv})  solves the problem at the global-scale.
The global scale solution ($\vm{U}_{g}$) is then transferred to the fine scale by the  UPDATE\_MICRO\_DOFS procedure (algorithm \ref{TS_algebra_micro_update}).
Compare to the classical \GL version, the scale loop is now controlled by the computation of a $resi$ criterion (COMPUTE\_RESIDUAL procedure described by the algorithm \ref{TS_algebra_residual}) which stops the iterations when it is lower than a threshold ( $\epsilon$) given by the user.
This  convergence check is an important contribution  because it allows to evaluate the \tS solver and to easily compare  different version of it.
This criterion is introduced in section \ref{TS_reffield} and its computation is given in section \ref{resi_computation}.
The following subsections present in more  detail the  algorithm \ref{TS_algebra} and the associated procedures without focusing too much on the parallelism  which will be discussed in section \ref{parallel_paradigm} (only the presence of the symbol $\triangleright$, on the right of a line of the algorithm, indicates that a certain communication between the processes exists for this step.).

\subsection{The reference field and its associated criterion} \label{TS_reffield}
The reference field $\vm{u}^{R}$ is  the  solution of the  problem defined by  \eqref{equilibrium} and \eqref{BC}, and  discretized by the reference mesh (figure \ref{ts_reference_field}).
Solving this reference problem  is to find $\vm{u}^{R}\in \mathcal{M}^{\mathsf{R}}\subset \mathcal{M}^{\mathsf{C}}$ such that:
	\begin{equation}
		\forall \vm{v}^* \in \mathcal{M}_0^{\mathsf{R}},~~ A\left( \vm{u}^{R},\vm{v}^* \right)_{\Omega} = B\left(\vm{v}^*\right)_{\Omega},
		~~\vm{u}^{R}=\overline{\vm{u}} ~\text{on} ~\partial\Omega^D,
		\label{equilibR}
	\end{equation}
	where: 
	\begin{itemize}
		\item $\mathcal{M}^{\mathsf{R}}=\left\{\vm{u}(\vm{x}):\vm{u}(\vm{x}) =\displaystyle\sum_{k\in I^l} \mathtt{N}^k(\vm{x})~\vm{u}^{k},
		~  \vm{u}(\vm{x})=\overline{\vm{u}}(\vm{x})~ \text{for}~\vm{x}\in  \partial\Omega^D
		\right\}$
		\item $I^l$ is the set of nodes associated with the reference mesh, discretization of $\Omega$
		\item $\mathtt{N}^k$ is either $n^k$ or $N^k$  corresponding respectively to the shape functions defined in \ref{equilibl} or in \ref{equilibg} depending on whether $\vm{x}$ is in a micro or macro element. 
		\item  $\vm{u}^{k}$ is the vector of classical dofs  for the node $\vm{x}^k$ 
\end{itemize}

This solution $\vm{u}^R$  represents the continuous field solution $\vm{u}^C$  only polluted by the discretization error. 
The standard Gauss quadrature  integration of \eqref{equilibR} gives the following linear system without Dirichlet elimination:
\begin{equation}	
	\vm{A}_{RR}\cdot \vm{u}_R^R=\vm{B}_R
	 ~\text{with $\vm{u}_R^R$ the vector of $card(R)$ values defining $\vm{u}^R$}
	\label{R_sys}
\end{equation}	
As already mentioned in  the section \ref{TSoverview}, in the reference mesh there are hanging nodes between the refined and unrefined areas (visible in figure \ref{ts_reference_field}) that need to  be fixed to avoid  discontinuity of the $\vm{u}^{R}$ field.
A simple approach is to eliminate them from the system \eqref{R_sys} by using a linear relationship to fix their value using the following formula: 
\begin{equation}	
	\text{For an hanging node}~\vm{x}^j,~ \vm{u}^{j} =\displaystyle\sum_{i\in I_j}N^i(\vm{x}^j)~ \vm{u}^{i}
	\label{hanging_eq}
\end{equation}   
where:
\begin{itemize}
	\item $\vm{u}^{j}$ is the vector of fixed values (L-set) for the node $\vm{x}^j$ 
	\item $\vm{u}^{i}$ is the vector of free values (f-set) for the node $\vm{x}^i$ 
	\item $I_j$ is the set of nodes corresponding to the vertices of the edge or face on which hanging node $\vm{x}^j$ is located. Thus $card\left( I_j \right) \in \left\lbrace 2,3\right\rbrace $.
	\item $N^i(\vm{x})$ is the standard first-order finite element approximation function associated with the  $\vm{x}^i$ node, related to the coarse element support that holds the edge or face where hanging node is present.
\end{itemize} 
Eliminating the linear relations \eqref{hanging_eq} and the Dirichlet boundary condition of \eqref{R_sys} leads to the system:
\begin{equation}	
	\vm{A}_{rr}\cdot \vm{u}_r^R=\vm{B}_r 
	\label{r_sys}
\end{equation} 
The solution given by a direct solver of this system is  $\vm{u}_r^R$ the vector of $card(r)$ dofs defining $\vm{u}^R$.
But in this work, this system is never solved directly.
Instead, the proposed \tS solver, at each scale iteration, computes a  solution $\vm{u}^{ts}$ ($\vm{u}^{ts} \in \mathcal{M}^{\mathsf{R}}$) such that it quickly converges  to the  reference solution $\vm{u}^R$.

In this work, a $resi$ criterion is proposed to control this convergence.
The $resi$ criterion is taken as the relative residual error of  system \eqref{r_sys}:
\begin{equation}	
	resi=\frac{\left \| \vm{A}_{rr}.\vm{u}_{r}^{ts}-\vm{B}_{r} \right \|}{\left \|\vm{B}_{r} \right \|} 
	\label{residual_error}
\end{equation}
where $\left \| . \right \|$ is the L2-norm and $\vm{u}_r^{ts}$ is the vector of $card(r)$ dofs   corresponding to the $\vm{u}^{ts}$ solution obtained by the \tS solver.
For simplicity, in what follows, $\vm{u}_r$ will represent the $\vm{u}_r^{ts}$ vector of the field $\vm{u}^{ts}$ where  superscripts $ts$ is made implicit.
The question of computing \eqref{residual_error} is exposed in the following section.

To  relate this  criterion to an energy error, an energy norm is defined  as follows using \eqref{bilinform}:
\begin{equation}
		\left\| \vm{u} \right\|_{E_{\Omega}} =\sqrt{A\left( \vm{u},\vm{u} \right)_{\Omega}}=\sqrt{\stretchint{5ex}_{\Omega}
			\tens{\epsilon} \left( \vm{u} \right):\tens{C}:\tens{\epsilon}\left( \vm{u}\right)\mathrm{d}\Omega}
		\label{energynorm}
\end{equation}
The error introduced by the \tS solver (mainly on the boundary conditions imposed on the fine scale problems) in the indirect resolution of system \ref{r_sys} compare to its direct resolution can thus be expressed by $\left\| \vm{u}^{ts}-\vm{u}^R \right\|_{E_{\Omega}}$.
It is related to $resi$ by $\vm{A}_{rr}\cdot \vm{u}_{r}-\vm{B}_{r}$ vector:
\begin{equation}	
	\left\| \vm{u}^{ts}-\vm{u}^R \right\|_{E_{\Omega}}=\sqrt{A\left( \vm{u}^{ts}-\vm{u}^R, \vm{u}^{ts}-\vm{u}^R\right)}=\sqrt{\left( \vm{A}_{rr}\left( \vm{u}_{r}-\vm{u}_r^R\right)\right) \cdot \left( \vm{u}_{r}-\vm{u}_r^R\right)}=\sqrt{\left(\vm{A}_{rr}.\vm{u}_{r}-\vm{B}_{r}\right)\cdot \left( \vm{u}_{r}-\vm{u}_r^R\right)}
	\label{lien_resi_energ}
\end{equation} 
So when $\vm{A}_{rr}\cdot \vm{u}_{r}-\vm{B}_{r}$ vector tends to the zero vector, both $resi$ and $\left\| \vm{u}^{ts}-\vm{u}^R \right\|_{E_{\Omega}}$ tend to zero.

It is also interesting to  compare the \tS solution to the continuous solution.
The associated energy error is  $\left\| \vm{u}^{ts}-\vm{u}^C \right\|_{E_{\Omega}}$  and can be expressed as:
 \begin{equation}
	\begin{array}{r}
		\left\| \vm{u}^{ts}-\vm{u}^C \right\|_{E_{\Omega}}^2 =\stretchint{5ex}_{\Omega}
		\tens{\epsilon} \left( \vm{u}^{ts} -\vm{u}^R + \vm{u}^R-\vm{u}^C \right):\tens{C}:\tens{\epsilon}\left( \vm{u}^{ts}-\vm{u}^R + \vm{u}^R-\vm{u}^C\right)\mathrm{d}\Omega =\\
		\left\| \vm{u}^{ts}-\vm{u}^R \right\|_{E_{\Omega}}^2+\left\| \vm{u}^R-\vm{u}^C \right\|_{E_{\Omega}}^2+2A\left( \vm{u}^R-\vm{u}^C,\vm{u}^{ts}-\vm{u}^R \right)_{\Omega} 	\end{array}
	\label{scale error}
\end{equation}
Subtracting \eqref{equilib} to \eqref{equilibR} gives as $\mathcal{M}^{\mathsf{R}}\subset \mathcal{M}^{\mathsf{C}}$: 
\begin{equation}
		\forall \vm{v}^* \in \mathcal{M}_0^{\mathsf{R}}, A\left( \vm{u}^R-\vm{u}^C,\vm{v}^* \right)_{\Omega} =B\left(\vm{v}^*\right)_{\Omega}-B\left(\vm{v}^*\right)_{\Omega}=0
		\label{equilibA-R}
\end{equation}
Based on \eqref{equilibA-R} and as $\vm{u}^{ts}-\vm{u}^R \in \mathcal{M}_0^{\mathsf{R}}$, then the extra term  of \eqref{scale error}  vanishes:
\begin{equation}
	A\left( \vm{u}^R-\vm{u}^C,\vm{u}^{ts}-\vm{u}^R \right)_{\Omega} =0
	\label{Nullterm}
\end{equation}
Leading to the relation:
\begin{equation}	
	\left\| \vm{u}^{ts}-\vm{u}^C \right\|_{E_{\Omega}}^2 =\left\| \vm{u}^{ts}-\vm{u}^R \right\|_{E_{\Omega}}^2+\left\| \vm{u}^R-\vm{u}^C \right\|_{E_{\Omega}}^2
	\label{equality}
\end{equation}

The second term on the right-hand side of \eqref{equality} represents the error introduced by the reference mesh discretization with respect to the continuous solution.
The first one is the error introduced by the \tS method.
Most publications on the subject do compute $\vm{u}^R$ with another solver and verify a posteriori that $\left\| \vm{u}^{ts}-\vm{u}^C \right\|_{E_{\Omega}}$ tends to $\left\| \vm{u}^R-\vm{u}^C \right\|_{E_{\Omega}}$ (i.e. that $\left\| \vm{u}^{ts}-\vm{u}^R \right\|_{E_{\Omega}}$ tends to zero) in few iterations.
In this work $resi$ carries out this verification during the iterations.
A numerical comparison of these errors is given in section \ref{UOL}.

\subsection{Specific global-scale system construction: the algebraic operator}\label{TS_integ_schem}
The $\vm{u}_{r}$  (as mentioned above, the simplified notation of $\vm{u}_r^{ts}$) vector will act as a bridge between computational tasks.
It stores  classical values at the nodes of $\vm{u}^{ts}$ field which will be used to compute $resi$ and impose boundary conditions on fine-scale problems.
The $\vm{u}_{r}$ dofs of this field can be related to global-scale field dofs via the kinematic equation \eqref{kinematic}.
With this equation, for any node $\vm{x}^j$ of the reference mesh, it is possible to write, without any approximation, that the \tS solution  $\vm{u}^{ts}$ is equal to the  global enriched discrete solution at this location: 
\begin{equation}
	\vm{u}^{ts}(\vm{x}^j)=\vm{U}^G(\vm{x}^j)=\displaystyle\sum_{i\in I^g}  N^i(\vm{x}^j)~ \vm{U}^{i} + \displaystyle\sum_{p\in I_{e}^g}   N^p(\vm{x}^j)~\vm{E}^{p}*\vm{F}^p(\vm{x}^j)
	\label{relation}
\end{equation}
This can be rewritten in algebraic form as:
\begin{equation}
	\vm{u}_M=\vm{T}_{MG}\cdot \vm{U}_G
	\label{u_ru_g}
\end{equation} 
where $\vm{T}_{MG}$ is a linear  interpolation operator whose structure and construction are given in \ref{Trg_anexe}.
In this work, the process of eliminating linear relations \eqref{hanging_eq} is made transparent by choosing,  from the R-set, to  first eliminate the L-set before applying  linear interpolation.
Thus, the operator in \eqref{u_ru_g} interpolates from the G-set to the M-set.

Considering the matrices of \eqref{R_sys} with the  L-set eliminated, and the equation \eqref{u_ru_g} to apply a change of variable,  the equality between the potential energy and the work of  external forces can be written as follows:
\begin{equation}
	\vm{U}_G^t\cdot \vm{T}_{MG}^t\cdot  \vm{A}_{MM} \cdot \vm{T}_{MG}\cdot \vm{U}_G=\vm{U}_G^t\cdot \vm{T}_{MG}^t\cdot \vm{B}_{M}
	\label{energy_g}
\end{equation} 
This gives the following system (see \ref{Agg_anexe} for the detailed construction):
\begin{equation}
	\vm{A}_{gg}\cdot \vm{U}_g=\vm{B}_g
	\label{g_sys}
\end{equation}

From a practical point of view, the sub-blocks of $\vm{T}_{MG}$, $\vm{A}_{MM}$ and  $\vm{B}_{M}$ are large sparse matrices, so they are never fully assembled in memory but stored in blocks using a per macro-element strategy.
This strategy, which gathers in memory all the matrices and vectors necessary for the calculation relative to a macro element, is favorable to the reduction of cache misses and to multi-threading computation.
Moreover, during the \tS loop, only the terms  related to  $\vm{F}^p$ (related to $\vm{T}_{Fe}$ operator given in \ref{Agg_anexe} ) have to be updated at each scale iteration. Thus some sub-block of $\vm{A}_{gg}$ and $\vm{B}_g$  are constant during the iteration and are computed and assembled only once at initialization time.

The INIT procedure (algorithm \ref{TS_algebra_init} detailed in \ref{anexe_algo} ) is in charge of this initialization.
Let $\omega^{e_{macro}}$ be the region of $\Omega$ covered by a macro element $e_{macro}$.
$\omega^{e_{macro}}$ is discretized by one or more  micro-elements of the reference mesh.
With two loops on all macro elements, the INIT procedure integrates and assembles \eqref{equilibR}  over $\omega^{e_{macro}}$ with the elimination of \eqref{hanging_eq} but keeping Dirichlet terms.
The resulting matrices and vectors  $\vm{A}_{FF}^{e_{macro}}$ and $\vm{B}_{F}^{e_{macro}}$  are stored for all SP elements.
From this first matrix assembly of $\vm{A}_{FF}$ and $\vm{B}_{F}$ per block and the  creation  of the  sub blocks of  $\vm{T}_{MG}$, all the constant terms of $\vm{A}_{gg}$ and $\vm{B}_g$ can be computed algebraically  and assembled in two dedicated memory areas, $\vm{A}_{gg}^{ini}$ and $\vm{B}_{g}^{ini}$.

\subsection{global to local transfer}\label{TS_global_to_local}
In this work the idea is to obtain the $\vm{u}^{ts}$ field not by solving directly \eqref{r_sys}  but by solving \eqref{g_sys} and use \eqref{u_ru_g} to go from $\vm{U}^G$ to $\vm{u}_F$.
This last point is achieved by the UPDATE\_MICRO\_DOFS procedure (algorithm \ref{TS_algebra_micro_update} detailed in \ref{anexe_algo}) which uses the same  per-macro-element storage strategy and stores the $\vm{u}^{ts}$ field in the $\vm{u}_F^{e_{macro}}$ vectors attached to SP macro-elements.
Note that this procedure does not loop over the NSP elements and that no information about the h-set is stored in memory because in the computational tasks related to $\vm{u}^{ts}$, this information can be ignored as we will see later.
\subsection{Local problem}\label{TS_local_prb}
In this work, as in the section \ref{TSMethodOrig}, the solution of the local problem  for a patch $p$ will be given by solving \eqref{equilibl} but using a different mechanism to create the linear system.
Consider that $\omega_p$ is equal to $\bigcup\limits_{e_{macro}\in J^p}\omega_{e_{macro}}$, with $J^p$ the set of macro-elements of the patch $p$.
As  the local problems are parts of the reference problem sharing with it its definition  and discretization, the integration of \eqref{equilibl} has already been partially done  when integrating \eqref{equilibR}  over $\omega^{e_{macro}}$.   
Thus, creating the matrix and vector for a patch $p$ is just  assembling  the blocks $\vm{A}_{FF}^{e_{macro}}$ and $\vm{B}_F^{e_{macro}}$ with $e_{macro}\in J^p$.

As for the Dirichlet boundary conditions on $\partial\omega_p$, they are directly obtained from the $\vm{u}^{ts}$ field ($\vm{u}_F^{e_{macro}}$ dofs).
It is strictly the same as using $\vm{U}^G$ because $\vm{u}_F$ which comes from \eqref{u_ru_g} is computed with $\vm{T}_{MG}$ constructed from \eqref{kinematic} as the space $\mathcal{M}^{\mathsf{G}}$ constructed from same kinematic equation.
The application of these boundary conditions gives the following set of local linear systems:
\begin{equation}
	\forall p\in I_e^g~\vm{A}_{qq}^p\cdot \vm{u}_q^p=\vm{B}_q^p
	\label{p_sys}
\end{equation}
Again, it is the INIT procedure, with a loop over the patches, that is responsible for  creating  $\vm{A}_{qq}^p$ which is constant during the \tS iteration.
Similarly, the constant part of the vector $\vm{B}_q^p$, called  $\vm{BI}_q^p$ and the constant coupling Dirichlet term $\vm{D}_{qd}^p$ are computed in this loop over the patches (see \ref{anexe_algo} for more details). 
Note that the use of Dirichlet boundary conditions for local problems  simplifies the software implementation and offers good performance for solving  \eqref{p_sys}: smaller system size, simple constant algebraic operator ($\vm{D}_{qd}^p$) for imposing the prescribed displacements, and the matrix is not polluted by additional parameters (springs) to be set.
But the consequence of this choice on the proposed \tS loop has not been studied and should be in a future work.

\subsection{The \tS loop}\label{TS_loop}
When the INIT procedure is completed, all  local structures associated with  macro-elements or patches are allocated and partially computed.
Before entering the \tS loop,  the algorithm \ref{TS_algebra}  initializes the field $\vm{u}^{ts}$  by calling  UPDATE\_MICRO\_DOFS with   the vector $\vm{U}_g$, solution of a specific macro-scale problem  (in general  the unenriched macro-scale problem but it can come from another computation).	
 
The algorithm \ref{TS_algebra}  then enters  the \tS loop and  first computes the solutions $\vm{u}_q^p$ of \eqref{p_sys} by calling the  MICRO-SCALE\_RESOLUTION procedure (algorithm \ref{TS_algebra_patch} detailed in \ref{anexe_algo}).
This procedure creates vectors $\vm{B}_q^p$   from the constant terms and the coupling Dirichlet terms multiplied by the current subset of $\vm{u}_F$.
It then solves the systems with a direct solver ($\ddagger$ of the algorithm \ref{TS_algebra_patch}),  factoring only once the matrices.

Having all the $\vm{u}_q^p$, the global to local step is done and the local to global step can start by calling the  UPDATE\_MACRO\_PRB procedure ( algorithm \ref{TS_algebra_macro_update} detailed in \ref{anexe_algo}).
This procedure computes the enrichment function from the solution $\vm{u}_q^p$  and using  $\vm{A}_{gg}^{ini}$ and $\vm{B}_{g}^{ini}$ finishes the algebraic computation of $\vm{A}_{gg}$ and $\vm{B}_{g}$.
The system \eqref{g_sys} is then solved ($\dagger$ of the algorithm \ref{TS_algebra}) providing $\vm{U}_g$ solution.

A call to  UPDATE\_MICRO\_DOFS with this new $\vm{U}_g$ solution updates $\vm{u}^{ts}$ field so that the $resi$ criterion   \eqref{residual_error} can be computed for the current iteration.
This is done by COMPUTE\_RESIDUAL procedure described in the  next section. 
Depending on whether $resi$ is greater than the given  $\epsilon$, the loop must continue or not.

\subsection{residual criterion computation}\label{resi_computation}
The numerator of the equation \eqref{residual_error} can be divided as follows:
\begin{equation}	
	resi=\frac{\left \| \begin{array}{c}
			\vm{A}_{hr}.\vm{u}_{r}-\vm{B}_{h}\\
			\vm{A}_{fr}.\vm{u}_{r}-\vm{B}_{f}
		\end{array} \right \|}{\left \|\vm{B}_{r} \right \|} 
	\label{residual_error_split}
\end{equation}
In this work, the h-set part of the residual vector is considered as zero.
This  is justified by the fact that the discretization of the NSP elements  at both scales is the same: the h-set part of the r-set is the same as the h-set part of the g-set.
Thus, since the  $\dagger$ resolution, when performed by a direct solver, provides a "zero" residual error of the g-set system (zero at the machine accuracy), all  this residual vector is zero and in particular its h-set rows.
The computation of $resi$ is therefore simplified to:  
\begin{equation}	
	resi=\frac{\left \| 
		\vm{A}_{fr}\cdot \vm{u}_{r}-\vm{B}_{f}
		\right \|}{\left \|\vm{B}_{r} \right \|} 
	\label{residual_error_simp}
\end{equation}   
The organization of the data by macro-element both simplifies and complicates the computation of  \eqref{residual_error_simp} performed by the COMPUTE\_RESIDUAL procedure  (algorithm \ref{TS_algebra_residual} detailed in \ref{anexe_algo}). 
The simplicity comes from the computation by $e_{macro}$  of the vector $\vm{A}_{FF}^{e_{macro}}\cdot\vm{u}_{F}^{e_{macro}}-\vm{B}_{F}^{e_{macro}}$ with all the  vectors and the matrix available as data attached  to the $e_{macro}$ element. 
The complexity comes from the accumulation of the special local scalar product of these vectors to obtain term $\left \| \vm{A}_{fr}\cdot \vm{u}_{r}-\vm{B}_{f} \right \|$ .

The $\left \|\vm{B}_{r} \right \|$ term of the equation \eqref{residual_error_simp} is computed by the procedure COMPUTE\_B\_NORM  described by the algorithm \ref{TS_algebra_B_norm} in \ref{anexe_algo}. 
This algorithm is very similar to the algorithm \ref{TS_algebra_residual} in that the  scalar product of the vector $\vm{B}_r$  by itself is transformed into  a sum of local dot product.
But as $\vm{B}_r$ remains unchanged during the  \tS loop, this procedure is called only  once during the initialization by the  INIT procedure.
The $\left \|\vm{B}_{r} \right \|$ term is therefore supplied as an argument to the COMPUTE\_RESIDUAL procedure  which computes the $resi$ numerator and just divide by this argument to obtains $resi$.

\subsection{\tS enrichment function}\label{TS_enrich_strategy}
\begin{figure}[h]
	\subfloat[$\vm{u}^{Q^p}(\vm{x})$]{\includegraphics[width=0.26\textwidth]{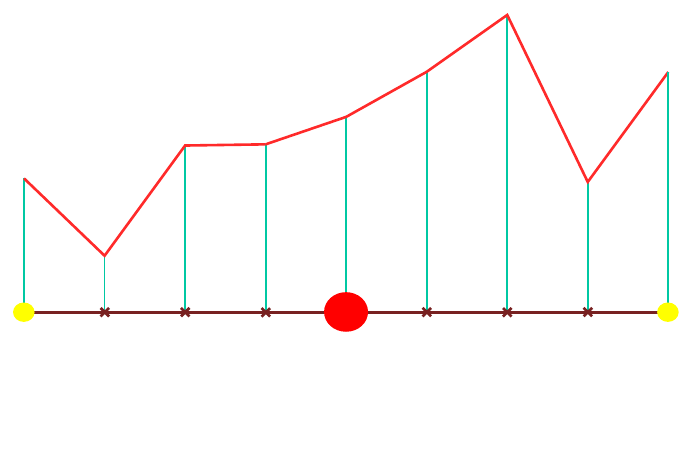}\label{enriche1D_construct_a}}
	\subfloat[$\vm{F}^p(\vm{x})$]{\includegraphics[width=0.26\textwidth]{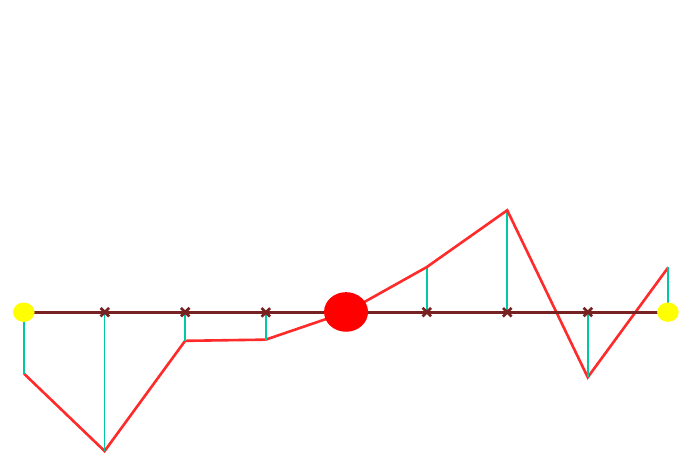}\label{enriche1D_construct_b}}
	\subfloat[$N^p(\vm{x})~\vm{F}^p(\vm{x})$]{\includegraphics[width=0.26\textwidth]{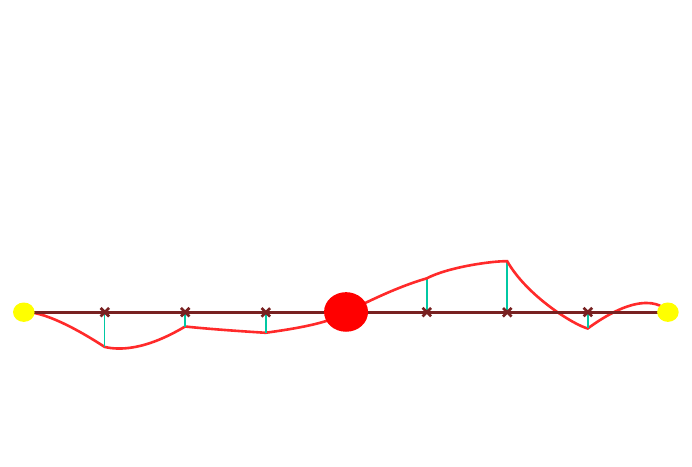}\label{enriche1D_construct_c}}
	\subfloat[$\displaystyle\sum_{m\in {\bar{I}_p^l}}  n_m(\vm{x})N^p(\vm{x}_m)~\vm{F}^p(\vm{x}_m)$]{\includegraphics[width=0.26\textwidth]{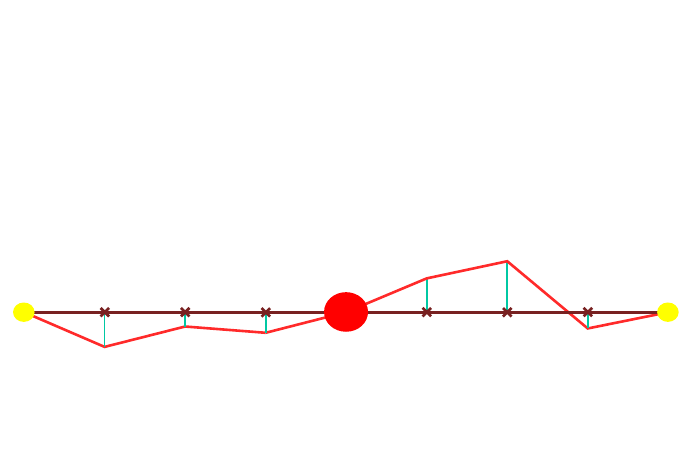}\label{enriche1D_construct_d}}
	\caption{Construction of an interpolated generalized approximation function  for a 1D patch $p$.
		For  visualization  purposes, the quantities   here are  artificial elevations of the patch mesh. 
		\label{enriche1D_construct}}
\end{figure}

\begin{figure}[h]
	\subfloat[$\vm{u}^{Q^p}(\vm{x})$]{
\includegraphics[width=0.26\textwidth]{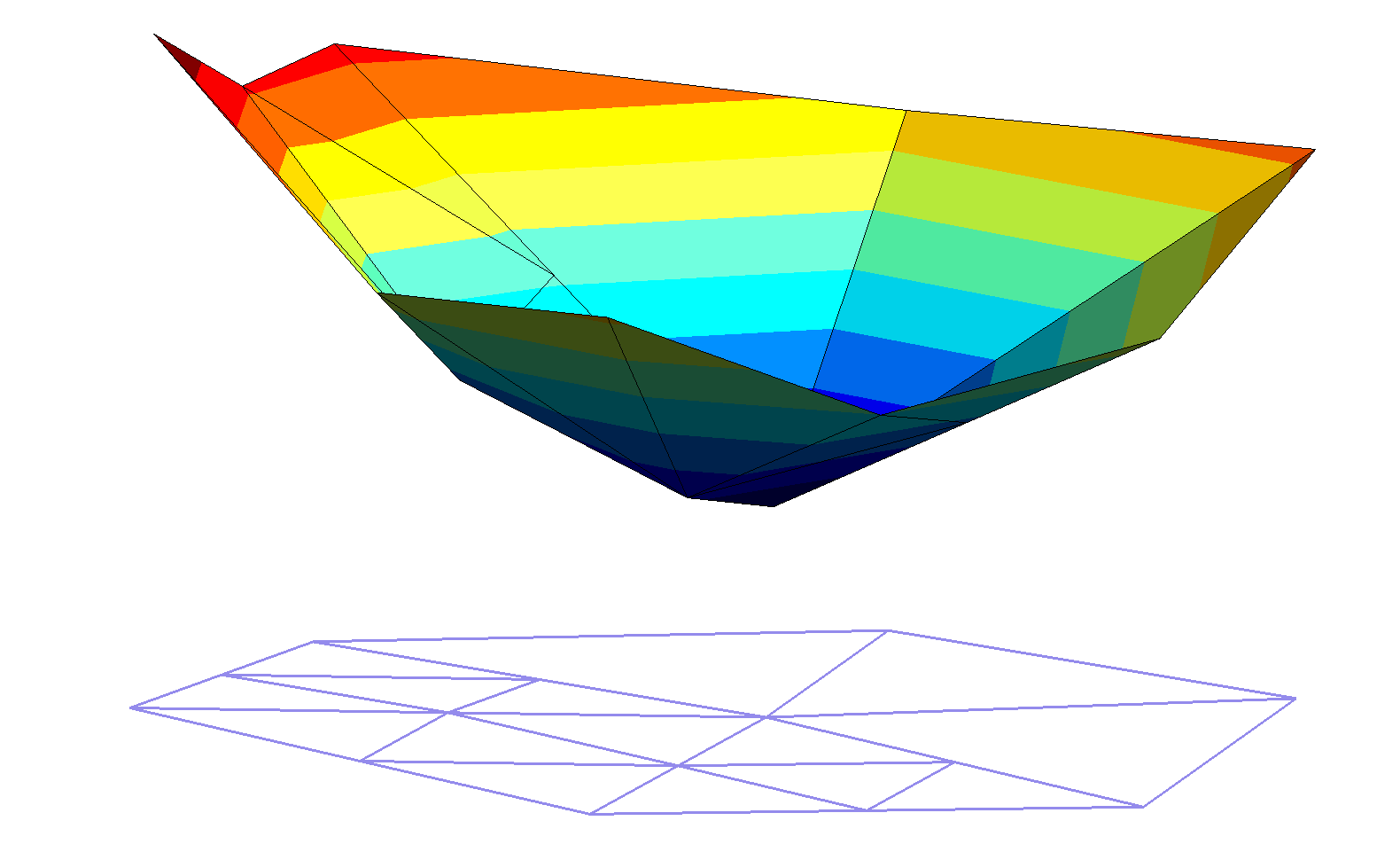}
		\label{enriche_construct_a}}
	\subfloat[$\vm{F}^p(\vm{x})$]{
		\includegraphics[width=0.26\textwidth]{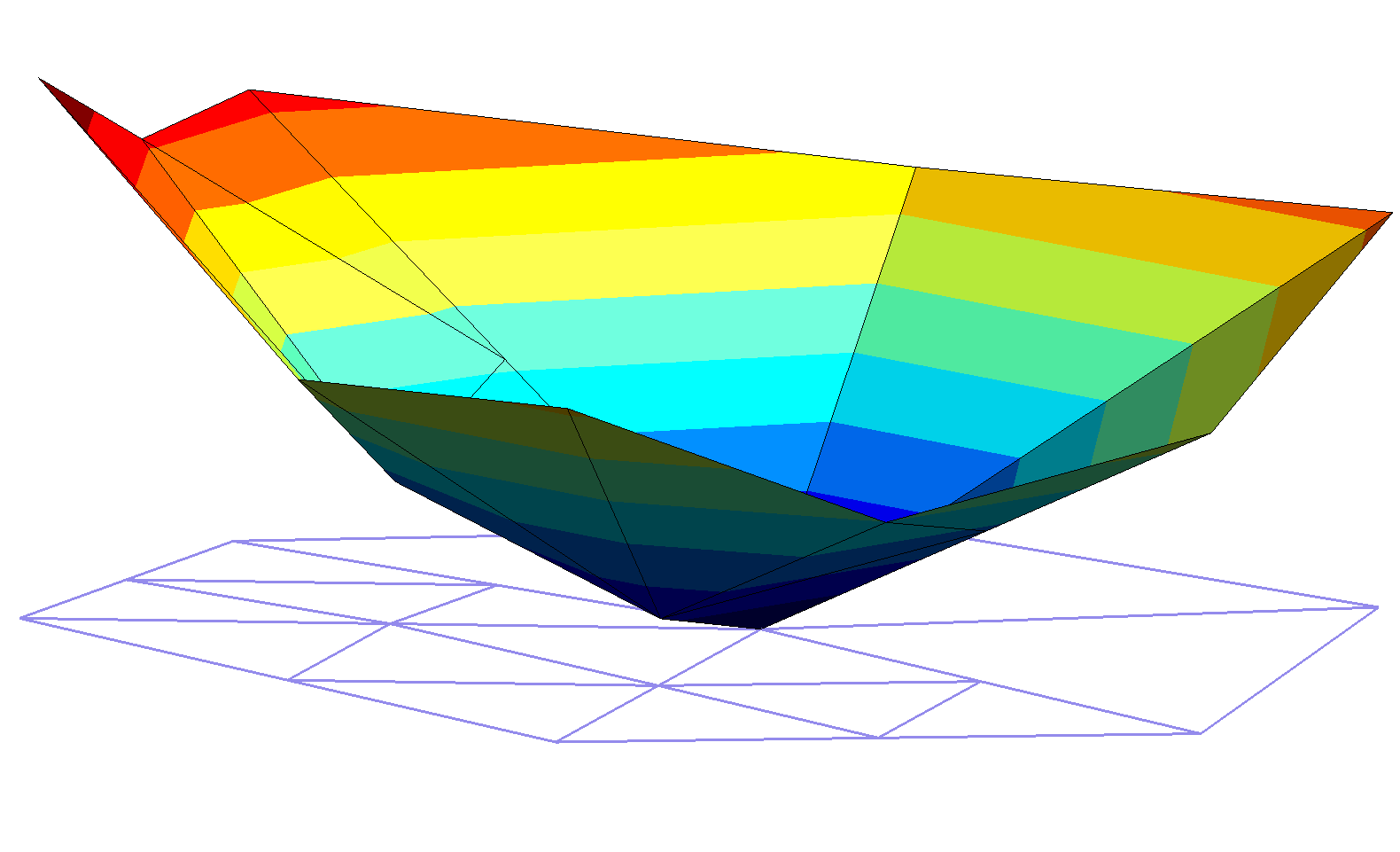}
		\label{enriche_construct_b}}
	\subfloat[$N^p(\vm{x})~\vm{F}^p(\vm{x})$]{
	    \includegraphics[width=0.26\textwidth]{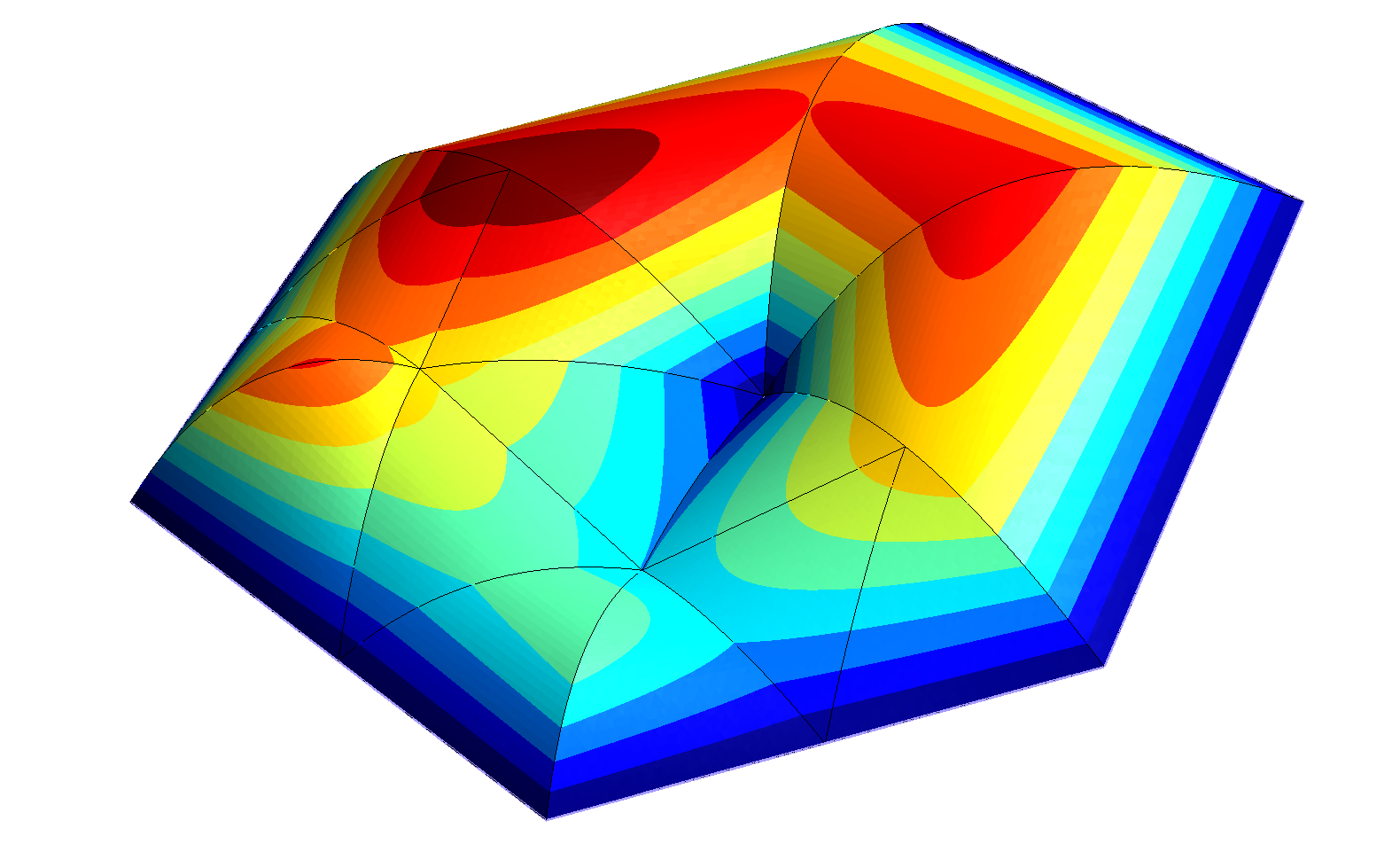}
        \label{enriche_construct_c}}
    \subfloat[$\displaystyle\sum_{m\in {\bar{I}_p^l}} n_m(\vm{x})N^p(\vm{x}_m)~\vm{F}^p(\vm{x}_m)$]{
		\includegraphics[width=0.26\textwidth]{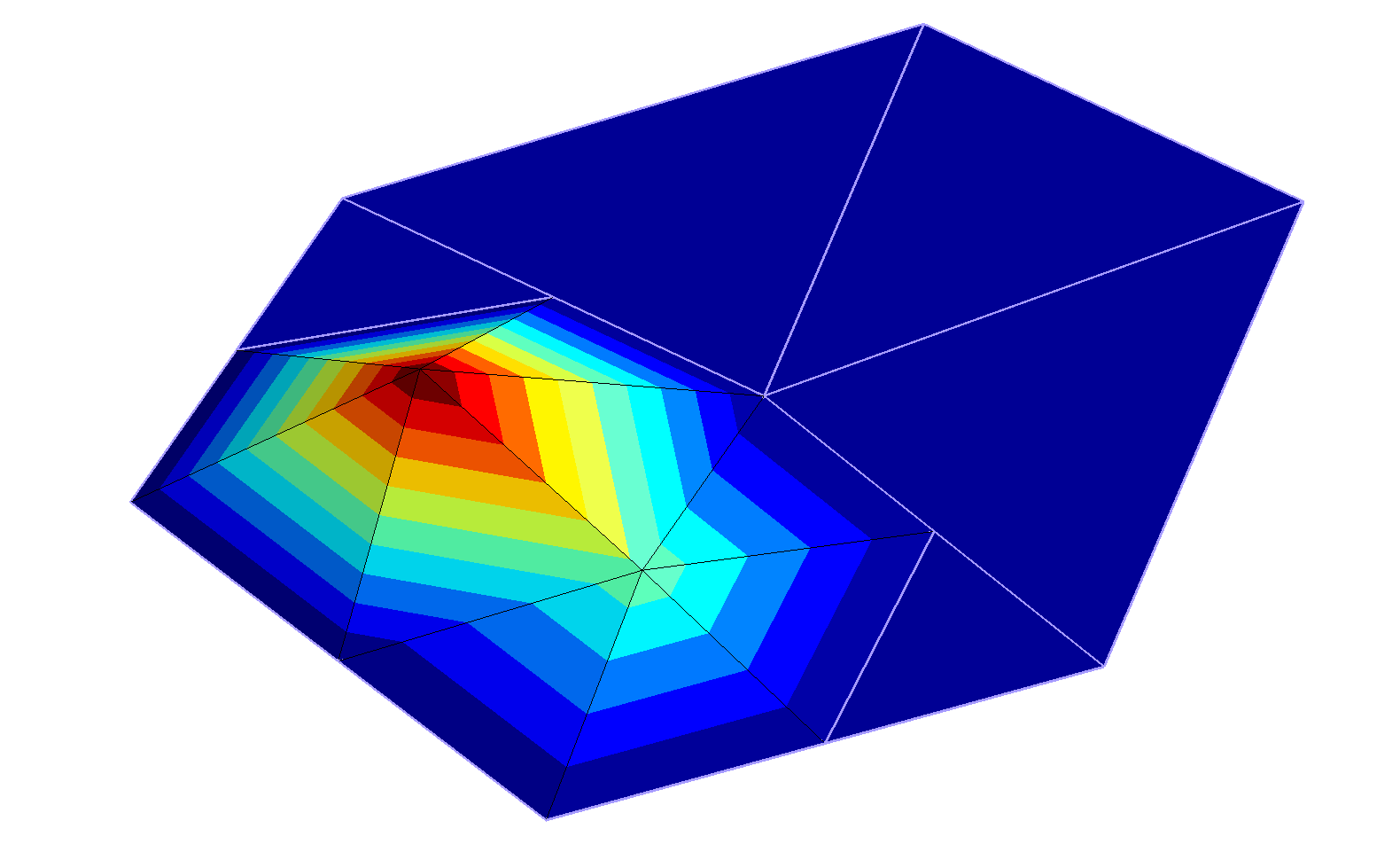}
		\label{enriche_construct_d}}
	\caption{Construction of an interpolated generalized approximation function  for a  $p$ patch (based on the example in figure  \ref{ts_fig} ).
   For visualization purposes, the quantities transformed into scalars are   drawn here as artificial elevations outside the plane formed by the patch mesh. (\protect\subref*{enriche_construct_a},\protect\subref*{enriche_construct_b}) and (\protect\subref*{enriche_construct_c},\protect\subref*{enriche_construct_d}) are represented from two different viewpoints in space for better visualization.
   \label{enriche_construct}}
\end{figure}

In this work we choose to approach blending element problem  by using a shift to cancel the enrichment function at the coarse-scale enriched node location and by forcing discretization transition to be embedded in patches, contrasting with SGFEM \cite{Babuska2012} presented in \eqref{sgfem}.
This last point is somewhat similar to the idea in \cite{Fries2008} where enrichment was added to nodes initially not enriched  in blending elements, with the use of a special ramp function to correct the original enrichment function.
Here, mixed discretization patches are added, so enrichments are added (red node in figure \ref{ts_sp})  on initially unenriched  blending elements nodes (unenriched  because their support did not cover the behavior being captured).
And the treatment of hanging nodes replaces the effect of the ramp function.
In this way, the original blending element (orange element connected to red node in figure \ref{ts_sp})  regain their partition of the unity property and the new blending elements (yellow element in figure \ref{ts_sp}) are not perturbed  thanks to the use of a null enrichment function.
Some tests not provided in this work show that the proposed enrichment give a good convergence compare to that given by the SGFEM enrichment.  
In this work, the  \tS enrichment function $\vm{F}^p(\vm{x})$ is :
\begin{equation}
	\forall p\in I_{e}^g,~\forall \vm{x}\in \omega_p ~\vm{F}^p(\vm{x})=\vm{u}^{Q^p}(\vm{x})-\vm{u}^{Q^p}(\vm{x}^p)~ \text{and}~ \forall \vm{x}\notin \omega_p~\vm{F}^p(\vm{x})=\vm{0}
	\label{shift_enrich}
\end{equation}
where:
\begin{itemize}
	\item $\vm{u}^{Q^p}(\vm{x})$ and $I_{e}^g$ are defined in section \ref{TSMethodOrig}
	\item $\vm{x}^p$ is the location of the enriched node $p$ 
\end{itemize}
Another specificity of this work is the wish that $\vm{u}^{ts}$ converges to  $\vm{u}^{R}$ which can only be achieved if the enrichment function is linearly interpolated on the reference mesh following  the construction and use of $\vm{T}_{MG}$.
Thus, the  kinematic equation actually used at the global scale, corresponding to the application of the $\vm{T}_{MG}$ operator, is: 
\begin{equation}
		\vm{U}(\vm{x})=\displaystyle\sum_{i\in I^g}  N^i(\vm{x})~ \vm{U}^{i} + \displaystyle\sum_{p\in I_{e}^g} \vm{E}^{p} \displaystyle\sum_{m\in {\bar{I}_p^l}} n^m(\vm{x})~  N^p(\vm{x}_m)~\vm{F}^p(\vm{x}_m)
		\label{kinematicnew}
\end{equation}
where $\bar{I}_p^l={I}_p^l\setminus \left\lbrace \text{hanging nodes} \right\rbrace$ since the L-set has also been eliminated in the patch problems.
The  construction of the second term of \eqref{kinematicnew} is first illustrated, in 1D, in the figure \ref{enriche1D_construct} where  artificial elevations show the different quantities of the 8 micro-edges of a  $p$ patch (composed of 2 macro-edges not shown here).
The same type of illustration of the 2D example in figure \ref{ts_fig} is shown  in figure \ref{enriche_construct}.
Again, artificial elevations  show the various quantities outside the  plane formed by a patch $p$.
The solutions $\vm{u}^{Q^p}(\vm{x})$  (figure \ref{enriche1D_construct_a} and \ref{enriche_construct_a}) are completely arbitrary. 
The enriched function $\vm{F}^p(\vm{x})$ corresponding to the equation \eqref{shift_enrich} is given in  figure \ref{enriche1D_construct_b} and \ref{enriche_construct_b}.
By construction, this function is null at the location of the enriched node.
When multiplied by the standard finite element approximation function to form the generalized finite element approximation function (figure \ref{enriche1D_construct_c} and \ref{enriche_construct_c}), all the boundary values of the  patches become zero and we obtain second order curves/surfaces.
The final linear interpolation is shown in the figure \ref{enriche1D_construct_d} and \ref{enriche_construct_d}. 
In the 2D example, where a mixture of refined and unrefined elements forms the patch, the interpolated generalized function is necessarily zero on the unrefined element (yellow blending  elements in the figure \ref{ts_sp}).
This property implies that blending elements have unenriched nodes  and enriched nodes  with zero contribution from enrichment function.
Thus, the enrichment contribution is zero for the blending element and the partition of unity is saved because only the classical terms matter.
For refined elements, $N^p(\vm{x}_m)~\vm{F}^p(\vm{x}_m)$ is not zero for any location inside the patch (if $\vm{F}^p(\vm{x}_m)\neq\vm{0}$).
For the  hanging nodes, since they are eliminated in $\bar{I}_p^l$ their values are therefore zero which is consistent with the values in the unrefined elements.

Note that if the refined elements correspond only to a terminal connection specific to the adaptation tool (element added in between  figure \ref{split_2} and figure \ref{split_2f} in \ref{rsplit_anexe}), then, despite  the refinement, the associated generalized finite element approximation function is zero on these elements.
To avoid this situation, which can lead to a singular problem, the refined SP elements are always divided at least once.

\section{Theoretical sequential algebra performance}\label{theoSeqAlgPerf}
Since the \tS method uses a resolution with sparse matrix of different patterns and densities, and since the number of iterations of the \tS loop is not known in advance, estimating the algebraic performance of the \tS method is difficult in a general context.
It is even more complex when using parallelism.
Thus, to obtain the  numerical complexity of the \tS solver, we will use a simple sequential numerical example to be able to count the number of flops \footnote{flop=floating-point operation: one addition, subtraction, division or multiplication of two floating-point numbers}  in each part of the algorithm.
The  number of flops used will then be expressed as  a polynomial function of the problem dimensions keeping only the leading terms to obtain the order of magnitude ($\mathcal{O}(n)$ notation).
If we take a meshed cube with a finite hexahedral element of order one in 3D, we can, using octree refinement, simply count the coarse scale dofs, the patches and the patches' dofs.

For the linear resolution steps of the \tS method, only direct sparse solver will be considered.
Thus, the  number of flops for the factorization, in terms of the $n$ dimension of the problem, is approximately between $n^2$  and $n^3$ (dense matrix).
And for backward and forward substitution the cost is in between $2.n$ and $2.n^2$.
Depending on the strategy of the solver (multi-frontal,left-looking,...), the reordering of the symbolic factorization (nested dissection, minimum degree, ...) and the density of the initial matrix, the number of flops varies greatly and depends on a variable number of parameters.
Thus, to estimate these quantities, a sparse ratio parameter $SR$ is introduced such that the dimension used to count the number of flops $n_f$ is:
\begin{equation}
	n_f(n,SR)=\sqrt{SR. n^2}
\end{equation}
with $SR\in ]0,1]$.
Then  the number of flops for factorization ($count_{fact}$), backward/forward computation  ($count_{bf}$) and resolution ($count_{resolv}$) are, using a dense estimation:
\begin{subequations}
	\label{TAP:costr}
	\begin{gather}\tag{\ref{TAP:costr}}
		count_{resolv}(n,SR)=count_{fact}(n,SR)+count_{bf}(n,SR) \\
		\begin{align}
			count_{fact}(n,SR)=(n_f(n,SR))^3 \\
			count_{bf}(n,SR)=2.(n_f(n,SR))^2 
		\end{align}
	\end{gather}
\end{subequations}
This $SR$ parameter  is completely arbitrary.
It is introduced to verify that the conclusions obtained in this section are not impacted by the lack of accurate estimation of the factorization and solving phase.

For discretization, the root of the octree  uses one element to represent the problem.
It corresponds to  zero level.
A level, noted $L$ in the following, describes the depth in the octree tree.
With $L=0$ being the starting level, to go from $L$ to $L+1$, the octree refinement adds a new node in the middle of all edges, faces and elements of level $L$.
Then, each element of  the level $L$ is replaced by  8 hexahedra encapsulated in it and connected to the old and new nodes.
The number of classical dofs for a level $L$ is easily deduced from this refinement strategy and is :
\begin{equation}nb_{dof}(L) = 3.(2^{L}+1)^{3}
	\label{TAP:nbdofL}
\end{equation}
Thus, the \tS coarse enriched problem defined at the $L_c$ level  will cost for $N_l$ \tS iteration with a sparse ratio parameter $SR_c$ using \eqref{TAP:costr} and \eqref{TAP:nbdofL}:
\begin{equation}
	\label{TAP:costc}
	cost_{coarse}(L_{c},N_l,SR_c)=N_l\times count_{resolv}(2\times nb_{dof}(L_{c}),SR_c)
\end{equation}

\begin{figure}[!htb]
	\centering
	\subfloat[][$N_l=30$,$SR_c=SR_p=0.017$]{
		\includegraphics[width=70mm,keepaspectratio]{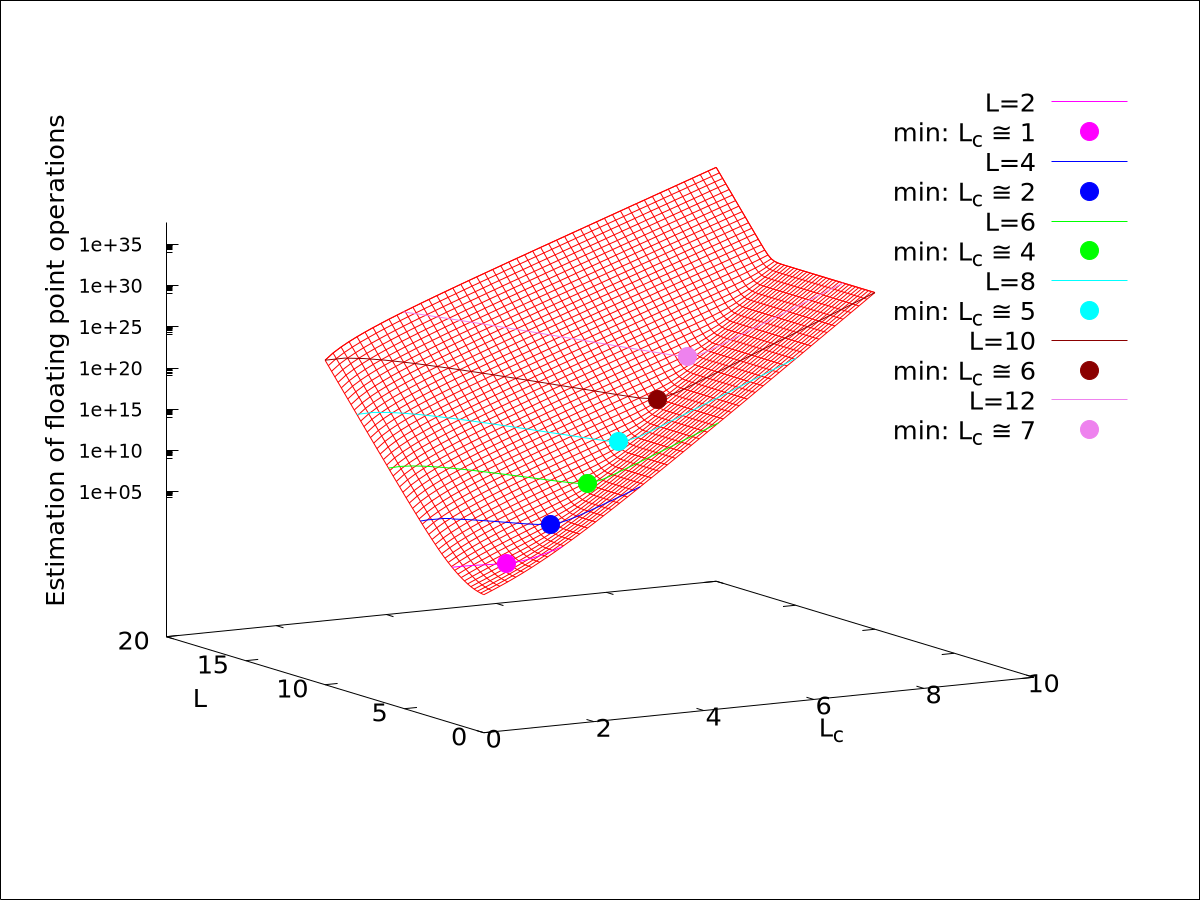}
		\label{TAP:f_1}
	}	
	\subfloat[][$L=8,N_L=10,SR_F=SR_C=0.017$]{
		\includegraphics[width=72mm,keepaspectratio]{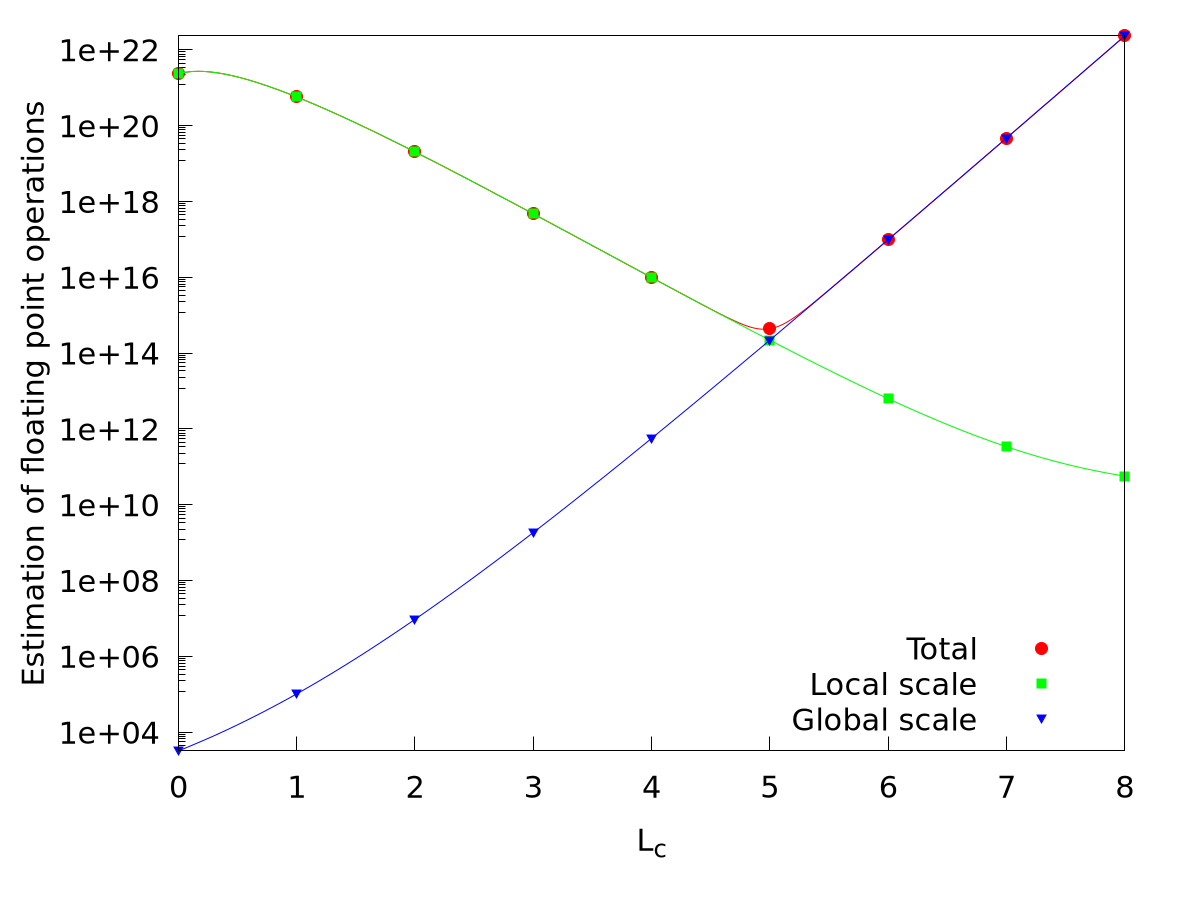}
		\label{TAP:f_3}
	}
\\
	\subfloat[][$L=8$]{
		\includegraphics[width=72mm,keepaspectratio]{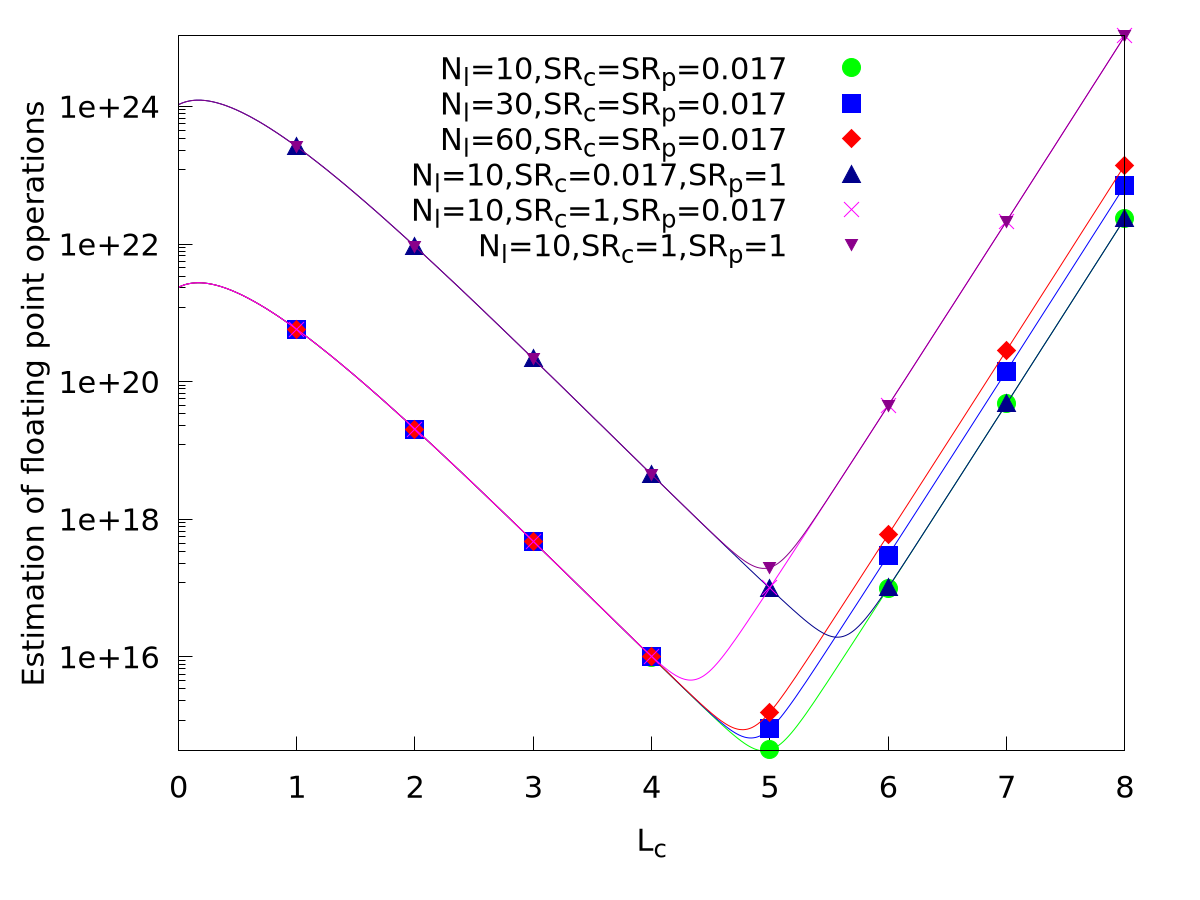}
		\label{TAP:f_2}}
	\subfloat[][ratio of floating point estimation (full rank/\tS resolution) with $N_l=30$,$SR_c=SR_p=0.017$]{
		\includegraphics[width=70mm,keepaspectratio]{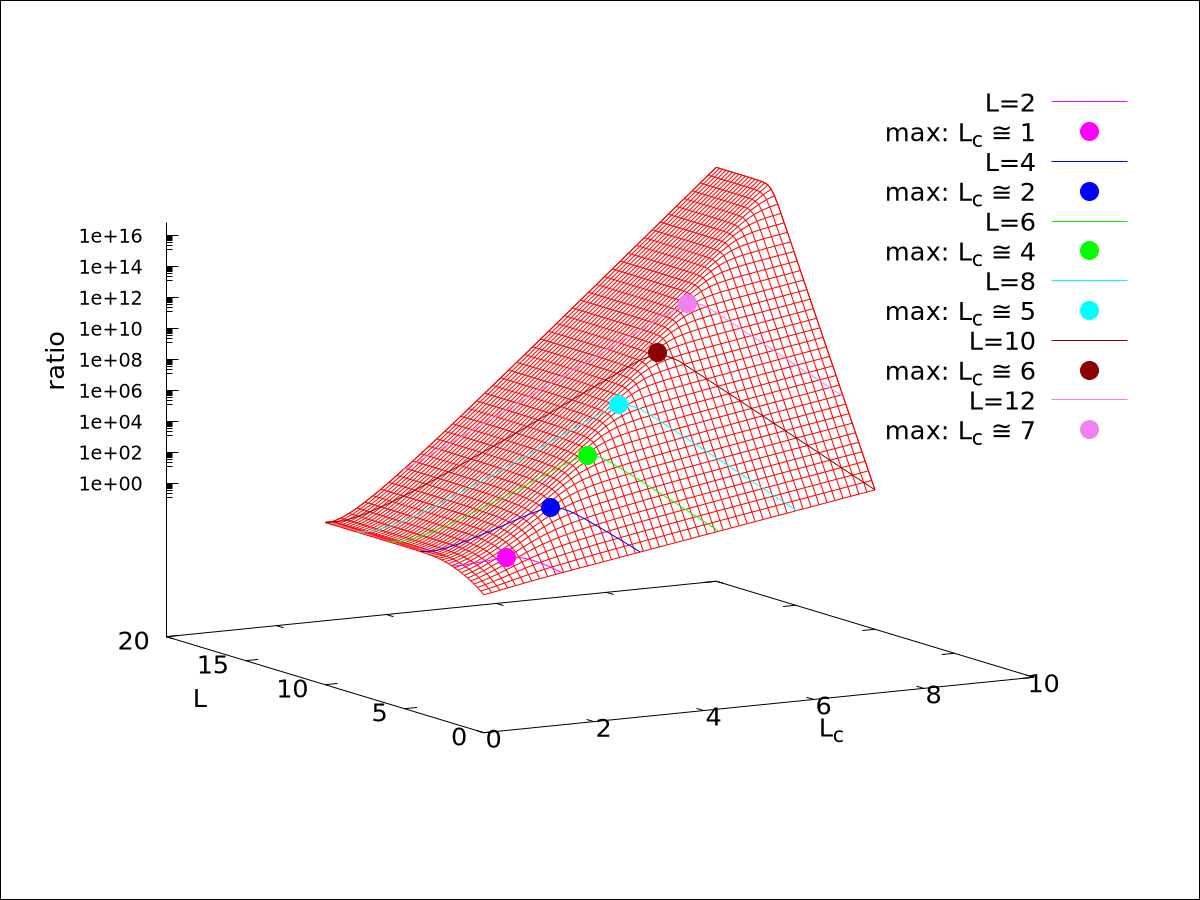}
		\label{TAP:f_4}		
	}
	\caption{floating point estimation for different value of $L_c$,L, $SR_C$ (sparse ratio for coarse problem), $SR_F$(sparse ratio for fine problems) and $N_l$ (number of  \tS iteration)}
	\label{theo_alg}
\end{figure}

\noindent For a given  coarse level $L_c$, it is possible to obtain in the same way the cost of computing patches (factorization and solving) for $N_l$ \tS iterations, a sparse ratio parameter $SR_p$ and a final refinement level $L$.
This is given by the function $cost_{patch}(L_c,L,N_l,SR_p)$ detailed in \ref{TAP_anexe}.
Using \eqref{TAP:costc} and \eqref{TAP_anexe_cp} we can state that the cost of the \tS method is:
\begin{equation}
	\label{TAP:cost_ts}
	cost_{\tS}(L_c,L,N_l,SR_c,SR_p)=cost_{patch}(L_c,L,N_l,SR_p)+\\cost_{coarse}(L_{c},N_l,SR_c)
\end{equation}
Fixing  $N_l,SR_c,SR_p$  the equation \eqref{TAP:cost_ts} for some values of $L_c$ and $L$ give the surface of figure \ref{TAP:f_1}.
We can see on this continuous representation of \eqref{TAP:cost_ts} that the iso-$L$ curves have a minimum for a specific $L_c$ value.
It means that \tS method has an optimal coarse level for a fixed target level that provides the best performance (i.e. the lowest number of flop consumption).
This can be understood by looking at figure \ref{TAP:f_3} for a specific target $L=8$  and fixed $N_l,SR_c,SR_p$ where the equation \eqref{TAP:cost_ts} is divided into patch contribution (equation \eqref{TAP_anexe_cp}) and coarse-level contribution (equation \eqref{TAP:costr}).
In this figure, it is clear that the antagonistic  evolution of the performance of patch solving and coarse-scale problem solving  lead to having a minimal value at the intersection of the curves.
In figure \ref{TAP:f_2} the variation of the parameters $N_l,SR_c,SR_p$ does not change the fact that a minimum exists. 
In this example where $L=8$, when the number of loop increases, the minimum remains the same.
When we consider the variation of the sparsity of the matrix (up to dense case), the minimum is moved around its original position depending on whether the coarse problem or the patch problems are denser or not.
With denser matrices for patches, it is naturally preferable to increase $l_c$ to have smaller patches.
And with a denser coarse matrix, it is better to decrease $l_c$ to have a smaller coarse problem to solve.
When both sparsities increase, the minimum remains the same but with a higher value.

Figure \ref{TAP:f_4} presents, as a ratio, an estimate  of the performance of a full rank solver (using the equation \eqref{TAP:costr}) compared to the \tS solver.
A ratio greater than, equal to, or less than 1 indicates that the full rank solver uses more, the same, or less flop than the \tS solver respectively.
The surface  naturally has a maximum for iso-$L$ curves with the same $L_c$ as the minimum of the figure \ref{TAP:f_1} (on iso-$L$ curves the cost of the full rank solver is constant).
In this example, the  theoretical performance of the \tS solver compared to the full rank solver is much better for $L>2$.
For $iso-L=2$ the optimal ratio at $L_c=1$ is only of 1.76.
And for $L=1$, as in this case the patch problems size are almost equivalent to the $L_c$ level of discretization  no gain is obtained.
This illustrates the fact that the \tS method, to be effective, must be used at a target level that allows to start from a coarse level that induces a scaling effect.

The conclusion drawn from the analysis of this octree refinement example were also observed in parallel with the different test cases studied bellow.
Parallelism, the nature of the unstructured mesh, and the cost of assembly also impact the choice of the optimal coarse mesh for a given fine-scale target discretization.
But the observation of the test case \ref{UOL} led to use the following formula to define the scale jump for most cases with regular refinement:
\begin{equation}
	L=L_c+3
	\label{ts_jump_level_eq}
\end{equation}

\begin{figure}[!htb]
	\subfloat[][Weight: micro embemded element count	]{
		\def\svgwidth{0.32\textwidth}
\begingroup%
  \makeatletter%
  \providecommand\color[2][]{%
    \errmessage{(Inkscape) Color is used for the text in Inkscape, but the package 'color.sty' is not loaded}%
    \renewcommand\color[2][]{}%
  }%
  \providecommand\transparent[1]{%
    \errmessage{(Inkscape) Transparency is used (non-zero) for the text in Inkscape, but the package 'transparent.sty' is not loaded}%
    \renewcommand\transparent[1]{}%
  }%
  \providecommand\rotatebox[2]{#2}%
  \ifx\svgwidth\undefined%
    \setlength{\unitlength}{342.09281158bp}%
    \ifx\svgscale\undefined%
      \relax%
    \else%
      \setlength{\unitlength}{\unitlength * \real{\svgscale}}%
    \fi%
  \else%
    \setlength{\unitlength}{\svgwidth}%
  \fi%
  \global\let\svgwidth\undefined%
  \global\let\svgscale\undefined%
  \makeatother%
  \begin{picture}(1,1.20965849)%
    \lineheight{1}%
    \setlength\tabcolsep{0pt}%
    \put(0,0){\includegraphics[width=\unitlength,page=1]{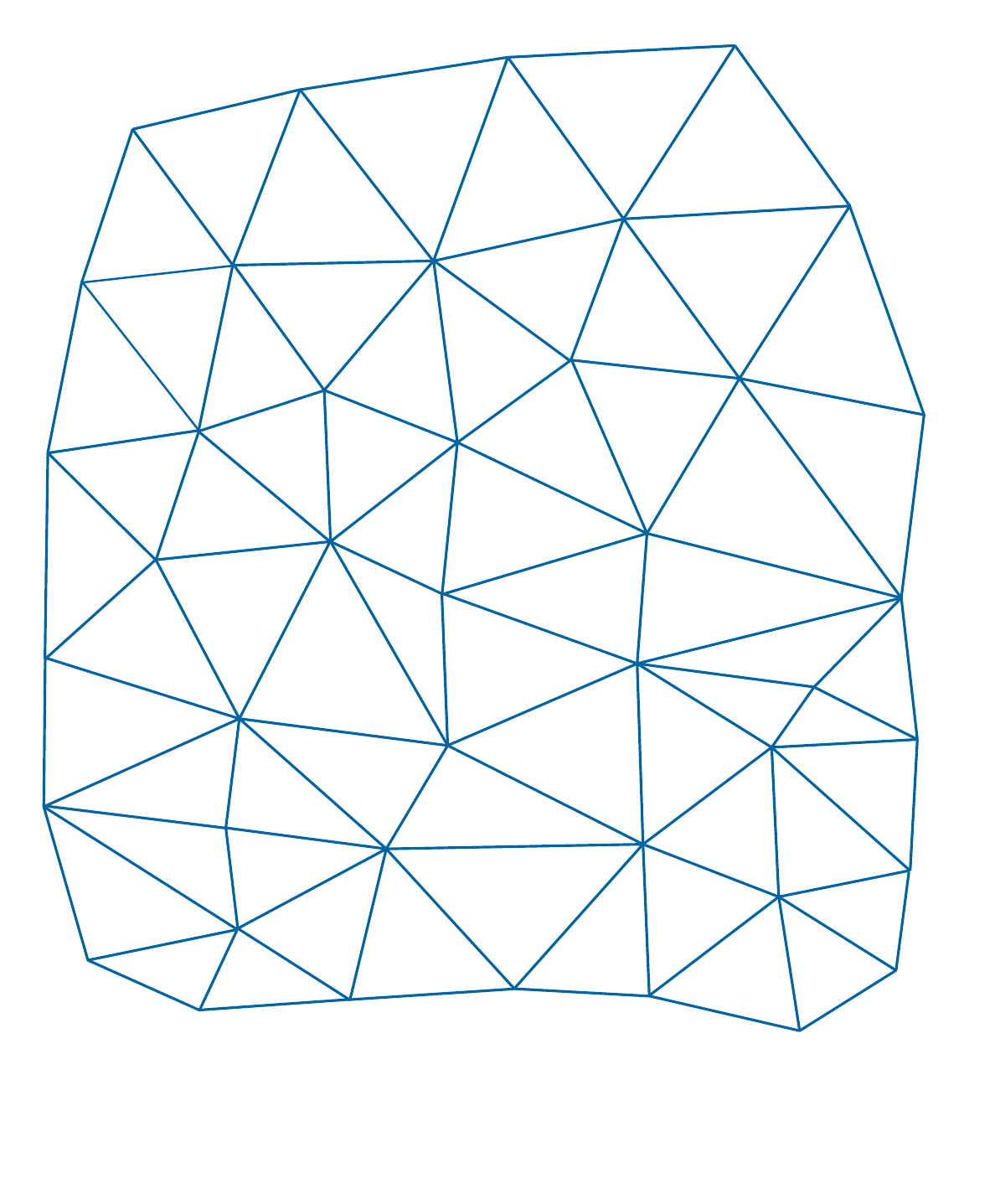}}%
    \put(0.48620807,1.01966624){\color[rgb]{0,0,0}\makebox(0,0)[lt]{\lineheight{0}\smash{\begin{tabular}[t]{l}\tiny{1}\end{tabular}}}}%
    \put(0.38960371,1.07364175){\color[rgb]{0,0,0}\makebox(0,0)[lt]{\lineheight{0}\smash{\begin{tabular}[t]{l}\tiny{1}\end{tabular}}}}%
    \put(0.28172572,1.0039157){\color[rgb]{0,0,0}\makebox(0,0)[lt]{\lineheight{0}\smash{\begin{tabular}[t]{l}\tiny{1}\end{tabular}}}}%
    \put(0.28231252,0.89935038){\color[rgb]{0,0,0}\makebox(0,0)[lt]{\lineheight{0}\smash{\begin{tabular}[t]{l}\tiny{1}\end{tabular}}}}%
    \put(0.24297861,0.83701503){\color[rgb]{0,0,0}\makebox(0,0)[lt]{\lineheight{0}\smash{\begin{tabular}[t]{l}\tiny{1}\end{tabular}}}}%
    \put(0.2496675,0.75714058){\color[rgb]{0,0,0}\makebox(0,0)[lt]{\lineheight{0}\smash{\begin{tabular}[t]{l}\tiny{4}\end{tabular}}}}%
    \put(0.58748812,0.25372639){\color[rgb]{0,0,0}\makebox(0,0)[lt]{\lineheight{0}\smash{\begin{tabular}[t]{l}\tiny{1}\end{tabular}}}}%
    \put(0.159648,0.87462083){\color[rgb]{0,0,0}\makebox(0,0)[lt]{\lineheight{0}\smash{\begin{tabular}[t]{l}\tiny{1}\end{tabular}}}}%
    \put(0.112509,0.71594782){\color[rgb]{0,0,0}\makebox(0,0)[lt]{\lineheight{0}\smash{\begin{tabular}[t]{l}\tiny{1}\end{tabular}}}}%
    \put(0.09064281,0.81277434){\color[rgb]{0,0,0}\makebox(0,0)[lt]{\lineheight{0}\smash{\begin{tabular}[t]{l}\tiny{1}\end{tabular}}}}%
    \put(0.07656086,0.64690988){\color[rgb]{0,0,0}\makebox(0,0)[lt]{\lineheight{0}\smash{\begin{tabular}[t]{l}\tiny{1}\end{tabular}}}}%
    \put(0.09512734,0.29753546){\color[rgb]{0,0,0}\makebox(0,0)[lt]{\lineheight{0}\smash{\begin{tabular}[t]{l}\tiny{1}\end{tabular}}}}%
    \put(0.16094553,0.33610485){\color[rgb]{0,0,0}\makebox(0,0)[lt]{\lineheight{0}\smash{\begin{tabular}[t]{l}\tiny{1}\end{tabular}}}}%
    \put(0.17300051,0.23268861){\color[rgb]{0,0,0}\makebox(0,0)[lt]{\lineheight{0}\smash{\begin{tabular}[t]{l}\tiny{1}\end{tabular}}}}%
    \put(0.24387448,0.22210127){\color[rgb]{0,0,0}\makebox(0,0)[lt]{\lineheight{0}\smash{\begin{tabular}[t]{l}\tiny{1}\end{tabular}}}}%
    \put(0.67312637,0.28561178){\color[rgb]{0,0,0}\makebox(0,0)[lt]{\lineheight{0}\smash{\begin{tabular}[t]{l}\tiny{1}\end{tabular}}}}%
    \put(0.70607184,0.36781382){\color[rgb]{0,0,0}\makebox(0,0)[lt]{\lineheight{0}\smash{\begin{tabular}[t]{l}\tiny{1}\end{tabular}}}}%
    \put(0.67158406,0.44832568){\color[rgb]{0,0,0}\makebox(0,0)[lt]{\lineheight{0}\smash{\begin{tabular}[t]{l}\tiny{1}\end{tabular}}}}%
    \put(0.75722783,0.49095477){\color[rgb]{0,0,0}\makebox(0,0)[lt]{\lineheight{0}\smash{\begin{tabular}[t]{l}\tiny{1}\end{tabular}}}}%
    \put(0.77621411,0.54152484){\color[rgb]{0,0,0}\makebox(0,0)[lt]{\lineheight{0}\smash{\begin{tabular}[t]{l}\tiny{1}\end{tabular}}}}%
    \put(0.86846602,0.51749065){\color[rgb]{0,0,0}\makebox(0,0)[lt]{\lineheight{0}\smash{\begin{tabular}[t]{l}\tiny{1}\end{tabular}}}}%
    \put(0.70233329,1.03744956){\color[rgb]{0,0,0}\makebox(0,0)[lt]{\lineheight{0}\smash{\begin{tabular}[t]{l}\tiny{1}\end{tabular}}}}%
    \put(0.57770659,1.10487054){\color[rgb]{0,0,0}\makebox(0,0)[lt]{\lineheight{0}\smash{\begin{tabular}[t]{l}\tiny{1}\end{tabular}}}}%
    \put(0.17828825,0.68604314){\color[rgb]{0,0,0}\makebox(0,0)[lt]{\lineheight{0}\smash{\begin{tabular}[t]{l}\tiny{4}\end{tabular}}}}%
    \put(0.12142873,0.55805318){\color[rgb]{0,0,0}\makebox(0,0)[lt]{\lineheight{0}\smash{\begin{tabular}[t]{l}\tiny{4}\end{tabular}}}}%
    \put(0.08933707,0.47867248){\color[rgb]{0,0,0}\makebox(0,0)[lt]{\lineheight{0}\smash{\begin{tabular}[t]{l}\tiny{4}\end{tabular}}}}%
    \put(0.71944226,0.92969867){\color[rgb]{0,0,0}\makebox(0,0)[lt]{\lineheight{0}\smash{\begin{tabular}[t]{l}\tiny{4}\end{tabular}}}}%
    \put(0.80988203,0.85701882){\color[rgb]{0,0,0}\makebox(0,0)[lt]{\lineheight{0}\smash{\begin{tabular}[t]{l}\tiny{4}\end{tabular}}}}%
    \put(0.82412991,0.75197841){\color[rgb]{0,0,0}\makebox(0,0)[lt]{\lineheight{0}\smash{\begin{tabular}[t]{l}\tiny{4}\end{tabular}}}}%
    \put(0.52165609,0.45319518){\color[rgb]{0,0,0}\makebox(0,0)[lt]{\lineheight{0}\smash{\begin{tabular}[t]{l}\tiny{4}\end{tabular}}}}%
    \put(0.46363682,0.30433885){\color[rgb]{0,0,0}\makebox(0,0)[lt]{\lineheight{0}\smash{\begin{tabular}[t]{l}\tiny{4}\end{tabular}}}}%
    \put(0.37510141,0.24550601){\color[rgb]{0,0,0}\makebox(0,0)[lt]{\lineheight{0}\smash{\begin{tabular}[t]{l}\tiny{4}\end{tabular}}}}%
    \put(0.29773097,0.26397338){\color[rgb]{0,0,0}\makebox(0,0)[lt]{\lineheight{0}\smash{\begin{tabular}[t]{l}\tiny{4}\end{tabular}}}}%
    \put(0.52493321,0.91469284){\color[rgb]{0,0,0}\makebox(0,0)[lt]{\lineheight{0}\smash{\begin{tabular}[t]{l}\tiny{4}\end{tabular}}}}%
    \put(0.37910998,0.82399794){\color[rgb]{0,0,0}\makebox(0,0)[lt]{\lineheight{0}\smash{\begin{tabular}[t]{l}\tiny{5}\end{tabular}}}}%
    \put(0.35506051,0.74733493){\color[rgb]{0,0,0}\makebox(0,0)[lt]{\lineheight{0}\smash{\begin{tabular}[t]{l}\tiny{6}\end{tabular}}}}%
    \put(0.19639562,0.60330428){\color[rgb]{0,0,0}\makebox(0,0)[lt]{\lineheight{0}\smash{\begin{tabular}[t]{l}\tiny{5}\end{tabular}}}}%
    \put(0.16788287,0.41153961){\color[rgb]{0,0,0}\makebox(0,0)[lt]{\lineheight{0}\smash{\begin{tabular}[t]{l}\tiny{5}\end{tabular}}}}%
    \put(0.26287349,0.32730654){\color[rgb]{0,0,0}\makebox(0,0)[lt]{\lineheight{0}\smash{\begin{tabular}[t]{l}\tiny{5}\end{tabular}}}}%
    \put(0.63255992,0.87263066){\color[rgb]{0,0,0}\makebox(0,0)[lt]{\lineheight{0}\smash{\begin{tabular}[t]{l}\tiny{5}\end{tabular}}}}%
    \put(0.4193694,0.38877781){\color[rgb]{0,0,0}\makebox(0,0)[lt]{\lineheight{0}\smash{\begin{tabular}[t]{l}\tiny{6}\end{tabular}}}}%
    \put(0.47688278,0.67339481){\color[rgb]{0,0,0}\makebox(0,0)[lt]{\lineheight{0}\smash{\begin{tabular}[t]{l}\tiny{35}\end{tabular}}}}%
    \put(0.38262946,0.67252342){\color[rgb]{0,0,0}\makebox(0,0)[lt]{\lineheight{0}\smash{\begin{tabular}[t]{l}\tiny{29}\end{tabular}}}}%
    \put(0.39384811,0.58204683){\color[rgb]{0,0,0}\makebox(0,0)[lt]{\lineheight{0}\smash{\begin{tabular}[t]{l}\tiny{46}\end{tabular}}}}%
    \put(0.31249198,0.533099){\color[rgb]{0,0,0}\makebox(0,0)[lt]{\lineheight{0}\smash{\begin{tabular}[t]{l}\tiny{31}\end{tabular}}}}%
    \put(0.32364596,0.43214955){\color[rgb]{0,0,0}\makebox(0,0)[lt]{\lineheight{0}\smash{\begin{tabular}[t]{l}\tiny{40}\end{tabular}}}}%
    \put(0.26685324,0.39580687){\color[rgb]{0,0,0}\makebox(0,0)[lt]{\lineheight{0}\smash{\begin{tabular}[t]{l}\tiny{19}\end{tabular}}}}%
    \put(0.49803322,0.53179206){\color[rgb]{0,0,0}\makebox(0,0)[lt]{\lineheight{0}\smash{\begin{tabular}[t]{l}\tiny{10}\end{tabular}}}}%
    \put(0.55323489,0.60423419){\color[rgb]{0,0,0}\makebox(0,0)[lt]{\lineheight{0}\smash{\begin{tabular}[t]{l}\tiny{14}\end{tabular}}}}%
    \put(0.71216092,0.60087089){\color[rgb]{0,0,0}\makebox(0,0)[lt]{\lineheight{0}\smash{\begin{tabular}[t]{l}\tiny{8}\end{tabular}}}}%
    \put(0.73617535,0.69727998){\color[rgb]{0,0,0}\makebox(0,0)[lt]{\lineheight{0}\smash{\begin{tabular}[t]{l}\tiny{9}\end{tabular}}}}%
    \put(0.52166126,0.76056055){\color[rgb]{0,0,0}\makebox(0,0)[lt]{\lineheight{0}\smash{\begin{tabular}[t]{l}\tiny{47}\end{tabular}}}}%
    \put(0.65062642,0.78155619){\color[rgb]{0,0,0}\makebox(0,0)[lt]{\lineheight{0}\smash{\begin{tabular}[t]{l}\tiny{24}\end{tabular}}}}%
    \put(0.47428471,0.84238185){\color[rgb]{0,0,0}\makebox(0,0)[lt]{\lineheight{0}\smash{\begin{tabular}[t]{l}\tiny{11}\end{tabular}}}}%
  \end{picture}%
\endgroup%
		\label{ts_dist_weigth_elem}
	}
	\subfloat[][Global element dispatching on 4 processes partially based on weight (\protect\subref*{ts_dist_weigth_elem})]{
	   \def\svgwidth{0.32\textwidth}
\begingroup%
  \makeatletter%
  \providecommand\color[2][]{%
    \errmessage{(Inkscape) Color is used for the text in Inkscape, but the package 'color.sty' is not loaded}%
    \renewcommand\color[2][]{}%
  }%
  \providecommand\transparent[1]{%
    \errmessage{(Inkscape) Transparency is used (non-zero) for the text in Inkscape, but the package 'transparent.sty' is not loaded}%
    \renewcommand\transparent[1]{}%
  }%
  \providecommand\rotatebox[2]{#2}%
  \ifx\svgwidth\undefined%
    \setlength{\unitlength}{342.09281158bp}%
    \ifx\svgscale\undefined%
      \relax%
    \else%
      \setlength{\unitlength}{\unitlength * \real{\svgscale}}%
    \fi%
  \else%
    \setlength{\unitlength}{\svgwidth}%
  \fi%
  \global\let\svgwidth\undefined%
  \global\let\svgscale\undefined%
  \makeatother%
  \begin{picture}(1,1.20965849)%
    \lineheight{1}%
    \setlength\tabcolsep{0pt}%
    \put(0,0){\includegraphics[width=\unitlength,page=1]{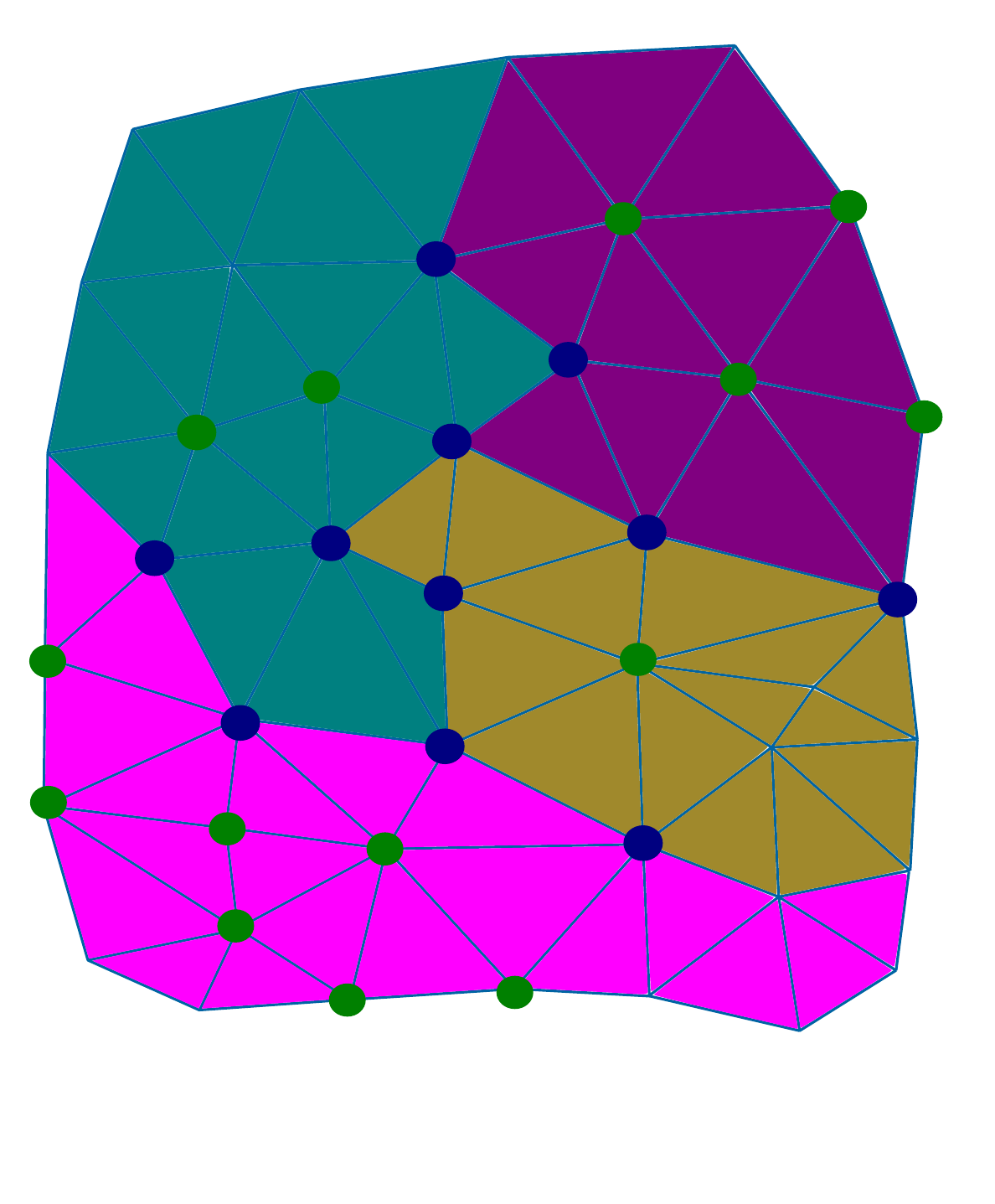}}%
    \put(0.89885423,0.88901305){\color[rgb]{0.50196078,0,0.50196078}\makebox(0,0)[lt]{\lineheight{0}\smash{\begin{tabular}[t]{l}\tiny{process 1}\end{tabular}}}}%
    \put(0.92777561,0.46741474){\color[rgb]{0.50196078,0.50196078,0}\makebox(0,0)[lt]{\lineheight{0}\smash{\begin{tabular}[t]{l}\tiny{process 2 }\end{tabular}}}}%
    \put(0.40797012,0.15228167){\color[rgb]{1,0,1}\makebox(0,0)[lt]{\lineheight{0}\smash{\begin{tabular}[t]{l}\tiny{process 3}\end{tabular}}}}%
    \put(-0.04530285,0.99875923){\color[rgb]{0,0.50196078,0.50196078}\makebox(0,0)[lt]{\lineheight{0}\smash{\begin{tabular}[t]{l}\tiny{process 0 }\end{tabular}}}}%
  \end{picture}%
\endgroup%
		\label{ts_dist_load_balance}
	}
	\subfloat[][Global element dispatching on 8 processes partially based on weight (\protect\subref*{ts_dist_weigth_elem})]{
	\def\svgwidth{0.32\textwidth}
\begingroup%
  \makeatletter%
  \providecommand\color[2][]{%
    \errmessage{(Inkscape) Color is used for the text in Inkscape, but the package 'color.sty' is not loaded}%
    \renewcommand\color[2][]{}%
  }%
  \providecommand\transparent[1]{%
    \errmessage{(Inkscape) Transparency is used (non-zero) for the text in Inkscape, but the package 'transparent.sty' is not loaded}%
    \renewcommand\transparent[1]{}%
  }%
  \providecommand\rotatebox[2]{#2}%
  \ifx\svgwidth\undefined%
    \setlength{\unitlength}{342.09281158bp}%
    \ifx\svgscale\undefined%
      \relax%
    \else%
      \setlength{\unitlength}{\unitlength * \real{\svgscale}}%
    \fi%
  \else%
    \setlength{\unitlength}{\svgwidth}%
  \fi%
  \global\let\svgwidth\undefined%
  \global\let\svgscale\undefined%
  \makeatother%
  \begin{picture}(1,1.20965849)%
    \lineheight{1}%
    \setlength\tabcolsep{0pt}%
    \put(0,0){\includegraphics[width=\unitlength,page=1]{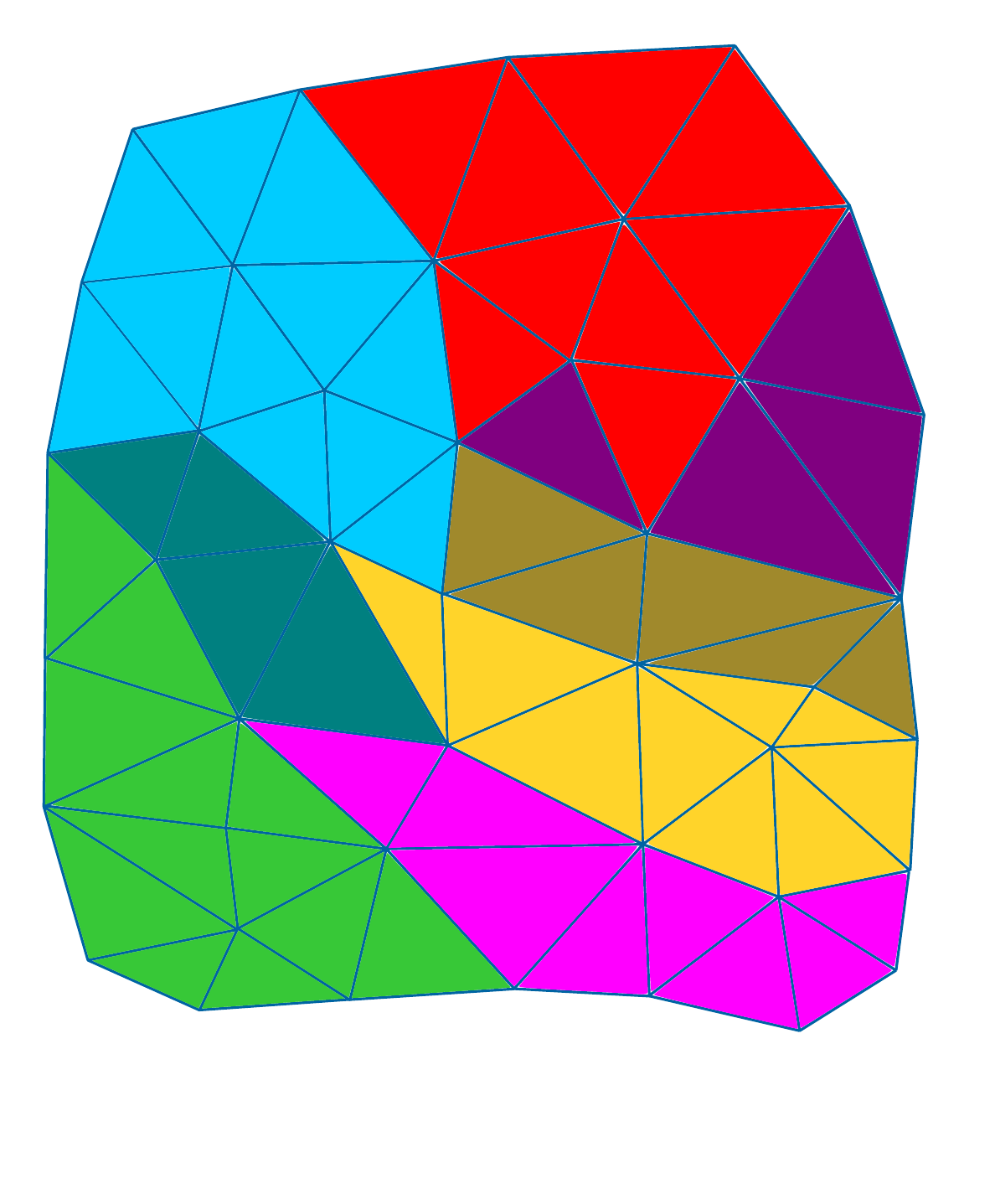}}%
    \put(0.59348497,0.8466883){\color[rgb]{0,0,0}\makebox(0,0)[lt]{\lineheight{0}\smash{\begin{tabular}[t]{l}\tiny{7}\end{tabular}}}}%
    \put(0.96309194,0.79464572){\color[rgb]{0,0,0}\makebox(0,0)[lt]{\lineheight{0}\smash{\begin{tabular}[t]{l}\tiny{4}\end{tabular}}}}%
    \put(0.65617276,0.98159756){\color[rgb]{0,0,0}\makebox(0,0)[lt]{\lineheight{0}\smash{\begin{tabular}[t]{l}\tiny{2}\end{tabular}}}}%
    \put(0.47072037,0.75909525){\color[rgb]{0,0,0}\makebox(0,0)[lt]{\lineheight{0}\smash{\begin{tabular}[t]{l}\tiny{10}\end{tabular}}}}%
    \put(0.25942799,0.27036035){\color[rgb]{0,0,0}\makebox(0,0)[lt]{\lineheight{0}\smash{\begin{tabular}[t]{l}\tiny{23}\end{tabular}}}}%
    \put(0.54620276,0.20435236){\color[rgb]{0,0,0}\makebox(0,0)[lt]{\lineheight{0}\smash{\begin{tabular}[t]{l}\tiny{25}\end{tabular}}}}%
    \put(0.77146822,0.81465496){\color[rgb]{0,0,0}\makebox(0,0)[lt]{\lineheight{0}\smash{\begin{tabular}[t]{l}\tiny{3}\end{tabular}}}}%
    \put(0.88940624,1.00073606){\color[rgb]{0,0,0}\makebox(0,0)[lt]{\lineheight{0}\smash{\begin{tabular}[t]{l}\tiny{1}\end{tabular}}}}%
    \put(0.46901196,0.45223444){\color[rgb]{0,0,0}\makebox(0,0)[lt]{\lineheight{0}\smash{\begin{tabular}[t]{l}\tiny{18}\end{tabular}}}}%
    \put(0.93363823,0.60079238){\color[rgb]{0,0,0}\makebox(0,0)[lt]{\lineheight{0}\smash{\begin{tabular}[t]{l}\tiny{5}\end{tabular}}}}%
    \put(0.45925715,0.94723531){\color[rgb]{0,0,0}\makebox(0,0)[lt]{\lineheight{0}\smash{\begin{tabular}[t]{l}\tiny{8}\end{tabular}}}}%
    \put(0.6633666,0.53758497){\color[rgb]{0,0,0}\makebox(0,0)[lt]{\lineheight{0}\smash{\begin{tabular}[t]{l}\tiny{11}\end{tabular}}}}%
    \put(0.47021805,0.60493473){\color[rgb]{0,0,0}\makebox(0,0)[lt]{\lineheight{0}\smash{\begin{tabular}[t]{l}\tiny{12}\end{tabular}}}}%
    \put(0.25061402,0.36555192){\color[rgb]{0,0,0}\makebox(0,0)[lt]{\lineheight{0}\smash{\begin{tabular}[t]{l}\tiny{21}\end{tabular}}}}%
    \put(0.26394086,0.4730861){\color[rgb]{0,0,0}\makebox(0,0)[lt]{\lineheight{0}\smash{\begin{tabular}[t]{l}\tiny{17}\end{tabular}}}}%
    \put(0.67006109,0.35231346){\color[rgb]{0,0,0}\makebox(0,0)[lt]{\lineheight{0}\smash{\begin{tabular}[t]{l}\tiny{19}\end{tabular}}}}%
    \put(0.35077369,0.65872394){\color[rgb]{0,0,0}\makebox(0,0)[lt]{\lineheight{0}\smash{\begin{tabular}[t]{l}\tiny{13}\end{tabular}}}}%
    \put(0.07001943,0.39062464){\color[rgb]{0,0,0}\makebox(0,0)[lt]{\lineheight{0}\smash{\begin{tabular}[t]{l}\tiny{22}\end{tabular}}}}%
    \put(0.68163635,0.66931691){\color[rgb]{0,0,0}\makebox(0,0)[lt]{\lineheight{0}\smash{\begin{tabular}[t]{l}\tiny{6}\end{tabular}}}}%
    \put(0.21926386,0.7677135){\color[rgb]{0,0,0}\makebox(0,0)[lt]{\lineheight{0}\smash{\begin{tabular}[t]{l}\tiny{14}\end{tabular}}}}%
    \put(0.17118822,0.63719051){\color[rgb]{0,0,0}\makebox(0,0)[lt]{\lineheight{0}\smash{\begin{tabular}[t]{l}\tiny{15}\end{tabular}}}}%
    \put(0.34572906,0.8138057){\color[rgb]{0,0,0}\makebox(0,0)[lt]{\lineheight{0}\smash{\begin{tabular}[t]{l}\tiny{9}\end{tabular}}}}%
    \put(0.40517341,0.35157653){\color[rgb]{0,0,0}\makebox(0,0)[lt]{\lineheight{0}\smash{\begin{tabular}[t]{l}\tiny{20}\end{tabular}}}}%
    \put(0.37525285,0.19908235){\color[rgb]{0,0,0}\makebox(0,0)[lt]{\lineheight{0}\smash{\begin{tabular}[t]{l}\tiny{24}\end{tabular}}}}%
    \put(0.07372234,0.52871972){\color[rgb]{0,0,0}\makebox(0,0)[lt]{\lineheight{0}\smash{\begin{tabular}[t]{l}\tiny{16}\end{tabular}}}}%
    \put(0,0){\includegraphics[width=\unitlength,page=2]{Fig05c.pdf}}%
    \put(-0.09201796,0.74600391){\color[rgb]{0,0.50196078,0.50196078}\makebox(0,0)[lt]{\lineheight{0}\smash{\begin{tabular}[t]{l}\tiny{process 7}\end{tabular}}}}%
    \put(0.89885423,0.88901305){\color[rgb]{0.50196078,0,0.50196078}\makebox(0,0)[lt]{\lineheight{0}\smash{\begin{tabular}[t]{l}\tiny{process 2}\end{tabular}}}}%
    \put(0.92021202,0.51974065){\color[rgb]{0.50196078,0.50196078,0}\makebox(0,0)[lt]{\lineheight{0}\smash{\begin{tabular}[t]{l}\tiny{process 3 }\end{tabular}}}}%
    \put(0.75579117,0.14654528){\color[rgb]{1,0,1}\makebox(0,0)[lt]{\lineheight{0}\smash{\begin{tabular}[t]{l}\tiny{process 5}\end{tabular}}}}%
    \put(-0.04530285,0.99875923){\color[rgb]{0,0.8,1}\makebox(0,0)[lt]{\lineheight{0}\smash{\begin{tabular}[t]{l}\tiny{process 0 }\end{tabular}}}}%
    \put(0.56634056,1.17854491){\color[rgb]{1,0,0}\makebox(0,0)[lt]{\lineheight{0}\smash{\begin{tabular}[t]{l}\tiny{process 1}\end{tabular}}}}%
    \put(0.91647254,0.38665167){\color[rgb]{1,0.83137255,0.16470588}\makebox(0,0)[lt]{\lineheight{0}\smash{\begin{tabular}[t]{l}\tiny{process 4}\end{tabular}}}}%
    \put(0.01965601,0.18307101){\color[rgb]{0.21568627,0.78431373,0.44313725}\makebox(0,0)[lt]{\lineheight{0}\smash{\begin{tabular}[t]{l}\tiny{process 6}\end{tabular}}}}%
  \end{picture}%
\endgroup%
    \label{ts_dist_load_balance_8p}
    }
	\\
	\subfloat[][Legend]{
		\def\svgwidth{0.37\textwidth}
\begingroup%
  \makeatletter%
  \providecommand\color[2][]{%
    \errmessage{(Inkscape) Color is used for the text in Inkscape, but the package 'color.sty' is not loaded}%
    \renewcommand\color[2][]{}%
  }%
  \providecommand\transparent[1]{%
    \errmessage{(Inkscape) Transparency is used (non-zero) for the text in Inkscape, but the package 'transparent.sty' is not loaded}%
    \renewcommand\transparent[1]{}%
  }%
  \providecommand\rotatebox[2]{#2}%
  \ifx\svgwidth\undefined%
    \setlength{\unitlength}{303.08874107bp}%
    \ifx\svgscale\undefined%
      \relax%
    \else%
      \setlength{\unitlength}{\unitlength * \real{\svgscale}}%
    \fi%
  \else%
    \setlength{\unitlength}{\svgwidth}%
  \fi%
  \global\let\svgwidth\undefined%
  \global\let\svgscale\undefined%
  \makeatother%
  \begin{picture}(1,0.71704398)%
    \lineheight{1}%
    \setlength\tabcolsep{0pt}%
    \put(0,0){\includegraphics[width=\unitlength,page=1]{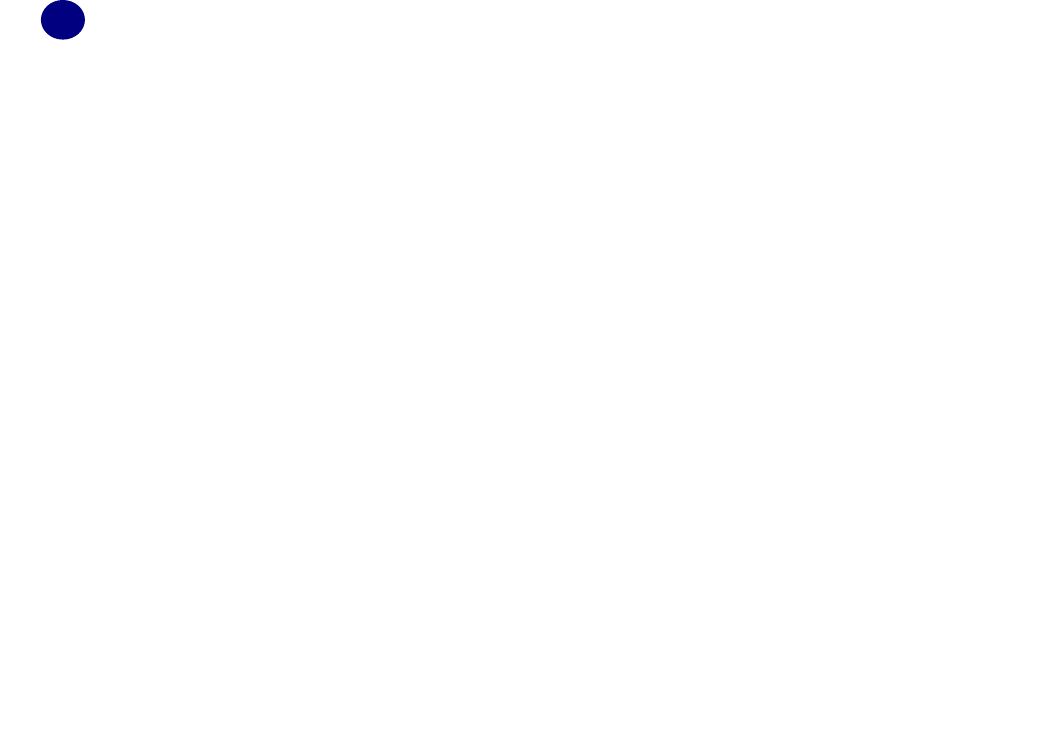}}%
    \put(0.24956891,0.69252891){\color[rgb]{0,0,0}\makebox(0,0)[lt]{\lineheight{0}\smash{\begin{tabular}[t]{l}\scriptsize{Global scale node with distributed support}\end{tabular}}}}%
    \put(0,0){\includegraphics[width=\unitlength,page=2]{Fig05d.pdf}}%
    \put(0.24961644,0.59985417){\color[rgb]{0,0,0}\makebox(0,0)[lt]{\lineheight{0}\smash{\begin{tabular}[t]{l}\scriptsize{Global scale node with non distributed support}\end{tabular}}}}%
    \put(0,0){\includegraphics[width=\unitlength,page=3]{Fig05d.pdf}}%
    \put(0.25355642,0.30656925){\color[rgb]{0,0,0}\makebox(0,0)[lt]{\lineheight{0}\smash{\begin{tabular}[t]{l}\scriptsize{Global-scale mesh}\end{tabular}}}}%
    \put(0,0){\includegraphics[width=\unitlength,page=4]{Fig05d.pdf}}%
    \put(0.40199406,0.04134013){\color[rgb]{0,0,0}\makebox(0,0)[lt]{\lineheight{0}\smash{\begin{tabular}[t]{l}\scriptsize{Process id color}\end{tabular}}}}%
    \put(0,0){\includegraphics[width=\unitlength,page=5]{Fig05d.pdf}}%
    \put(0.31214749,0.51877158){\color[rgb]{0,0,0}\makebox(0,0)[lt]{\lineheight{0}\smash{\begin{tabular}[t]{l}\scriptsize{Interprocess interface}\end{tabular}}}}%
    \put(0,0){\includegraphics[width=\unitlength,page=6]{Fig05d.pdf}}%
    \put(0.26072418,0.18327801){\color[rgb]{0,0,0}\makebox(0,0)[lt]{\lineheight{0}\smash{\begin{tabular}[t]{l}\scriptsize{Fine-scale mesh}\end{tabular}}}}%
    \put(0,0){\includegraphics[width=\unitlength,page=7]{Fig05d.pdf}}%
    \put(0.31484353,0.43445197){\color[rgb]{0,0,0}\makebox(0,0)[lt]{\lineheight{0}\smash{\begin{tabular}[t]{l}\scriptsize{Process boundary}\end{tabular}}}}%
  \end{picture}%
\endgroup%
		\label{ts_dist_legend}
	}
	\subfloat[][Fine scale problems (2 out 25) following global element dispatching on 4 processes (\protect\subref*{ts_dist_load_balance})]{
	\def\svgwidth{0.43\textwidth}
\begingroup%
  \makeatletter%
  \providecommand\color[2][]{%
    \errmessage{(Inkscape) Color is used for the text in Inkscape, but the package 'color.sty' is not loaded}%
    \renewcommand\color[2][]{}%
  }%
  \providecommand\transparent[1]{%
    \errmessage{(Inkscape) Transparency is used (non-zero) for the text in Inkscape, but the package 'transparent.sty' is not loaded}%
    \renewcommand\transparent[1]{}%
  }%
  \providecommand\rotatebox[2]{#2}%
  \ifx\svgwidth\undefined%
    \setlength{\unitlength}{457.82834852bp}%
    \ifx\svgscale\undefined%
      \relax%
    \else%
      \setlength{\unitlength}{\unitlength * \real{\svgscale}}%
    \fi%
  \else%
    \setlength{\unitlength}{\svgwidth}%
  \fi%
  \global\let\svgwidth\undefined%
  \global\let\svgscale\undefined%
  \makeatother%
  \begin{picture}(1,0.73166042)%
    \lineheight{1}%
    \setlength\tabcolsep{0pt}%
    \put(0,0){\includegraphics[width=\unitlength,page=1]{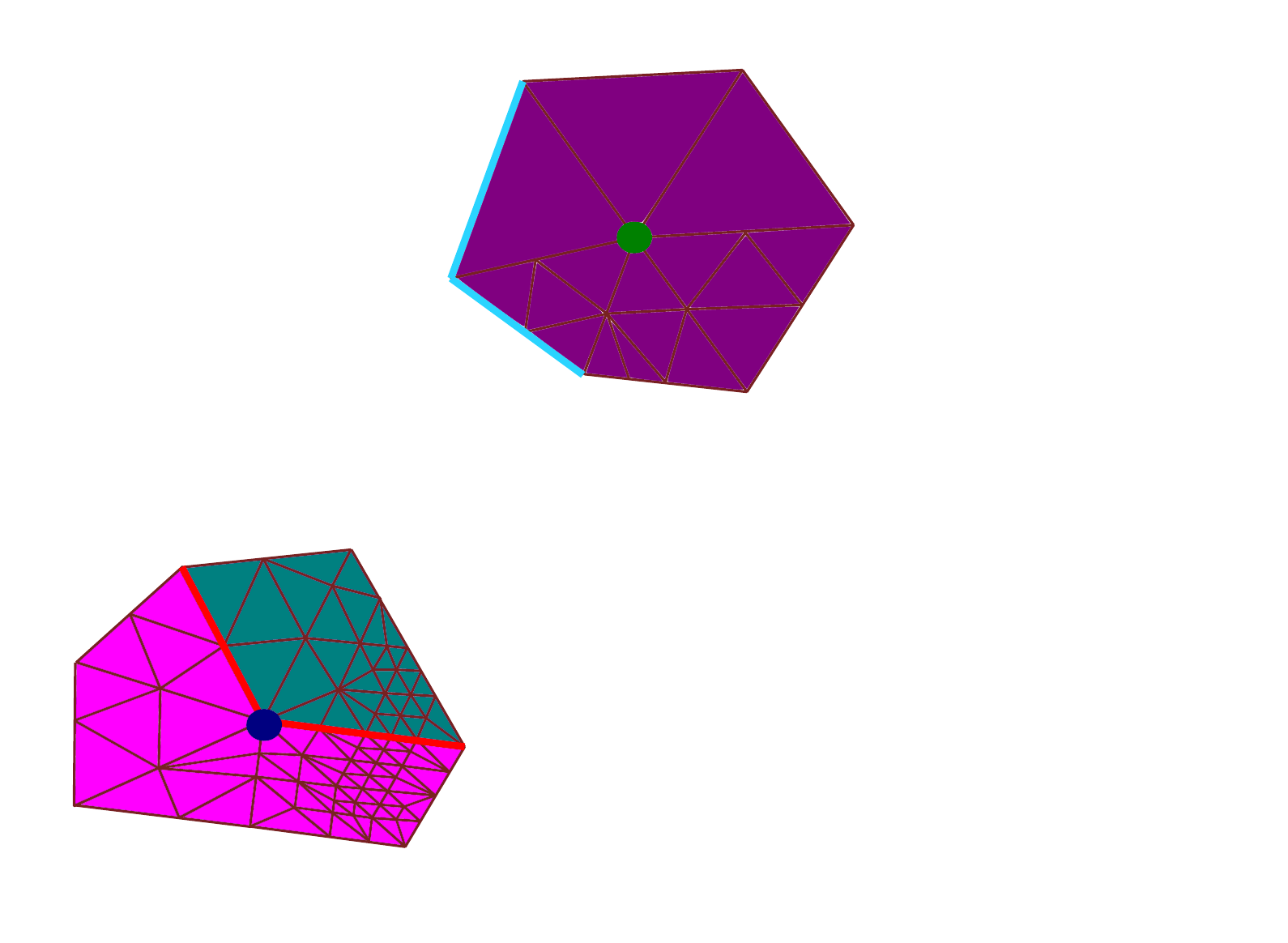}}%
    \put(0.03114623,0.68638138){\color[rgb]{0.50196078,0,0.50196078}\makebox(0,0)[lt]{\lineheight{0}\smash{\begin{tabular}[t]{l}\scriptsize{Global comm: process 1} / \scriptsize{Local comm: process 0}\end{tabular}}}}%
    \put(0.07733094,0.03550769){\color[rgb]{1,0,1}\makebox(0,0)[lt]{\lineheight{0}\smash{\begin{tabular}[t]{l}\scriptsize{Global comm: process 3} / \scriptsize{Local comm: process 1}\end{tabular}}}}%
    \put(0.33973626,0.2837736){\color[rgb]{0,0.50196078,0.50196078}\makebox(0,0)[lt]{\lineheight{0}\smash{\begin{tabular}[t]{l}\scriptsize{Global comm: process 0} / \scriptsize{Local comm: process 0 }\end{tabular}}}}%
  \end{picture}%
\endgroup%
	\label{ts_dist_fine}
}
	\caption{\tS fictitious 2D problem of figure\ref{ts_fig} distributed on processes}
	\label{ts_dist_fig}
\end{figure}

\section{Parallel paradigm}\label{parallel_paradigm}
This work does not use shared-memory multithreaded parallelism and uses only the distributed-memory message passing paradigm \footnote{In this work we use the MPI-3.0 (or higher) message passing standard.}.
This choice comes from the thread-safety constraint which would have generated too much implementation work for this first study (for example the in-house library used to code the  \tS method is not thread-safe).
As a result, multi-threading is never used in this study and to be fair when comparing with other methods, their multi-threading capability was also disabled.
Thus, in this section, parallelism will always concern the message passing paradigm.

Parallelism is introduced at all scales in all \tS frameworks, but with a quite different granularity.
At the global level, the mesh is distributed over processes (see  figures \ref{ts_dist_load_balance} and \ref{ts_dist_load_balance_8p} for the example of figure \ref{ts_fig} distributed over 4 and 8 processes respectively), based on a macro-element partitioning,  in order to allow a good load balancing.
This partitioning (not evaluated in depth in this paper), given by ParMetis (\cite{Karypis1997}) in this work, is supposed to balance the cost of creating linear systems, assembling them and solving them at both scales.
To do this, a cost per macro-element is obtained  directly  from  the number of its encapsulated micro-elements (see figure \ref{ts_dist_weigth_elem}) and the number of enriched nodes.
If the entire global-scale mesh is refined (as in the examples in the section \ref{UOL} and \ref{MS}), all macro elements have roughly the same number of encapsulated micro elements and enriched nodes.
Thus, the cost per element can be considered as a constant.
Now, if only a specific region is involved (see the 2D example in Figure \ref{ts_coarse_scale} or the example in section \ref{PO}), the number of encapsulated micro-elements and enriched nodes provides weights to the ParaMetis multi-objective partitioning optimization algorithm.
In the examples studied, this partitioning provides great scalability for the assembly task, as all processes receive roughly the same amount of micro-elements to integrate and assemble (see in particular the table \ref{UOL_table_ass}  in section \ref{UOL}).
It also allows the use of a fairly well balanced parallel resolution for the linear system related to the global-scale problem ($\dagger$ in the algorithm \ref{TS_algebra}).
The section \ref{gsolv} gives more details on the choice of solver to use at this level.

At the fine scale, patches follow the distribution of the global scale mesh (see figure \ref{ts_dist_fine}) and are therefore naturally distributed in most processes.
This results, as in \cite{Kim2011}, in a high level of coarse-grained parallelism because many fine-scale problems are solved independently (without communication) at the same time  by many processes.
But some patches are split between the computational units due to the global mesh distribution (see the left patch in the figure \ref{ts_dist_fine}).
In this work the distributed patch problems  are solved entirely in parallel without any overlapping mechanism (phantom variables, ...).
The  assembly of the problems ($\vm{A}_{qq}^{p}$ and $\vm{BI}_{q}^{p}$ of algorithm \ref{TS_algebra_init} and $\vm{B}_{q}^{p}$ of algorithm \ref{TS_algebra_patch}) is thus done in parallel without communication.
Then, the matrix contributions of the linear system are provided to a parallel solver ($\ddagger$ of the algorithm \ref{TS_algebra_patch}) and treated with a mid-grain parallelism.
But this introduces some constraints on the overall process of fine-scale parallel resolution.
The distributed patches are no longer independent.
Typically, if two distributed patches share the same process, it will be impossible to run the resolution of both in that process.
The solver is not able to treat two independent matrices at the same time.
This implies that the computation of distributed patches  must be scheduled as explained in \ref{scheduling}.
These distributed patches have an ambivalent impact on performance.
On the one hand, they reduce the time to solve a rather large fine-scale problem.
And on the other hand, if their number increases too much, it expands the amount of communication and can also have an impact on sequencing.

\subsection{global scale problem parallel resolution }\label{gsolv}
In this work, we first chose an asynchronous parallel MPI  sparse  direct  solver (MUMPS \cite{Amestoy2001,Amestoy2006}) for the resolution of  the linear system associated with the global problem.
It was a safe choice (even if the matrix becomes badly conditioned) for this first implementation but by no means a requirement.
The symbolic step is performed only once at the beginning of the \tS loop.
Factorization and solving step are repeated for each iteration with the updated $\vm{A}_{gg}$ and $\vm{B}_{g}$.
With many patches, as the fine scale resolution has good scalability (see the results presented in section \ref{NSOLV}), many processes can be used well beyond the scalability of the parallel direct solver which solves at the coarse-scale, a problem of rather smaller size.
To limit this effect, in all simulations tested,  a subset of the available processes (taking at most $nbp_{max}$ processes, $nbp_{max}$ arbitrarily chosen to stay within the range of good performance of the direct solver) was used to solve this coarse-scale problem (gathering $\vm{A}_{gg}$ contribution of the processes eliminated from the resolution in the $nbp_{max}$ retained processes).
This version is referred to as \TSD below.
This reduction had a positive effect, but the overall speed-up\footnote{the speed-up is the ratio between the time used by an application in single-process mode and the time used in multi-process mode.
	Ideally, it should give the number of processes used. } of the \tS resolution was still deteriorated by the coarse-scale efficiency.
To reduce this effect, as the \tS loop iterates over a coarse-scale problem without changing its discretization, it is possible to use a parallel preconditioned iterative solver (conjugate gradient) that starts with the solution of the previous iteration.
This choice  intrinsically reduces the coarse-scale resolution by considering that the successive solutions are not too far apart from each other, so that the number of iterations in the iterative solver should be small.
And this choice  adds scalability because the conjugate gradient scales well (to the point where the scalar product becomes a problem).
A natural choice is to use for  the preconditioner  a previous factorization performed by the direct solver.
The preconditioning step  then only costs a forward/backward resolution, performed in parallel, following the distribution of the direct solver.
If not done to many times, it is cheap compared to a full factorization.
Note that this technique does not take into account the stable nature of the classical dofs of the coarse-scale problem as  is the case in \cite{Kim2015}.
This would have added additional complexity regarding load balancing and resolution.
Nevertheless, the proposed approach offers apriori a slightly richer preconditioner  compared to the block Jacobi proposed in this quote because it incorporates all terms (including the classic x enriched coupling blocks).
The use of the iterative solver and the quality of the preconditioner  are controlled using the following arbitrary rules:
\begin{enumerate}
	\item The iterative resolution is used instead of the direct resolution as soon as the criterion "$resi$"  becomes smaller than $\epsilon\times 10000$ ($\epsilon$ being the precision to be reach by "$resi$" in algorithm \ref{TS_algebra}  ) and as soon as 2 \tS iterations have been performed
	\item When the iterative resolution takes more than 13 iterations to converge, the preconditioner is recomputed: at the next \tS iteration, the direct solver is used again to provide a new solution and a new factorization.
	The iterative solver is then reused in the next \tS iteration.
\end{enumerate} 
The second rule avoids spending too much time in the iterative solver when the  preconditioner is not of high quality.
The first rule requires  waiting for the \tS solution to reach a minimum quality so that the variation between the \tS iteration is small enough for the preconditioner to efficiently reduce the number of iterations of the iterative solver.
This version is hereafter referred to  as \TSIn.
Also note that in this version, for the direct resolution phase,  low-rank resolution (see \ref{BLRRES}) can also be used instead of full rank resolution.
This is left as future work.

At the coarse scale, another type of solver ( a domain decomposition solver presented in \ref{DDRES}) was also tested during this study in order to gain scalability as illustrated in section \ref{PO}.

\begin{figure}[!htb]
	\subfloat[][weight per patch]{
		\def\svgwidth{30mm}
\begingroup%
  \makeatletter%
  \providecommand\color[2][]{%
    \errmessage{(Inkscape) Color is used for the text in Inkscape, but the package 'color.sty' is not loaded}%
    \renewcommand\color[2][]{}%
  }%
  \providecommand\transparent[1]{%
    \errmessage{(Inkscape) Transparency is used (non-zero) for the text in Inkscape, but the package 'transparent.sty' is not loaded}%
    \renewcommand\transparent[1]{}%
  }%
  \providecommand\rotatebox[2]{#2}%
  \ifx\svgwidth\undefined%
    \setlength{\unitlength}{342.09281158bp}%
    \ifx\svgscale\undefined%
      \relax%
    \else%
      \setlength{\unitlength}{\unitlength * \real{\svgscale}}%
    \fi%
  \else%
    \setlength{\unitlength}{\svgwidth}%
  \fi%
  \global\let\svgwidth\undefined%
  \global\let\svgscale\undefined%
  \makeatother%
  \begin{picture}(1,1.20965849)%
    \lineheight{1}%
    \setlength\tabcolsep{0pt}%
    \put(0.59348497,0.8466883){\color[rgb]{0,0,0}\makebox(0,0)[lt]{\lineheight{0}\smash{\begin{tabular}[t]{l}\tiny{87}\end{tabular}}}}%
    \put(0.96309194,0.79464572){\color[rgb]{0,0,0}\makebox(0,0)[lt]{\lineheight{0}\smash{\begin{tabular}[t]{l}\tiny{8}\end{tabular}}}}%
    \put(0.65617276,0.98159756){\color[rgb]{0,0,0}\makebox(0,0)[lt]{\lineheight{0}\smash{\begin{tabular}[t]{l}\tiny{16}\end{tabular}}}}%
    \put(0.47072037,0.75909525){\color[rgb]{0,0,0}\makebox(0,0)[lt]{\lineheight{0}\smash{\begin{tabular}[t]{l}\tiny{133}\end{tabular}}}}%
    \put(0.25942799,0.27036035){\color[rgb]{0,0,0}\makebox(0,0)[lt]{\lineheight{0}\smash{\begin{tabular}[t]{l}\tiny{13}\end{tabular}}}}%
    \put(0.54620276,0.20435236){\color[rgb]{0,0,0}\makebox(0,0)[lt]{\lineheight{0}\smash{\begin{tabular}[t]{l}\tiny{9}\end{tabular}}}}%
    \put(0.77146822,0.81465496){\color[rgb]{0,0,0}\makebox(0,0)[lt]{\lineheight{0}\smash{\begin{tabular}[t]{l}\tiny{50}\end{tabular}}}}%
    \put(0.88940624,1.00073606){\color[rgb]{0,0,0}\makebox(0,0)[lt]{\lineheight{0}\smash{\begin{tabular}[t]{l}\tiny{9}\end{tabular}}}}%
    \put(0.46901196,0.45223444){\color[rgb]{0,0,0}\makebox(0,0)[lt]{\lineheight{0}\smash{\begin{tabular}[t]{l}\tiny{137}\end{tabular}}}}%
    \put(0.93363823,0.60079238){\color[rgb]{0,0,0}\makebox(0,0)[lt]{\lineheight{0}\smash{\begin{tabular}[t]{l}\tiny{23}\end{tabular}}}}%
    \put(0.45925715,0.94723531){\color[rgb]{0,0,0}\makebox(0,0)[lt]{\lineheight{0}\smash{\begin{tabular}[t]{l}\tiny{22}\end{tabular}}}}%
    \put(0.6633666,0.53758497){\color[rgb]{0,0,0}\makebox(0,0)[lt]{\lineheight{0}\smash{\begin{tabular}[t]{l}\tiny{39}\end{tabular}}}}%
    \put(0.47021805,0.60493473){\color[rgb]{0,0,0}\makebox(0,0)[lt]{\lineheight{0}\smash{\begin{tabular}[t]{l}\tiny{134}\end{tabular}}}}%
    \put(0.25061402,0.36555192){\color[rgb]{0,0,0}\makebox(0,0)[lt]{\lineheight{0}\smash{\begin{tabular}[t]{l}\tiny{30}\end{tabular}}}}%
    \put(0.26394086,0.4730861){\color[rgb]{0,0,0}\makebox(0,0)[lt]{\lineheight{0}\smash{\begin{tabular}[t]{l}\tiny{108}\end{tabular}}}}%
    \put(0.67006109,0.35231346){\color[rgb]{0,0,0}\makebox(0,0)[lt]{\lineheight{0}\smash{\begin{tabular}[t]{l}\tiny{18}\end{tabular}}}}%
    \put(0.35077369,0.65872394){\color[rgb]{0,0,0}\makebox(0,0)[lt]{\lineheight{0}\smash{\begin{tabular}[t]{l}\tiny{125}\end{tabular}}}}%
    \put(0.07001943,0.39062464){\color[rgb]{0,0,0}\makebox(0,0)[lt]{\lineheight{0}\smash{\begin{tabular}[t]{l}\tiny{11}\end{tabular}}}}%
    \put(0.68163635,0.66931691){\color[rgb]{0,0,0}\makebox(0,0)[lt]{\lineheight{0}\smash{\begin{tabular}[t]{l}\tiny{137}\end{tabular}}}}%
    \put(0.21926386,0.7677135){\color[rgb]{0,0,0}\makebox(0,0)[lt]{\lineheight{0}\smash{\begin{tabular}[t]{l}\tiny{12}\end{tabular}}}}%
    \put(0.17118822,0.63719051){\color[rgb]{0,0,0}\makebox(0,0)[lt]{\lineheight{0}\smash{\begin{tabular}[t]{l}\tiny{15}\end{tabular}}}}%
    \put(0.34572906,0.8138057){\color[rgb]{0,0,0}\makebox(0,0)[lt]{\lineheight{0}\smash{\begin{tabular}[t]{l}\tiny{17}\end{tabular}}}}%
    \put(0.40517341,0.35157653){\color[rgb]{0,0,0}\makebox(0,0)[lt]{\lineheight{0}\smash{\begin{tabular}[t]{l}\tiny{82}\end{tabular}}}}%
    \put(0.37525285,0.19908235){\color[rgb]{0,0,0}\makebox(0,0)[lt]{\lineheight{0}\smash{\begin{tabular}[t]{l}\tiny{9}\end{tabular}}}}%
    \put(0.07372234,0.52871972){\color[rgb]{0,0,0}\makebox(0,0)[lt]{\lineheight{0}\smash{\begin{tabular}[t]{l}\tiny{9}\end{tabular}}}}%
    \put(0,0){\includegraphics[width=\unitlength,page=1]{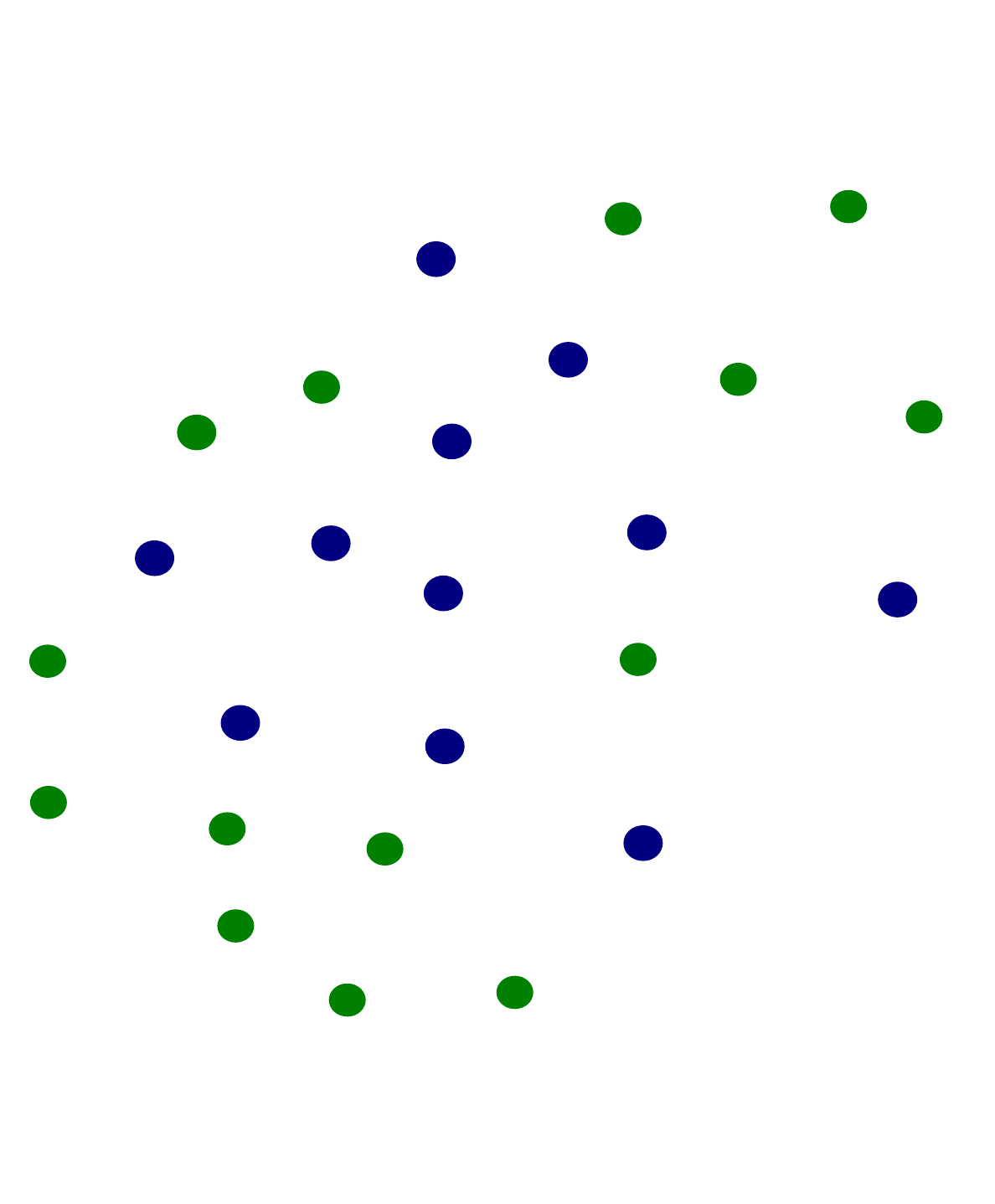}}%
  \end{picture}%
\endgroup%
		\label{seq:4weigt}
	}
	\subfloat[][Initial graph with arbitrary node id]{
		\def\svgwidth{30mm}
\begingroup%
  \makeatletter%
  \providecommand\color[2][]{%
    \errmessage{(Inkscape) Color is used for the text in Inkscape, but the package 'color.sty' is not loaded}%
    \renewcommand\color[2][]{}%
  }%
  \providecommand\transparent[1]{%
    \errmessage{(Inkscape) Transparency is used (non-zero) for the text in Inkscape, but the package 'transparent.sty' is not loaded}%
    \renewcommand\transparent[1]{}%
  }%
  \providecommand\rotatebox[2]{#2}%
  \ifx\svgwidth\undefined%
    \setlength{\unitlength}{342.09281158bp}%
    \ifx\svgscale\undefined%
      \relax%
    \else%
      \setlength{\unitlength}{\unitlength * \real{\svgscale}}%
    \fi%
  \else%
    \setlength{\unitlength}{\svgwidth}%
  \fi%
  \global\let\svgwidth\undefined%
  \global\let\svgscale\undefined%
  \makeatother%
  \begin{picture}(1,1.20965849)%
    \lineheight{1}%
    \setlength\tabcolsep{0pt}%
    \put(0,0){\includegraphics[width=\unitlength,page=1]{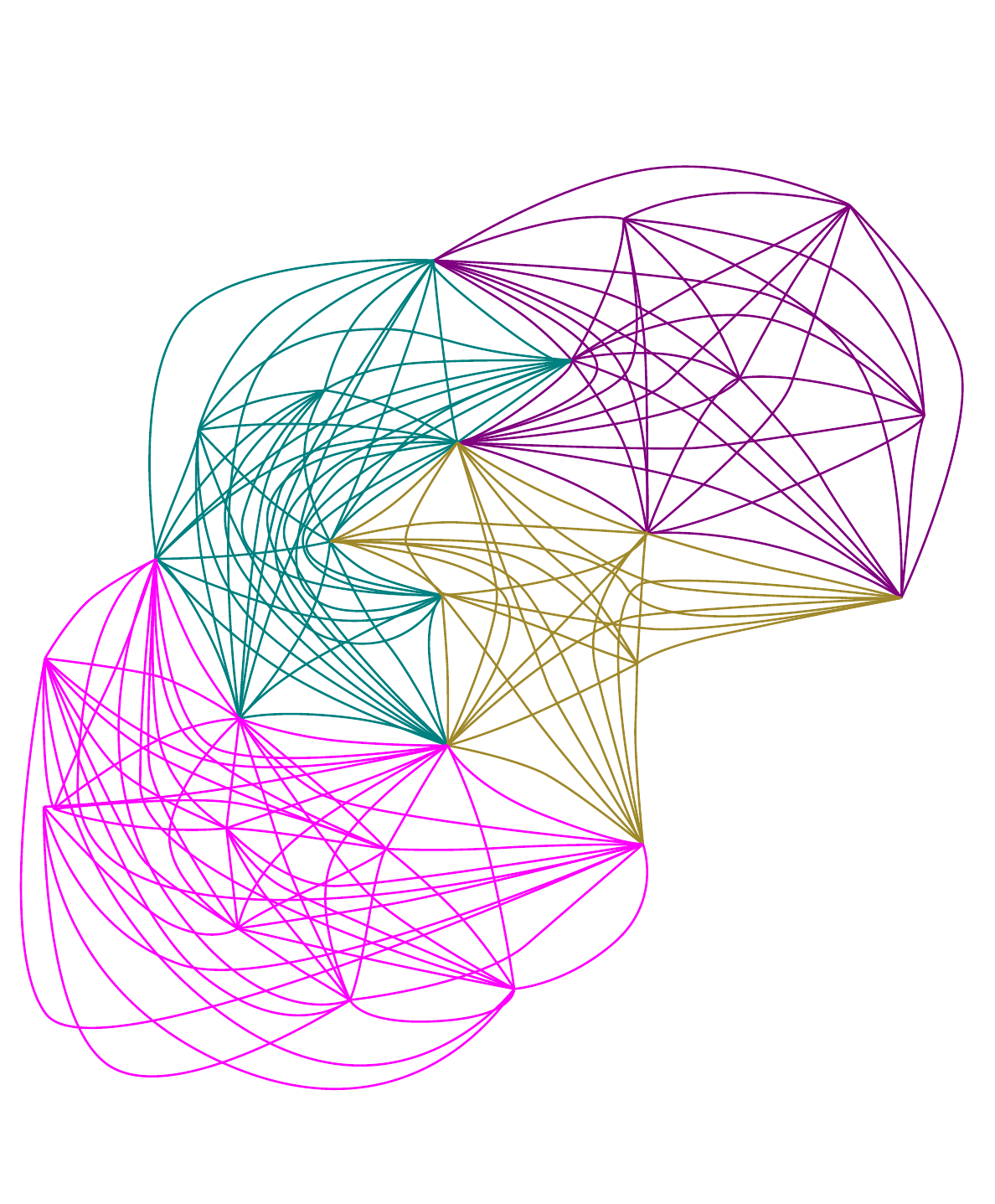}}%
    \put(0.59348497,0.8466883){\color[rgb]{0,0,0}\makebox(0,0)[lt]{\lineheight{0}\smash{\begin{tabular}[t]{l}\tiny{7}\end{tabular}}}}%
    \put(0.96309194,0.79464572){\color[rgb]{0,0,0}\makebox(0,0)[lt]{\lineheight{0}\smash{\begin{tabular}[t]{l}\tiny{4}\end{tabular}}}}%
    \put(0.65617276,0.98159756){\color[rgb]{0,0,0}\makebox(0,0)[lt]{\lineheight{0}\smash{\begin{tabular}[t]{l}\tiny{2}\end{tabular}}}}%
    \put(0.47072037,0.75909525){\color[rgb]{0,0,0}\makebox(0,0)[lt]{\lineheight{0}\smash{\begin{tabular}[t]{l}\tiny{10}\end{tabular}}}}%
    \put(0.25942799,0.27036035){\color[rgb]{0,0,0}\makebox(0,0)[lt]{\lineheight{0}\smash{\begin{tabular}[t]{l}\tiny{23}\end{tabular}}}}%
    \put(0.54620276,0.20435236){\color[rgb]{0,0,0}\makebox(0,0)[lt]{\lineheight{0}\smash{\begin{tabular}[t]{l}\tiny{25}\end{tabular}}}}%
    \put(0.77146822,0.81465496){\color[rgb]{0,0,0}\makebox(0,0)[lt]{\lineheight{0}\smash{\begin{tabular}[t]{l}\tiny{3}\end{tabular}}}}%
    \put(0.88940624,1.00073606){\color[rgb]{0,0,0}\makebox(0,0)[lt]{\lineheight{0}\smash{\begin{tabular}[t]{l}\tiny{1}\end{tabular}}}}%
    \put(0.46901196,0.45223444){\color[rgb]{0,0,0}\makebox(0,0)[lt]{\lineheight{0}\smash{\begin{tabular}[t]{l}\tiny{18}\end{tabular}}}}%
    \put(0.93363823,0.60079238){\color[rgb]{0,0,0}\makebox(0,0)[lt]{\lineheight{0}\smash{\begin{tabular}[t]{l}\tiny{5}\end{tabular}}}}%
    \put(0.45925715,0.94723531){\color[rgb]{0,0,0}\makebox(0,0)[lt]{\lineheight{0}\smash{\begin{tabular}[t]{l}\tiny{8}\end{tabular}}}}%
    \put(0.6633666,0.53758497){\color[rgb]{0,0,0}\makebox(0,0)[lt]{\lineheight{0}\smash{\begin{tabular}[t]{l}\tiny{11}\end{tabular}}}}%
    \put(0.47021805,0.60493473){\color[rgb]{0,0,0}\makebox(0,0)[lt]{\lineheight{0}\smash{\begin{tabular}[t]{l}\tiny{12}\end{tabular}}}}%
    \put(0.25061402,0.36555192){\color[rgb]{0,0,0}\makebox(0,0)[lt]{\lineheight{0}\smash{\begin{tabular}[t]{l}\tiny{21}\end{tabular}}}}%
    \put(0.26394086,0.4730861){\color[rgb]{0,0,0}\makebox(0,0)[lt]{\lineheight{0}\smash{\begin{tabular}[t]{l}\tiny{17}\end{tabular}}}}%
    \put(0.67006109,0.35231346){\color[rgb]{0,0,0}\makebox(0,0)[lt]{\lineheight{0}\smash{\begin{tabular}[t]{l}\tiny{19}\end{tabular}}}}%
    \put(0.35077369,0.65872394){\color[rgb]{0,0,0}\makebox(0,0)[lt]{\lineheight{0}\smash{\begin{tabular}[t]{l}\tiny{13}\end{tabular}}}}%
    \put(0.07001943,0.39062464){\color[rgb]{0,0,0}\makebox(0,0)[lt]{\lineheight{0}\smash{\begin{tabular}[t]{l}\tiny{22}\end{tabular}}}}%
    \put(0.68163635,0.66931691){\color[rgb]{0,0,0}\makebox(0,0)[lt]{\lineheight{0}\smash{\begin{tabular}[t]{l}\tiny{6}\end{tabular}}}}%
    \put(0.21926386,0.7677135){\color[rgb]{0,0,0}\makebox(0,0)[lt]{\lineheight{0}\smash{\begin{tabular}[t]{l}\tiny{14}\end{tabular}}}}%
    \put(0.17118822,0.63719051){\color[rgb]{0,0,0}\makebox(0,0)[lt]{\lineheight{0}\smash{\begin{tabular}[t]{l}\tiny{15}\end{tabular}}}}%
    \put(0.34572906,0.8138057){\color[rgb]{0,0,0}\makebox(0,0)[lt]{\lineheight{0}\smash{\begin{tabular}[t]{l}\tiny{9}\end{tabular}}}}%
    \put(0.40517341,0.35157653){\color[rgb]{0,0,0}\makebox(0,0)[lt]{\lineheight{0}\smash{\begin{tabular}[t]{l}\tiny{20}\end{tabular}}}}%
    \put(0.37525285,0.19908235){\color[rgb]{0,0,0}\makebox(0,0)[lt]{\lineheight{0}\smash{\begin{tabular}[t]{l}\tiny{24}\end{tabular}}}}%
    \put(0.07372234,0.52871972){\color[rgb]{0,0,0}\makebox(0,0)[lt]{\lineheight{0}\smash{\begin{tabular}[t]{l}\tiny{16}\end{tabular}}}}%
    \put(0,0){\includegraphics[width=\unitlength,page=2]{Fig06b.pdf}}%
  \end{picture}%
\endgroup%
		\label{seq:4pini}
	}
	\subfloat[][sequence 1]{
		\def\svgwidth{30mm}
\begingroup%
  \makeatletter%
  \providecommand\color[2][]{%
    \errmessage{(Inkscape) Color is used for the text in Inkscape, but the package 'color.sty' is not loaded}%
    \renewcommand\color[2][]{}%
  }%
  \providecommand\transparent[1]{%
    \errmessage{(Inkscape) Transparency is used (non-zero) for the text in Inkscape, but the package 'transparent.sty' is not loaded}%
    \renewcommand\transparent[1]{}%
  }%
  \providecommand\rotatebox[2]{#2}%
  \ifx\svgwidth\undefined%
    \setlength{\unitlength}{342.09281158bp}%
    \ifx\svgscale\undefined%
      \relax%
    \else%
      \setlength{\unitlength}{\unitlength * \real{\svgscale}}%
    \fi%
  \else%
    \setlength{\unitlength}{\svgwidth}%
  \fi%
  \global\let\svgwidth\undefined%
  \global\let\svgscale\undefined%
  \makeatother%
  \begin{picture}(1,1.20965849)%
    \lineheight{1}%
    \setlength\tabcolsep{0pt}%
    \put(0,0){\includegraphics[width=\unitlength,page=1]{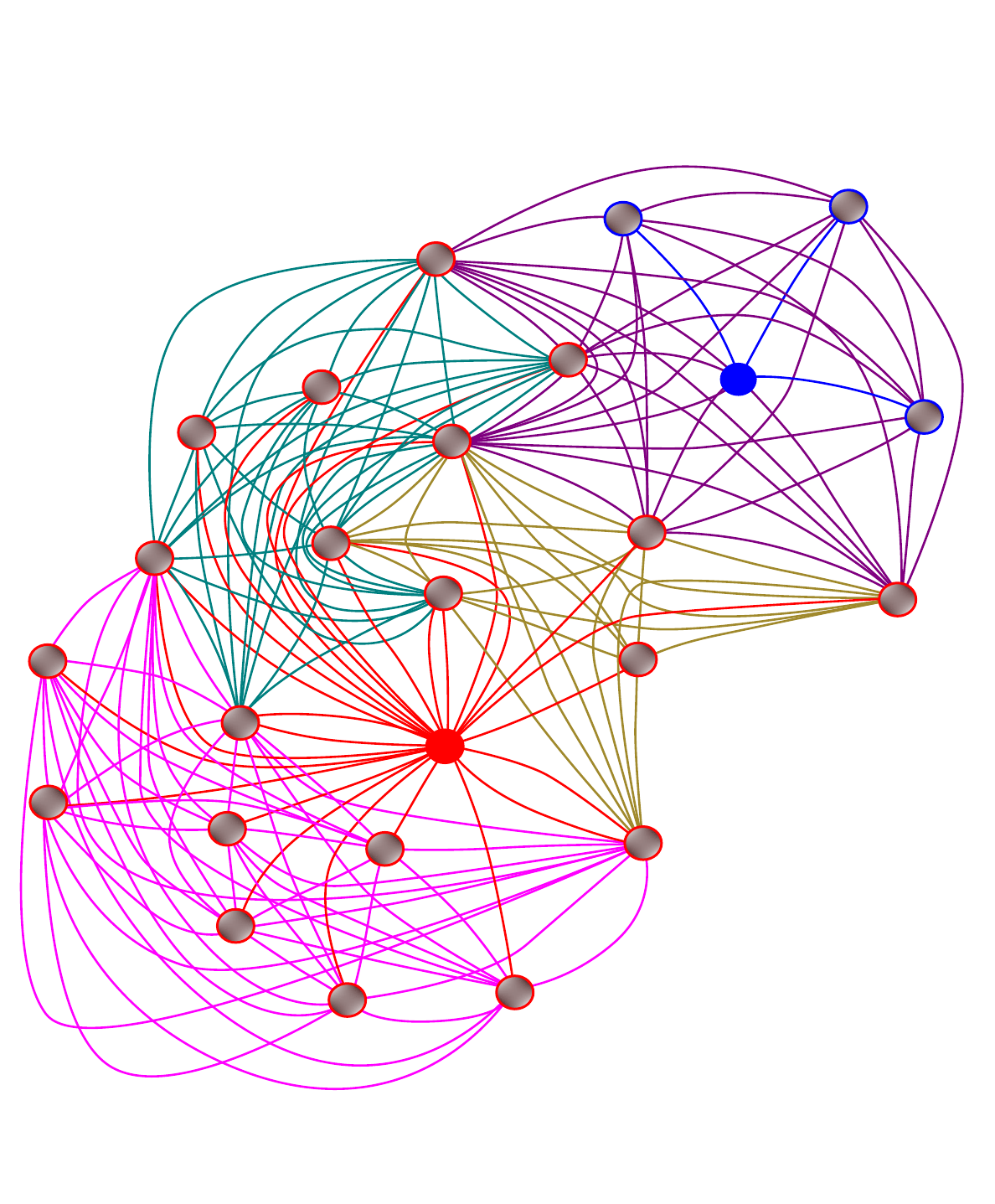}}%
    \put(0.59348497,0.8466883){\color[rgb]{0,0,0}\makebox(0,0)[lt]{\lineheight{0}\smash{\begin{tabular}[t]{l}\tiny{7}\end{tabular}}}}%
    \put(0.96309194,0.79464572){\color[rgb]{0,0,0}\makebox(0,0)[lt]{\lineheight{0}\smash{\begin{tabular}[t]{l}\tiny{4}\end{tabular}}}}%
    \put(0.65617276,0.98159756){\color[rgb]{0,0,0}\makebox(0,0)[lt]{\lineheight{0}\smash{\begin{tabular}[t]{l}\tiny{2}\end{tabular}}}}%
    \put(0.47072037,0.75909525){\color[rgb]{0,0,0}\makebox(0,0)[lt]{\lineheight{0}\smash{\begin{tabular}[t]{l}\tiny{10}\end{tabular}}}}%
    \put(0.25942799,0.27036035){\color[rgb]{0,0,0}\makebox(0,0)[lt]{\lineheight{0}\smash{\begin{tabular}[t]{l}\tiny{23}\end{tabular}}}}%
    \put(0.54620276,0.20435236){\color[rgb]{0,0,0}\makebox(0,0)[lt]{\lineheight{0}\smash{\begin{tabular}[t]{l}\tiny{25}\end{tabular}}}}%
    \put(0.77146822,0.81465496){\color[rgb]{0,0,0}\makebox(0,0)[lt]{\lineheight{0}\smash{\begin{tabular}[t]{l}\tiny{3}\end{tabular}}}}%
    \put(0.88940624,1.00073606){\color[rgb]{0,0,0}\makebox(0,0)[lt]{\lineheight{0}\smash{\begin{tabular}[t]{l}\tiny{1}\end{tabular}}}}%
    \put(0.46901196,0.45223444){\color[rgb]{0,0,0}\makebox(0,0)[lt]{\lineheight{0}\smash{\begin{tabular}[t]{l}\tiny{18}\end{tabular}}}}%
    \put(0.93363823,0.60079238){\color[rgb]{0,0,0}\makebox(0,0)[lt]{\lineheight{0}\smash{\begin{tabular}[t]{l}\tiny{5}\end{tabular}}}}%
    \put(0.45925715,0.94723531){\color[rgb]{0,0,0}\makebox(0,0)[lt]{\lineheight{0}\smash{\begin{tabular}[t]{l}\tiny{8}\end{tabular}}}}%
    \put(0.6633666,0.53758497){\color[rgb]{0,0,0}\makebox(0,0)[lt]{\lineheight{0}\smash{\begin{tabular}[t]{l}\tiny{11}\end{tabular}}}}%
    \put(0.47021805,0.60493473){\color[rgb]{0,0,0}\makebox(0,0)[lt]{\lineheight{0}\smash{\begin{tabular}[t]{l}\tiny{12}\end{tabular}}}}%
    \put(0.25061402,0.36555192){\color[rgb]{0,0,0}\makebox(0,0)[lt]{\lineheight{0}\smash{\begin{tabular}[t]{l}\tiny{21}\end{tabular}}}}%
    \put(0.26394086,0.4730861){\color[rgb]{0,0,0}\makebox(0,0)[lt]{\lineheight{0}\smash{\begin{tabular}[t]{l}\tiny{17}\end{tabular}}}}%
    \put(0.67006109,0.35231346){\color[rgb]{0,0,0}\makebox(0,0)[lt]{\lineheight{0}\smash{\begin{tabular}[t]{l}\tiny{19}\end{tabular}}}}%
    \put(0.35077369,0.65872394){\color[rgb]{0,0,0}\makebox(0,0)[lt]{\lineheight{0}\smash{\begin{tabular}[t]{l}\tiny{13}\end{tabular}}}}%
    \put(0.07001943,0.39062464){\color[rgb]{0,0,0}\makebox(0,0)[lt]{\lineheight{0}\smash{\begin{tabular}[t]{l}\tiny{22}\end{tabular}}}}%
    \put(0.68163635,0.66931691){\color[rgb]{0,0,0}\makebox(0,0)[lt]{\lineheight{0}\smash{\begin{tabular}[t]{l}\tiny{6}\end{tabular}}}}%
    \put(0.21926386,0.7677135){\color[rgb]{0,0,0}\makebox(0,0)[lt]{\lineheight{0}\smash{\begin{tabular}[t]{l}\tiny{14}\end{tabular}}}}%
    \put(0.17118822,0.63719051){\color[rgb]{0,0,0}\makebox(0,0)[lt]{\lineheight{0}\smash{\begin{tabular}[t]{l}\tiny{15}\end{tabular}}}}%
    \put(0.34572906,0.8138057){\color[rgb]{0,0,0}\makebox(0,0)[lt]{\lineheight{0}\smash{\begin{tabular}[t]{l}\tiny{9}\end{tabular}}}}%
    \put(0.40517341,0.35157653){\color[rgb]{0,0,0}\makebox(0,0)[lt]{\lineheight{0}\smash{\begin{tabular}[t]{l}\tiny{20}\end{tabular}}}}%
    \put(0.37525285,0.19908235){\color[rgb]{0,0,0}\makebox(0,0)[lt]{\lineheight{0}\smash{\begin{tabular}[t]{l}\tiny{24}\end{tabular}}}}%
    \put(0.07372234,0.52871972){\color[rgb]{0,0,0}\makebox(0,0)[lt]{\lineheight{0}\smash{\begin{tabular}[t]{l}\tiny{16}\end{tabular}}}}%
  \end{picture}%
\endgroup%
		\label{seq:4pseq1}
	}
	\subfloat[][sequence 2 ]{
		\def\svgwidth{30mm}
\begingroup%
  \makeatletter%
  \providecommand\color[2][]{%
    \errmessage{(Inkscape) Color is used for the text in Inkscape, but the package 'color.sty' is not loaded}%
    \renewcommand\color[2][]{}%
  }%
  \providecommand\transparent[1]{%
    \errmessage{(Inkscape) Transparency is used (non-zero) for the text in Inkscape, but the package 'transparent.sty' is not loaded}%
    \renewcommand\transparent[1]{}%
  }%
  \providecommand\rotatebox[2]{#2}%
  \ifx\svgwidth\undefined%
    \setlength{\unitlength}{342.09281158bp}%
    \ifx\svgscale\undefined%
      \relax%
    \else%
      \setlength{\unitlength}{\unitlength * \real{\svgscale}}%
    \fi%
  \else%
    \setlength{\unitlength}{\svgwidth}%
  \fi%
  \global\let\svgwidth\undefined%
  \global\let\svgscale\undefined%
  \makeatother%
  \begin{picture}(1,1.20965849)%
    \lineheight{1}%
    \setlength\tabcolsep{0pt}%
    \put(0,0){\includegraphics[width=\unitlength,page=1]{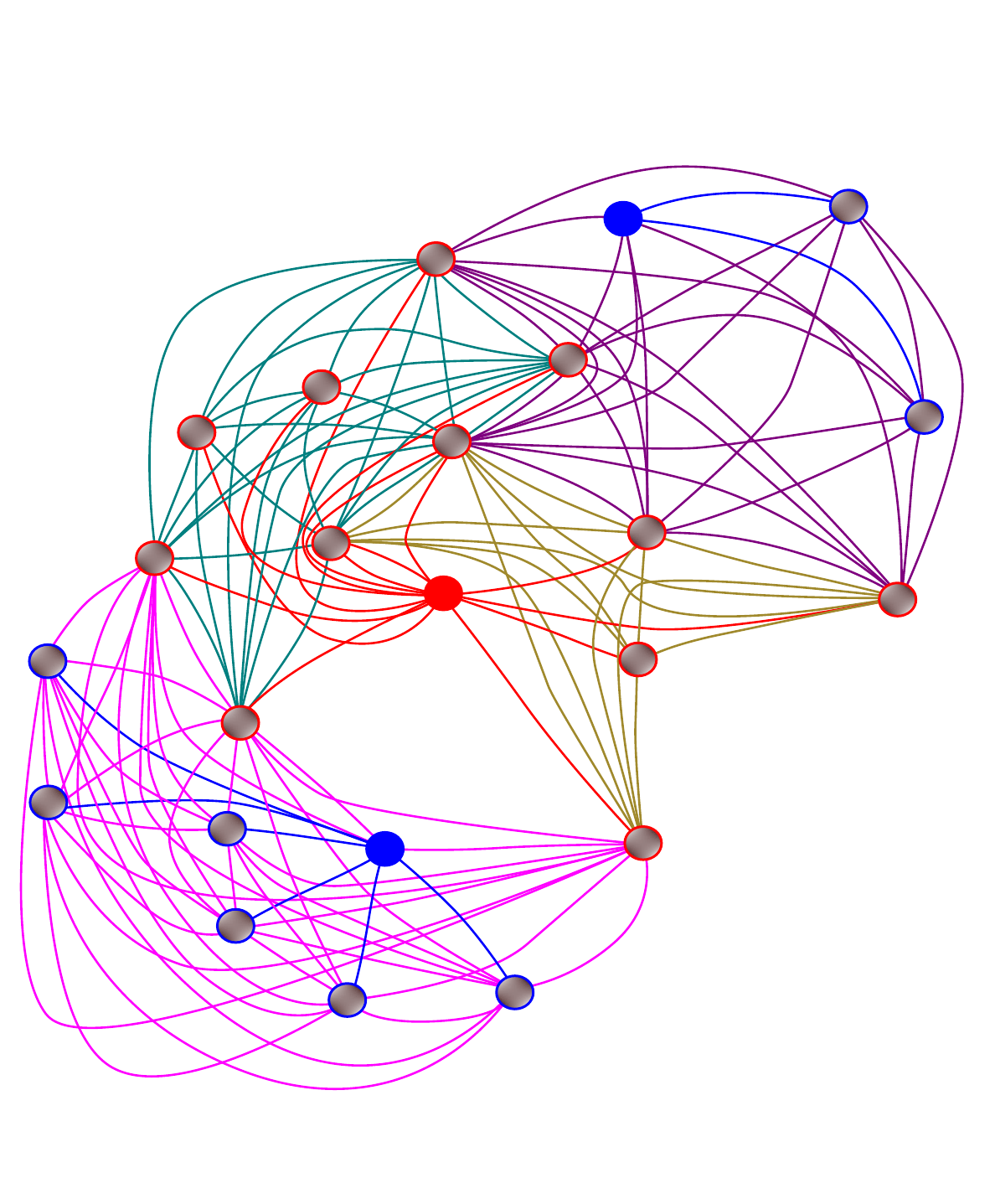}}%
    \put(0.59348497,0.8466883){\color[rgb]{0,0,0}\makebox(0,0)[lt]{\lineheight{0}\smash{\begin{tabular}[t]{l}\tiny{7}\end{tabular}}}}%
    \put(0.96309194,0.79464572){\color[rgb]{0,0,0}\makebox(0,0)[lt]{\lineheight{0}\smash{\begin{tabular}[t]{l}\tiny{4}\end{tabular}}}}%
    \put(0.65617276,0.98159756){\color[rgb]{0,0,0}\makebox(0,0)[lt]{\lineheight{0}\smash{\begin{tabular}[t]{l}\tiny{2}\end{tabular}}}}%
    \put(0.47072037,0.75909525){\color[rgb]{0,0,0}\makebox(0,0)[lt]{\lineheight{0}\smash{\begin{tabular}[t]{l}\tiny{10}\end{tabular}}}}%
    \put(0.25942799,0.27036035){\color[rgb]{0,0,0}\makebox(0,0)[lt]{\lineheight{0}\smash{\begin{tabular}[t]{l}\tiny{23}\end{tabular}}}}%
    \put(0.54620276,0.20435236){\color[rgb]{0,0,0}\makebox(0,0)[lt]{\lineheight{0}\smash{\begin{tabular}[t]{l}\tiny{25}\end{tabular}}}}%
    \put(0.88940624,1.00073606){\color[rgb]{0,0,0}\makebox(0,0)[lt]{\lineheight{0}\smash{\begin{tabular}[t]{l}\tiny{1}\end{tabular}}}}%
    \put(0.93363823,0.60079238){\color[rgb]{0,0,0}\makebox(0,0)[lt]{\lineheight{0}\smash{\begin{tabular}[t]{l}\tiny{5}\end{tabular}}}}%
    \put(0.45925715,0.94723531){\color[rgb]{0,0,0}\makebox(0,0)[lt]{\lineheight{0}\smash{\begin{tabular}[t]{l}\tiny{8}\end{tabular}}}}%
    \put(0.6633666,0.53758497){\color[rgb]{0,0,0}\makebox(0,0)[lt]{\lineheight{0}\smash{\begin{tabular}[t]{l}\tiny{11}\end{tabular}}}}%
    \put(0.47021805,0.60493473){\color[rgb]{0,0,0}\makebox(0,0)[lt]{\lineheight{0}\smash{\begin{tabular}[t]{l}\tiny{12}\end{tabular}}}}%
    \put(0.25061402,0.36555192){\color[rgb]{0,0,0}\makebox(0,0)[lt]{\lineheight{0}\smash{\begin{tabular}[t]{l}\tiny{21}\end{tabular}}}}%
    \put(0.26394086,0.4730861){\color[rgb]{0,0,0}\makebox(0,0)[lt]{\lineheight{0}\smash{\begin{tabular}[t]{l}\tiny{17}\end{tabular}}}}%
    \put(0.67006109,0.35231346){\color[rgb]{0,0,0}\makebox(0,0)[lt]{\lineheight{0}\smash{\begin{tabular}[t]{l}\tiny{19}\end{tabular}}}}%
    \put(0.35077369,0.65872394){\color[rgb]{0,0,0}\makebox(0,0)[lt]{\lineheight{0}\smash{\begin{tabular}[t]{l}\tiny{13}\end{tabular}}}}%
    \put(0.07001943,0.39062464){\color[rgb]{0,0,0}\makebox(0,0)[lt]{\lineheight{0}\smash{\begin{tabular}[t]{l}\tiny{22}\end{tabular}}}}%
    \put(0.68163635,0.66931691){\color[rgb]{0,0,0}\makebox(0,0)[lt]{\lineheight{0}\smash{\begin{tabular}[t]{l}\tiny{6}\end{tabular}}}}%
    \put(0.21926386,0.7677135){\color[rgb]{0,0,0}\makebox(0,0)[lt]{\lineheight{0}\smash{\begin{tabular}[t]{l}\tiny{14}\end{tabular}}}}%
    \put(0.17118822,0.63719051){\color[rgb]{0,0,0}\makebox(0,0)[lt]{\lineheight{0}\smash{\begin{tabular}[t]{l}\tiny{15}\end{tabular}}}}%
    \put(0.34572906,0.8138057){\color[rgb]{0,0,0}\makebox(0,0)[lt]{\lineheight{0}\smash{\begin{tabular}[t]{l}\tiny{9}\end{tabular}}}}%
    \put(0.40517341,0.35157653){\color[rgb]{0,0,0}\makebox(0,0)[lt]{\lineheight{0}\smash{\begin{tabular}[t]{l}\tiny{20}\end{tabular}}}}%
    \put(0.37525285,0.19908235){\color[rgb]{0,0,0}\makebox(0,0)[lt]{\lineheight{0}\smash{\begin{tabular}[t]{l}\tiny{24}\end{tabular}}}}%
    \put(0.07372234,0.52871972){\color[rgb]{0,0,0}\makebox(0,0)[lt]{\lineheight{0}\smash{\begin{tabular}[t]{l}\tiny{16}\end{tabular}}}}%
  \end{picture}%
\endgroup%
		\label{seq:4pseq2}
	}
	\subfloat[][sequence 3]{
		\def\svgwidth{30mm}
\begingroup%
  \makeatletter%
  \providecommand\color[2][]{%
    \errmessage{(Inkscape) Color is used for the text in Inkscape, but the package 'color.sty' is not loaded}%
    \renewcommand\color[2][]{}%
  }%
  \providecommand\transparent[1]{%
    \errmessage{(Inkscape) Transparency is used (non-zero) for the text in Inkscape, but the package 'transparent.sty' is not loaded}%
    \renewcommand\transparent[1]{}%
  }%
  \providecommand\rotatebox[2]{#2}%
  \ifx\svgwidth\undefined%
    \setlength{\unitlength}{342.09281158bp}%
    \ifx\svgscale\undefined%
      \relax%
    \else%
      \setlength{\unitlength}{\unitlength * \real{\svgscale}}%
    \fi%
  \else%
    \setlength{\unitlength}{\svgwidth}%
  \fi%
  \global\let\svgwidth\undefined%
  \global\let\svgscale\undefined%
  \makeatother%
  \begin{picture}(1,1.20965849)%
    \lineheight{1}%
    \setlength\tabcolsep{0pt}%
    \put(0,0){\includegraphics[width=\unitlength,page=1]{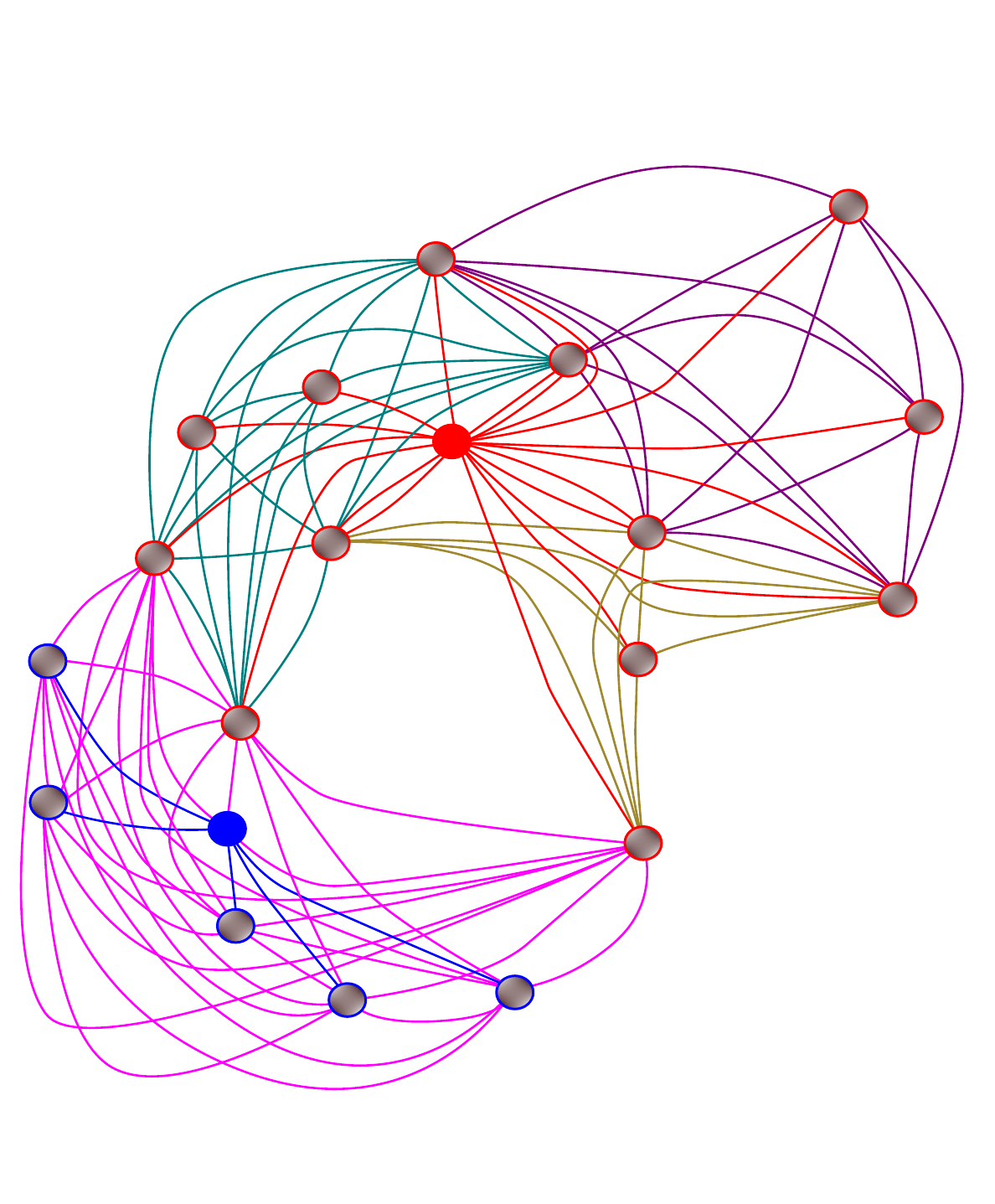}}%
    \put(0.59348497,0.8466883){\color[rgb]{0,0,0}\makebox(0,0)[lt]{\lineheight{0}\smash{\begin{tabular}[t]{l}\tiny{7}\end{tabular}}}}%
    \put(0.96309194,0.79464572){\color[rgb]{0,0,0}\makebox(0,0)[lt]{\lineheight{0}\smash{\begin{tabular}[t]{l}\tiny{4}\end{tabular}}}}%
    \put(0.47072037,0.75909525){\color[rgb]{0,0,0}\makebox(0,0)[lt]{\lineheight{0}\smash{\begin{tabular}[t]{l}\tiny{10}\end{tabular}}}}%
    \put(0.25942799,0.27036035){\color[rgb]{0,0,0}\makebox(0,0)[lt]{\lineheight{0}\smash{\begin{tabular}[t]{l}\tiny{23}\end{tabular}}}}%
    \put(0.54620276,0.20435236){\color[rgb]{0,0,0}\makebox(0,0)[lt]{\lineheight{0}\smash{\begin{tabular}[t]{l}\tiny{25}\end{tabular}}}}%
    \put(0.88940624,1.00073606){\color[rgb]{0,0,0}\makebox(0,0)[lt]{\lineheight{0}\smash{\begin{tabular}[t]{l}\tiny{1}\end{tabular}}}}%
    \put(0.93363823,0.60079238){\color[rgb]{0,0,0}\makebox(0,0)[lt]{\lineheight{0}\smash{\begin{tabular}[t]{l}\tiny{5}\end{tabular}}}}%
    \put(0.45925715,0.94723531){\color[rgb]{0,0,0}\makebox(0,0)[lt]{\lineheight{0}\smash{\begin{tabular}[t]{l}\tiny{8}\end{tabular}}}}%
    \put(0.6633666,0.53758497){\color[rgb]{0,0,0}\makebox(0,0)[lt]{\lineheight{0}\smash{\begin{tabular}[t]{l}\tiny{11}\end{tabular}}}}%
    \put(0.25061402,0.36555192){\color[rgb]{0,0,0}\makebox(0,0)[lt]{\lineheight{0}\smash{\begin{tabular}[t]{l}\tiny{21}\end{tabular}}}}%
    \put(0.26394086,0.4730861){\color[rgb]{0,0,0}\makebox(0,0)[lt]{\lineheight{0}\smash{\begin{tabular}[t]{l}\tiny{17}\end{tabular}}}}%
    \put(0.67006109,0.35231346){\color[rgb]{0,0,0}\makebox(0,0)[lt]{\lineheight{0}\smash{\begin{tabular}[t]{l}\tiny{19}\end{tabular}}}}%
    \put(0.35077369,0.65872394){\color[rgb]{0,0,0}\makebox(0,0)[lt]{\lineheight{0}\smash{\begin{tabular}[t]{l}\tiny{13}\end{tabular}}}}%
    \put(0.07001943,0.39062464){\color[rgb]{0,0,0}\makebox(0,0)[lt]{\lineheight{0}\smash{\begin{tabular}[t]{l}\tiny{22}\end{tabular}}}}%
    \put(0.68163635,0.66931691){\color[rgb]{0,0,0}\makebox(0,0)[lt]{\lineheight{0}\smash{\begin{tabular}[t]{l}\tiny{6}\end{tabular}}}}%
    \put(0.21926386,0.7677135){\color[rgb]{0,0,0}\makebox(0,0)[lt]{\lineheight{0}\smash{\begin{tabular}[t]{l}\tiny{14}\end{tabular}}}}%
    \put(0.17118822,0.63719051){\color[rgb]{0,0,0}\makebox(0,0)[lt]{\lineheight{0}\smash{\begin{tabular}[t]{l}\tiny{15}\end{tabular}}}}%
    \put(0.34572906,0.8138057){\color[rgb]{0,0,0}\makebox(0,0)[lt]{\lineheight{0}\smash{\begin{tabular}[t]{l}\tiny{9}\end{tabular}}}}%
    \put(0.37525285,0.19908235){\color[rgb]{0,0,0}\makebox(0,0)[lt]{\lineheight{0}\smash{\begin{tabular}[t]{l}\tiny{24}\end{tabular}}}}%
    \put(0.07372234,0.52871972){\color[rgb]{0,0,0}\makebox(0,0)[lt]{\lineheight{0}\smash{\begin{tabular}[t]{l}\tiny{16}\end{tabular}}}}%
  \end{picture}%
\endgroup%
		\label{seq:4pseq3}
	}\\
	\subfloat[][sequence 4]{
		\def\svgwidth{30mm}
\begingroup%
  \makeatletter%
  \providecommand\color[2][]{%
    \errmessage{(Inkscape) Color is used for the text in Inkscape, but the package 'color.sty' is not loaded}%
    \renewcommand\color[2][]{}%
  }%
  \providecommand\transparent[1]{%
    \errmessage{(Inkscape) Transparency is used (non-zero) for the text in Inkscape, but the package 'transparent.sty' is not loaded}%
    \renewcommand\transparent[1]{}%
  }%
  \providecommand\rotatebox[2]{#2}%
  \ifx\svgwidth\undefined%
    \setlength{\unitlength}{342.09281158bp}%
    \ifx\svgscale\undefined%
      \relax%
    \else%
      \setlength{\unitlength}{\unitlength * \real{\svgscale}}%
    \fi%
  \else%
    \setlength{\unitlength}{\svgwidth}%
  \fi%
  \global\let\svgwidth\undefined%
  \global\let\svgscale\undefined%
  \makeatother%
  \begin{picture}(1,1.20965849)%
    \lineheight{1}%
    \setlength\tabcolsep{0pt}%
    \put(0,0){\includegraphics[width=\unitlength,page=1]{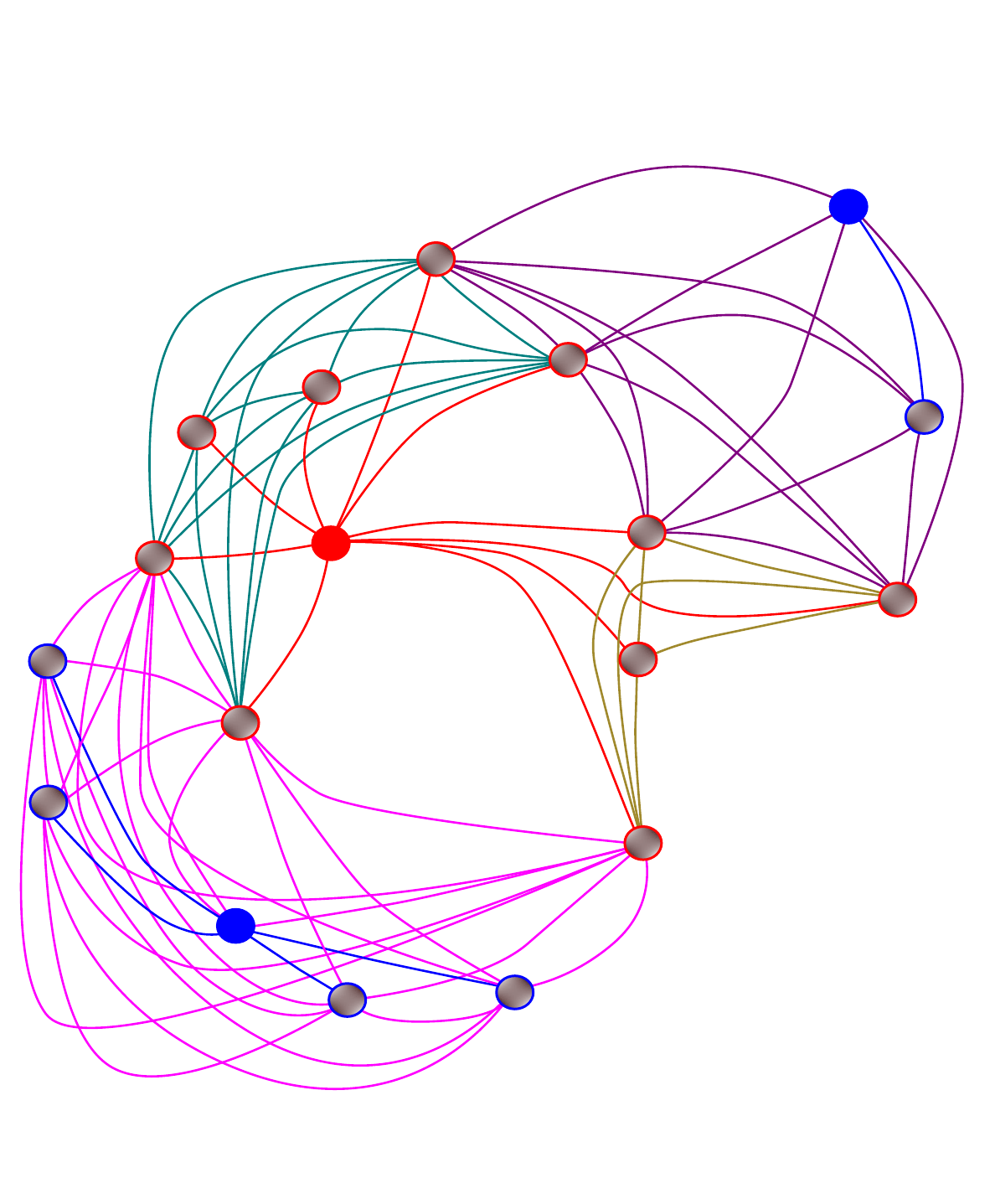}}%
    \put(0.59348497,0.8466883){\color[rgb]{0,0,0}\makebox(0,0)[lt]{\lineheight{0}\smash{\begin{tabular}[t]{l}\tiny{7}\end{tabular}}}}%
    \put(0.96309194,0.79464572){\color[rgb]{0,0,0}\makebox(0,0)[lt]{\lineheight{0}\smash{\begin{tabular}[t]{l}\tiny{4}\end{tabular}}}}%
    \put(0.25942799,0.27036035){\color[rgb]{0,0,0}\makebox(0,0)[lt]{\lineheight{0}\smash{\begin{tabular}[t]{l}\tiny{23}\end{tabular}}}}%
    \put(0.54620276,0.20435236){\color[rgb]{0,0,0}\makebox(0,0)[lt]{\lineheight{0}\smash{\begin{tabular}[t]{l}\tiny{25}\end{tabular}}}}%
    \put(0.88940624,1.00073606){\color[rgb]{0,0,0}\makebox(0,0)[lt]{\lineheight{0}\smash{\begin{tabular}[t]{l}\tiny{1}\end{tabular}}}}%
    \put(0.93363823,0.60079238){\color[rgb]{0,0,0}\makebox(0,0)[lt]{\lineheight{0}\smash{\begin{tabular}[t]{l}\tiny{5}\end{tabular}}}}%
    \put(0.45925715,0.94723531){\color[rgb]{0,0,0}\makebox(0,0)[lt]{\lineheight{0}\smash{\begin{tabular}[t]{l}\tiny{8}\end{tabular}}}}%
    \put(0.6633666,0.53758497){\color[rgb]{0,0,0}\makebox(0,0)[lt]{\lineheight{0}\smash{\begin{tabular}[t]{l}\tiny{11}\end{tabular}}}}%
    \put(0.26394086,0.4730861){\color[rgb]{0,0,0}\makebox(0,0)[lt]{\lineheight{0}\smash{\begin{tabular}[t]{l}\tiny{17}\end{tabular}}}}%
    \put(0.67006109,0.35231346){\color[rgb]{0,0,0}\makebox(0,0)[lt]{\lineheight{0}\smash{\begin{tabular}[t]{l}\tiny{19}\end{tabular}}}}%
    \put(0.35077369,0.65872394){\color[rgb]{0,0,0}\makebox(0,0)[lt]{\lineheight{0}\smash{\begin{tabular}[t]{l}\tiny{13}\end{tabular}}}}%
    \put(0.07001943,0.39062464){\color[rgb]{0,0,0}\makebox(0,0)[lt]{\lineheight{0}\smash{\begin{tabular}[t]{l}\tiny{22}\end{tabular}}}}%
    \put(0.68163635,0.66931691){\color[rgb]{0,0,0}\makebox(0,0)[lt]{\lineheight{0}\smash{\begin{tabular}[t]{l}\tiny{6}\end{tabular}}}}%
    \put(0.21926386,0.7677135){\color[rgb]{0,0,0}\makebox(0,0)[lt]{\lineheight{0}\smash{\begin{tabular}[t]{l}\tiny{14}\end{tabular}}}}%
    \put(0.17118822,0.63719051){\color[rgb]{0,0,0}\makebox(0,0)[lt]{\lineheight{0}\smash{\begin{tabular}[t]{l}\tiny{15}\end{tabular}}}}%
    \put(0.34572906,0.8138057){\color[rgb]{0,0,0}\makebox(0,0)[lt]{\lineheight{0}\smash{\begin{tabular}[t]{l}\tiny{9}\end{tabular}}}}%
    \put(0.37525285,0.19908235){\color[rgb]{0,0,0}\makebox(0,0)[lt]{\lineheight{0}\smash{\begin{tabular}[t]{l}\tiny{24}\end{tabular}}}}%
    \put(0.07372234,0.52871972){\color[rgb]{0,0,0}\makebox(0,0)[lt]{\lineheight{0}\smash{\begin{tabular}[t]{l}\tiny{16}\end{tabular}}}}%
  \end{picture}%
\endgroup%
		\label{seq:4pseq4}		}
	\subfloat[][Sequence 5]{
		\def\svgwidth{30mm}
\begingroup%
  \makeatletter%
  \providecommand\color[2][]{%
    \errmessage{(Inkscape) Color is used for the text in Inkscape, but the package 'color.sty' is not loaded}%
    \renewcommand\color[2][]{}%
  }%
  \providecommand\transparent[1]{%
    \errmessage{(Inkscape) Transparency is used (non-zero) for the text in Inkscape, but the package 'transparent.sty' is not loaded}%
    \renewcommand\transparent[1]{}%
  }%
  \providecommand\rotatebox[2]{#2}%
  \ifx\svgwidth\undefined%
    \setlength{\unitlength}{342.09281158bp}%
    \ifx\svgscale\undefined%
      \relax%
    \else%
      \setlength{\unitlength}{\unitlength * \real{\svgscale}}%
    \fi%
  \else%
    \setlength{\unitlength}{\svgwidth}%
  \fi%
  \global\let\svgwidth\undefined%
  \global\let\svgscale\undefined%
  \makeatother%
  \begin{picture}(1,1.20965849)%
    \lineheight{1}%
    \setlength\tabcolsep{0pt}%
    \put(0,0){\includegraphics[width=\unitlength,page=1]{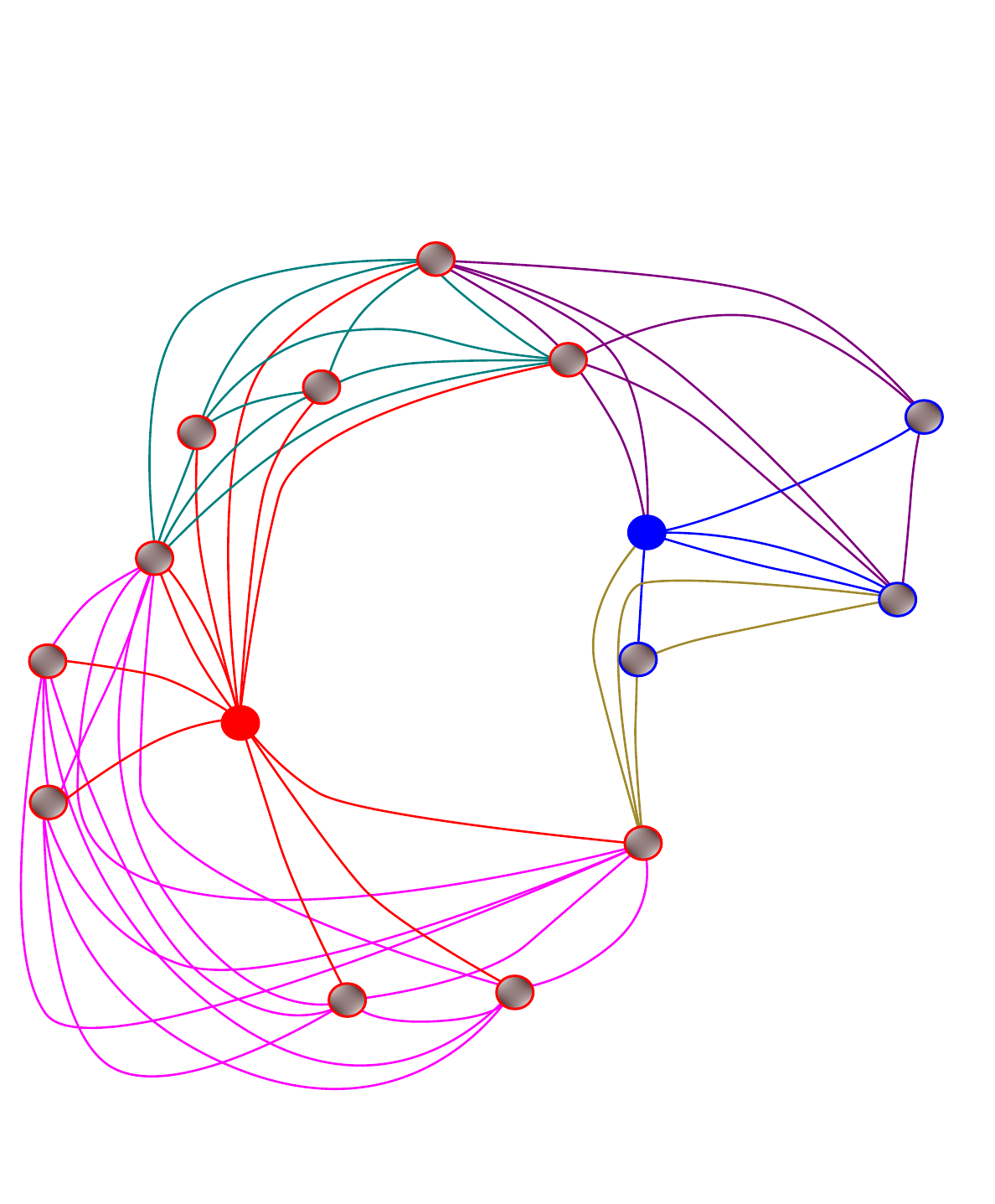}}%
    \put(0.59348497,0.8466883){\color[rgb]{0,0,0}\makebox(0,0)[lt]{\lineheight{0}\smash{\begin{tabular}[t]{l}\tiny{7}\end{tabular}}}}%
    \put(0.96309194,0.79464572){\color[rgb]{0,0,0}\makebox(0,0)[lt]{\lineheight{0}\smash{\begin{tabular}[t]{l}\tiny{4}\end{tabular}}}}%
    \put(0.54620276,0.20435236){\color[rgb]{0,0,0}\makebox(0,0)[lt]{\lineheight{0}\smash{\begin{tabular}[t]{l}\tiny{25}\end{tabular}}}}%
    \put(0.93363823,0.60079238){\color[rgb]{0,0,0}\makebox(0,0)[lt]{\lineheight{0}\smash{\begin{tabular}[t]{l}\tiny{5}\end{tabular}}}}%
    \put(0.45925715,0.94723531){\color[rgb]{0,0,0}\makebox(0,0)[lt]{\lineheight{0}\smash{\begin{tabular}[t]{l}\tiny{8}\end{tabular}}}}%
    \put(0.6633666,0.53758497){\color[rgb]{0,0,0}\makebox(0,0)[lt]{\lineheight{0}\smash{\begin{tabular}[t]{l}\tiny{11}\end{tabular}}}}%
    \put(0.26394086,0.4730861){\color[rgb]{0,0,0}\makebox(0,0)[lt]{\lineheight{0}\smash{\begin{tabular}[t]{l}\tiny{17}\end{tabular}}}}%
    \put(0.67006109,0.35231346){\color[rgb]{0,0,0}\makebox(0,0)[lt]{\lineheight{0}\smash{\begin{tabular}[t]{l}\tiny{19}\end{tabular}}}}%
    \put(0.07001943,0.39062464){\color[rgb]{0,0,0}\makebox(0,0)[lt]{\lineheight{0}\smash{\begin{tabular}[t]{l}\tiny{22}\end{tabular}}}}%
    \put(0.68163635,0.66931691){\color[rgb]{0,0,0}\makebox(0,0)[lt]{\lineheight{0}\smash{\begin{tabular}[t]{l}\tiny{6}\end{tabular}}}}%
    \put(0.21926386,0.7677135){\color[rgb]{0,0,0}\makebox(0,0)[lt]{\lineheight{0}\smash{\begin{tabular}[t]{l}\tiny{14}\end{tabular}}}}%
    \put(0.17118822,0.63719051){\color[rgb]{0,0,0}\makebox(0,0)[lt]{\lineheight{0}\smash{\begin{tabular}[t]{l}\tiny{15}\end{tabular}}}}%
    \put(0.34572906,0.8138057){\color[rgb]{0,0,0}\makebox(0,0)[lt]{\lineheight{0}\smash{\begin{tabular}[t]{l}\tiny{9}\end{tabular}}}}%
    \put(0.37525285,0.19908235){\color[rgb]{0,0,0}\makebox(0,0)[lt]{\lineheight{0}\smash{\begin{tabular}[t]{l}\tiny{24}\end{tabular}}}}%
    \put(0.07372234,0.52871972){\color[rgb]{0,0,0}\makebox(0,0)[lt]{\lineheight{0}\smash{\begin{tabular}[t]{l}\tiny{16}\end{tabular}}}}%
  \end{picture}%
\endgroup%
		\label{seq:4pseq5}		}
	\subfloat[][Sequence 6]{
		\def\svgwidth{30mm}
\begingroup%
  \makeatletter%
  \providecommand\color[2][]{%
    \errmessage{(Inkscape) Color is used for the text in Inkscape, but the package 'color.sty' is not loaded}%
    \renewcommand\color[2][]{}%
  }%
  \providecommand\transparent[1]{%
    \errmessage{(Inkscape) Transparency is used (non-zero) for the text in Inkscape, but the package 'transparent.sty' is not loaded}%
    \renewcommand\transparent[1]{}%
  }%
  \providecommand\rotatebox[2]{#2}%
  \ifx\svgwidth\undefined%
    \setlength{\unitlength}{342.09281158bp}%
    \ifx\svgscale\undefined%
      \relax%
    \else%
      \setlength{\unitlength}{\unitlength * \real{\svgscale}}%
    \fi%
  \else%
    \setlength{\unitlength}{\svgwidth}%
  \fi%
  \global\let\svgwidth\undefined%
  \global\let\svgscale\undefined%
  \makeatother%
  \begin{picture}(1,1.20965849)%
    \lineheight{1}%
    \setlength\tabcolsep{0pt}%
    \put(0,0){\includegraphics[width=\unitlength,page=1]{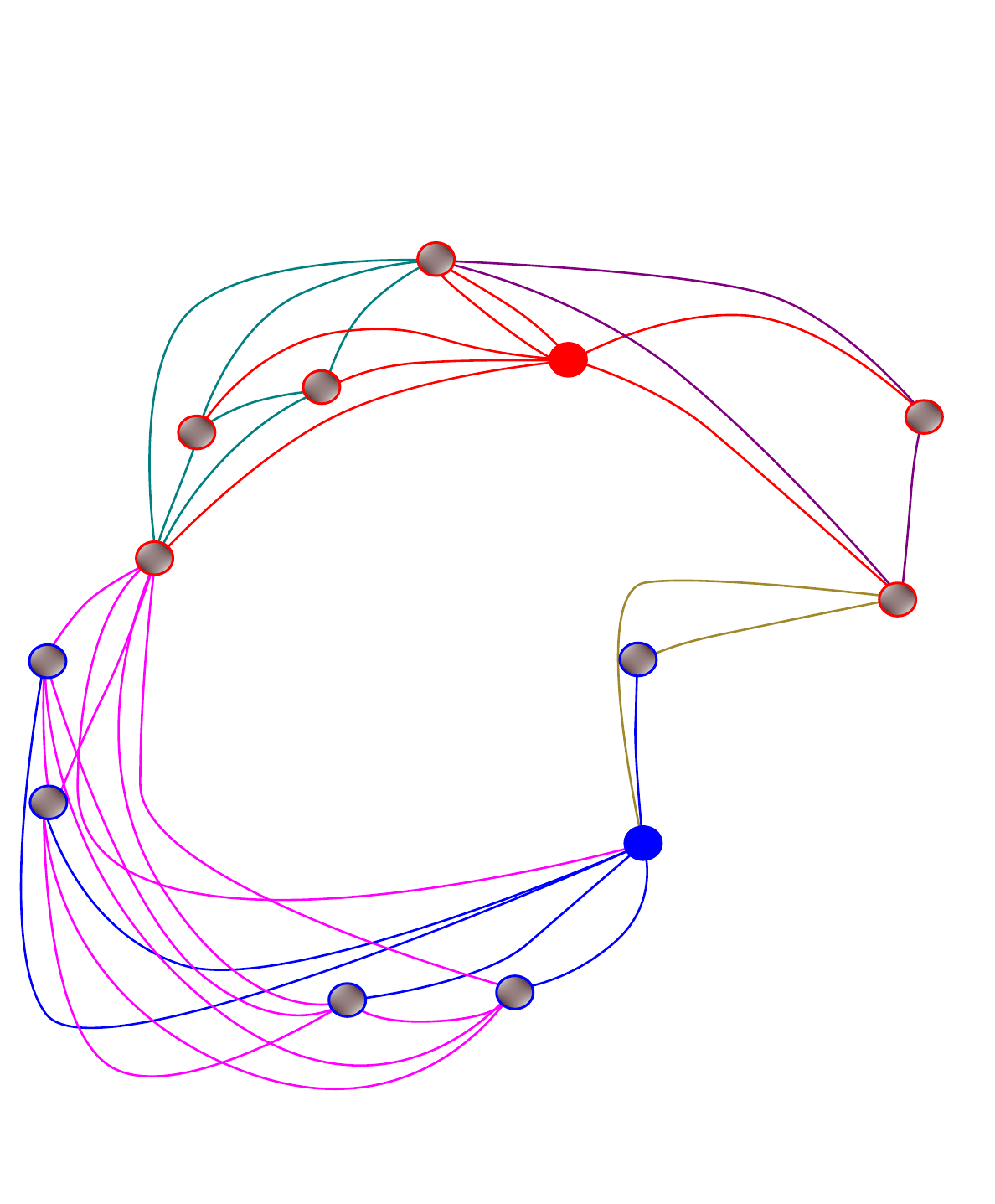}}%
    \put(0.59348497,0.8466883){\color[rgb]{0,0,0}\makebox(0,0)[lt]{\lineheight{0}\smash{\begin{tabular}[t]{l}\tiny{7}\end{tabular}}}}%
    \put(0.96309194,0.79464572){\color[rgb]{0,0,0}\makebox(0,0)[lt]{\lineheight{0}\smash{\begin{tabular}[t]{l}\tiny{4}\end{tabular}}}}%
    \put(0.54620276,0.20435236){\color[rgb]{0,0,0}\makebox(0,0)[lt]{\lineheight{0}\smash{\begin{tabular}[t]{l}\tiny{25}\end{tabular}}}}%
    \put(0.93363823,0.60079238){\color[rgb]{0,0,0}\makebox(0,0)[lt]{\lineheight{0}\smash{\begin{tabular}[t]{l}\tiny{5}\end{tabular}}}}%
    \put(0.45925715,0.94723531){\color[rgb]{0,0,0}\makebox(0,0)[lt]{\lineheight{0}\smash{\begin{tabular}[t]{l}\tiny{8}\end{tabular}}}}%
    \put(0.6633666,0.53758497){\color[rgb]{0,0,0}\makebox(0,0)[lt]{\lineheight{0}\smash{\begin{tabular}[t]{l}\tiny{11}\end{tabular}}}}%
    \put(0.67006109,0.35231346){\color[rgb]{0,0,0}\makebox(0,0)[lt]{\lineheight{0}\smash{\begin{tabular}[t]{l}\tiny{19}\end{tabular}}}}%
    \put(0.07001943,0.39062464){\color[rgb]{0,0,0}\makebox(0,0)[lt]{\lineheight{0}\smash{\begin{tabular}[t]{l}\tiny{22}\end{tabular}}}}%
    \put(0.21926386,0.7677135){\color[rgb]{0,0,0}\makebox(0,0)[lt]{\lineheight{0}\smash{\begin{tabular}[t]{l}\tiny{14}\end{tabular}}}}%
    \put(0.17118822,0.63719051){\color[rgb]{0,0,0}\makebox(0,0)[lt]{\lineheight{0}\smash{\begin{tabular}[t]{l}\tiny{15}\end{tabular}}}}%
    \put(0.34572906,0.8138057){\color[rgb]{0,0,0}\makebox(0,0)[lt]{\lineheight{0}\smash{\begin{tabular}[t]{l}\tiny{9}\end{tabular}}}}%
    \put(0.37525285,0.19908235){\color[rgb]{0,0,0}\makebox(0,0)[lt]{\lineheight{0}\smash{\begin{tabular}[t]{l}\tiny{24}\end{tabular}}}}%
    \put(0.07372234,0.52871972){\color[rgb]{0,0,0}\makebox(0,0)[lt]{\lineheight{0}\smash{\begin{tabular}[t]{l}\tiny{16}\end{tabular}}}}%
  \end{picture}%
\endgroup%
		\label{seq:4pseq6}		}
	\subfloat[][Sequence 7]{
		\def\svgwidth{30mm}
\begingroup%
  \makeatletter%
  \providecommand\color[2][]{%
    \errmessage{(Inkscape) Color is used for the text in Inkscape, but the package 'color.sty' is not loaded}%
    \renewcommand\color[2][]{}%
  }%
  \providecommand\transparent[1]{%
    \errmessage{(Inkscape) Transparency is used (non-zero) for the text in Inkscape, but the package 'transparent.sty' is not loaded}%
    \renewcommand\transparent[1]{}%
  }%
  \providecommand\rotatebox[2]{#2}%
  \ifx\svgwidth\undefined%
    \setlength{\unitlength}{342.09281158bp}%
    \ifx\svgscale\undefined%
      \relax%
    \else%
      \setlength{\unitlength}{\unitlength * \real{\svgscale}}%
    \fi%
  \else%
    \setlength{\unitlength}{\svgwidth}%
  \fi%
  \global\let\svgwidth\undefined%
  \global\let\svgscale\undefined%
  \makeatother%
  \begin{picture}(1,1.20965849)%
    \lineheight{1}%
    \setlength\tabcolsep{0pt}%
    \put(0,0){\includegraphics[width=\unitlength,page=1]{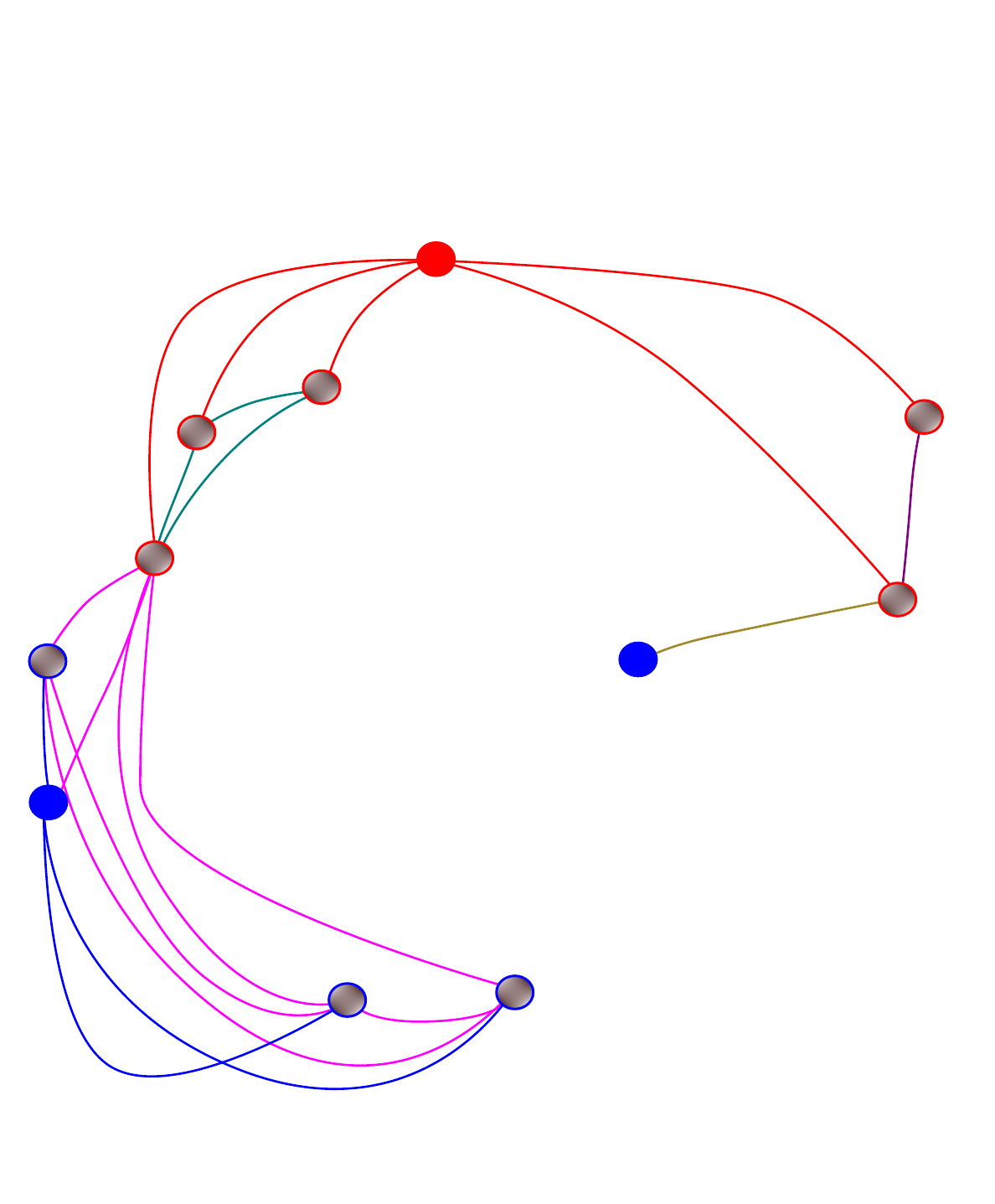}}%
    \put(0.96309194,0.79464572){\color[rgb]{0,0,0}\makebox(0,0)[lt]{\lineheight{0}\smash{\begin{tabular}[t]{l}\tiny{4}\end{tabular}}}}%
    \put(0.54620276,0.20435236){\color[rgb]{0,0,0}\makebox(0,0)[lt]{\lineheight{0}\smash{\begin{tabular}[t]{l}\tiny{25}\end{tabular}}}}%
    \put(0.93363823,0.60079238){\color[rgb]{0,0,0}\makebox(0,0)[lt]{\lineheight{0}\smash{\begin{tabular}[t]{l}\tiny{5}\end{tabular}}}}%
    \put(0.45925715,0.94723531){\color[rgb]{0,0,0}\makebox(0,0)[lt]{\lineheight{0}\smash{\begin{tabular}[t]{l}\tiny{8}\end{tabular}}}}%
    \put(0.6633666,0.53758497){\color[rgb]{0,0,0}\makebox(0,0)[lt]{\lineheight{0}\smash{\begin{tabular}[t]{l}\tiny{11}\end{tabular}}}}%
    \put(0.07001943,0.39062464){\color[rgb]{0,0,0}\makebox(0,0)[lt]{\lineheight{0}\smash{\begin{tabular}[t]{l}\tiny{22}\end{tabular}}}}%
    \put(0.21926386,0.7677135){\color[rgb]{0,0,0}\makebox(0,0)[lt]{\lineheight{0}\smash{\begin{tabular}[t]{l}\tiny{14}\end{tabular}}}}%
    \put(0.17118822,0.63719051){\color[rgb]{0,0,0}\makebox(0,0)[lt]{\lineheight{0}\smash{\begin{tabular}[t]{l}\tiny{15}\end{tabular}}}}%
    \put(0.34572906,0.8138057){\color[rgb]{0,0,0}\makebox(0,0)[lt]{\lineheight{0}\smash{\begin{tabular}[t]{l}\tiny{9}\end{tabular}}}}%
    \put(0.37525285,0.19908235){\color[rgb]{0,0,0}\makebox(0,0)[lt]{\lineheight{0}\smash{\begin{tabular}[t]{l}\tiny{24}\end{tabular}}}}%
    \put(0.07372234,0.52871972){\color[rgb]{0,0,0}\makebox(0,0)[lt]{\lineheight{0}\smash{\begin{tabular}[t]{l}\tiny{16}\end{tabular}}}}%
  \end{picture}%
\endgroup%
		\label{seq:4pseq7}		}
	\subfloat[][Sequence 8]{
		\def\svgwidth{30mm}
\begingroup%
  \makeatletter%
  \providecommand\color[2][]{%
    \errmessage{(Inkscape) Color is used for the text in Inkscape, but the package 'color.sty' is not loaded}%
    \renewcommand\color[2][]{}%
  }%
  \providecommand\transparent[1]{%
    \errmessage{(Inkscape) Transparency is used (non-zero) for the text in Inkscape, but the package 'transparent.sty' is not loaded}%
    \renewcommand\transparent[1]{}%
  }%
  \providecommand\rotatebox[2]{#2}%
  \ifx\svgwidth\undefined%
    \setlength{\unitlength}{342.09281158bp}%
    \ifx\svgscale\undefined%
      \relax%
    \else%
      \setlength{\unitlength}{\unitlength * \real{\svgscale}}%
    \fi%
  \else%
    \setlength{\unitlength}{\svgwidth}%
  \fi%
  \global\let\svgwidth\undefined%
  \global\let\svgscale\undefined%
  \makeatother%
  \begin{picture}(1,1.20965849)%
    \lineheight{1}%
    \setlength\tabcolsep{0pt}%
    \put(0,0){\includegraphics[width=\unitlength,page=1]{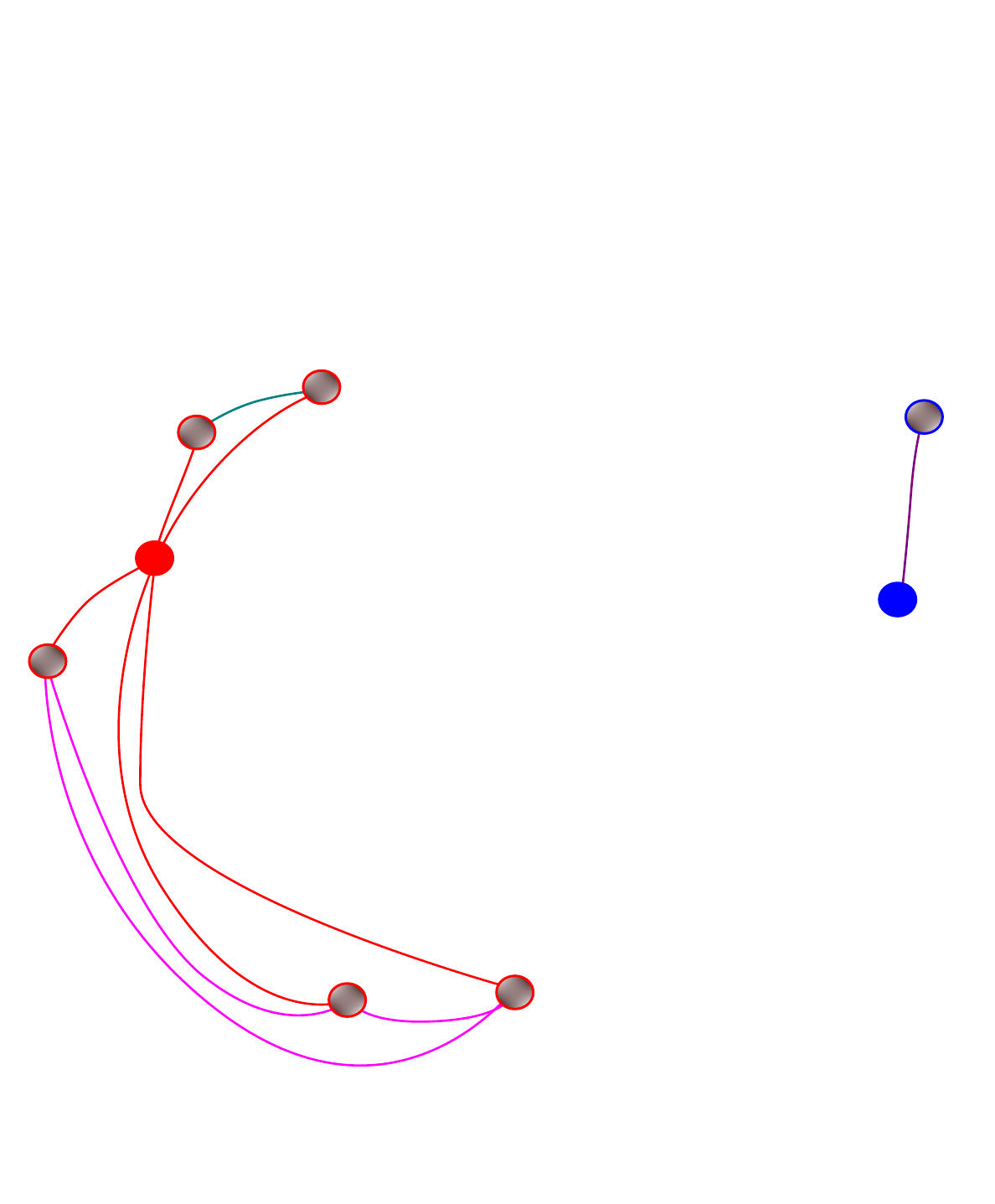}}%
    \put(0.96309194,0.79464572){\color[rgb]{0,0,0}\makebox(0,0)[lt]{\lineheight{0}\smash{\begin{tabular}[t]{l}\tiny{4}\end{tabular}}}}%
    \put(0.54620276,0.20435236){\color[rgb]{0,0,0}\makebox(0,0)[lt]{\lineheight{0}\smash{\begin{tabular}[t]{l}\tiny{25}\end{tabular}}}}%
    \put(0.93363823,0.60079238){\color[rgb]{0,0,0}\makebox(0,0)[lt]{\lineheight{0}\smash{\begin{tabular}[t]{l}\tiny{5}\end{tabular}}}}%
    \put(0.21926386,0.7677135){\color[rgb]{0,0,0}\makebox(0,0)[lt]{\lineheight{0}\smash{\begin{tabular}[t]{l}\tiny{14}\end{tabular}}}}%
    \put(0.17118822,0.63719051){\color[rgb]{0,0,0}\makebox(0,0)[lt]{\lineheight{0}\smash{\begin{tabular}[t]{l}\tiny{15}\end{tabular}}}}%
    \put(0.34572906,0.8138057){\color[rgb]{0,0,0}\makebox(0,0)[lt]{\lineheight{0}\smash{\begin{tabular}[t]{l}\tiny{9}\end{tabular}}}}%
    \put(0.37525285,0.19908235){\color[rgb]{0,0,0}\makebox(0,0)[lt]{\lineheight{0}\smash{\begin{tabular}[t]{l}\tiny{24}\end{tabular}}}}%
    \put(0.07372234,0.52871972){\color[rgb]{0,0,0}\makebox(0,0)[lt]{\lineheight{0}\smash{\begin{tabular}[t]{l}\tiny{16}\end{tabular}}}}%
  \end{picture}%
\endgroup%
		\label{seq:4pseq8}		}\\
	\subfloat[][Remaining local patch]{
		\def\svgwidth{30mm}
\begingroup%
  \makeatletter%
  \providecommand\color[2][]{%
    \errmessage{(Inkscape) Color is used for the text in Inkscape, but the package 'color.sty' is not loaded}%
    \renewcommand\color[2][]{}%
  }%
  \providecommand\transparent[1]{%
    \errmessage{(Inkscape) Transparency is used (non-zero) for the text in Inkscape, but the package 'transparent.sty' is not loaded}%
    \renewcommand\transparent[1]{}%
  }%
  \providecommand\rotatebox[2]{#2}%
  \ifx\svgwidth\undefined%
    \setlength{\unitlength}{342.09281158bp}%
    \ifx\svgscale\undefined%
      \relax%
    \else%
      \setlength{\unitlength}{\unitlength * \real{\svgscale}}%
    \fi%
  \else%
    \setlength{\unitlength}{\svgwidth}%
  \fi%
  \global\let\svgwidth\undefined%
  \global\let\svgscale\undefined%
  \makeatother%
  \begin{picture}(1,1.20965849)%
    \lineheight{1}%
    \setlength\tabcolsep{0pt}%
    \put(0,0){\includegraphics[width=\unitlength,page=1]{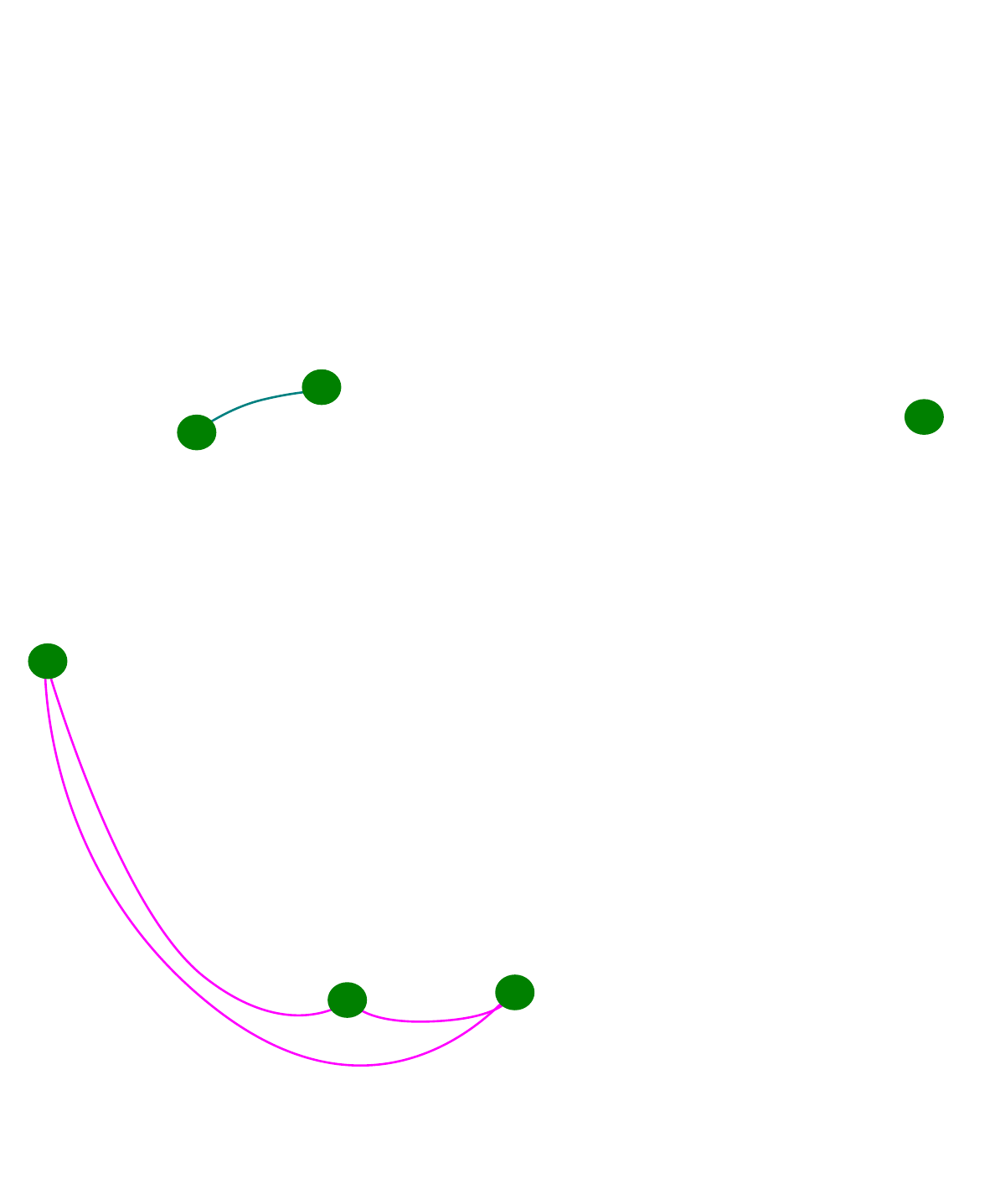}}%
    \put(0.96309194,0.79464572){\color[rgb]{0,0,0}\makebox(0,0)[lt]{\lineheight{0}\smash{\begin{tabular}[t]{l}\tiny{4}\end{tabular}}}}%
    \put(0.54620276,0.20435236){\color[rgb]{0,0,0}\makebox(0,0)[lt]{\lineheight{0}\smash{\begin{tabular}[t]{l}\tiny{25}\end{tabular}}}}%
    \put(0.21926386,0.7677135){\color[rgb]{0,0,0}\makebox(0,0)[lt]{\lineheight{0}\smash{\begin{tabular}[t]{l}\tiny{14}\end{tabular}}}}%
    \put(0.34572906,0.8138057){\color[rgb]{0,0,0}\makebox(0,0)[lt]{\lineheight{0}\smash{\begin{tabular}[t]{l}\tiny{9}\end{tabular}}}}%
    \put(0.37525285,0.19908235){\color[rgb]{0,0,0}\makebox(0,0)[lt]{\lineheight{0}\smash{\begin{tabular}[t]{l}\tiny{24}\end{tabular}}}}%
    \put(0.07372234,0.52871972){\color[rgb]{0,0,0}\makebox(0,0)[lt]{\lineheight{0}\smash{\begin{tabular}[t]{l}\tiny{16}\end{tabular}}}}%
  \end{picture}%
\endgroup%
		\label{seq:4pseq9}		}
	\subfloat[][Legend]{
		\def\svgwidth{30mm}
\begingroup%
  \makeatletter%
  \providecommand\color[2][]{%
    \errmessage{(Inkscape) Color is used for the text in Inkscape, but the package 'color.sty' is not loaded}%
    \renewcommand\color[2][]{}%
  }%
  \providecommand\transparent[1]{%
    \errmessage{(Inkscape) Transparency is used (non-zero) for the text in Inkscape, but the package 'transparent.sty' is not loaded}%
    \renewcommand\transparent[1]{}%
  }%
  \providecommand\rotatebox[2]{#2}%
  \ifx\svgwidth\undefined%
    \setlength{\unitlength}{321.66164485bp}%
    \ifx\svgscale\undefined%
      \relax%
    \else%
      \setlength{\unitlength}{\unitlength * \real{\svgscale}}%
    \fi%
  \else%
    \setlength{\unitlength}{\svgwidth}%
  \fi%
  \global\let\svgwidth\undefined%
  \global\let\svgscale\undefined%
  \makeatother%
  \begin{picture}(1,1.26668155)%
    \lineheight{1}%
    \setlength\tabcolsep{0pt}%
    \put(0,0){\includegraphics[width=\unitlength,page=1]{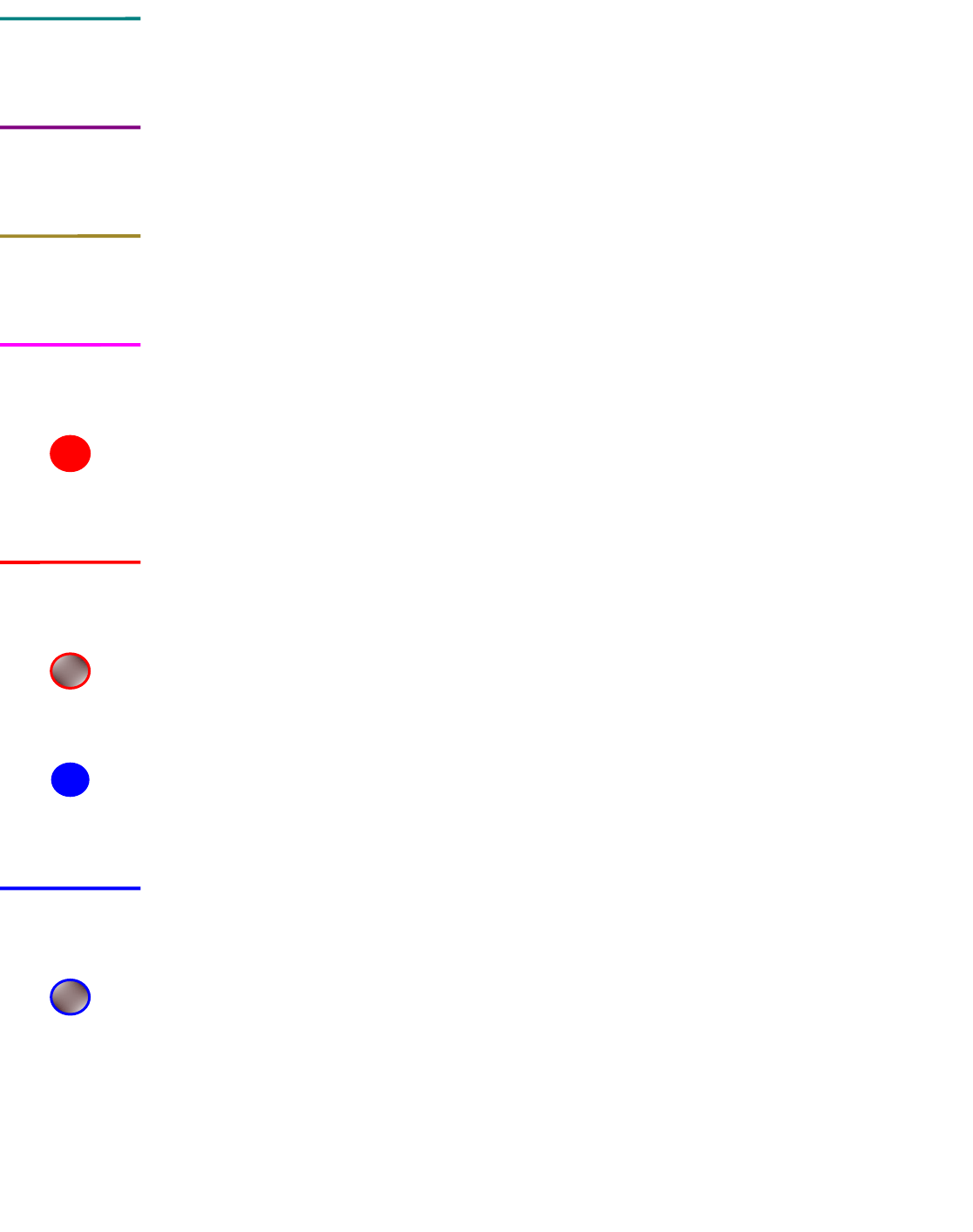}}%
    \put(0.20066397,1.2462974){\color[rgb]{0,0,0}\makebox(0,0)[lt]{\lineheight{0}\smash{\begin{tabular}[t]{l}\scriptsize{edge related to process 0 }\end{tabular}}}}%
    \put(0.20066397,1.13389424){\color[rgb]{0,0,0}\makebox(0,0)[lt]{\lineheight{0}\smash{\begin{tabular}[t]{l}\scriptsize{edge related to process 1}\end{tabular}}}}%
    \put(0.20066397,1.02149109){\color[rgb]{0,0,0}\makebox(0,0)[lt]{\lineheight{0}\smash{\begin{tabular}[t]{l}\scriptsize{edge related to process 2 }\end{tabular}}}}%
    \put(0.20066397,0.90908791){\color[rgb]{0,0,0}\makebox(0,0)[lt]{\lineheight{0}\smash{\begin{tabular}[t]{l}\scriptsize{edge related to process 3 }\end{tabular}}}}%
    \put(0.20066397,0.79668479){\color[rgb]{0,0,0}\makebox(0,0)[lt]{\lineheight{0}\smash{\begin{tabular}[t]{l}\scriptsize{first node selected for this sequence }\end{tabular}}}}%
    \put(0.20066397,0.68428159){\color[rgb]{0,0,0}\makebox(0,0)[lt]{\lineheight{0}\smash{\begin{tabular}[t]{l}\scriptsize{edge folowed to find impacted node by first choice }\end{tabular}}}}%
    \put(0.20066397,0.57187843){\color[rgb]{0,0,0}\makebox(0,0)[lt]{\lineheight{0}\smash{\begin{tabular}[t]{l}\scriptsize{blocked node by first choice }\end{tabular}}}}%
    \put(0.20066397,0.4594752){\color[rgb]{0,0,0}\makebox(0,0)[lt]{\lineheight{0}\smash{\begin{tabular}[t]{l}\scriptsize{other node(s) selected for this sequence }\end{tabular}}}}%
    \put(0.20066397,0.34707211){\color[rgb]{0,0,0}\makebox(0,0)[lt]{\lineheight{0}\smash{\begin{tabular}[t]{l}\scriptsize{edge folowed to find impacted node by other choice(s) }\end{tabular}}}}%
    \put(0.20066397,0.23466895){\color[rgb]{0,0,0}\makebox(0,0)[lt]{\lineheight{0}\smash{\begin{tabular}[t]{l}\scriptsize{blocked node by other choice(s)}\end{tabular}}}}%
    \put(0,0){\includegraphics[width=\unitlength,page=2]{Fig06l.pdf}}%
    \put(0.20061921,0.1234222){\color[rgb]{0,0,0}\makebox(0,0)[lt]{\lineheight{0}\smash{\begin{tabular}[t]{l}\scriptsize{node with distributed support}\end{tabular}}}}%
    \put(0,0){\includegraphics[width=\unitlength,page=3]{Fig06l.pdf}}%
    \put(0.20066397,0.00986271){\color[rgb]{0,0,0}\makebox(0,0)[lt]{\lineheight{0}\smash{\begin{tabular}[t]{l}\scriptsize{Node with non distributed support}\end{tabular}}}}%
  \end{picture}%
\endgroup%
		\label{seq:4leg}
	}
	\caption{Patch sequencing for figure \ref{ts_dist_load_balance} example. At each sequence, a set $\mathcal{S}$ is chosen so that as many distributed patches as possible are solved at the same time and during the same time period.  A sequence is first built  by the lowest process which chooses  one of its distributed nodes  (red node) in the graph. Then, all the other processes that are not already blocked  by this first choice do the same (blue node) in ascending order of the process identification numbers. See \ref{static_scheduling} for a detailed description and the table \ref{tabsequencing2dexample4p} for each $\mathcal{S}$. After a sequence is constructed, its nodes and connected edges are removed from the graph to give a new sub-graph used to create the next sequence until no more distributed nodes are available. }
	\label{seq:4}
\end{figure}

\subsection{local scale problem scheduling }\label{scheduling}
The 2D example in figure \ref{ts_fig} will be used to illustrate the static scheduling algorithm.
As already mentioned, two distributed patches that share the same process cannot start solving their linear system at the same time in that process.
To avoid this situation, a simple idea is to force distributed patches that share the same process to be computed at different times.
To do this, a static\footnote{Avoiding dynamic scheduling is normally more efficient} scheduling is created that  requires  all distributed patches over all processes to be computed in a specific order during the \tS loop.
The general idea is to solve as many distributed patches as possible at once, thus maximizing the use of processes and reducing the number of sequences.

To achieve this goal, the first mandatory task is to identify the dependencies between the distributed patches so that they can be sequenced.
A distributed patch depends on another if it shares a process with it.
This dependency can be extended to local patches to create an undirected graph where the vertices are all patches and the edges represent a process connection (i.e. two patches are connected by an edge if they belong even partially to the same process).
Using the mesh distribution in figure \ref{ts_dist_load_balance} this graph concept is illustrated in figure \ref{seq:4pini}: the vertices of this graph are patches (here, for a quick visual recognition of patches, the vertices of the graph are located at the same place as the enriched nodes of the mesh and  each of them has been assigned an arbitrary identifier), the edges of this graph represent a process connection (in this figure, the edges are colored the same way as macro-element of figure \ref{ts_dist_load_balance} to clearly identify the process identifier given in figure \ref{seq:4leg}).

From this undirected graph, it is possible to create a set $\mathcal{S}$ of distributed vertices such that for all vertices in $\mathcal{S}$, no edge connect them.
Thus, these vertices  can be computed at the same time because they do not share a common process (no edges connect them).
Finding this  set $\mathcal{S}$ (not unique)  corresponds in graph theory to  the "independent set problem" and will correspond to finding a sequence of the scheduling. 
It is interesting to maximize the size of $\mathcal{S}$ (which corresponds to the NP-hard  "maximum independent set problem" in graph theory) because on the one hand it maximizes a priori the number of used processes  and therefore the load balancing is good. On the other hand, the distributed vertices are consumed faster and therefore the number of sequences can be minimized.
When  the set $\mathcal{S}$  is assigned to a sequence (mainly by creating a colored  MPI communicator associated with it),  its vertices (and their connecting edges) can be removed from the graph.
The maximum independent set problem is then solved on the remaining subgraph to find the next sequence.
All sequences are found this way  until all distributed patches are removed from the subgraph.	
Then the remaining local patches can be processed in any order as they are independent. 

\begin{table}[h]
	\footnotesize
	\subfloat[distribution of figure \ref{ts_dist_load_balance}]{\begin{tabular}{|c|c|c|c|c|c|c|c|}
			\hline
			\multicolumn{8}{|c|}{Total number of distributed patches: 11}\\
			\hline
			\multicolumn{8}{|c|}{Maximum number of distributed patches per process: 8}\\
			\hline
			\multicolumn{8}{|c|}{Minimum number of distributed patches per process: 4}\\
			\hline
			\multicolumn{8}{|c|}{Average number of distributed patches per process: 6}   \\ 
			\hline
			\rotatebox[origin=c]{90}{sequence}& \rotatebox[origin=c]{90}{Active patches}& \rotatebox[origin=c]{90}{Active distributed} \rotatebox[origin=c]{90}{patches}& \rotatebox[origin=c]{90}{\% Active distributed} \rotatebox[origin=c]{90}{ patches}	& \rotatebox[origin=c]{90}{Active processes}  &
			\rotatebox[origin=c]{90}{\% Active processes} &
			\rotatebox[origin=c]{90}{$\mathcal{S}$} 	&	 \rotatebox[origin=c]{90}{Figure}\\
			\hline
			1&2&1&50&4&100&18,3&\ref{seq:4pseq1}\\
			\hline
			2&3&1&33.3&4&100&12,2,20&\ref{seq:4pseq2}\\
			\hline
			3&2&1&50&4&100&10,21&\ref{seq:4pseq3}\\
			\hline
			4&3&1&33.3&4&100&13,1,23&\ref{seq:4pseq4}\\
			\hline
			5&2&2&100&4&100&17,6&\ref{seq:4pseq5}\\
			\hline
			6&2&2&100&4&100&7,19&\ref{seq:4pseq6}\\
			\hline
			7&3&1&33.3&4&100&8,11,22&\ref{seq:4pseq7}\\
			\hline
			8&2&2&100&4&100&15,5&\ref{seq:4pseq8}\\
			\hline
			\hline
			average&2.4&1.4&62.5&4&100&\multicolumn{2}{c}{~}\\
			\cline{1-6}
		\end{tabular}\label{tabsequencing2dexample4p}}
	\hspace{2em}
	\subfloat[distribution of figure  \ref{ts_dist_load_balance_8p}]{
				\raisebox{2.6ex}{
			\begin{tabular}{|c|c|c|c|c|c|c|}
			\hline
			\multicolumn{7}{|c|}{Total number of distributed patches: 17}\\
			\hline
			\multicolumn{7}{|c|}{Maximum number of distributed patches per process: 6}\\
			\hline
			\multicolumn{7}{|c|}{Minimum number of distributed patches per process: 4}\\
			\hline
			\multicolumn{7}{|c|}{Average number of distributed patches per process: 5.1}\\
			\hline
			\rotatebox[origin=c]{90}{sequence}& \rotatebox[origin=c]{90}{Active patches}& \rotatebox[origin=c]{90}{Active distributed} \rotatebox[origin=c]{90}{patches}& \rotatebox[origin=c]{90}{\% Active distributed} \rotatebox[origin=c]{90}{patches}	& \rotatebox[origin=c]{90}{Active processes}  &
			\rotatebox[origin=c]{90}{\% Active processes} &
			\rotatebox[origin=c]{90}{$\mathcal{S}$} 	\\
			\hline
			1&3&3&100&8&100&12,7,17\\
			\hline
			2&3&2&66.7&8&100&10,18,21\\
			\hline
			3&3&3&100&8&100&13,6,20\\
			\hline
			4&4&4&100&8&100&8,5,19,15\\
			\hline
			5&4&4&100&8&100&14,3,11,25\\
			\hline
			6&3&1&33.3&4&50&9,1,23\\
			\hline
			\hline
			average&3.3&2.8&83.3&7.3&91.7&\multicolumn{1}{c}{~}\\
			\cline{1-6}
		\end{tabular}
\label{tabsequencing2dexample8p}	
    }	
	}
	\caption{Distributed patches sequencing with the algorithm \ref{sequencing_algo:general} for the   fictitious 2D problem in figure\ref{ts_fig} distributed over 4 and 8 processes (figure \ref{ts_dist_fig}). It does not correspond  the calculation of all patches  because the average does not take into account the computation of the non-constrained local patches. 	
	}
	\label{tabsequencing2dexample}
\end{table}

The condition of maximizing the size of $\mathcal{S}$  is not sufficient to obtain a balanced load.
The vertices of $\mathcal{S}$ must have the same "load" (in a sequence, if a vertex $v_i$ of $\mathcal{S}$ takes longer to compute than the others, most of the process will be idle while waiting for $v_i$ to finish).
Thus, with distributed patches  having different computational consumption (called weights hereafter), the choice of $\mathcal{S}$ must  be oriented by these weights in order to have approximately the same "load" for all the vertices of  $\mathcal{S}$.
In this work, these weights are simply the number of micro-elements in a patch which reflects its number of dofs and thus its computational cost (see figure \ref{seq:4weigt} computed from the weights in figure \ref{ts_dist_weigth_elem}).
The NP-hard nature of the maximum independent set problem, this weighting constraint, and the a priori distributed nature of the undirected graph lead to the choice of a basic sub-optimal  heuristic to quickly find a static scheduling.

More precisely, the chosen algorithm, formalized in the algorithms \ref{sequencing_algo:general}, \ref{sequencing_algo:pick_first} and \ref{sequencing_algo:pick}  detailed  in \ref{static_scheduling}, is played in the table \ref{tabsequencing2dexample4p} and figure \ref{seq:4} for the example of the figure \ref{ts_dist_load_balance}.
In this 2D example, 8 sequences are required, which corresponds to the maximum number of distributed patches per process: in the process with 8 distributed patches, there cannot be less than 8 sequences.
With the example of figure \ref{ts_fig}, distributed on  8 processes (see figure \ref{ts_dist_load_balance_8p}), the table \ref{tabsequencing2dexample8p} shows that in this case 6 sequences must be computed.
But in this case, the last sequence is less efficient because process 3, 4, 5 and 7 are note used.
This shows that the scalability of the distributed \tS method will be influenced by the quality of the algorithm \ref{sequencing_algo:general}.
In the section \ref{PO}, and in particular in the table \ref{tabsortingversion} of this section, we study in a real example the heuristic effect of this algorithm.

\section{Numerical simulation}\label{NSOLV}
All simulation were run on the Liger cluster (see  \ref{cluster} for  details) in the exclusive node condition (no other users at the node level but  shared resources at the network and I/O level). 
In many cases, the \tS solver  is  compared  to the full rank direct solver MUMPS (denoted by "fr" hereafter), to a low-rank direct solver based on MUMPS (see \ref{BLRRES}, denoted by "blr" hereafter) and to a domain decomposition solver (see \ref{DDRES}, denoted by "dd" hereafter). 
Hereinafter, the  \TSD and \TSI versions are referred to as  "ts" and "tsi" respectively.
In terms of notation E will be the Young's modulus and $\nu$ the Poisson's ratio.
\subsection{Analytic solution}\label{UOL}
This test is intended to validate the proposed work using a given cubic displacement field over a plate:

\noindent\begin{minipage}{0.35\textwidth}
	\def\svgwidth{1\textwidth}
\begingroup%
  \makeatletter%
  \providecommand\color[2][]{%
    \errmessage{(Inkscape) Color is used for the text in Inkscape, but the package 'color.sty' is not loaded}%
    \renewcommand\color[2][]{}%
  }%
  \providecommand\transparent[1]{%
    \errmessage{(Inkscape) Transparency is used (non-zero) for the text in Inkscape, but the package 'transparent.sty' is not loaded}%
    \renewcommand\transparent[1]{}%
  }%
  \providecommand\rotatebox[2]{#2}%
  \ifx\svgwidth\undefined%
    \setlength{\unitlength}{314.96022949bp}%
    \ifx\svgscale\undefined%
      \relax%
    \else%
      \setlength{\unitlength}{\unitlength * \real{\svgscale}}%
    \fi%
  \else%
    \setlength{\unitlength}{\svgwidth}%
  \fi%
  \global\let\svgwidth\undefined%
  \global\let\svgscale\undefined%
  \makeatother%
  \begin{picture}(1,0.56063393)%
    \put(0,0){\includegraphics[width=\unitlength]{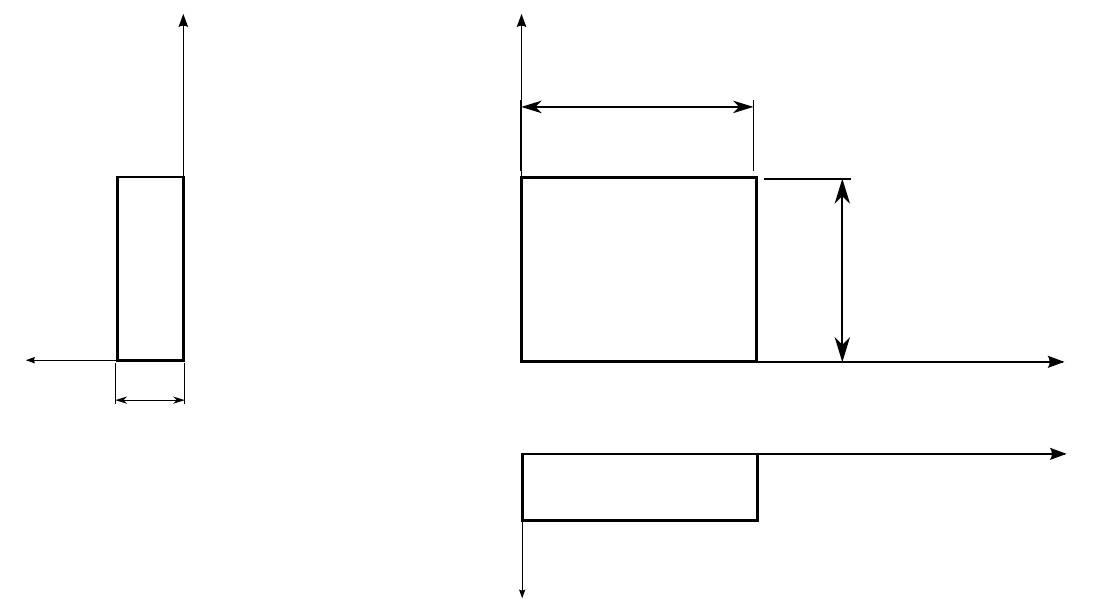}}%
    \put(0.46884916,0.19704033){\color[rgb]{0,0,0}\makebox(0,0)[lb]{\smash{\tiny{A}}}}%
    \put(0.67043044,0.19855609){\color[rgb]{0,0,0}\makebox(0,0)[lb]{\smash{\tiny{B}}}}%
    \put(0.17319528,0.23412865){\color[rgb]{0,0,0}\makebox(0,0)[lb]{\smash{\tiny{A}}}}%
    \put(0.69587481,0.12385927){\color[rgb]{0,0,0}\makebox(0,0)[lb]{\smash{\tiny{B}}}}%
    \put(0.48365742,0.12242525){\color[rgb]{0,0,0}\makebox(0,0)[lb]{\smash{\tiny{A}}}}%
    \put(0.17511695,0.21236751){\color[rgb]{0,0,0}\makebox(0,0)[lb]{\smash{\tiny{B}}}}%
    \put(0.56908905,0.47847215){\color[rgb]{0,0,0}\makebox(0,0)[lb]{\smash{$K$}}}%
    \put(0.78681724,0.31719665){\color[rgb]{0,0,0}\makebox(0,0)[lb]{\smash{$H$}}}%
    \put(0.95109502,0.19){\color[rgb]{0,0,0}\makebox(0,0)[lb]{\smash{$x$}}}%
    \put(0.11908153,0.1450){\color[rgb]{0,0,0}\makebox(0,0)[lb]{\smash{$T$}}}%
    \put(0.42877775,0.52476912){\color[rgb]{0,0,0}\makebox(0,0)[lb]{\smash{$y$}}}%
    \put(0.12466061,0.52933283){\color[rgb]{0,0,0}\makebox(0,0)[lb]{\smash{$y$}}}%
    \put(0.0110133,0.19){\color[rgb]{0,0,0}\makebox(0,0)[lb]{\smash{$z$}}}%
    \put(0.4463329,0.4155655){\color[rgb]{0,0,0}\makebox(0,0)[lb]{\smash{\tiny{C}}}}%
    \put(0.17405057,0.38717304){\color[rgb]{0,0,0}\makebox(0,0)[lb]{\smash{\tiny{C}}}}%
    \put(0.94349493,0.1){\color[rgb]{0,0,0}\makebox(0,0)[lb]{\smash{$x$}}}%
    \put(0.48527821,0.01568651){\color[rgb]{0,0,0}\makebox(0,0)[lb]{\smash{$z$}}}%
    \put(0.48063054,0.14967192){\color[rgb]{0,0,0}\makebox(0,0)[lb]{\smash{\tiny{C}}}}%
  \end{picture}%
\endgroup%
	\captionof{figure}{Analytical test:  plate under volume and surface loading. A, B, and C are the corner points used to fix the rigid body modes.
		\label{UOL_geo}}
\end{minipage}\hfill
\begin{minipage}{0.62\textwidth}
\begin{equation}
	\Vec{u}^C(x,y,z)= \frac{(\nu+1)(1-2\nu)F}{E}\left[ \begin{array}{l} (x^2(\frac{K}{2}-\frac{x}{3})+2\nu(y^2-z^2))\\
		-4\nu xy\\
		4\nu xz
	\end{array}
	\right]
	\label{eqUC}
\end{equation}
where : 
\begin{itemize}
	\item $E=36.5$GPa
	\item $\nu=0.2$
	\item $F$ is a scalar corresponding to a force
	\item $K$ is the plate length  (figure \ref{UOL_geo} )
\end{itemize}
\end{minipage}
\vspace{1em}

The imposed loads are obtained from the equation \eqref{eqUC} (writing the equilibrium of the system) and the plate is made of a  homogeneous isotropic elastic material.
The plate is discretized with a coarse mesh (595 tetrahedrons, 188 vertices) which corresponds, hereafter, to the level 0 (L0).
Up to 4 other discretizations (L1,L2,L3,L4) derived from this level are used as the global-scale problem.
Each of them is built using the adaptation strategy described in \ref{rsplit_anexe}  (considering the whole mesh as an area of interest: no hanging nodes).
Each level corresponds to the previous one with all its elements divided once (L2 is L1 divided once). 
This corresponds exactly to the refinement process used for scaling.
Thus, since all nodes are enriched, the 4 meshes are identical to  any refined SP targeting the same level (L2 mesh will be the same as the refined \SP created by refining L0 mesh 2 times).

We evaluate the performance of the proposed parallel \tS solver  and  compare it to the solvers "fr", "blr" and "dd".
In all figures, for the \tS curves,  "L$Y$ from L$X$" or "from L$X$" (when the fine scale-level "$Y$" is implicit) means to use a "$X$" level  for the discretization of the global level and to use  a "$Y$" level  for the discretization of the \SP ($Y>X$).
Since the full rank direct solver provides a residual error close to the machine accuracy, a relatively fair comparison requires to adopt a value of $1.e-7$  for $\epsilon$ in the algorithm \ref{TS_algebra},\ref{domain_decomposition_algo} and \ref{blr_algo}.
As noted in section \ref{parallel_paradigm}, for comparison purposes, the multithreading capability of MUMPS and underlying BLAS is not used.
Thus, the number of cores (abscissa in many figures) directly represents  the number of MPI process dispatched on one or more nodes of the cluster.
\begin{figure}[h]
	\subfloat[L2]{

\label{UOL_curve_times_L7}}
	\caption{Cubic field: elapsed time in second versus the number of cores. Log/log scale. The results for a fine-scale refinement at  L2 , L3, L4, L5, L6  and L7 level are presented  in 
	\protect\subref*{UOL_curve_times_L2},
	\protect\subref*{UOL_curve_times_L3},
	\protect\subref*{UOL_curve_times_L4},
	\protect\subref*{UOL_curve_times_L5},
	\protect\subref*{UOL_curve_times_L6} and \protect\subref*{UOL_curve_times_L7} respectively.
	The ideal curve shows the  slope when the elapsed time of a single process is perfectly divided by the number of processes used for the calculation (this corresponds to the ideal speed-up).
	This curve is always shifted so that it starts at the same point as the best \tS solution for the smallest number of cores. 
		\label{UOL_curve_times}}
\end{figure}
Figure \ref{UOL_curve_times} presents, for different SP refinement, an analysis of the elapsed time used by all solvers.
The points not present on these curves are mainly due to the fact that, with this cluster, it was impossible to compute them: not enough memory, too much computation time, too many processes for the solver or integer overflow.
These limitations only apply to the actual cluster installation and are in no way an absolute judgment on the software used, especially in a non-multi-threaded context.

\begin{figure}[h]
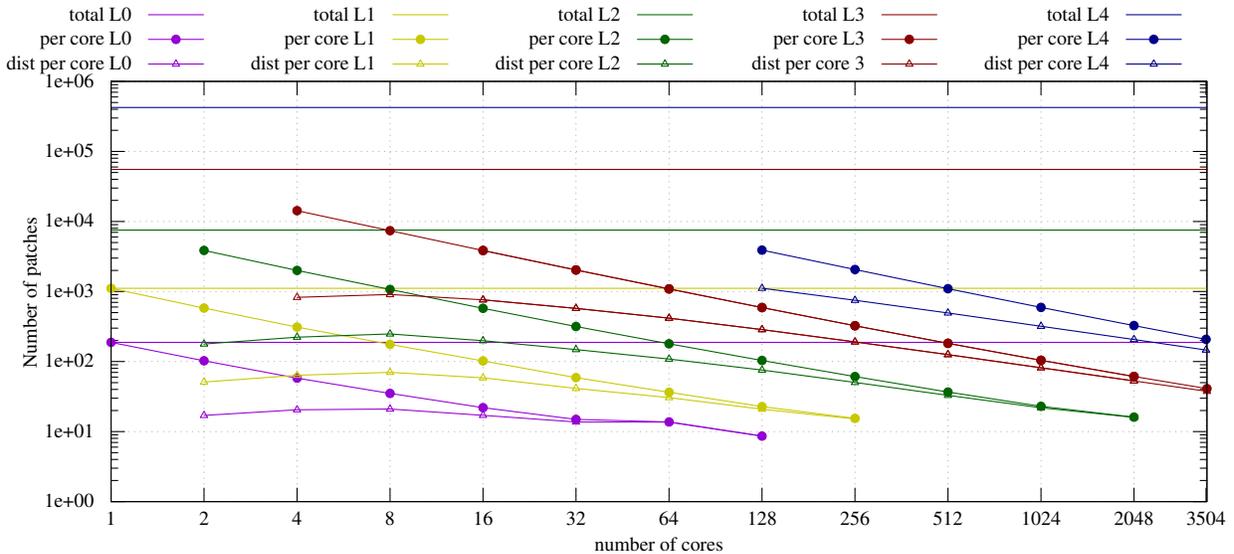

	\centering

	\caption{Cubic field: number of patches. Log/log scale.\label{UOL_curve_patch}.}
\end{figure}

As discussed  in section \ref{theoSeqAlgPerf}, the \tS solver has an ideal scale ratio.
Thus, for the \SP refinement of levels L3, L4, L5 and L6 different coarse-scale levels were tested.
The observed  optimal scale  ratio  is 3 (best elapsed times curves are "tsi from L1", "tsi from L2" and  "tsi from L3" respectively for the \SP refinement of levels L4 (\ref{UOL_curve_times_L4}), L5 (\ref{UOL_curve_times_L5}) and L6 (\ref{UOL_curve_times_L6}): 1+3=4, 2+3=5 and 3+3=6), which is in good agreement with the analysis in section \ref{theoSeqAlgPerf}  and gives  equation  \eqref{ts_jump_level_eq}.
For L7, only this ratio was used to save computing resources.
Regarding the number of dofs, the figure \ref{UOL_curve_dofs_ts} shows that this ratio corresponds to a fixed number of dofs per patch appearing visually as a horizontal line (circle between 3435 and 4928).
\begin{figure}
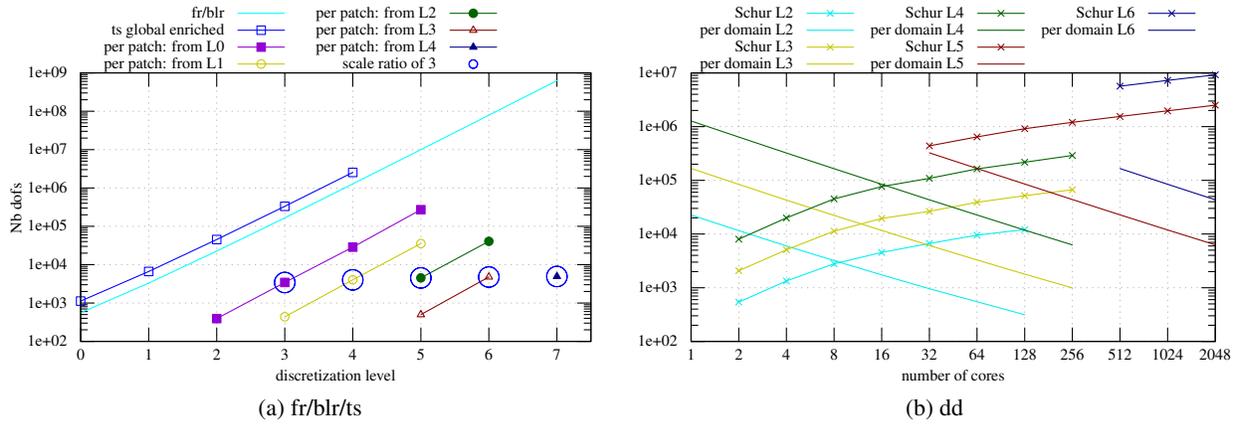

	\subfloat[fr/blr/ts]{

\label{UOL_curve_dofs_dd}}
	\caption{Cubic field: degree of freedom numbers \label{UOL_curve_dofs}.}
\end{figure}
Compare to the number of "fr/br"   dofs  that increases with discretization, the \tS solver can keep the size of the fine-scale resolutions  within a specific range as long as the scale ratio is maintained.
Keeping the same ratio implies increasing the discretization level for the global scale problem.
In the figure \ref{UOL_curve_dofs_ts} the "ts global enriched" curve is simply 2 times the "fr/blr" curve (all nodes are enriched).
This already highlights the issue of having  a global problem with a not so small size.
This point is addressed by the \TSI version as verified below.
Another impact can be seen in figure \ref{UOL_curve_patch} where increasing the discretization level for a global scale problem implies an increase in the "total" number of patches involved in the computation.
This last point is naturally treated by the parallelism as can be seen on the same figure.
The number of patches calculated per process decreases when using an increasing number of cores.

For  the \SP refinement of level L2  (figure \ref{UOL_curve_times_L2})   the \tS solver does not provide better performance compared to the tested solvers (for "fr" and "blr", it is only with a high number of cores that the \tS solver reaches the same level of performance of these solvers before their efficiency loss).
One of the reasons is certainly that the optimal ratio of 3 cannot be reached and thus the cost of overlapping patches is not compensated by an improvement in algorithmic complexity (fine-scale problems size to small).
When this ratio of 3 is reached (figure \ref{UOL_curve_times_L3}), the \tS solver starts to provide performances similar to those of other solutions.
From  Level 4 and above,  the \tS solver  performs better than all the other  solutions tested.

In figure \ref{UOL_curve_times_L5} (L5), corresponding to "fr/blr" $\sim1.e^{+7}$ dofs (figure \ref{UOL_curve_dofs_ts}),  the potentiality of the other solver is expressed, as expected, in different ranges of number of cores.
The "fr" solver provides an exact result but has a limited scalability and needs memory (the computation was not possible with less than 16 cores distributed over 16 nodes : 16x128Go$\approx$2To).
The "blr" solver reduces the memory footprint and elapsed time compared to the full rank version with  a controlled error ( $\epsilon=1.e^{-7}$).
The "dd" solver offers good scalability but fails for a small number of processes.
This is related to the choice of implementation which uses only one domain per process.
Thus, for a small number of cores, the domains become very large and their resolution cannot fit in a single node of the cluster: the figure \ref{UOL_curve_dofs_dd} shows that the average size of the domain increases to the  "fr/blr" size when the number of cores decreases to 1 process.
A classical alternative, not tested in this work, is to consider that the domain size is independent of the process (i.e. a process can compute more than one domain whose size is controlled by an arbitrary threshold).
But the Schur complement is then potentially handled by fewer processes compared to the proposed "dd" implementation, and thus would have added  additional complexity in the analyses.
However, it would have been easy to compare the case where the number of dofs for domains and patches is the same.
Nevertheless, with 2048 processes, the "dd" domains represent about the same number of dofs as the  patches (figure \ref{UOL_curve_dofs}) and the \tS solver  shows better performance in this condition.
This comparison should be addressed in more detail and is left as a future prospect.
As for the other solvers,  with a small number of cores, the \tS memory consumption becomes an problem (impossible computation  with less than 2 process for the "from L2" version)  but clearly at a lower level.
This validates the choice made in the implementation regarding the per-macro-element  storage (see section \ref{TS_integ_schem}) and shows the low memory consumption impact of the fines scale problems.
Also in figure \ref{UOL_curve_times_L5}, the curves show that the proposed \tS solver, even with a non optimal scale ratio, consistently gives  better performance from 2 to 2048 processes compared to other solvers.
Note also that another interesting reading of this figure \ref{UOL_curve_times_L5} gives the following conclusion:  with a small number of processes, the \tS solver gives the same level of performance (i.e. elapsed time) compared to other solvers using larger number of processes.

With L6 (figure \ref{UOL_curve_times_L6}, 78 851 421  fine scale dofs) and L7 (figure \ref{UOL_curve_times_L7}, 627 350 205 fine scale dofs), the \tS solver confirms its low resource requirements.
The memory consumption allows to launch calculations with "only" 16 and 128 processes.
For L7, no other solver passes with the chosen hardware and software configuration.
And for L6 only "dd" starts to provide results with 512 cores that the  \tS solver provides with only 16 processes in the  same time frame.

More detailed information about the  "L6 from L3" test is given in the table \ref{UOL_table_ass}.
\begin{table}[h]
	\begin{tabular}{r|c|c|c|c|c|c|}
		\cline{2-7}
		&\multicolumn{6}{c|}{Task: elapsed time in s and \% of total resolution }\\
		\hline
		\multicolumn{1}{|r|}{\TS resolution}&1314.4&100&67.3&100&29.6&100\\
		\hline
		\multicolumn{1}{|r|}{$resi$ computation}&34.4&2.6&1.8&2.6&0.6&2.2\\
        \hline
        \multicolumn{1}{|r|}{$\vm{u}_r$ update}&15.9&1.2&0.6&0.9&0.2&0.6\\
        \hline
		\multicolumn{1}{|r|}{$\vm{A}_{gg}$ (product,operator construction and 		assembly of equation \ref{Aggs} )}&261.6&19.9&8.2&12.1&2.0&6.9\\
		\multicolumn{1}{|r|}{$\vm{B}_{g}$ (product and assembly of equation \ref{Bg})}&3.4&0.3&0.1&0.2&0.03&0.1\\
		\multicolumn{1}{|r|}{$\vm{A}_{FF}^{e_{macro}}$ (computation and assembly per macro-element)}&153.3&11.7&4.9&7.3&1.2&4.1\\
		\multicolumn{1}{|r|}{$\vm{B}_{F}^{e_{macro}}$ (computation and assembly per macro-element)}&36.1&2.7&1.1&1.7&0.3&0.9\\
		
		\multicolumn{1}{|r|}{System of equation \ref{g_sys} construction (sum of above)}&454.4&34.6&14.3&21.3&3.5&12\\
		\hline
		nb of cores&\multicolumn{2}{c|}{16}&\multicolumn{2}{c|}{512}&\multicolumn{2}{c|}{2048}\\
		\cline{2-7}
	\end{tabular}
\caption{Detailed elapsed times  of the global matrix construction, $resi$ computation and $\vm{u}_r$ update  for some points of "tsi" L6 from L3 curve (figure \ref{UOL_curve_times_L6}) \label{UOL_table_ass} }
\end{table}
This table shows that the construction time of the global system (equation \ref{g_sys}) represents less than 35\% of the \tS solver resolution with 16 processes.
It is only  12\% with 2048 processes, which can be explained by the perfect scalability of this task (454.4/128$\approx$3.5) compared to the scalability of the \tS solver  (1314.4/128$\approx$10.3$<$29.6).
The different calculation($\vm{A}_{gg}$,$\vm{B}_{g}$,$\vm{A}_{FF}^{e_{macro}}$ and $\vm{B}_{F}^{e_{macro}}$) have the same perfect scalability and contribute with the same ratio to the construction task.
The most expensive calculation is related to the construction of $\vm{A}_{gg}$  (19.9\%, 12.1\% and 6.9\%) and in particular the product $\vm{A}_{FF}^{e_{macro}}\cdot \vm{T}_{Fe}^{e_{macro}}$  which represents $\sim$56\% of the construction of $\vm{A}_{gg}$  (result not given in the table \ref{UOL_table_ass}).
It is a sparse matrix by sparse matrix product (not using the  level 3 BLAS routine) performed on all SP elements at each \tS iteration.
Perhaps an effort on its implementation  would allow to obtain better performances but it would be anyway only on a small part of the \tS resolution.   
The performance of the global scale system construction is therefore no longer analyzed in this work (but included in the results) because it is  second order in the calculations and  scales perfectly.
Regarding the computation of $resi$,  it consistently represents less than 3\% of the \tS resolution and scales relatively (few communication) well (34.4/128$\approx$0.3$<$0.6).
This validates the proposed algorithm \ref{TS_algebra_residual}.
The same conclusion applies to the update of $\vm{u}_r$ (algorithm \ref{TS_algebra_micro_update}) which costs less than 1,3\% and scales well (15,9/128$\approx$0,12$\lessapprox$0,2).

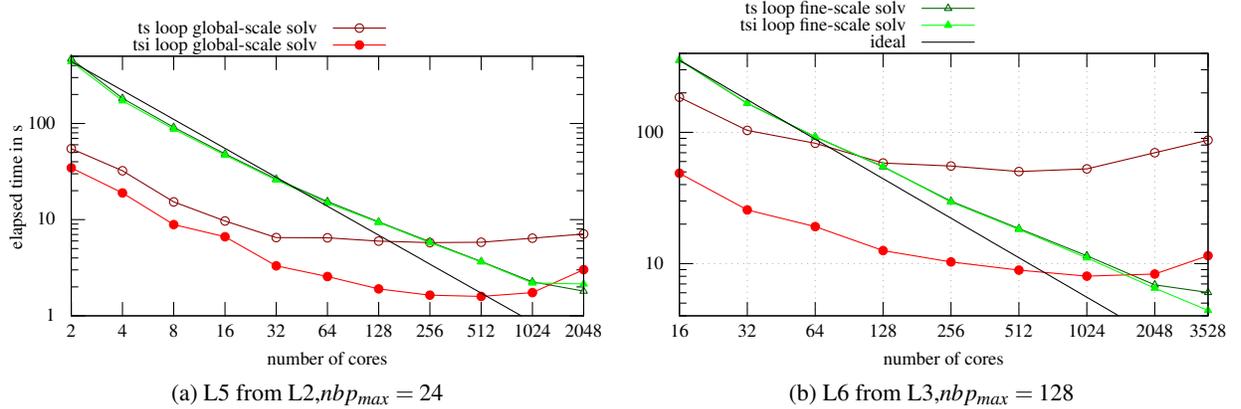
\begin{figure}
	\subfloat[L5 from L2,$nbp_{max}=24$]{
\begin{tikzpicture}[gnuplot]
\tikzset{every node/.append style={scale=0.70}}
\path (0.000,0.000) rectangle (8.125,4.812);
\gpcolor{color=gp lt color border}
\gpsetlinetype{gp lt border}
\gpsetdashtype{gp dt solid}
\gpsetlinewidth{1.00}
\draw[gp path] (0.925,0.691)--(1.105,0.691);
\draw[gp path] (7.737,0.691)--(7.557,0.691);
\node[gp node right] at (0.796,0.691) {$1$};
\draw[gp path] (0.925,1.076)--(1.015,1.076);
\draw[gp path] (7.737,1.076)--(7.647,1.076);
\draw[gp path] (0.925,1.302)--(1.015,1.302);
\draw[gp path] (7.737,1.302)--(7.647,1.302);
\draw[gp path] (0.925,1.461)--(1.015,1.461);
\draw[gp path] (7.737,1.461)--(7.647,1.461);
\draw[gp path] (0.925,1.586)--(1.015,1.586);
\draw[gp path] (7.737,1.586)--(7.647,1.586);
\draw[gp path] (0.925,1.687)--(1.015,1.687);
\draw[gp path] (7.737,1.687)--(7.647,1.687);
\draw[gp path] (0.925,1.773)--(1.015,1.773);
\draw[gp path] (7.737,1.773)--(7.647,1.773);
\draw[gp path] (0.925,1.847)--(1.015,1.847);
\draw[gp path] (7.737,1.847)--(7.647,1.847);
\draw[gp path] (0.925,1.912)--(1.015,1.912);
\draw[gp path] (7.737,1.912)--(7.647,1.912);
\draw[gp path] (0.925,1.971)--(1.105,1.971);
\draw[gp path] (7.737,1.971)--(7.557,1.971);
\node[gp node right] at (0.796,1.971) {$10$};
\draw[gp path] (0.925,2.356)--(1.015,2.356);
\draw[gp path] (7.737,2.356)--(7.647,2.356);
\draw[gp path] (0.925,2.581)--(1.015,2.581);
\draw[gp path] (7.737,2.581)--(7.647,2.581);
\draw[gp path] (0.925,2.741)--(1.015,2.741);
\draw[gp path] (7.737,2.741)--(7.647,2.741);
\draw[gp path] (0.925,2.865)--(1.015,2.865);
\draw[gp path] (7.737,2.865)--(7.647,2.865);
\draw[gp path] (0.925,2.967)--(1.015,2.967);
\draw[gp path] (7.737,2.967)--(7.647,2.967);
\draw[gp path] (0.925,3.052)--(1.015,3.052);
\draw[gp path] (7.737,3.052)--(7.647,3.052);
\draw[gp path] (0.925,3.126)--(1.015,3.126);
\draw[gp path] (7.737,3.126)--(7.647,3.126);
\draw[gp path] (0.925,3.192)--(1.015,3.192);
\draw[gp path] (7.737,3.192)--(7.647,3.192);
\draw[gp path] (0.925,3.250)--(1.105,3.250);
\draw[gp path] (7.737,3.250)--(7.557,3.250);
\node[gp node right] at (0.796,3.250) {$100$};
\draw[gp path] (0.925,3.636)--(1.015,3.636);
\draw[gp path] (7.737,3.636)--(7.647,3.636);
\draw[gp path] (0.925,3.861)--(1.015,3.861);
\draw[gp path] (7.737,3.861)--(7.647,3.861);
\draw[gp path] (0.925,4.021)--(1.015,4.021);
\draw[gp path] (7.737,4.021)--(7.647,4.021);
\draw[gp path] (0.925,4.145)--(1.015,4.145);
\draw[gp path] (7.737,4.145)--(7.647,4.145);
\draw[gp path] (0.925,0.691)--(0.925,0.871);
\draw[gp path] (0.925,4.145)--(0.925,3.965);
\node[gp node center,font={\fontsize{10.0pt}{12.0pt}\selectfont}] at (0.925,0.475) {2};
\draw[gp path] (1.606,0.691)--(1.606,0.871);
\draw[gp path] (1.606,4.145)--(1.606,3.965);
\node[gp node center,font={\fontsize{10.0pt}{12.0pt}\selectfont}] at (1.606,0.475) {4};
\draw[gp path] (2.287,0.691)--(2.287,0.871);
\draw[gp path] (2.287,4.145)--(2.287,3.965);
\node[gp node center,font={\fontsize{10.0pt}{12.0pt}\selectfont}] at (2.287,0.475) {8};
\draw[gp path] (2.969,0.691)--(2.969,0.871);
\draw[gp path] (2.969,4.145)--(2.969,3.965);
\node[gp node center,font={\fontsize{10.0pt}{12.0pt}\selectfont}] at (2.969,0.475) {16};
\draw[gp path] (3.650,0.691)--(3.650,0.871);
\draw[gp path] (3.650,4.145)--(3.650,3.965);
\node[gp node center,font={\fontsize{10.0pt}{12.0pt}\selectfont}] at (3.650,0.475) {32};
\draw[gp path] (4.331,0.691)--(4.331,0.871);
\draw[gp path] (4.331,4.145)--(4.331,3.965);
\node[gp node center,font={\fontsize{10.0pt}{12.0pt}\selectfont}] at (4.331,0.475) {64};
\draw[gp path] (5.012,0.691)--(5.012,0.871);
\draw[gp path] (5.012,4.145)--(5.012,3.965);
\node[gp node center,font={\fontsize{10.0pt}{12.0pt}\selectfont}] at (5.012,0.475) {128};
\draw[gp path] (5.693,0.691)--(5.693,0.871);
\draw[gp path] (5.693,4.145)--(5.693,3.965);
\node[gp node center,font={\fontsize{10.0pt}{12.0pt}\selectfont}] at (5.693,0.475) {256};
\draw[gp path] (6.375,0.691)--(6.375,0.871);
\draw[gp path] (6.375,4.145)--(6.375,3.965);
\node[gp node center,font={\fontsize{10.0pt}{12.0pt}\selectfont}] at (6.375,0.475) {512};
\draw[gp path] (7.056,0.691)--(7.056,0.871);
\draw[gp path] (7.056,4.145)--(7.056,3.965);
\node[gp node center,font={\fontsize{10.0pt}{12.0pt}\selectfont}] at (7.056,0.475) {1024};
\draw[gp path] (7.737,0.691)--(7.737,0.871);
\draw[gp path] (7.737,4.145)--(7.737,3.965);
\node[gp node center,font={\fontsize{10.0pt}{12.0pt}\selectfont}] at (7.737,0.475) {2048};
\draw[gp path] (0.925,4.145)--(0.925,0.691)--(7.737,0.691)--(7.737,4.145)--cycle;
\node[gp node center,rotate=-270] at (0.218,2.418) {elapsed time in s};
\node[gp node center] at (4.331,0.105) {number of cores};
\node[gp node right] at (4.279,4.519) {ts loop global-scale solv};
\gpcolor{rgb color={0.545,0.000,0.000}}
\draw[gp path] (4.408,4.519)--(5.104,4.519);
\draw[gp path] (0.925,2.915)--(1.606,2.621)--(2.287,2.207)--(2.969,1.953)--(3.650,1.732)%
  --(4.331,1.729)--(5.012,1.685)--(5.693,1.665)--(6.375,1.670)--(7.056,1.724)--(7.737,1.780);
\gpsetpointsize{4.00}
\gp3point{gp mark 6}{}{(0.925,2.915)}
\gp3point{gp mark 6}{}{(1.606,2.621)}
\gp3point{gp mark 6}{}{(2.287,2.207)}
\gp3point{gp mark 6}{}{(2.969,1.953)}
\gp3point{gp mark 6}{}{(3.650,1.732)}
\gp3point{gp mark 6}{}{(4.331,1.729)}
\gp3point{gp mark 6}{}{(5.012,1.685)}
\gp3point{gp mark 6}{}{(5.693,1.665)}
\gp3point{gp mark 6}{}{(6.375,1.670)}
\gp3point{gp mark 6}{}{(7.056,1.724)}
\gp3point{gp mark 6}{}{(7.737,1.780)}
\gp3point{gp mark 6}{}{(4.756,4.519)}
\gpcolor{rgb color={0.000,0.392,0.000}}
\draw[gp path] (0.925,4.110)--(1.606,3.585)--(2.287,3.199)--(2.969,2.846)--(3.650,2.511)%
  --(4.331,2.212)--(5.012,1.941)--(5.693,1.676)--(6.375,1.414)--(7.056,1.143)--(7.737,1.020);
\gp3point{gp mark 8}{}{(0.925,4.110)}
\gp3point{gp mark 8}{}{(1.606,3.585)}
\gp3point{gp mark 8}{}{(2.287,3.199)}
\gp3point{gp mark 8}{}{(2.969,2.846)}
\gp3point{gp mark 8}{}{(3.650,2.511)}
\gp3point{gp mark 8}{}{(4.331,2.212)}
\gp3point{gp mark 8}{}{(5.012,1.941)}
\gp3point{gp mark 8}{}{(5.693,1.676)}
\gp3point{gp mark 8}{}{(6.375,1.414)}
\gp3point{gp mark 8}{}{(7.056,1.143)}
\gp3point{gp mark 8}{}{(7.737,1.020)}
\gpcolor{color=gp lt color border}
\node[gp node right] at (4.279,4.294) {tsi loop global-scale solv};
\gpcolor{rgb color={1.000,0.000,0.000}}
\draw[gp path] (4.408,4.294)--(5.104,4.294);
\draw[gp path] (0.925,2.660)--(1.606,2.326)--(2.287,1.905)--(2.969,1.743)--(3.650,1.356)%
  --(4.331,1.214)--(5.012,1.049)--(5.693,0.966)--(6.375,0.949)--(7.056,1.000)--(7.737,1.305);
\gp3point{gp mark 7}{}{(0.925,2.660)}
\gp3point{gp mark 7}{}{(1.606,2.326)}
\gp3point{gp mark 7}{}{(2.287,1.905)}
\gp3point{gp mark 7}{}{(2.969,1.743)}
\gp3point{gp mark 7}{}{(3.650,1.356)}
\gp3point{gp mark 7}{}{(4.331,1.214)}
\gp3point{gp mark 7}{}{(5.012,1.049)}
\gp3point{gp mark 7}{}{(5.693,0.966)}
\gp3point{gp mark 7}{}{(6.375,0.949)}
\gp3point{gp mark 7}{}{(7.056,1.000)}
\gp3point{gp mark 7}{}{(7.737,1.305)}
\gp3point{gp mark 7}{}{(4.756,4.294)}
\gpcolor{rgb color={0.000,1.000,0.000}}
\draw[gp path] (0.925,4.075)--(1.606,3.552)--(2.287,3.175)--(2.969,2.829)--(3.650,2.498)%
  --(4.331,2.191)--(5.012,1.933)--(5.693,1.665)--(6.375,1.408)--(7.056,1.128)--(7.737,1.115);
\gp3point{gp mark 9}{}{(0.925,4.075)}
\gp3point{gp mark 9}{}{(1.606,3.552)}
\gp3point{gp mark 9}{}{(2.287,3.175)}
\gp3point{gp mark 9}{}{(2.969,2.829)}
\gp3point{gp mark 9}{}{(3.650,2.498)}
\gp3point{gp mark 9}{}{(4.331,2.191)}
\gp3point{gp mark 9}{}{(5.012,1.933)}
\gp3point{gp mark 9}{}{(5.693,1.665)}
\gp3point{gp mark 9}{}{(6.375,1.408)}
\gp3point{gp mark 9}{}{(7.056,1.128)}
\gp3point{gp mark 9}{}{(7.737,1.115)}
\gpcolor{rgb color={0.000,0.000,0.000}}
\draw[gp path] (0.925,4.075)--(0.987,4.040)--(1.049,4.005)--(1.111,3.970)--(1.173,3.935)%
  --(1.235,3.900)--(1.297,3.865)--(1.358,3.830)--(1.420,3.795)--(1.482,3.760)--(1.544,3.725)%
  --(1.606,3.690)--(1.668,3.655)--(1.730,3.620)--(1.792,3.585)--(1.854,3.550)--(1.916,3.515)%
  --(1.978,3.480)--(2.040,3.445)--(2.102,3.410)--(2.164,3.375)--(2.225,3.340)--(2.287,3.305)%
  --(2.349,3.270)--(2.411,3.235)--(2.473,3.200)--(2.535,3.165)--(2.597,3.130)--(2.659,3.095)%
  --(2.721,3.060)--(2.783,3.024)--(2.845,2.989)--(2.907,2.954)--(2.969,2.919)--(3.031,2.884)%
  --(3.092,2.849)--(3.154,2.814)--(3.216,2.779)--(3.278,2.744)--(3.340,2.709)--(3.402,2.674)%
  --(3.464,2.639)--(3.526,2.604)--(3.588,2.569)--(3.650,2.534)--(3.712,2.499)--(3.774,2.464)%
  --(3.836,2.429)--(3.897,2.394)--(3.959,2.359)--(4.021,2.324)--(4.083,2.289)--(4.145,2.254)%
  --(4.207,2.219)--(4.269,2.184)--(4.331,2.149)--(4.393,2.114)--(4.455,2.079)--(4.517,2.044)%
  --(4.579,2.009)--(4.641,1.974)--(4.703,1.939)--(4.764,1.904)--(4.826,1.869)--(4.888,1.834)%
  --(4.950,1.799)--(5.012,1.764)--(5.074,1.729)--(5.136,1.694)--(5.198,1.659)--(5.260,1.624)%
  --(5.322,1.589)--(5.384,1.554)--(5.446,1.519)--(5.508,1.484)--(5.570,1.449)--(5.631,1.413)%
  --(5.693,1.378)--(5.755,1.343)--(5.817,1.308)--(5.879,1.273)--(5.941,1.238)--(6.003,1.203)%
  --(6.065,1.168)--(6.127,1.133)--(6.189,1.098)--(6.251,1.063)--(6.313,1.028)--(6.375,0.993)%
  --(6.436,0.958)--(6.498,0.923)--(6.560,0.888)--(6.622,0.853)--(6.684,0.818)--(6.746,0.783)%
  --(6.808,0.748)--(6.870,0.713)--(6.908,0.691);
\gpcolor{color=gp lt color border}
\draw[gp path] (0.925,4.145)--(0.925,0.691)--(7.737,0.691)--(7.737,4.145)--cycle;
\gpdefrectangularnode{gp plot 1}{\pgfpoint{0.925cm}{0.691cm}}{\pgfpoint{7.737cm}{4.145cm}}
\end{tikzpicture}
   \label{UOL_curve_percent_L5}}
\subfloat[L6 from L3,$nbp_{max}=128$]{
\begin{tikzpicture}[gnuplot]
\tikzset{every node/.append style={scale=0.70}}
\path (0.000,0.000) rectangle (8.125,5.075);
\gpcolor{color=gp lt color border}
\gpsetlinetype{gp lt border}
\gpsetdashtype{gp dt solid}
\gpsetlinewidth{1.00}
\draw[gp path] (0.709,0.691)--(0.799,0.691);
\draw[gp path] (7.737,0.691)--(7.647,0.691);
\draw[gp path] (0.709,0.860)--(0.799,0.860);
\draw[gp path] (7.737,0.860)--(7.647,0.860);
\draw[gp path] (0.709,0.998)--(0.799,0.998);
\draw[gp path] (7.737,0.998)--(7.647,0.998);
\draw[gp path] (0.709,1.115)--(0.799,1.115);
\draw[gp path] (7.737,1.115)--(7.647,1.115);
\draw[gp path] (0.709,1.217)--(0.799,1.217);
\draw[gp path] (7.737,1.217)--(7.647,1.217);
\draw[gp path] (0.709,1.306)--(0.799,1.306);
\draw[gp path] (7.737,1.306)--(7.647,1.306);
\gpcolor{color=gp lt color axes}
\gpsetlinetype{gp lt axes}
\gpsetdashtype{gp dt axes}
\gpsetlinewidth{0.50}
\draw[gp path] (0.709,1.386)--(7.737,1.386);
\gpcolor{color=gp lt color border}
\gpsetlinetype{gp lt border}
\gpsetdashtype{gp dt solid}
\gpsetlinewidth{1.00}
\draw[gp path] (0.709,1.386)--(0.889,1.386);
\draw[gp path] (7.737,1.386)--(7.557,1.386);
\node[gp node right] at (0.580,1.386) {$10$};
\draw[gp path] (0.709,1.911)--(0.799,1.911);
\draw[gp path] (7.737,1.911)--(7.647,1.911);
\draw[gp path] (0.709,2.219)--(0.799,2.219);
\draw[gp path] (7.737,2.219)--(7.647,2.219);
\draw[gp path] (0.709,2.437)--(0.799,2.437);
\draw[gp path] (7.737,2.437)--(7.647,2.437);
\draw[gp path] (0.709,2.606)--(0.799,2.606);
\draw[gp path] (7.737,2.606)--(7.647,2.606);
\draw[gp path] (0.709,2.744)--(0.799,2.744);
\draw[gp path] (7.737,2.744)--(7.647,2.744);
\draw[gp path] (0.709,2.861)--(0.799,2.861);
\draw[gp path] (7.737,2.861)--(7.647,2.861);
\draw[gp path] (0.709,2.963)--(0.799,2.963);
\draw[gp path] (7.737,2.963)--(7.647,2.963);
\draw[gp path] (0.709,3.052)--(0.799,3.052);
\draw[gp path] (7.737,3.052)--(7.647,3.052);
\gpcolor{color=gp lt color axes}
\gpsetlinetype{gp lt axes}
\gpsetdashtype{gp dt axes}
\gpsetlinewidth{0.50}
\draw[gp path] (0.709,3.132)--(7.737,3.132);
\gpcolor{color=gp lt color border}
\gpsetlinetype{gp lt border}
\gpsetdashtype{gp dt solid}
\gpsetlinewidth{1.00}
\draw[gp path] (0.709,3.132)--(0.889,3.132);
\draw[gp path] (7.737,3.132)--(7.557,3.132);
\node[gp node right] at (0.580,3.132) {$100$};
\draw[gp path] (0.709,3.657)--(0.799,3.657);
\draw[gp path] (7.737,3.657)--(7.647,3.657);
\draw[gp path] (0.709,3.965)--(0.799,3.965);
\draw[gp path] (7.737,3.965)--(7.647,3.965);
\draw[gp path] (0.709,4.183)--(0.799,4.183);
\draw[gp path] (7.737,4.183)--(7.647,4.183);
\gpcolor{color=gp lt color axes}
\gpsetlinetype{gp lt axes}
\gpsetdashtype{gp dt axes}
\gpsetlinewidth{0.50}
\draw[gp path] (0.709,0.691)--(0.709,4.183);
\gpcolor{color=gp lt color border}
\gpsetlinetype{gp lt border}
\gpsetdashtype{gp dt solid}
\gpsetlinewidth{1.00}
\draw[gp path] (0.709,0.691)--(0.709,0.871);
\draw[gp path] (0.709,4.183)--(0.709,4.003);
\node[gp node center,font={\fontsize{10.0pt}{12.0pt}\selectfont}] at (0.709,0.475) {16};
\gpcolor{color=gp lt color axes}
\gpsetlinetype{gp lt axes}
\gpsetdashtype{gp dt axes}
\gpsetlinewidth{0.50}
\draw[gp path] (1.612,0.691)--(1.612,4.183);
\gpcolor{color=gp lt color border}
\gpsetlinetype{gp lt border}
\gpsetdashtype{gp dt solid}
\gpsetlinewidth{1.00}
\draw[gp path] (1.612,0.691)--(1.612,0.871);
\draw[gp path] (1.612,4.183)--(1.612,4.003);
\node[gp node center,font={\fontsize{10.0pt}{12.0pt}\selectfont}] at (1.612,0.475) {32};
\gpcolor{color=gp lt color axes}
\gpsetlinetype{gp lt axes}
\gpsetdashtype{gp dt axes}
\gpsetlinewidth{0.50}
\draw[gp path] (2.515,0.691)--(2.515,4.183);
\gpcolor{color=gp lt color border}
\gpsetlinetype{gp lt border}
\gpsetdashtype{gp dt solid}
\gpsetlinewidth{1.00}
\draw[gp path] (2.515,0.691)--(2.515,0.871);
\draw[gp path] (2.515,4.183)--(2.515,4.003);
\node[gp node center,font={\fontsize{10.0pt}{12.0pt}\selectfont}] at (2.515,0.475) {64};
\gpcolor{color=gp lt color axes}
\gpsetlinetype{gp lt axes}
\gpsetdashtype{gp dt axes}
\gpsetlinewidth{0.50}
\draw[gp path] (3.417,0.691)--(3.417,4.183);
\gpcolor{color=gp lt color border}
\gpsetlinetype{gp lt border}
\gpsetdashtype{gp dt solid}
\gpsetlinewidth{1.00}
\draw[gp path] (3.417,0.691)--(3.417,0.871);
\draw[gp path] (3.417,4.183)--(3.417,4.003);
\node[gp node center,font={\fontsize{10.0pt}{12.0pt}\selectfont}] at (3.417,0.475) {128};
\gpcolor{color=gp lt color axes}
\gpsetlinetype{gp lt axes}
\gpsetdashtype{gp dt axes}
\gpsetlinewidth{0.50}
\draw[gp path] (4.320,0.691)--(4.320,4.183);
\gpcolor{color=gp lt color border}
\gpsetlinetype{gp lt border}
\gpsetdashtype{gp dt solid}
\gpsetlinewidth{1.00}
\draw[gp path] (4.320,0.691)--(4.320,0.871);
\draw[gp path] (4.320,4.183)--(4.320,4.003);
\node[gp node center,font={\fontsize{10.0pt}{12.0pt}\selectfont}] at (4.320,0.475) {256};
\gpcolor{color=gp lt color axes}
\gpsetlinetype{gp lt axes}
\gpsetdashtype{gp dt axes}
\gpsetlinewidth{0.50}
\draw[gp path] (5.223,0.691)--(5.223,4.183);
\gpcolor{color=gp lt color border}
\gpsetlinetype{gp lt border}
\gpsetdashtype{gp dt solid}
\gpsetlinewidth{1.00}
\draw[gp path] (5.223,0.691)--(5.223,0.871);
\draw[gp path] (5.223,4.183)--(5.223,4.003);
\node[gp node center,font={\fontsize{10.0pt}{12.0pt}\selectfont}] at (5.223,0.475) {512};
\gpcolor{color=gp lt color axes}
\gpsetlinetype{gp lt axes}
\gpsetdashtype{gp dt axes}
\gpsetlinewidth{0.50}
\draw[gp path] (6.126,0.691)--(6.126,4.183);
\gpcolor{color=gp lt color border}
\gpsetlinetype{gp lt border}
\gpsetdashtype{gp dt solid}
\gpsetlinewidth{1.00}
\draw[gp path] (6.126,0.691)--(6.126,0.871);
\draw[gp path] (6.126,4.183)--(6.126,4.003);
\node[gp node center,font={\fontsize{10.0pt}{12.0pt}\selectfont}] at (6.126,0.475) {1024};
\gpcolor{color=gp lt color axes}
\gpsetlinetype{gp lt axes}
\gpsetdashtype{gp dt axes}
\gpsetlinewidth{0.50}
\draw[gp path] (7.029,0.691)--(7.029,4.183);
\gpcolor{color=gp lt color border}
\gpsetlinetype{gp lt border}
\gpsetdashtype{gp dt solid}
\gpsetlinewidth{1.00}
\draw[gp path] (7.029,0.691)--(7.029,0.871);
\draw[gp path] (7.029,4.183)--(7.029,4.003);
\node[gp node center,font={\fontsize{10.0pt}{12.0pt}\selectfont}] at (7.029,0.475) {2048};
\gpcolor{color=gp lt color axes}
\gpsetlinetype{gp lt axes}
\gpsetdashtype{gp dt axes}
\gpsetlinewidth{0.50}
\draw[gp path] (7.737,0.691)--(7.737,4.183);
\gpcolor{color=gp lt color border}
\gpsetlinetype{gp lt border}
\gpsetdashtype{gp dt solid}
\gpsetlinewidth{1.00}
\draw[gp path] (7.737,0.691)--(7.737,0.871);
\draw[gp path] (7.737,4.183)--(7.737,4.003);
\node[gp node center,font={\fontsize{10.0pt}{12.0pt}\selectfont}] at (7.737,0.475) {3528};
\draw[gp path] (0.709,4.183)--(0.709,0.691)--(7.737,0.691)--(7.737,4.183)--cycle;
\node[gp node center] at (4.223,0.105) {number of cores};
\gpcolor{rgb color={0.545,0.000,0.000}}
\draw[gp path] (0.709,3.599)--(1.612,3.158)--(2.515,2.987)--(3.417,2.724)--(4.320,2.684)%
  --(5.223,2.611)--(6.126,2.646)--(7.029,2.861)--(7.737,3.028);
\gpsetpointsize{4.00}
\gp3point{gp mark 6}{}{(0.709,3.599)}
\gp3point{gp mark 6}{}{(1.612,3.158)}
\gp3point{gp mark 6}{}{(2.515,2.987)}
\gp3point{gp mark 6}{}{(3.417,2.724)}
\gp3point{gp mark 6}{}{(4.320,2.684)}
\gp3point{gp mark 6}{}{(5.223,2.611)}
\gp3point{gp mark 6}{}{(6.126,2.646)}
\gp3point{gp mark 6}{}{(7.029,2.861)}
\gp3point{gp mark 6}{}{(7.737,3.028)}
\gpcolor{color=gp lt color border}
\node[gp node right] at (3.805,4.782) {ts loop fine-scale solv};
\gpcolor{rgb color={0.000,0.392,0.000}}
\draw[gp path] (3.934,4.782)--(4.630,4.782);
\draw[gp path] (0.709,4.093)--(1.612,3.518)--(2.515,3.072)--(3.417,2.674)--(4.320,2.216)%
  --(5.223,1.852)--(6.126,1.490)--(7.029,1.103)--(7.737,1.002);
\gp3point{gp mark 8}{}{(0.709,4.093)}
\gp3point{gp mark 8}{}{(1.612,3.518)}
\gp3point{gp mark 8}{}{(2.515,3.072)}
\gp3point{gp mark 8}{}{(3.417,2.674)}
\gp3point{gp mark 8}{}{(4.320,2.216)}
\gp3point{gp mark 8}{}{(5.223,1.852)}
\gp3point{gp mark 8}{}{(6.126,1.490)}
\gp3point{gp mark 8}{}{(7.029,1.103)}
\gp3point{gp mark 8}{}{(7.737,1.002)}
\gp3point{gp mark 8}{}{(4.282,4.782)}
\gpcolor{rgb color={1.000,0.000,0.000}}
\draw[gp path] (0.709,2.588)--(1.612,2.100)--(2.515,1.878)--(3.417,1.559)--(4.320,1.408)%
  --(5.223,1.299)--(6.126,1.219)--(7.029,1.247)--(7.737,1.491);
\gp3point{gp mark 7}{}{(0.709,2.588)}
\gp3point{gp mark 7}{}{(1.612,2.100)}
\gp3point{gp mark 7}{}{(2.515,1.878)}
\gp3point{gp mark 7}{}{(3.417,1.559)}
\gp3point{gp mark 7}{}{(4.320,1.408)}
\gp3point{gp mark 7}{}{(5.223,1.299)}
\gp3point{gp mark 7}{}{(6.126,1.219)}
\gp3point{gp mark 7}{}{(7.029,1.247)}
\gp3point{gp mark 7}{}{(7.737,1.491)}
\gpcolor{color=gp lt color border}
\node[gp node right] at (3.805,4.557) {tsi loop fine-scale solv};
\gpcolor{rgb color={0.000,1.000,0.000}}
\draw[gp path] (3.934,4.557)--(4.630,4.557);
\draw[gp path] (0.709,4.092)--(1.612,3.518)--(2.515,3.073)--(3.417,2.679)--(4.320,2.204)%
  --(5.223,1.840)--(6.126,1.462)--(7.029,1.059)--(7.737,0.763);
\gp3point{gp mark 9}{}{(0.709,4.092)}
\gp3point{gp mark 9}{}{(1.612,3.518)}
\gp3point{gp mark 9}{}{(2.515,3.073)}
\gp3point{gp mark 9}{}{(3.417,2.679)}
\gp3point{gp mark 9}{}{(4.320,2.204)}
\gp3point{gp mark 9}{}{(5.223,1.840)}
\gp3point{gp mark 9}{}{(6.126,1.462)}
\gp3point{gp mark 9}{}{(7.029,1.059)}
\gp3point{gp mark 9}{}{(7.737,0.763)}
\gp3point{gp mark 9}{}{(4.282,4.557)}
\gpcolor{color=gp lt color border}
\node[gp node right] at (3.805,4.332) {ideal};
\gpcolor{rgb color={0.000,0.000,0.000}}
\draw[gp path] (3.934,4.332)--(4.630,4.332);
\draw[gp path] (0.709,4.092)--(0.773,4.055)--(0.837,4.018)--(0.901,3.981)--(0.964,3.943)%
  --(1.028,3.906)--(1.092,3.869)--(1.156,3.832)--(1.220,3.795)--(1.284,3.758)--(1.347,3.720)%
  --(1.411,3.683)--(1.475,3.646)--(1.539,3.609)--(1.603,3.572)--(1.667,3.535)--(1.730,3.497)%
  --(1.794,3.460)--(1.858,3.423)--(1.922,3.386)--(1.986,3.349)--(2.050,3.312)--(2.113,3.275)%
  --(2.177,3.237)--(2.241,3.200)--(2.305,3.163)--(2.369,3.126)--(2.433,3.089)--(2.496,3.052)%
  --(2.560,3.014)--(2.624,2.977)--(2.688,2.940)--(2.752,2.903)--(2.816,2.866)--(2.879,2.829)%
  --(2.943,2.791)--(3.007,2.754)--(3.071,2.717)--(3.135,2.680)--(3.199,2.643)--(3.262,2.606)%
  --(3.326,2.568)--(3.390,2.531)--(3.454,2.494)--(3.518,2.457)--(3.582,2.420)--(3.645,2.383)%
  --(3.709,2.345)--(3.773,2.308)--(3.837,2.271)--(3.901,2.234)--(3.965,2.197)--(4.028,2.160)%
  --(4.092,2.122)--(4.156,2.085)--(4.220,2.048)--(4.284,2.011)--(4.348,1.974)--(4.411,1.937)%
  --(4.475,1.899)--(4.539,1.862)--(4.603,1.825)--(4.667,1.788)--(4.731,1.751)--(4.794,1.714)%
  --(4.858,1.676)--(4.922,1.639)--(4.986,1.602)--(5.050,1.565)--(5.114,1.528)--(5.177,1.491)%
  --(5.241,1.454)--(5.305,1.416)--(5.369,1.379)--(5.433,1.342)--(5.497,1.305)--(5.560,1.268)%
  --(5.624,1.231)--(5.688,1.193)--(5.752,1.156)--(5.816,1.119)--(5.880,1.082)--(5.943,1.045)%
  --(6.007,1.008)--(6.071,0.970)--(6.135,0.933)--(6.199,0.896)--(6.263,0.859)--(6.326,0.822)%
  --(6.390,0.785)--(6.454,0.747)--(6.518,0.710)--(6.550,0.691);
\gpcolor{color=gp lt color border}
\draw[gp path] (0.709,4.183)--(0.709,0.691)--(7.737,0.691)--(7.737,4.183)--cycle;
\gpdefrectangularnode{gp plot 1}{\pgfpoint{0.709cm}{0.691cm}}{\pgfpoint{7.737cm}{4.183cm}}
\end{tikzpicture}
\label{UOL_curve_percent_L6}}
	\caption{Cubic field: elapsed time in second of the \tS loop resolutions (ts and tsi versions) divided in two parts: the global-scale resolution and the local-scale resolutions. Log/log scale.\label{UOL_curve_percent}.}
\end{figure}

The curve "fr" in figure \ref{UOL_curve_times} for L2,L3,L4 and L5 shows that beyond a limit, the elapsed times increase when the number of used cores increases.
This justifies why, in the section \ref{gsolv}, the version \TSD takes at most $nbp_{max}$ processes to compute the problem at the global scale.
But this has an impact on scalability.
In the figure \ref{UOL_curve_percent}, for two configurations ("L5 from L2" and "L6 from L3"), the elapsed times of the scale loop resolutions are isolated and divided into global scale contribution(resolution $\dagger$ in algorithm \ref{TS_algebra}) and fine scale  contribution (MICRO-SCALE\_RESOLUTION call  of algorithm \ref{TS_algebra}).
It shows that the elapsed time of the global-scale resolution  corresponding to the \TSD version ("ts loop global-scale solv") decreases until $nbp_{max}$ is reached.
For a larger number of processes, the elapsed time remains stable and even increases.
This last point can be explained by the way the subset of $nbp_{max}$ processes is chosen.
With the current implementation, this subset follows a uniform distribution with respect to the mapping of all the processes.
This ensures that problem solving spans many nodes to ensure balanced memory consumption.
But this leads, with a high number of processes, to force almost all the communications of factoring and  forward/backward solving, to go through the network  between nodes that may be relatively far apart (in the sens of the network topology).
Another parameter explaining this increase, is related to the gathering of $\vm{A}_{gg}$ and $\vm{B}_g$, and the scattering of $\vm{U}_g$, from/to a large number of processes that add communication.
This poor performance for a number of processes higher than $nbp_{max}$  justifies the introduction of the \TSI version.
Its effect is first to reduce the global-scale consumption: the "tsi loop global-scale solv" is shifted down compared  to the "ts loop global-scale solv" (a facorization is replaced by some forward/backward resolution and sparse matrix multiplication per vector  at some iterations).
This is even more true in the case of  figure \ref{UOL_curve_percent_L6} where  the global scale level is larger  than in the case of figure \ref{UOL_curve_percent_L5}.  
The second effect of this version is to  maintain a decreasing slop of the "tsi loop global-scale solv" curves when using more than $nbp_{max}$ processes.
The iterative solver adds scalability.
With more than 1024 processes, the elapsed time increases again for "tsi loop global-scale solv".
The non-perfect scalability and this increase are again related  to the 
choice of the $nbp_{max}$ cores.
The preconditioner uses the factorization of the direct solver to perform  forward/backward resolution but only on $nbp_{max}$ processes spread over many nodes.
This limits the scalability because this part will not use more than $nbp_{max}$ processes and, as said above, network will be more stressed.
Note that with Mumps, we choose to have the solution and  the right-hand side  centralized on process 0 to simplify the exchange of these vectors between the  $nbp_{max}$ processes and all processes.
Perhaps using the ability of MUMPS to spread the solution vector and the right-hand side vector across processes would improve the situation.
But it's not clear because the vector gathering/scattering operation   should still be  performed, now directly, between  $nbp_{max}$ processes and all  processes, but not longer collectively.
In comparison, the elapsed time of the fine-scale resolution   (curves "ts loop fine-scale solv" and "tsi loop fine-scale solv"  in figure \ref{UOL_curve_percent} ) decrease steadily and confirms a fairly  good scalability of the fine-scale resolution. 

Now, regarding the  strong scaling efficiency \footnote{Strong scaling efficiency is the ratio of  speed-up to the number of processes used.
It is multiplied by 100 to obtain a \%. 100\% implies a perfect use of resource(cores)}  of the \tS solver, the figure  \ref{UOL_curve_seff2} shows a decrease in performance with increasing number of cores.
This deterioration is mainly related to the poor performance of the global-scale  as mentioned above in the analysis of the figure \ref{UOL_curve_percent}.
But it is not the only reason.
The strong scaling efficiency (not shown here) of "ts loop fine-scale solv" and "tsi loop fine-scale solv" in figure \ref{UOL_curve_percent}  shows that the performance remains above 50\% over a wide range, but decreases.
This decrease in fine-scale problem solving performance  can be understood by analyzing the evolution of the number of distributed patches per process.
The figure \ref{UOL_curve_patch} shows that the number of distributed patches per process reaches the number of  patches per process when the number of processes increases.
This implies that the amount of communication increases as  more patches need to communicate during their resolution.
This in itself has a significant impact on overall efficiency.
But this also increases the number of sequences that the algorithm \ref{sequencing_algo:general} must process.
And, because this algorithm uses a heuristic, perfect sequencing may be more difficult to achieve when almost all patches are distributed (the performance of the sequencing algorithm will be discussed in the section \ref{PO}).

\begin{figure}[h]
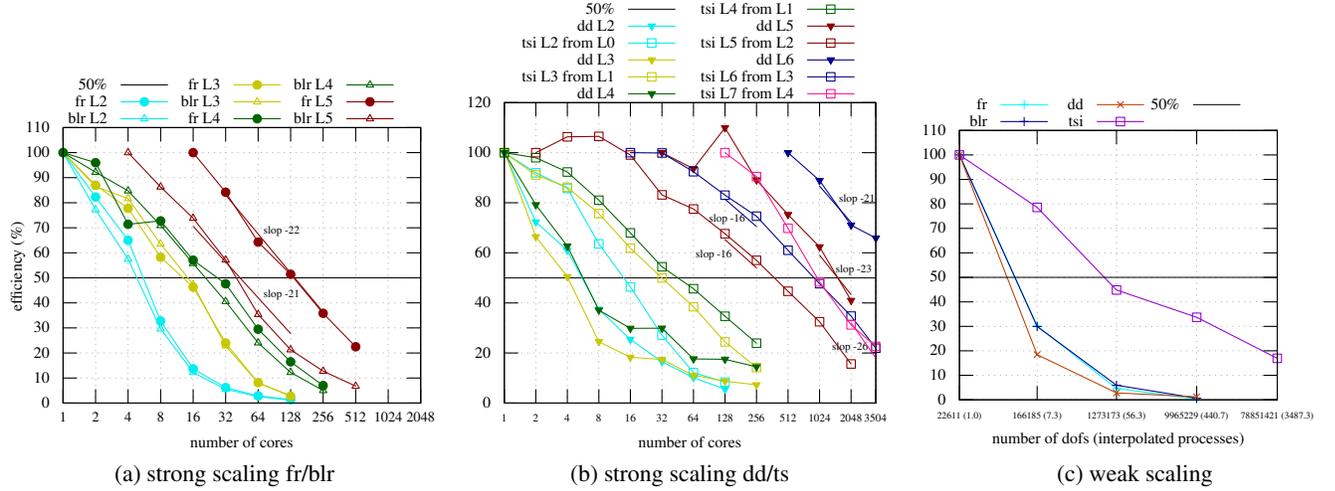

	\subfloat[strong scaling fr/blr]{

\label{UOL_curve_seff1}}
	\subfloat[strong scaling dd/ts]{
%
\label{UOL_curve_weff}}
	\caption{Cubic field: scaling efficiency. log/lin scale.\label{UOL_curve_eff}.}
\end{figure}

A direct comparison of the strong scaling efficiency with other solvers is possible for small problems (L2, L3, L4) where the data for one process is available for all solvers.
In figure  \ref{UOL_curve_seff1} and  \ref{UOL_curve_seff2}, the efficiency drops bellow 50\% when using more then 8 (L2) or 32 (L3 and L4)  processes with the \tS solver.
In comparison, this 50\% drop appears earlier for "dd" (4 processes) and "fr","blr' (4 processes for L2, 8 processes for L3 and 16 processes for L4).
For the other discretization level  (L5,L6 and L7), the data for "one process" are missing and the first available point is therefore arbitrarily considered as perfect (100\%).
In these cases, only the slopes can be compared.
But by using a linear fit of the curves  on their regularly decreasing part, the approximate slope can be extracted ("slop" in figure \ref{UOL_curve_seff1} and  \ref{UOL_curve_seff2}).
This shows that the \tS solver has lower slopes compared to "dd","fr" and "blr".
In conclusion, in all cases, the  \tS solver presents slightly better strong scaling performances than the other "solvers".

With respect to weak scaling efficiency\footnote{For a problem of size $N$, the weak scaling efficiency is the ratio of the time used by an application to solve a problem in single-process mode divided by the time used in multi-process mode (with $nbp$ processes) to solve a problem of size $nbp\times N$.
It is multiply by 100 to obtain a \%. 100\% implies perfect use of the resource(cores) } a full simulation campaign has not been conducted.
Interpolation  between the data in figures \ref{UOL_curve_seff1} and  \ref{UOL_curve_seff2} is used to produce the curves in figure \ref{UOL_curve_weff} where "tsi" uses a  different coarse-scale problem size.
The efficiency remain above 50\% up to $1.e^{+6}$ dofs and drops to 17\% for $7.8e^{+7}$ dofs, which in itself is positive (i.e. the efficiency does not drop below  5\% immediately).
Compared to other solvers, the \tS solver  offers a better efficiency.
Note that for the "dd" solver, the elapsed time reference for a single process is that of "fr" because the implementation does not handle the single domain case.
Thus, the weak efficiency in this case should certainly be moved up.

\vspace{1em}

\begin{figure}[h]
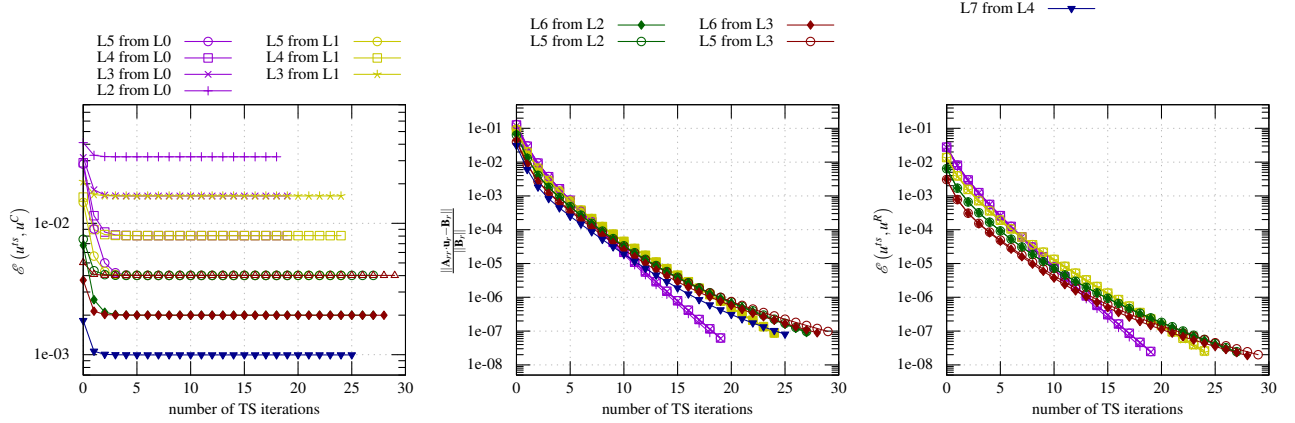


	\caption{Cubic field: errors \label{UOL_curve_error}.}
\end{figure}

 We now analyze the quality of the \tS solver using the following relative errors: 
 \begin{equation}
 	\mathcal{E}\left( \vm{u}^{ts},\vm{u}^C\right) =\frac{\left\| \vm{u}^{ts}-\vm{u}^C \right\|_{E_{\Omega}}}{\left\| \vm{u}^C \right\|_{E_{\Omega}}} ~\text{and}~ \mathcal{E}\left( \vm{u}^{ts},\vm{u}^R\right) =\frac{\left\| \vm{u}^{ts}-\vm{u}^R \right\|_{E_{\Omega}}}{\left\| \vm{u}^R \right\|_{E_{\Omega}}}
 	\label{TSvsA}
 \end{equation}
where $\vm{u}^C$ is the analytical solution given by \eqref{eqUC} (equivalent to the continuous solution of the section \ref{TS_reffield} ) and $\vm{u}^R$ is the reference solution given by solving of \eqref{r_sys}  with an auxiliary computation ("fr" or "dd" with  $\epsilon=10^{-13}$).
Figure \ref{UOL_curve_error} shows in the left graph $\mathcal{E}\left( \vm{u}^{ts},\vm{u}^C\right)$, in the right graph  $\mathcal{E}\left( \vm{u}^{ts},\vm{u}^R\right)$   and in the central graph the residual error $\frac{\left \| \vm{A}_{rr}\cdot \vm{u}_{r}-\vm{B}_{r} \right \|}{\left \|\vm{B}_{r} \right \|}$ ($resi$) provided by the algorithm \ref{TS_algebra_residual}.
These graphs show the evolution of the error with respect to the  \tS iterations.
The equation \eqref{equality} is verified numerically with these curves.
This explains  why, as $\mathcal{E}\left( \vm{u}^{ts},\vm{u}^R\right)$ tends towards zero in all cases, $\mathcal{E}\left( \vm{u}^{ts},\vm{u}^C\right)$ decreases to a constant value which corresponds to the discretization error   $\left\| \vm{u}^R-\vm{u}^C \right\|_{E_{\Omega}}/\left\| \vm{u}^C \right\|_{E_{\Omega}}$.
As the fine-scale discretization increases,  this constant naturally decrease.
But in all cases, only a few iterations make $\mathcal{E}\left( \vm{u}^{ts},\vm{u}^C\right)$ close to the discretization error.
This confirms the same observation made in the literature that the error on the boundary conditions of fine scale problems decreases rapidly with respect to $\left\| \vm{u}^R-\vm{u}^C \right\|_{E_{\Omega}}$.
The proposed \tS solver behaves in the same way, regardless of the number of processes used (there is no side effect introduced by parallelism).
In these graphs, we can see that the  $\mathcal{E}\left( \vm{u}^{ts},\vm{u}^R\right)$ error, corresponding mainly to the error on the boundary condition of fine-scale problems, evolves in a similar way to the residual error.
For $\epsilon=1.e^{-7}$, the $\mathcal{E}\left( \vm{u}^{ts},\vm{u}^R\right)$ error is bellow this threshold a few iterations before the residual error.
Thus, this last criterion is conservative: when it is respected, it ensures that $\mathcal{E}\left( \vm{u}^{ts},\vm{u}^R\right)$ is respected.
We can also observe that the curves, whether for the residual or for $\mathcal{E}\left( \vm{u}^{ts},\vm{u}^R\right)$ error, are almost the same for a given coarse-mesh discretization.
In particular L2, L3, L4 and L5 "from L0" are very close.
This is also the case for the curves "from L1" and "from L2".
For the "from L3", a small gap appears  between L5 and L6.
After a few iterations of the scale loop, the curves behave almost linearly (in linear/logarithmic scale) with  a softening of the slopes when at the global-scale the discretization grows.
These observations, which still need to be analyzed theoretically, show that  $\mathcal{E}\left( \vm{u}^{ts},\vm{u}^R\right)$and $resi$ are influenced  by the level of discretization on the global-scale.
But the next test case will already give us  a more refined interpretation of these observations.  

\vspace{1em}

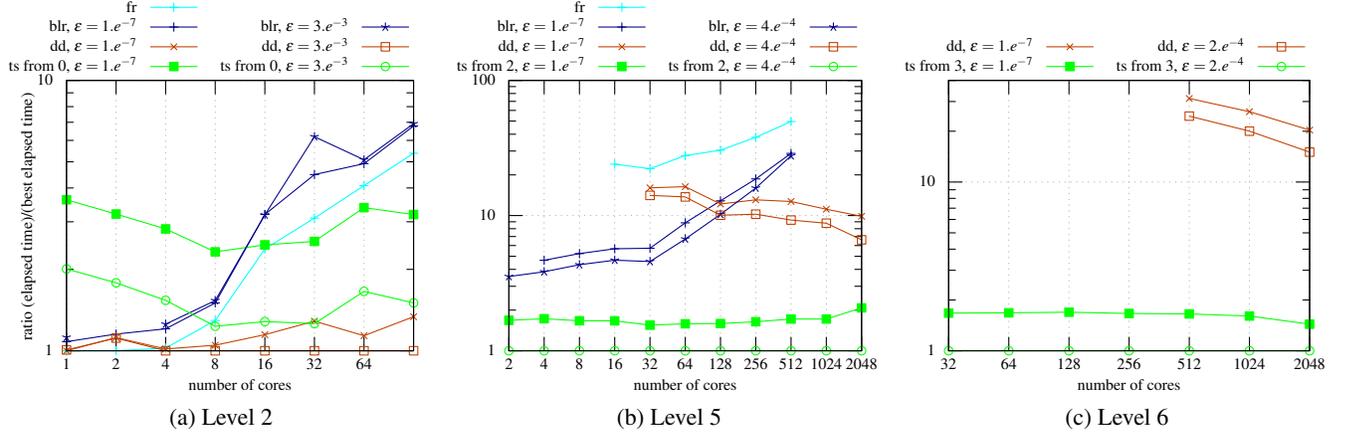
\begin{figure}[h]
 	\subfloat[Level 2]{
\begin{tikzpicture}[gnuplot]
\tikzset{every node/.append style={scale=0.60}}
\path (0.000,0.000) rectangle (5.625,4.375);
\gpcolor{color=gp lt color axes}
\gpsetlinetype{gp lt axes}
\gpsetdashtype{gp dt axes}
\gpsetlinewidth{0.50}
\draw[gp path] (0.680,0.592)--(5.294,0.592);
\gpcolor{color=gp lt color border}
\gpsetlinetype{gp lt border}
\gpsetdashtype{gp dt solid}
\gpsetlinewidth{1.00}
\draw[gp path] (0.680,0.592)--(0.860,0.592);
\draw[gp path] (5.294,0.592)--(5.114,0.592);
\node[gp node right] at (0.570,0.592) {$1$};
\draw[gp path] (0.680,1.675)--(0.770,1.675);
\draw[gp path] (5.294,1.675)--(5.204,1.675);
\draw[gp path] (0.680,2.308)--(0.770,2.308);
\draw[gp path] (5.294,2.308)--(5.204,2.308);
\draw[gp path] (0.680,2.758)--(0.770,2.758);
\draw[gp path] (5.294,2.758)--(5.204,2.758);
\draw[gp path] (0.680,3.106)--(0.770,3.106);
\draw[gp path] (5.294,3.106)--(5.204,3.106);
\draw[gp path] (0.680,3.391)--(0.770,3.391);
\draw[gp path] (5.294,3.391)--(5.204,3.391);
\draw[gp path] (0.680,3.632)--(0.770,3.632);
\draw[gp path] (5.294,3.632)--(5.204,3.632);
\draw[gp path] (0.680,3.840)--(0.770,3.840);
\draw[gp path] (5.294,3.840)--(5.204,3.840);
\draw[gp path] (0.680,4.024)--(0.770,4.024);
\draw[gp path] (5.294,4.024)--(5.204,4.024);
\gpcolor{color=gp lt color axes}
\gpsetlinetype{gp lt axes}
\gpsetdashtype{gp dt axes}
\gpsetlinewidth{0.50}
\draw[gp path] (0.680,4.189)--(5.294,4.189);
\gpcolor{color=gp lt color border}
\gpsetlinetype{gp lt border}
\gpsetdashtype{gp dt solid}
\gpsetlinewidth{1.00}
\draw[gp path] (0.680,4.189)--(0.860,4.189);
\draw[gp path] (5.294,4.189)--(5.114,4.189);
\node[gp node right] at (0.570,4.189) {$10$};
\gpcolor{color=gp lt color axes}
\gpsetlinetype{gp lt axes}
\gpsetdashtype{gp dt axes}
\gpsetlinewidth{0.50}
\draw[gp path] (0.680,0.592)--(0.680,4.189);
\gpcolor{color=gp lt color border}
\gpsetlinetype{gp lt border}
\gpsetdashtype{gp dt solid}
\gpsetlinewidth{1.00}
\draw[gp path] (0.680,0.592)--(0.680,0.772);
\draw[gp path] (0.680,4.189)--(0.680,4.009);
\node[gp node center] at (0.680,0.407) {1};
\gpcolor{color=gp lt color axes}
\gpsetlinetype{gp lt axes}
\gpsetdashtype{gp dt axes}
\gpsetlinewidth{0.50}
\draw[gp path] (1.339,0.592)--(1.339,4.189);
\gpcolor{color=gp lt color border}
\gpsetlinetype{gp lt border}
\gpsetdashtype{gp dt solid}
\gpsetlinewidth{1.00}
\draw[gp path] (1.339,0.592)--(1.339,0.772);
\draw[gp path] (1.339,4.189)--(1.339,4.009);
\node[gp node center] at (1.339,0.407) {2};
\gpcolor{color=gp lt color axes}
\gpsetlinetype{gp lt axes}
\gpsetdashtype{gp dt axes}
\gpsetlinewidth{0.50}
\draw[gp path] (1.998,0.592)--(1.998,4.189);
\gpcolor{color=gp lt color border}
\gpsetlinetype{gp lt border}
\gpsetdashtype{gp dt solid}
\gpsetlinewidth{1.00}
\draw[gp path] (1.998,0.592)--(1.998,0.772);
\draw[gp path] (1.998,4.189)--(1.998,4.009);
\node[gp node center] at (1.998,0.407) {4};
\gpcolor{color=gp lt color axes}
\gpsetlinetype{gp lt axes}
\gpsetdashtype{gp dt axes}
\gpsetlinewidth{0.50}
\draw[gp path] (2.657,0.592)--(2.657,4.189);
\gpcolor{color=gp lt color border}
\gpsetlinetype{gp lt border}
\gpsetdashtype{gp dt solid}
\gpsetlinewidth{1.00}
\draw[gp path] (2.657,0.592)--(2.657,0.772);
\draw[gp path] (2.657,4.189)--(2.657,4.009);
\node[gp node center] at (2.657,0.407) {8};
\gpcolor{color=gp lt color axes}
\gpsetlinetype{gp lt axes}
\gpsetdashtype{gp dt axes}
\gpsetlinewidth{0.50}
\draw[gp path] (3.317,0.592)--(3.317,4.189);
\gpcolor{color=gp lt color border}
\gpsetlinetype{gp lt border}
\gpsetdashtype{gp dt solid}
\gpsetlinewidth{1.00}
\draw[gp path] (3.317,0.592)--(3.317,0.772);
\draw[gp path] (3.317,4.189)--(3.317,4.009);
\node[gp node center] at (3.317,0.407) {16};
\gpcolor{color=gp lt color axes}
\gpsetlinetype{gp lt axes}
\gpsetdashtype{gp dt axes}
\gpsetlinewidth{0.50}
\draw[gp path] (3.976,0.592)--(3.976,4.189);
\gpcolor{color=gp lt color border}
\gpsetlinetype{gp lt border}
\gpsetdashtype{gp dt solid}
\gpsetlinewidth{1.00}
\draw[gp path] (3.976,0.592)--(3.976,0.772);
\draw[gp path] (3.976,4.189)--(3.976,4.009);
\node[gp node center] at (3.976,0.407) {32};
\gpcolor{color=gp lt color axes}
\gpsetlinetype{gp lt axes}
\gpsetdashtype{gp dt axes}
\gpsetlinewidth{0.50}
\draw[gp path] (4.635,0.592)--(4.635,4.189);
\gpcolor{color=gp lt color border}
\gpsetlinetype{gp lt border}
\gpsetdashtype{gp dt solid}
\gpsetlinewidth{1.00}
\draw[gp path] (4.635,0.592)--(4.635,0.772);
\draw[gp path] (4.635,4.189)--(4.635,4.009);
\node[gp node center] at (4.635,0.407) {64};
\draw[gp path] (0.680,4.189)--(0.680,0.592)--(5.294,0.592)--(5.294,4.189)--cycle;
\node[gp node center,rotate=-270] at (0.175,2.390) {ratio (elapsed time)/(best elapsed time)};
\node[gp node center] at (2.987,0.130) {number of cores};
\node[gp node right] at (1.687,5.161) {fr};
\gpcolor{rgb color={0.000,1.000,1.000}}
\draw[gp path] (1.797,5.161)--(2.417,5.161);
\draw[gp path] (0.680,0.592)--(1.339,0.592)--(1.998,0.627)--(2.657,0.995)--(3.317,1.947)%
  --(3.976,2.357)--(4.635,2.789)--(5.294,3.220);
\gpsetpointsize{4.00}
\gp3point{gp mark 1}{}{(0.680,0.592)}
\gp3point{gp mark 1}{}{(1.339,0.592)}
\gp3point{gp mark 1}{}{(1.998,0.627)}
\gp3point{gp mark 1}{}{(2.657,0.995)}
\gp3point{gp mark 1}{}{(3.317,1.947)}
\gp3point{gp mark 1}{}{(3.976,2.357)}
\gp3point{gp mark 1}{}{(4.635,2.789)}
\gp3point{gp mark 1}{}{(5.294,3.220)}
\gp3point{gp mark 1}{}{(2.107,5.161)}
\gpcolor{color=gp lt color border}
\node[gp node right] at (1.687,4.899) {blr, $\epsilon =1.e^{-7}$};
\gpcolor{rgb color={0.000,0.000,0.545}}
\draw[gp path] (1.797,4.899)--(2.417,4.899);
\draw[gp path] (0.680,0.713)--(1.339,0.814)--(1.998,0.884)--(2.657,1.228)--(3.317,2.405)%
  --(3.976,2.936)--(4.635,3.081)--(5.294,3.582);
\gp3point{gp mark 1}{}{(0.680,0.713)}
\gp3point{gp mark 1}{}{(1.339,0.814)}
\gp3point{gp mark 1}{}{(1.998,0.884)}
\gp3point{gp mark 1}{}{(2.657,1.228)}
\gp3point{gp mark 1}{}{(3.317,2.405)}
\gp3point{gp mark 1}{}{(3.976,2.936)}
\gp3point{gp mark 1}{}{(4.635,3.081)}
\gp3point{gp mark 1}{}{(5.294,3.582)}
\gp3point{gp mark 1}{}{(2.107,4.899)}
\gpcolor{color=gp lt color border}
\node[gp node right] at (1.687,4.636) {dd, $\epsilon=1.e^{-7}$};
\gpcolor{rgb color={0.753,0.251,0.000}}
\draw[gp path] (1.797,4.636)--(2.417,4.636);
\draw[gp path] (0.680,0.592)--(1.339,0.767)--(1.998,0.616)--(2.657,0.666)--(3.317,0.806)%
  --(3.976,0.986)--(4.635,0.792)--(5.294,1.045);
\gp3point{gp mark 2}{}{(0.680,0.592)}
\gp3point{gp mark 2}{}{(1.339,0.767)}
\gp3point{gp mark 2}{}{(1.998,0.616)}
\gp3point{gp mark 2}{}{(2.657,0.666)}
\gp3point{gp mark 2}{}{(3.317,0.806)}
\gp3point{gp mark 2}{}{(3.976,0.986)}
\gp3point{gp mark 2}{}{(4.635,0.792)}
\gp3point{gp mark 2}{}{(5.294,1.045)}
\gp3point{gp mark 2}{}{(2.107,4.636)}
\gpcolor{color=gp lt color border}
\node[gp node right] at (1.687,4.374) {ts from 0, $\epsilon=1.e^{-7}$};
\gpcolor{rgb color={0.000,1.000,0.000}}
\draw[gp path] (1.797,4.374)--(2.417,4.374);
\draw[gp path] (0.680,2.601)--(1.339,2.411)--(1.998,2.213)--(2.657,1.908)--(3.317,2.001)%
  --(3.976,2.047)--(4.635,2.496)--(5.294,2.405);
\gp3point{gp mark 5}{}{(0.680,2.601)}
\gp3point{gp mark 5}{}{(1.339,2.411)}
\gp3point{gp mark 5}{}{(1.998,2.213)}
\gp3point{gp mark 5}{}{(2.657,1.908)}
\gp3point{gp mark 5}{}{(3.317,2.001)}
\gp3point{gp mark 5}{}{(3.976,2.047)}
\gp3point{gp mark 5}{}{(4.635,2.496)}
\gp3point{gp mark 5}{}{(5.294,2.405)}
\gp3point{gp mark 5}{}{(2.107,4.374)}
\gpcolor{color=gp lt color border}
\node[gp node right] at (4.499,4.899) {blr, $\epsilon=3.e^{-3}$};
\gpcolor{rgb color={0.000,0.000,0.545}}
\draw[gp path] (4.609,4.899)--(5.229,4.899);
\draw[gp path] (1.998,0.946)--(2.657,1.260)--(3.317,2.400)--(3.976,3.440)--(4.635,3.133)%
  --(5.294,3.613);
\gp3point{gp mark 3}{}{(0.680,0.771)}
\gp3point{gp mark 3}{}{(1.998,0.946)}
\gp3point{gp mark 3}{}{(2.657,1.260)}
\gp3point{gp mark 3}{}{(3.317,2.400)}
\gp3point{gp mark 3}{}{(3.976,3.440)}
\gp3point{gp mark 3}{}{(4.635,3.133)}
\gp3point{gp mark 3}{}{(5.294,3.613)}
\gp3point{gp mark 3}{}{(4.919,4.899)}
\gpcolor{color=gp lt color border}
\node[gp node right] at (4.499,4.636) {dd, $\epsilon=3.e^{-3}$};
\gpcolor{rgb color={0.753,0.251,0.000}}
\draw[gp path] (4.609,4.636)--(5.229,4.636);
\draw[gp path] (0.680,0.603)--(1.339,0.761)--(1.998,0.592)--(2.657,0.592)--(3.317,0.592)%
  --(3.976,0.592)--(4.635,0.592)--(5.294,0.592);
\gp3point{gp mark 4}{}{(0.680,0.603)}
\gp3point{gp mark 4}{}{(1.339,0.761)}
\gp3point{gp mark 4}{}{(1.998,0.592)}
\gp3point{gp mark 4}{}{(2.657,0.592)}
\gp3point{gp mark 4}{}{(3.317,0.592)}
\gp3point{gp mark 4}{}{(3.976,0.592)}
\gp3point{gp mark 4}{}{(4.635,0.592)}
\gp3point{gp mark 4}{}{(5.294,0.592)}
\gp3point{gp mark 4}{}{(4.919,4.636)}
\gpcolor{color=gp lt color border}
\node[gp node right] at (4.499,4.374) {ts from 0, $\epsilon=3.e^{-3}$};
\gpcolor{rgb color={0.000,1.000,0.000}}
\draw[gp path] (4.609,4.374)--(5.229,4.374);
\draw[gp path] (0.680,1.681)--(1.339,1.496)--(1.998,1.263)--(2.657,0.918)--(3.317,0.980)%
  --(3.976,0.953)--(4.635,1.382)--(5.294,1.229);
\gp3point{gp mark 6}{}{(0.680,1.681)}
\gp3point{gp mark 6}{}{(1.339,1.496)}
\gp3point{gp mark 6}{}{(1.998,1.263)}
\gp3point{gp mark 6}{}{(2.657,0.918)}
\gp3point{gp mark 6}{}{(3.317,0.980)}
\gp3point{gp mark 6}{}{(3.976,0.953)}
\gp3point{gp mark 6}{}{(4.635,1.382)}
\gp3point{gp mark 6}{}{(5.294,1.229)}
\gp3point{gp mark 6}{}{(4.919,4.374)}
\gpcolor{color=gp lt color border}
\draw[gp path] (0.680,4.189)--(0.680,0.592)--(5.294,0.592)--(5.294,4.189)--cycle;
\gpdefrectangularnode{gp plot 1}{\pgfpoint{0.680cm}{0.592cm}}{\pgfpoint{5.294cm}{4.189cm}}
\end{tikzpicture}
\label{UOL_precision_impact_L2}}
 	\subfloat[Level 5]{
\begin{tikzpicture}[gnuplot]
\tikzset{every node/.append style={scale=0.60}}
\path (0.000,0.000) rectangle (5.625,4.375);
\gpcolor{color=gp lt color axes}
\gpsetlinetype{gp lt axes}
\gpsetdashtype{gp dt axes}
\gpsetlinewidth{0.50}
\draw[gp path] (0.605,0.592)--(5.294,0.592);
\gpcolor{color=gp lt color border}
\gpsetlinetype{gp lt border}
\gpsetdashtype{gp dt solid}
\gpsetlinewidth{1.00}
\draw[gp path] (0.605,0.592)--(0.785,0.592);
\draw[gp path] (5.294,0.592)--(5.114,0.592);
\node[gp node right] at (0.495,0.592) {$1$};
\draw[gp path] (0.605,1.133)--(0.695,1.133);
\draw[gp path] (5.294,1.133)--(5.204,1.133);
\draw[gp path] (0.605,1.450)--(0.695,1.450);
\draw[gp path] (5.294,1.450)--(5.204,1.450);
\draw[gp path] (0.605,1.675)--(0.695,1.675);
\draw[gp path] (5.294,1.675)--(5.204,1.675);
\draw[gp path] (0.605,1.849)--(0.695,1.849);
\draw[gp path] (5.294,1.849)--(5.204,1.849);
\draw[gp path] (0.605,1.992)--(0.695,1.992);
\draw[gp path] (5.294,1.992)--(5.204,1.992);
\draw[gp path] (0.605,2.112)--(0.695,2.112);
\draw[gp path] (5.294,2.112)--(5.204,2.112);
\draw[gp path] (0.605,2.216)--(0.695,2.216);
\draw[gp path] (5.294,2.216)--(5.204,2.216);
\draw[gp path] (0.605,2.308)--(0.695,2.308);
\draw[gp path] (5.294,2.308)--(5.204,2.308);
\gpcolor{color=gp lt color axes}
\gpsetlinetype{gp lt axes}
\gpsetdashtype{gp dt axes}
\gpsetlinewidth{0.50}
\draw[gp path] (0.605,2.391)--(5.294,2.391);
\gpcolor{color=gp lt color border}
\gpsetlinetype{gp lt border}
\gpsetdashtype{gp dt solid}
\gpsetlinewidth{1.00}
\draw[gp path] (0.605,2.391)--(0.785,2.391);
\draw[gp path] (5.294,2.391)--(5.114,2.391);
\node[gp node right] at (0.495,2.391) {$10$};
\draw[gp path] (0.605,2.932)--(0.695,2.932);
\draw[gp path] (5.294,2.932)--(5.204,2.932);
\draw[gp path] (0.605,3.249)--(0.695,3.249);
\draw[gp path] (5.294,3.249)--(5.204,3.249);
\draw[gp path] (0.605,3.473)--(0.695,3.473);
\draw[gp path] (5.294,3.473)--(5.204,3.473);
\draw[gp path] (0.605,3.648)--(0.695,3.648);
\draw[gp path] (5.294,3.648)--(5.204,3.648);
\draw[gp path] (0.605,3.790)--(0.695,3.790);
\draw[gp path] (5.294,3.790)--(5.204,3.790);
\draw[gp path] (0.605,3.910)--(0.695,3.910);
\draw[gp path] (5.294,3.910)--(5.204,3.910);
\draw[gp path] (0.605,4.015)--(0.695,4.015);
\draw[gp path] (5.294,4.015)--(5.204,4.015);
\draw[gp path] (0.605,4.107)--(0.695,4.107);
\draw[gp path] (5.294,4.107)--(5.204,4.107);
\gpcolor{color=gp lt color axes}
\gpsetlinetype{gp lt axes}
\gpsetdashtype{gp dt axes}
\gpsetlinewidth{0.50}
\draw[gp path] (0.605,4.189)--(5.294,4.189);
\gpcolor{color=gp lt color border}
\gpsetlinetype{gp lt border}
\gpsetdashtype{gp dt solid}
\gpsetlinewidth{1.00}
\draw[gp path] (0.605,4.189)--(0.785,4.189);
\draw[gp path] (5.294,4.189)--(5.114,4.189);
\node[gp node right] at (0.495,4.189) {$100$};
\gpcolor{color=gp lt color axes}
\gpsetlinetype{gp lt axes}
\gpsetdashtype{gp dt axes}
\gpsetlinewidth{0.50}
\draw[gp path] (0.605,0.592)--(0.605,4.189);
\gpcolor{color=gp lt color border}
\gpsetlinetype{gp lt border}
\gpsetdashtype{gp dt solid}
\gpsetlinewidth{1.00}
\draw[gp path] (0.605,0.592)--(0.605,0.772);
\draw[gp path] (0.605,4.189)--(0.605,4.009);
\node[gp node center] at (0.605,0.407) {2};
\gpcolor{color=gp lt color axes}
\gpsetlinetype{gp lt axes}
\gpsetdashtype{gp dt axes}
\gpsetlinewidth{0.50}
\draw[gp path] (1.074,0.592)--(1.074,4.189);
\gpcolor{color=gp lt color border}
\gpsetlinetype{gp lt border}
\gpsetdashtype{gp dt solid}
\gpsetlinewidth{1.00}
\draw[gp path] (1.074,0.592)--(1.074,0.772);
\draw[gp path] (1.074,4.189)--(1.074,4.009);
\node[gp node center] at (1.074,0.407) {4};
\gpcolor{color=gp lt color axes}
\gpsetlinetype{gp lt axes}
\gpsetdashtype{gp dt axes}
\gpsetlinewidth{0.50}
\draw[gp path] (1.543,0.592)--(1.543,4.189);
\gpcolor{color=gp lt color border}
\gpsetlinetype{gp lt border}
\gpsetdashtype{gp dt solid}
\gpsetlinewidth{1.00}
\draw[gp path] (1.543,0.592)--(1.543,0.772);
\draw[gp path] (1.543,4.189)--(1.543,4.009);
\node[gp node center] at (1.543,0.407) {8};
\gpcolor{color=gp lt color axes}
\gpsetlinetype{gp lt axes}
\gpsetdashtype{gp dt axes}
\gpsetlinewidth{0.50}
\draw[gp path] (2.012,0.592)--(2.012,4.189);
\gpcolor{color=gp lt color border}
\gpsetlinetype{gp lt border}
\gpsetdashtype{gp dt solid}
\gpsetlinewidth{1.00}
\draw[gp path] (2.012,0.592)--(2.012,0.772);
\draw[gp path] (2.012,4.189)--(2.012,4.009);
\node[gp node center] at (2.012,0.407) {16};
\gpcolor{color=gp lt color axes}
\gpsetlinetype{gp lt axes}
\gpsetdashtype{gp dt axes}
\gpsetlinewidth{0.50}
\draw[gp path] (2.481,0.592)--(2.481,4.189);
\gpcolor{color=gp lt color border}
\gpsetlinetype{gp lt border}
\gpsetdashtype{gp dt solid}
\gpsetlinewidth{1.00}
\draw[gp path] (2.481,0.592)--(2.481,0.772);
\draw[gp path] (2.481,4.189)--(2.481,4.009);
\node[gp node center] at (2.481,0.407) {32};
\gpcolor{color=gp lt color axes}
\gpsetlinetype{gp lt axes}
\gpsetdashtype{gp dt axes}
\gpsetlinewidth{0.50}
\draw[gp path] (2.950,0.592)--(2.950,4.189);
\gpcolor{color=gp lt color border}
\gpsetlinetype{gp lt border}
\gpsetdashtype{gp dt solid}
\gpsetlinewidth{1.00}
\draw[gp path] (2.950,0.592)--(2.950,0.772);
\draw[gp path] (2.950,4.189)--(2.950,4.009);
\node[gp node center] at (2.950,0.407) {64};
\gpcolor{color=gp lt color axes}
\gpsetlinetype{gp lt axes}
\gpsetdashtype{gp dt axes}
\gpsetlinewidth{0.50}
\draw[gp path] (3.418,0.592)--(3.418,4.189);
\gpcolor{color=gp lt color border}
\gpsetlinetype{gp lt border}
\gpsetdashtype{gp dt solid}
\gpsetlinewidth{1.00}
\draw[gp path] (3.418,0.592)--(3.418,0.772);
\draw[gp path] (3.418,4.189)--(3.418,4.009);
\node[gp node center] at (3.418,0.407) {128};
\gpcolor{color=gp lt color axes}
\gpsetlinetype{gp lt axes}
\gpsetdashtype{gp dt axes}
\gpsetlinewidth{0.50}
\draw[gp path] (3.887,0.592)--(3.887,4.189);
\gpcolor{color=gp lt color border}
\gpsetlinetype{gp lt border}
\gpsetdashtype{gp dt solid}
\gpsetlinewidth{1.00}
\draw[gp path] (3.887,0.592)--(3.887,0.772);
\draw[gp path] (3.887,4.189)--(3.887,4.009);
\node[gp node center] at (3.887,0.407) {256};
\gpcolor{color=gp lt color axes}
\gpsetlinetype{gp lt axes}
\gpsetdashtype{gp dt axes}
\gpsetlinewidth{0.50}
\draw[gp path] (4.356,0.592)--(4.356,4.189);
\gpcolor{color=gp lt color border}
\gpsetlinetype{gp lt border}
\gpsetdashtype{gp dt solid}
\gpsetlinewidth{1.00}
\draw[gp path] (4.356,0.592)--(4.356,0.772);
\draw[gp path] (4.356,4.189)--(4.356,4.009);
\node[gp node center] at (4.356,0.407) {512};
\gpcolor{color=gp lt color axes}
\gpsetlinetype{gp lt axes}
\gpsetdashtype{gp dt axes}
\gpsetlinewidth{0.50}
\draw[gp path] (4.825,0.592)--(4.825,4.189);
\gpcolor{color=gp lt color border}
\gpsetlinetype{gp lt border}
\gpsetdashtype{gp dt solid}
\gpsetlinewidth{1.00}
\draw[gp path] (4.825,0.592)--(4.825,0.772);
\draw[gp path] (4.825,4.189)--(4.825,4.009);
\node[gp node center] at (4.825,0.407) {1024};
\gpcolor{color=gp lt color axes}
\gpsetlinetype{gp lt axes}
\gpsetdashtype{gp dt axes}
\gpsetlinewidth{0.50}
\draw[gp path] (5.294,0.592)--(5.294,4.189);
\gpcolor{color=gp lt color border}
\gpsetlinetype{gp lt border}
\gpsetdashtype{gp dt solid}
\gpsetlinewidth{1.00}
\draw[gp path] (5.294,0.592)--(5.294,0.772);
\draw[gp path] (5.294,4.189)--(5.294,4.009);
\node[gp node center] at (5.294,0.407) {2048};
\draw[gp path] (0.605,4.189)--(0.605,0.592)--(5.294,0.592)--(5.294,4.189)--cycle;
\node[gp node center] at (2.949,0.130) {number of cores};
\node[gp node right] at (1.687,5.161) {fr};
\gpcolor{rgb color={0.000,1.000,1.000}}
\draw[gp path] (1.797,5.161)--(2.417,5.161);
\draw[gp path] (2.012,3.074)--(2.481,3.014)--(2.950,3.189)--(3.418,3.260)--(3.887,3.432)%
  --(4.356,3.642);
\gpsetpointsize{4.00}
\gp3point{gp mark 1}{}{(2.012,3.074)}
\gp3point{gp mark 1}{}{(2.481,3.014)}
\gp3point{gp mark 1}{}{(2.950,3.189)}
\gp3point{gp mark 1}{}{(3.418,3.260)}
\gp3point{gp mark 1}{}{(3.887,3.432)}
\gp3point{gp mark 1}{}{(4.356,3.642)}
\gp3point{gp mark 1}{}{(2.107,5.161)}
\gpcolor{color=gp lt color border}
\node[gp node right] at (1.687,4.899) {blr, $\epsilon =1.e^{-7}$};
\gpcolor{rgb color={0.000,0.000,0.545}}
\draw[gp path] (1.797,4.899)--(2.417,4.899);
\draw[gp path] (1.074,1.796)--(1.543,1.885)--(2.012,1.949)--(2.481,1.956)--(2.950,2.294)%
  --(3.418,2.587)--(3.887,2.878)--(4.356,3.218);
\gp3point{gp mark 1}{}{(1.074,1.796)}
\gp3point{gp mark 1}{}{(1.543,1.885)}
\gp3point{gp mark 1}{}{(2.012,1.949)}
\gp3point{gp mark 1}{}{(2.481,1.956)}
\gp3point{gp mark 1}{}{(2.950,2.294)}
\gp3point{gp mark 1}{}{(3.418,2.587)}
\gp3point{gp mark 1}{}{(3.887,2.878)}
\gp3point{gp mark 1}{}{(4.356,3.218)}
\gp3point{gp mark 1}{}{(2.107,4.899)}
\gpcolor{color=gp lt color border}
\node[gp node right] at (1.687,4.636) {dd, $\epsilon=1.e^{-7}$};
\gpcolor{rgb color={0.753,0.251,0.000}}
\draw[gp path] (1.797,4.636)--(2.417,4.636);
\draw[gp path] (2.481,2.760)--(2.950,2.775)--(3.418,2.545)--(3.887,2.600)--(4.356,2.578)%
  --(4.825,2.475)--(5.294,2.381);
\gp3point{gp mark 2}{}{(2.481,2.760)}
\gp3point{gp mark 2}{}{(2.950,2.775)}
\gp3point{gp mark 2}{}{(3.418,2.545)}
\gp3point{gp mark 2}{}{(3.887,2.600)}
\gp3point{gp mark 2}{}{(4.356,2.578)}
\gp3point{gp mark 2}{}{(4.825,2.475)}
\gp3point{gp mark 2}{}{(5.294,2.381)}
\gp3point{gp mark 2}{}{(2.107,4.636)}
\gpcolor{color=gp lt color border}
\node[gp node right] at (1.687,4.374) {ts from 2, $\epsilon=1.e^{-7}$};
\gpcolor{rgb color={0.000,1.000,0.000}}
\draw[gp path] (1.797,4.374)--(2.417,4.374);
\draw[gp path] (0.605,0.997)--(1.074,1.018)--(1.543,0.991)--(2.012,0.990)--(2.481,0.933)%
  --(2.950,0.951)--(3.418,0.954)--(3.887,0.978)--(4.356,1.014)--(4.825,1.014)--(5.294,1.163);
\gp3point{gp mark 5}{}{(0.605,0.997)}
\gp3point{gp mark 5}{}{(1.074,1.018)}
\gp3point{gp mark 5}{}{(1.543,0.991)}
\gp3point{gp mark 5}{}{(2.012,0.990)}
\gp3point{gp mark 5}{}{(2.481,0.933)}
\gp3point{gp mark 5}{}{(2.950,0.951)}
\gp3point{gp mark 5}{}{(3.418,0.954)}
\gp3point{gp mark 5}{}{(3.887,0.978)}
\gp3point{gp mark 5}{}{(4.356,1.014)}
\gp3point{gp mark 5}{}{(4.825,1.014)}
\gp3point{gp mark 5}{}{(5.294,1.163)}
\gp3point{gp mark 5}{}{(2.107,4.374)}
\gpcolor{color=gp lt color border}
\node[gp node right] at (4.499,4.899) {blr, $\epsilon=4.e^{-4}$};
\gpcolor{rgb color={0.000,0.000,0.545}}
\draw[gp path] (4.609,4.899)--(5.229,4.899);
\draw[gp path] (0.605,1.580)--(1.074,1.643)--(1.543,1.737)--(2.012,1.797)--(2.481,1.779)%
  --(2.950,2.079)--(3.418,2.407)--(3.887,2.757)--(4.356,3.186);
\gp3point{gp mark 3}{}{(0.605,1.580)}
\gp3point{gp mark 3}{}{(1.074,1.643)}
\gp3point{gp mark 3}{}{(1.543,1.737)}
\gp3point{gp mark 3}{}{(2.012,1.797)}
\gp3point{gp mark 3}{}{(2.481,1.779)}
\gp3point{gp mark 3}{}{(2.950,2.079)}
\gp3point{gp mark 3}{}{(3.418,2.407)}
\gp3point{gp mark 3}{}{(3.887,2.757)}
\gp3point{gp mark 3}{}{(4.356,3.186)}
\gp3point{gp mark 3}{}{(4.919,4.899)}
\gpcolor{color=gp lt color border}
\node[gp node right] at (4.499,4.636) {dd, $\epsilon=4.e^{-4}$};
\gpcolor{rgb color={0.753,0.251,0.000}}
\draw[gp path] (4.609,4.636)--(5.229,4.636);
\draw[gp path] (2.481,2.659)--(2.950,2.639)--(3.418,2.395)--(3.887,2.409)--(4.356,2.331)%
  --(4.825,2.289)--(5.294,2.068);
\gp3point{gp mark 4}{}{(2.481,2.659)}
\gp3point{gp mark 4}{}{(2.950,2.639)}
\gp3point{gp mark 4}{}{(3.418,2.395)}
\gp3point{gp mark 4}{}{(3.887,2.409)}
\gp3point{gp mark 4}{}{(4.356,2.331)}
\gp3point{gp mark 4}{}{(4.825,2.289)}
\gp3point{gp mark 4}{}{(5.294,2.068)}
\gp3point{gp mark 4}{}{(4.919,4.636)}
\gpcolor{color=gp lt color border}
\node[gp node right] at (4.499,4.374) {ts from 2, $\epsilon=4.e^{-4}$};
\gpcolor{rgb color={0.000,1.000,0.000}}
\draw[gp path] (4.609,4.374)--(5.229,4.374);
\draw[gp path] (0.605,0.592)--(1.074,0.592)--(1.543,0.592)--(2.012,0.592)--(2.481,0.592)%
  --(2.950,0.592)--(3.418,0.592)--(3.887,0.592)--(4.356,0.592)--(4.825,0.592)--(5.294,0.592);
\gp3point{gp mark 6}{}{(0.605,0.592)}
\gp3point{gp mark 6}{}{(1.074,0.592)}
\gp3point{gp mark 6}{}{(1.543,0.592)}
\gp3point{gp mark 6}{}{(2.012,0.592)}
\gp3point{gp mark 6}{}{(2.481,0.592)}
\gp3point{gp mark 6}{}{(2.950,0.592)}
\gp3point{gp mark 6}{}{(3.418,0.592)}
\gp3point{gp mark 6}{}{(3.887,0.592)}
\gp3point{gp mark 6}{}{(4.356,0.592)}
\gp3point{gp mark 6}{}{(4.825,0.592)}
\gp3point{gp mark 6}{}{(5.294,0.592)}
\gp3point{gp mark 6}{}{(4.919,4.374)}
\gpcolor{color=gp lt color border}
\draw[gp path] (0.605,4.189)--(0.605,0.592)--(5.294,0.592)--(5.294,4.189)--cycle;
\gpdefrectangularnode{gp plot 1}{\pgfpoint{0.605cm}{0.592cm}}{\pgfpoint{5.294cm}{4.189cm}}
\end{tikzpicture}
\label{UOL_precision_impact_L5}}
 	\subfloat[Level 6]{
\begin{tikzpicture}[gnuplot]
\tikzset{every node/.append style={scale=0.60}}
\path (0.000,0.000) rectangle (5.625,4.375);
\gpcolor{color=gp lt color axes}
\gpsetlinetype{gp lt axes}
\gpsetdashtype{gp dt axes}
\gpsetlinewidth{0.50}
\draw[gp path] (0.495,0.592)--(5.294,0.592);
\gpcolor{color=gp lt color border}
\gpsetlinetype{gp lt border}
\gpsetdashtype{gp dt solid}
\gpsetlinewidth{1.00}
\draw[gp path] (0.495,0.592)--(0.675,0.592);
\draw[gp path] (5.294,0.592)--(5.114,0.592);
\node[gp node right] at (0.385,0.592) {$1$};
\draw[gp path] (0.495,1.268)--(0.585,1.268);
\draw[gp path] (5.294,1.268)--(5.204,1.268);
\draw[gp path] (0.495,1.663)--(0.585,1.663);
\draw[gp path] (5.294,1.663)--(5.204,1.663);
\draw[gp path] (0.495,1.944)--(0.585,1.944);
\draw[gp path] (5.294,1.944)--(5.204,1.944);
\draw[gp path] (0.495,2.161)--(0.585,2.161);
\draw[gp path] (5.294,2.161)--(5.204,2.161);
\draw[gp path] (0.495,2.339)--(0.585,2.339);
\draw[gp path] (5.294,2.339)--(5.204,2.339);
\draw[gp path] (0.495,2.489)--(0.585,2.489);
\draw[gp path] (5.294,2.489)--(5.204,2.489);
\draw[gp path] (0.495,2.620)--(0.585,2.620);
\draw[gp path] (5.294,2.620)--(5.204,2.620);
\draw[gp path] (0.495,2.734)--(0.585,2.734);
\draw[gp path] (5.294,2.734)--(5.204,2.734);
\gpcolor{color=gp lt color axes}
\gpsetlinetype{gp lt axes}
\gpsetdashtype{gp dt axes}
\gpsetlinewidth{0.50}
\draw[gp path] (0.495,2.837)--(5.294,2.837);
\gpcolor{color=gp lt color border}
\gpsetlinetype{gp lt border}
\gpsetdashtype{gp dt solid}
\gpsetlinewidth{1.00}
\draw[gp path] (0.495,2.837)--(0.675,2.837);
\draw[gp path] (5.294,2.837)--(5.114,2.837);
\node[gp node right] at (0.385,2.837) {$10$};
\draw[gp path] (0.495,3.513)--(0.585,3.513);
\draw[gp path] (5.294,3.513)--(5.204,3.513);
\draw[gp path] (0.495,3.908)--(0.585,3.908);
\draw[gp path] (5.294,3.908)--(5.204,3.908);
\draw[gp path] (0.495,4.189)--(0.585,4.189);
\draw[gp path] (5.294,4.189)--(5.204,4.189);
\gpcolor{color=gp lt color axes}
\gpsetlinetype{gp lt axes}
\gpsetdashtype{gp dt axes}
\gpsetlinewidth{0.50}
\draw[gp path] (0.495,0.592)--(0.495,4.189);
\gpcolor{color=gp lt color border}
\gpsetlinetype{gp lt border}
\gpsetdashtype{gp dt solid}
\gpsetlinewidth{1.00}
\draw[gp path] (0.495,0.592)--(0.495,0.772);
\draw[gp path] (0.495,4.189)--(0.495,4.009);
\node[gp node center] at (0.495,0.407) {32};
\gpcolor{color=gp lt color axes}
\gpsetlinetype{gp lt axes}
\gpsetdashtype{gp dt axes}
\gpsetlinewidth{0.50}
\draw[gp path] (1.295,0.592)--(1.295,4.189);
\gpcolor{color=gp lt color border}
\gpsetlinetype{gp lt border}
\gpsetdashtype{gp dt solid}
\gpsetlinewidth{1.00}
\draw[gp path] (1.295,0.592)--(1.295,0.772);
\draw[gp path] (1.295,4.189)--(1.295,4.009);
\node[gp node center] at (1.295,0.407) {64};
\gpcolor{color=gp lt color axes}
\gpsetlinetype{gp lt axes}
\gpsetdashtype{gp dt axes}
\gpsetlinewidth{0.50}
\draw[gp path] (2.095,0.592)--(2.095,4.189);
\gpcolor{color=gp lt color border}
\gpsetlinetype{gp lt border}
\gpsetdashtype{gp dt solid}
\gpsetlinewidth{1.00}
\draw[gp path] (2.095,0.592)--(2.095,0.772);
\draw[gp path] (2.095,4.189)--(2.095,4.009);
\node[gp node center] at (2.095,0.407) {128};
\gpcolor{color=gp lt color axes}
\gpsetlinetype{gp lt axes}
\gpsetdashtype{gp dt axes}
\gpsetlinewidth{0.50}
\draw[gp path] (2.894,0.592)--(2.894,4.189);
\gpcolor{color=gp lt color border}
\gpsetlinetype{gp lt border}
\gpsetdashtype{gp dt solid}
\gpsetlinewidth{1.00}
\draw[gp path] (2.894,0.592)--(2.894,0.772);
\draw[gp path] (2.894,4.189)--(2.894,4.009);
\node[gp node center] at (2.894,0.407) {256};
\gpcolor{color=gp lt color axes}
\gpsetlinetype{gp lt axes}
\gpsetdashtype{gp dt axes}
\gpsetlinewidth{0.50}
\draw[gp path] (3.694,0.592)--(3.694,4.189);
\gpcolor{color=gp lt color border}
\gpsetlinetype{gp lt border}
\gpsetdashtype{gp dt solid}
\gpsetlinewidth{1.00}
\draw[gp path] (3.694,0.592)--(3.694,0.772);
\draw[gp path] (3.694,4.189)--(3.694,4.009);
\node[gp node center] at (3.694,0.407) {512};
\gpcolor{color=gp lt color axes}
\gpsetlinetype{gp lt axes}
\gpsetdashtype{gp dt axes}
\gpsetlinewidth{0.50}
\draw[gp path] (4.494,0.592)--(4.494,4.189);
\gpcolor{color=gp lt color border}
\gpsetlinetype{gp lt border}
\gpsetdashtype{gp dt solid}
\gpsetlinewidth{1.00}
\draw[gp path] (4.494,0.592)--(4.494,0.772);
\draw[gp path] (4.494,4.189)--(4.494,4.009);
\node[gp node center] at (4.494,0.407) {1024};
\gpcolor{color=gp lt color axes}
\gpsetlinetype{gp lt axes}
\gpsetdashtype{gp dt axes}
\gpsetlinewidth{0.50}
\draw[gp path] (5.294,0.592)--(5.294,4.189);
\gpcolor{color=gp lt color border}
\gpsetlinetype{gp lt border}
\gpsetdashtype{gp dt solid}
\gpsetlinewidth{1.00}
\draw[gp path] (5.294,0.592)--(5.294,0.772);
\draw[gp path] (5.294,4.189)--(5.294,4.009);
\node[gp node center] at (5.294,0.407) {2048};
\draw[gp path] (0.495,4.189)--(0.495,0.592)--(5.294,0.592)--(5.294,4.189)--cycle;
\node[gp node center] at (2.894,0.130) {number of cores};
\node[gp node right] at (1.687,4.374) {ts from 3, $\epsilon=1.e^{-7}$};
\gpcolor{rgb color={0.000,1.000,0.000}}
\draw[gp path] (1.797,4.374)--(2.417,4.374);
\draw[gp path] (0.495,1.093)--(1.295,1.097)--(2.095,1.105)--(2.894,1.088)--(3.694,1.082)%
  --(4.494,1.054)--(5.294,0.946);
\gpsetpointsize{4.00}
\gp3point{gp mark 5}{}{(0.495,1.093)}
\gp3point{gp mark 5}{}{(1.295,1.097)}
\gp3point{gp mark 5}{}{(2.095,1.105)}
\gp3point{gp mark 5}{}{(2.894,1.088)}
\gp3point{gp mark 5}{}{(3.694,1.082)}
\gp3point{gp mark 5}{}{(4.494,1.054)}
\gp3point{gp mark 5}{}{(5.294,0.946)}
\gp3point{gp mark 5}{}{(2.107,4.374)}
\gpcolor{color=gp lt color border}
\node[gp node right] at (4.499,4.374) {ts from 3, $\epsilon=2.e^{-4}$};
\gpcolor{rgb color={0.000,1.000,0.000}}
\draw[gp path] (4.609,4.374)--(5.229,4.374);
\draw[gp path] (0.495,0.592)--(1.295,0.592)--(2.095,0.592)--(2.894,0.592)--(3.694,0.592)%
  --(4.494,0.592)--(5.294,0.592);
\gp3point{gp mark 6}{}{(0.495,0.592)}
\gp3point{gp mark 6}{}{(1.295,0.592)}
\gp3point{gp mark 6}{}{(2.095,0.592)}
\gp3point{gp mark 6}{}{(2.894,0.592)}
\gp3point{gp mark 6}{}{(3.694,0.592)}
\gp3point{gp mark 6}{}{(4.494,0.592)}
\gp3point{gp mark 6}{}{(5.294,0.592)}
\gp3point{gp mark 6}{}{(4.919,4.374)}
\gpcolor{color=gp lt color border}
\node[gp node right] at (1.687,4.636) {dd, $\epsilon=1.e^{-7}$};
\gpcolor{rgb color={0.753,0.251,0.000}}
\draw[gp path] (1.797,4.636)--(2.417,4.636);
\draw[gp path] (3.694,3.949)--(4.494,3.773)--(5.294,3.528);
\gp3point{gp mark 2}{}{(3.694,3.949)}
\gp3point{gp mark 2}{}{(4.494,3.773)}
\gp3point{gp mark 2}{}{(5.294,3.528)}
\gp3point{gp mark 2}{}{(2.107,4.636)}
\gpcolor{color=gp lt color border}
\node[gp node right] at (4.499,4.636) {dd, $\epsilon=2.e^{-4}$};
\gpcolor{rgb color={0.753,0.251,0.000}}
\draw[gp path] (4.609,4.636)--(5.229,4.636);
\draw[gp path] (3.694,3.714)--(4.494,3.514)--(5.294,3.235);
\gp3point{gp mark 4}{}{(3.694,3.714)}
\gp3point{gp mark 4}{}{(4.494,3.514)}
\gp3point{gp mark 4}{}{(5.294,3.235)}
\gp3point{gp mark 4}{}{(4.919,4.636)}
\gpcolor{color=gp lt color border}
\draw[gp path] (0.495,4.189)--(0.495,0.592)--(5.294,0.592)--(5.294,4.189)--cycle;
\gpdefrectangularnode{gp plot 1}{\pgfpoint{0.495cm}{0.592cm}}{\pgfpoint{5.294cm}{4.189cm}}
\end{tikzpicture}
\label{UOL_precision_impact_L6}}
 	\caption{Accuracy impact: ratio indicating how much slower the solver associate with a curve is compared to the best solver (based on elapsed time)  \label{UOL_precision_impact}.}
 \end{figure}

This last part studies the choice of $\epsilon=1.e^{-7}$.
From the graph of $\mathcal{E}\left( \vm{u}^{ts},\vm{u}^C\right)$,  after a few iterations there is no real interest in continuing the scale loop because the  precision is clearly imposed by  $\left\| \vm{u}^R-\vm{u}^C \right\|_{E_{\Omega}}$.
A relative error with respect to this quantity may be satisfactory in some calculation contexts.
To evaluate the impact of such a looser error condition, some tests are performed with an arbitrary residual error precision: $\epsilon=\frac{\mathcal{E}\left( \vm{u}^R,\vm{u}^C\right)}{10}$.
As mentioned above, the residual error is conservative so if it is less than $\epsilon$ then $\mathcal{E}\left( \vm{u}^{ts},\vm{u}^R\right)$ is too.
This leads to the following condition:
\begin{equation}	
	\mathcal{E}\left( \vm{u}^{ts},\vm{u}^R\right)\leqslant \frac{\mathcal{E}\left( \vm{u}^R,\vm{u}^C\right)}{10}
	\label{inequal1}
\end{equation}
Passing at power of 2 and using equation \eqref{equality} the overall error in such condition will respect the following condition:
\begin{equation}	
	\mathcal{E}\left( \vm{u}^{ts},\vm{u}^C\right)^2 \leqslant \left( 1 + \frac{ \left\| \vm{u}^R \right\|_{E_{\Omega}}^2}{100\times \left\| \vm{u}^C \right\|_{E_{\Omega}}^2} \right) 	\times \mathcal{E}\left( \vm{u}^R,\vm{u}^C\right)^2  
	\label{inequal2}
\end{equation}
 And the relative error  between $\mathcal{E}\left( \vm{u}^{ts},\vm{u}^C\right)$ and $\mathcal{E}\left( \vm{u}^R,\vm{u}^C\right)$ is then limited by:
 \begin{equation}	
 	\frac{\mathcal{E}\left( \vm{u}^{ts},\vm{u}^C\right)-\mathcal{E}\left( \vm{u}^R,\vm{u}^C\right)}{\mathcal{E}\left( \vm{u}^R,\vm{u}^C\right)} \leqslant  \sqrt{1+\frac{ \left\| \vm{u}^R \right\|_{E_{\Omega}}^2}{100\times \left\| \vm{u}^C \right\|_{E_{\Omega}}^2}}-1  
 	\label{inequal3}
 \end{equation}
Thus, numerically, the  arbitrarily chosen $\epsilon$ leads to less than 0.5\% relative   error  between $\mathcal{E}\left( \vm{u}^{ts},\vm{u}^C\right)$ and $\mathcal{E}\left( \vm{u}^R,\vm{u}^C\right)$ using $0.9998$ for $\left\| \vm{u}^R \right\|_{E_{\Omega}}^2/\left\| \vm{u}^C \right\|_{E_{\Omega}}^2$ (average of data L3, L4 and L5).
The elapsed time performance for  $\epsilon=\frac{\mathcal{E}\left( \vm{u}^R,\vm{u}^C\right)}{10}$ are presented in figure \ref{UOL_precision_impact} as a ratio of the elapsed time for the current solver to  the best elapsed time among all solutions and all $\epsilon$.
 For the three fine-scale discretizations tested, L2, L5 and L6, in addition to $\epsilon=1.e^{-7}$, the tested values of $\epsilon$ are $3.e^{-3}$, $4.e^{-4}$ and $2.e^{-4}$ respectively.
 For all discretizations, a larger $\epsilon$ reduces the computation time of all solutions depending on it.
 
For L5 and L6 (figure \ref{UOL_precision_impact_L5} and \ref{UOL_precision_impact_L6}), the best solution is the \tS solver with a larger $\epsilon$.
The other solvers, even if some gains are observed with larger $\epsilon$, do not provide better performances than \tS solver.
For the "blr" solver, the gain in the factorization step due to the decrease in accuracy is partially lost by a larger number of iterations in the resolution of the conjugate gradient.
For the "dd" solver, the gain due to decrease in accuracy  is in the conjugate gradient resolution of the global boundary problem.
With large domains, this task represents about 50\% of the "dd" time consumed, so the gains are moderate.
As the number of cores increases, this task takes most of the time ( 97\% for L5 with 2048 processes) and the gains increase slightly.
For the \tS solver, the decrease in accuracy  reduces the number of iterations of the scale loop.
For L5 with 32, 128 and 2048 processes this task represents respectively about 46\%, 51\% and 66\% of the time consumed which explains a moderate gain (less than 2).
For L2 (figure \ref{UOL_precision_impact_L2}), the gain obtained with the \tS solver is more consistent compared to the other solvers but is not sufficient for this solver to show the best performance.

\subsection{Micro structural test case}\label{MS}
\begin{figure}[h!]
	\centering
	\subfloat[Fine scale]{\includegraphics[width=35mm]{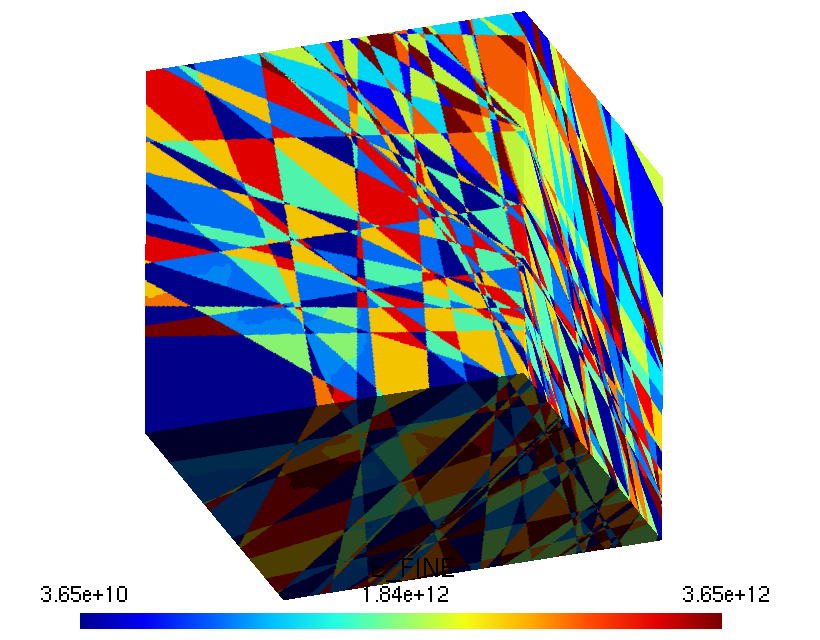}}
	\subfloat[Fine scale: cut at x=0.5]{\includegraphics[width=35mm]{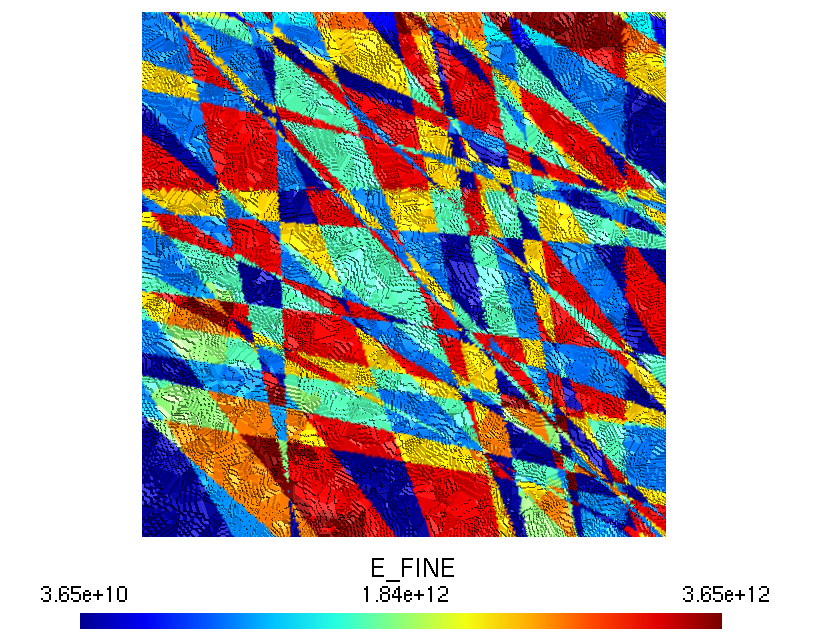}}
	\subfloat[Fine scale: discretization detail (cut at x=0.5)]{\includegraphics[width=35mm]{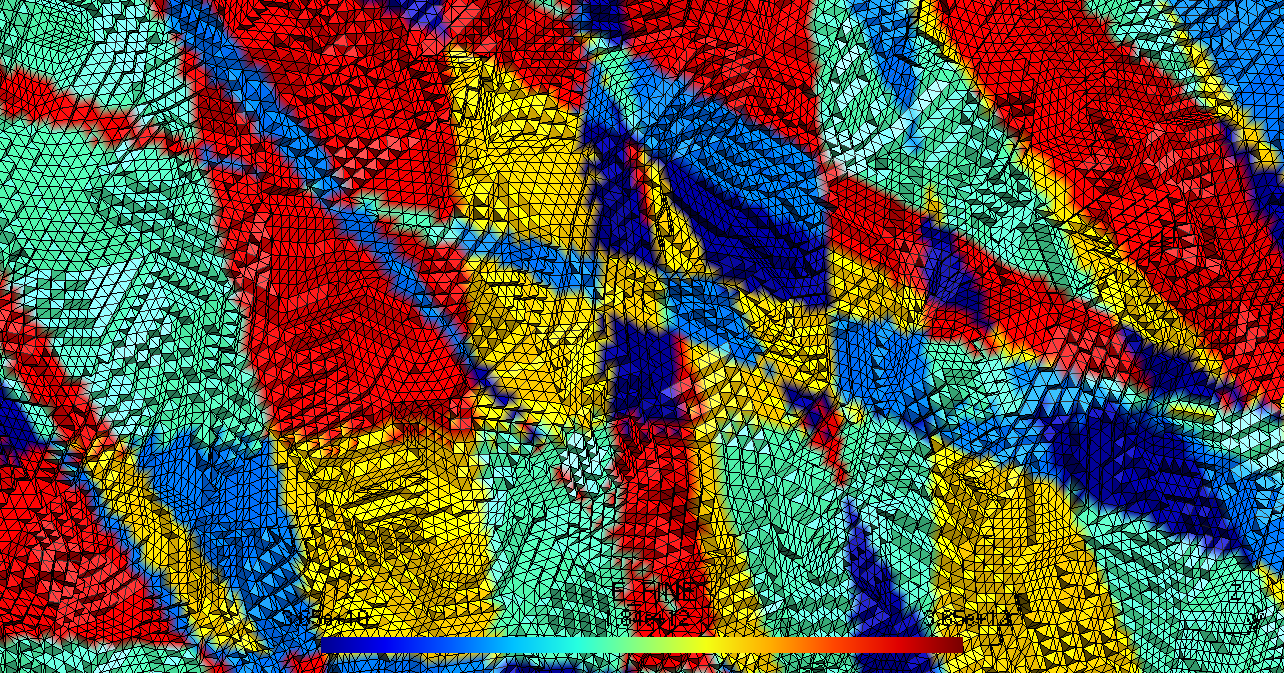}}\\
		\subfloat[Coarse scale]{\includegraphics[width=35mm]{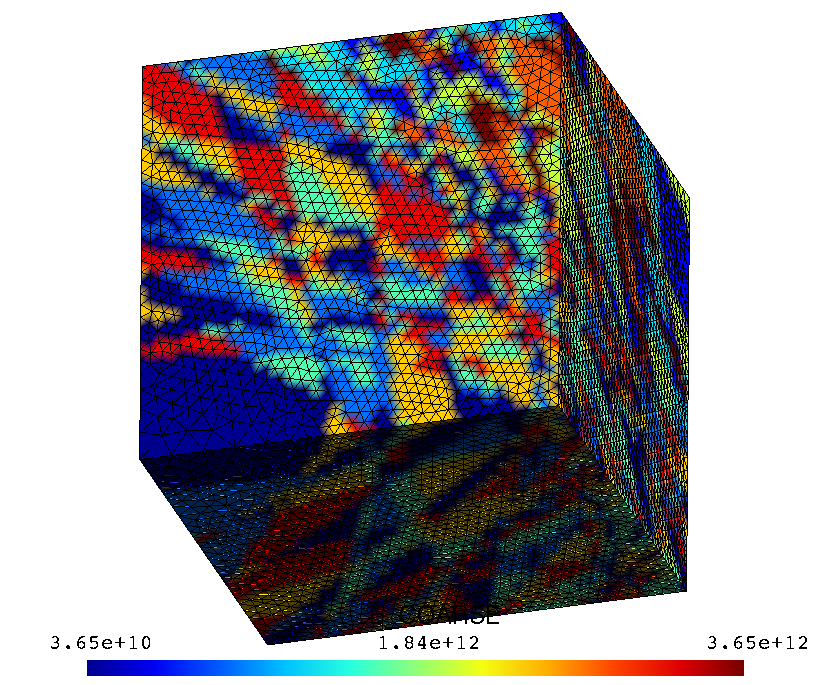}\label{MS_discret_coarse}}
	\subfloat[Coarse scale: cut at x=0.5]{\includegraphics[width=30mm]{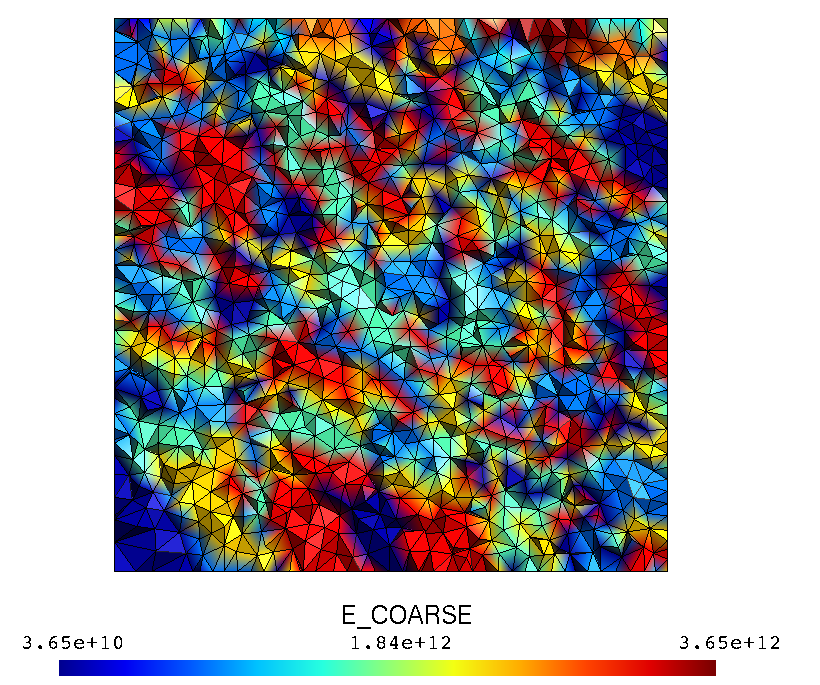}}
	\caption{Micro-structure discretization and Young's modulus values (in Pa).\label{MS_discret}}
\end{figure}

This example illustrates the cases where matrix conditioning can be a problem for a conventional iterative solver which then needs  good preconditionners.
A cube (with opposite diagonal  corners at the following coordinates: (0,0,0) and (2,2,2); unit: m) is randomly divided  by 64 plans of different orientation and origin (see \ref{MS_plan_anexe}).
In the cube, each sub-region delimited by plans  has a specific and uniform E.
These Young's moduli, which follow an arbitrary internal encoding, range from 36.5GPa to 3650GPa.
The coding makes the distribution of materials  relatively random, as illustrated in figure \ref{MS_discret}.

\begin{figure}[h]
	\subfloat[Domain decomposition]{\includegraphics[width=40mm]{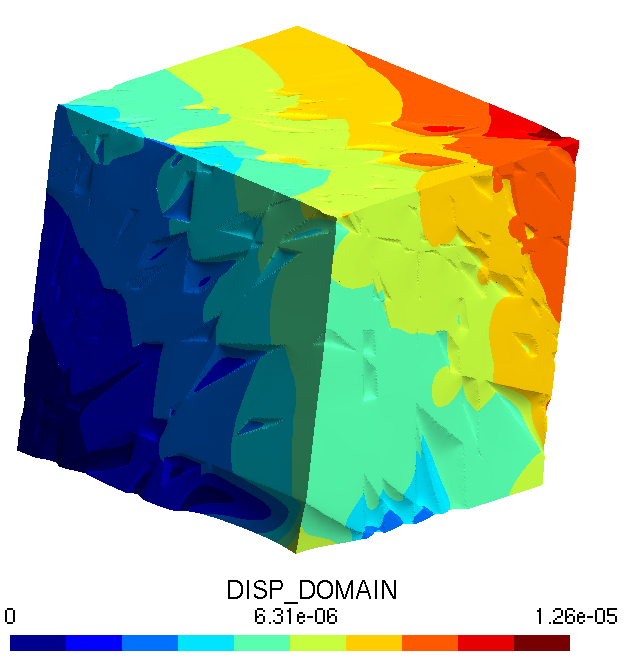}\label{MS_disp_dd}}
	\subfloat[Two-scale, fine scale]{\includegraphics[width=40mm]{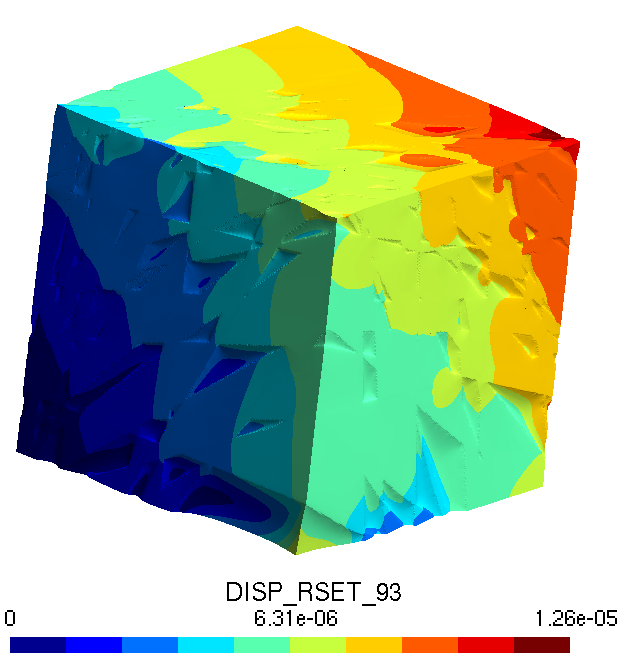}\label{MS_disp_rset}}
	\subfloat[Two-scale, coarse scale]{\includegraphics[width=40mm]{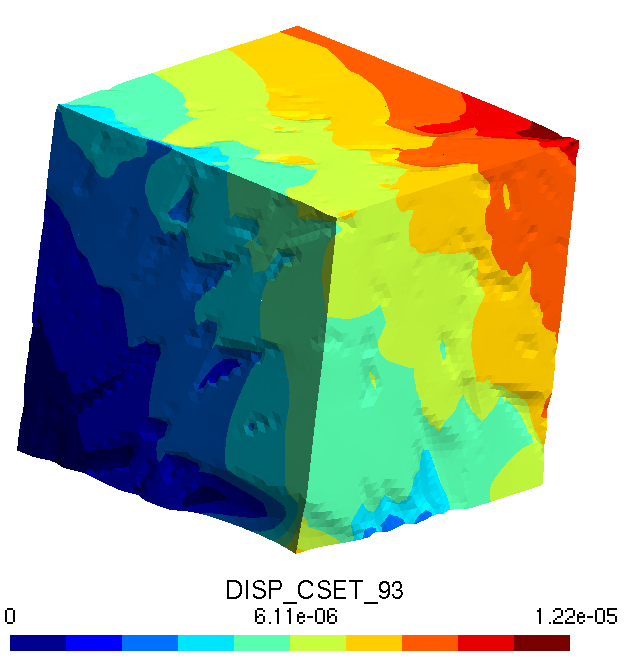}\label{MS_disp_cset}}
	\subfloat[Two-scale, fine scale]{\includegraphics[width=40mm]{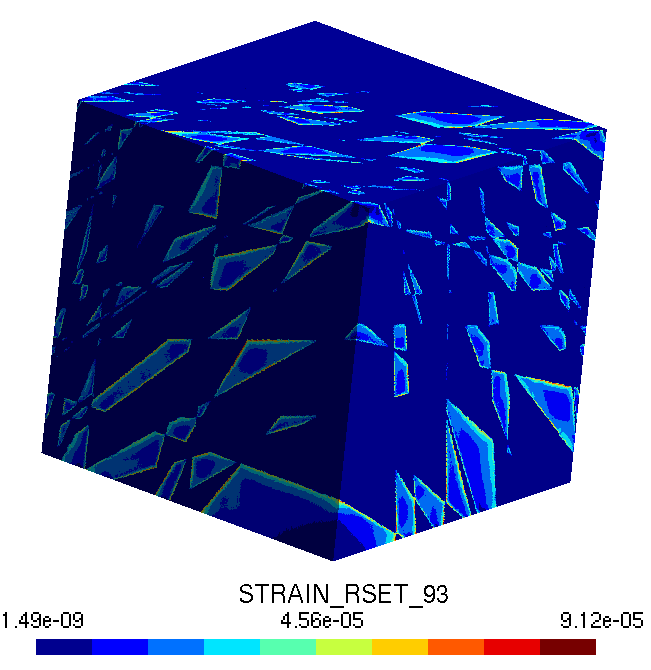}\label{MS_strain}}
	\caption{Micro-structure solutions: Displacement field with the "dd" solver \protect\subref{MS_disp_dd}. Displacement (\protect\subref*{MS_disp_rset},\protect\subref*{MS_disp_cset}) / strain \protect\subref{MS_strain} field at iteration 93 with the \tS solver. 
    The \tS solver can give the  displacement field  at both scales.
    At the coarse-scale \protect\subref{MS_disp_cset}, only the dofs in the C-set are used for post-processing, giving a crude visualization on the macro-element.
    At the fine-scale \protect\subref{MS_disp_rset}, using the R-set dofs, all details are given and can be compared to the solution obtained with the "dd" solver \protect\subref{MS_disp_dd}.
    \label{MS_disp}}
\end{figure}

In this test case, with the \tS solver,  all coarse nodes are enriched ($card\left( e\right) =card\left( c\right) $).
Thus, after choosing a fine-scale discretization  arbitrarily considered to meet our micro-structure representation needs, we consider a coarse-scale discretization that follows the equation \eqref{ts_jump_level_eq} (note that the coarse mesh is slightly adapted to the plane positions.
This can be seen in the figure \ref{MS_discret_coarse} where the visible left corner has larger element sizes, because no plane cuts this region).
The resulting fine-scale discretization has 60 278 925 dofs and corresponds to a section \ref{UOL} problem  between L5 and L6  (figure \ref{UOL_curve_dofs_ts}). 
This cube is compressed in all directions by the same constant load (4MN) applied on all surfaces.
As in the section \ref{UOL}, the \tS loop  starts with the displacement field of the macro problem without enrichment and the stopping criterion for the iterations is $\epsilon=1e^{-7}$.
The \tS solver (\TSI version) is again compared with the "dd"  solver.
This first simulation gives after 93 \tS iterations the fields presented in the figure \ref{MS_disp} with the  solution of "dd".
Note that the creation of post-processing files based on R-set dofs, obtained in parallel using MPI-IO functions, uses the per-macro-element organization of data  described in the section \ref{TS_integ_schem} and thus can be further parallelized using multi-threading.
As with the displacement field, we can obtain the strains from the R-set solution, as shown in figure \ref{MS_strain}, again by looping over the macro-elements using $\vm{u}_F^{e_{macro}}$.
 
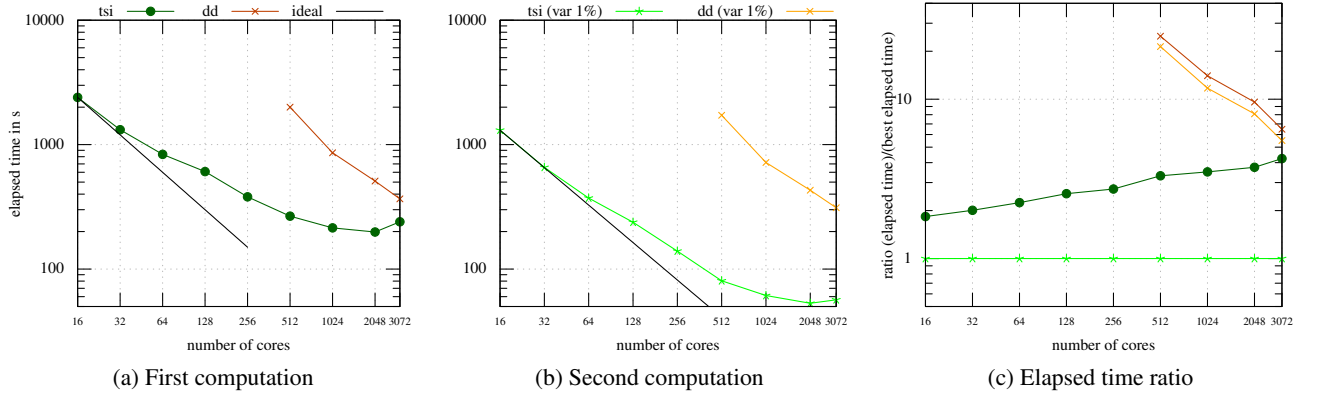
\begin{figure}[h]
	\subfloat[First computation]{
\begin{tikzpicture}[gnuplot]
\tikzset{every node/.append style={scale=0.60}}
\path (0.000,0.000) rectangle (5.625,4.812);
\gpcolor{color=gp lt color border}
\gpsetlinetype{gp lt border}
\gpsetdashtype{gp dt solid}
\gpsetlinewidth{1.00}
\draw[gp path] (1.010,0.592)--(1.100,0.592);
\draw[gp path] (5.294,0.592)--(5.204,0.592);
\draw[gp path] (1.010,0.723)--(1.100,0.723);
\draw[gp path] (5.294,0.723)--(5.204,0.723);
\draw[gp path] (1.010,0.834)--(1.100,0.834);
\draw[gp path] (5.294,0.834)--(5.204,0.834);
\draw[gp path] (1.010,0.930)--(1.100,0.930);
\draw[gp path] (5.294,0.930)--(5.204,0.930);
\draw[gp path] (1.010,1.015)--(1.100,1.015);
\draw[gp path] (5.294,1.015)--(5.204,1.015);
\gpcolor{color=gp lt color axes}
\gpsetlinetype{gp lt axes}
\gpsetdashtype{gp dt axes}
\gpsetlinewidth{0.50}
\draw[gp path] (1.010,1.090)--(5.294,1.090);
\gpcolor{color=gp lt color border}
\gpsetlinetype{gp lt border}
\gpsetdashtype{gp dt solid}
\gpsetlinewidth{1.00}
\draw[gp path] (1.010,1.090)--(1.190,1.090);
\draw[gp path] (5.294,1.090)--(5.114,1.090);
\node[gp node right] at (0.900,1.090) {$100$};
\draw[gp path] (1.010,1.589)--(1.100,1.589);
\draw[gp path] (5.294,1.589)--(5.204,1.589);
\draw[gp path] (1.010,1.880)--(1.100,1.880);
\draw[gp path] (5.294,1.880)--(5.204,1.880);
\draw[gp path] (1.010,2.087)--(1.100,2.087);
\draw[gp path] (5.294,2.087)--(5.204,2.087);
\draw[gp path] (1.010,2.247)--(1.100,2.247);
\draw[gp path] (5.294,2.247)--(5.204,2.247);
\draw[gp path] (1.010,2.378)--(1.100,2.378);
\draw[gp path] (5.294,2.378)--(5.204,2.378);
\draw[gp path] (1.010,2.489)--(1.100,2.489);
\draw[gp path] (5.294,2.489)--(5.204,2.489);
\draw[gp path] (1.010,2.585)--(1.100,2.585);
\draw[gp path] (5.294,2.585)--(5.204,2.585);
\draw[gp path] (1.010,2.670)--(1.100,2.670);
\draw[gp path] (5.294,2.670)--(5.204,2.670);
\gpcolor{color=gp lt color axes}
\gpsetlinetype{gp lt axes}
\gpsetdashtype{gp dt axes}
\gpsetlinewidth{0.50}
\draw[gp path] (1.010,2.746)--(5.294,2.746);
\gpcolor{color=gp lt color border}
\gpsetlinetype{gp lt border}
\gpsetdashtype{gp dt solid}
\gpsetlinewidth{1.00}
\draw[gp path] (1.010,2.746)--(1.190,2.746);
\draw[gp path] (5.294,2.746)--(5.114,2.746);
\node[gp node right] at (0.900,2.746) {$1000$};
\draw[gp path] (1.010,3.244)--(1.100,3.244);
\draw[gp path] (5.294,3.244)--(5.204,3.244);
\draw[gp path] (1.010,3.535)--(1.100,3.535);
\draw[gp path] (5.294,3.535)--(5.204,3.535);
\draw[gp path] (1.010,3.742)--(1.100,3.742);
\draw[gp path] (5.294,3.742)--(5.204,3.742);
\draw[gp path] (1.010,3.903)--(1.100,3.903);
\draw[gp path] (5.294,3.903)--(5.204,3.903);
\draw[gp path] (1.010,4.034)--(1.100,4.034);
\draw[gp path] (5.294,4.034)--(5.204,4.034);
\draw[gp path] (1.010,4.145)--(1.100,4.145);
\draw[gp path] (5.294,4.145)--(5.204,4.145);
\draw[gp path] (1.010,4.241)--(1.100,4.241);
\draw[gp path] (5.294,4.241)--(5.204,4.241);
\draw[gp path] (1.010,4.325)--(1.100,4.325);
\draw[gp path] (5.294,4.325)--(5.204,4.325);
\gpcolor{color=gp lt color axes}
\gpsetlinetype{gp lt axes}
\gpsetdashtype{gp dt axes}
\gpsetlinewidth{0.50}
\draw[gp path] (1.010,4.401)--(5.294,4.401);
\gpcolor{color=gp lt color border}
\gpsetlinetype{gp lt border}
\gpsetdashtype{gp dt solid}
\gpsetlinewidth{1.00}
\draw[gp path] (1.010,4.401)--(1.190,4.401);
\draw[gp path] (5.294,4.401)--(5.114,4.401);
\node[gp node right] at (0.900,4.401) {$10000$};
\gpcolor{color=gp lt color axes}
\gpsetlinetype{gp lt axes}
\gpsetdashtype{gp dt axes}
\gpsetlinewidth{0.50}
\draw[gp path] (1.010,0.592)--(1.010,4.401);
\gpcolor{color=gp lt color border}
\gpsetlinetype{gp lt border}
\gpsetdashtype{gp dt solid}
\gpsetlinewidth{1.00}
\draw[gp path] (1.010,0.592)--(1.010,0.772);
\draw[gp path] (1.010,4.401)--(1.010,4.221);
\node[gp node center,font={\fontsize{7.0pt}{8.4pt}\selectfont}] at (1.010,0.407) {16};
\gpcolor{color=gp lt color axes}
\gpsetlinetype{gp lt axes}
\gpsetdashtype{gp dt axes}
\gpsetlinewidth{0.50}
\draw[gp path] (1.575,0.592)--(1.575,4.401);
\gpcolor{color=gp lt color border}
\gpsetlinetype{gp lt border}
\gpsetdashtype{gp dt solid}
\gpsetlinewidth{1.00}
\draw[gp path] (1.575,0.592)--(1.575,0.772);
\draw[gp path] (1.575,4.401)--(1.575,4.221);
\node[gp node center,font={\fontsize{7.0pt}{8.4pt}\selectfont}] at (1.575,0.407) {32};
\gpcolor{color=gp lt color axes}
\gpsetlinetype{gp lt axes}
\gpsetdashtype{gp dt axes}
\gpsetlinewidth{0.50}
\draw[gp path] (2.140,0.592)--(2.140,4.401);
\gpcolor{color=gp lt color border}
\gpsetlinetype{gp lt border}
\gpsetdashtype{gp dt solid}
\gpsetlinewidth{1.00}
\draw[gp path] (2.140,0.592)--(2.140,0.772);
\draw[gp path] (2.140,4.401)--(2.140,4.221);
\node[gp node center,font={\fontsize{7.0pt}{8.4pt}\selectfont}] at (2.140,0.407) {64};
\gpcolor{color=gp lt color axes}
\gpsetlinetype{gp lt axes}
\gpsetdashtype{gp dt axes}
\gpsetlinewidth{0.50}
\draw[gp path] (2.704,0.592)--(2.704,4.401);
\gpcolor{color=gp lt color border}
\gpsetlinetype{gp lt border}
\gpsetdashtype{gp dt solid}
\gpsetlinewidth{1.00}
\draw[gp path] (2.704,0.592)--(2.704,0.772);
\draw[gp path] (2.704,4.401)--(2.704,4.221);
\node[gp node center,font={\fontsize{7.0pt}{8.4pt}\selectfont}] at (2.704,0.407) {128};
\gpcolor{color=gp lt color axes}
\gpsetlinetype{gp lt axes}
\gpsetdashtype{gp dt axes}
\gpsetlinewidth{0.50}
\draw[gp path] (3.269,0.592)--(3.269,4.401);
\gpcolor{color=gp lt color border}
\gpsetlinetype{gp lt border}
\gpsetdashtype{gp dt solid}
\gpsetlinewidth{1.00}
\draw[gp path] (3.269,0.592)--(3.269,0.772);
\draw[gp path] (3.269,4.401)--(3.269,4.221);
\node[gp node center,font={\fontsize{7.0pt}{8.4pt}\selectfont}] at (3.269,0.407) {256};
\gpcolor{color=gp lt color axes}
\gpsetlinetype{gp lt axes}
\gpsetdashtype{gp dt axes}
\gpsetlinewidth{0.50}
\draw[gp path] (3.834,0.592)--(3.834,4.401);
\gpcolor{color=gp lt color border}
\gpsetlinetype{gp lt border}
\gpsetdashtype{gp dt solid}
\gpsetlinewidth{1.00}
\draw[gp path] (3.834,0.592)--(3.834,0.772);
\draw[gp path] (3.834,4.401)--(3.834,4.221);
\node[gp node center,font={\fontsize{7.0pt}{8.4pt}\selectfont}] at (3.834,0.407) {512};
\gpcolor{color=gp lt color axes}
\gpsetlinetype{gp lt axes}
\gpsetdashtype{gp dt axes}
\gpsetlinewidth{0.50}
\draw[gp path] (4.399,0.592)--(4.399,4.401);
\gpcolor{color=gp lt color border}
\gpsetlinetype{gp lt border}
\gpsetdashtype{gp dt solid}
\gpsetlinewidth{1.00}
\draw[gp path] (4.399,0.592)--(4.399,0.772);
\draw[gp path] (4.399,4.401)--(4.399,4.221);
\node[gp node center,font={\fontsize{7.0pt}{8.4pt}\selectfont}] at (4.399,0.407) {1024};
\gpcolor{color=gp lt color axes}
\gpsetlinetype{gp lt axes}
\gpsetdashtype{gp dt axes}
\gpsetlinewidth{0.50}
\draw[gp path] (4.964,0.592)--(4.964,4.401);
\gpcolor{color=gp lt color border}
\gpsetlinetype{gp lt border}
\gpsetdashtype{gp dt solid}
\gpsetlinewidth{1.00}
\draw[gp path] (4.964,0.592)--(4.964,0.772);
\draw[gp path] (4.964,4.401)--(4.964,4.221);
\node[gp node center,font={\fontsize{7.0pt}{8.4pt}\selectfont}] at (4.964,0.407) {2048};
\gpcolor{color=gp lt color axes}
\gpsetlinetype{gp lt axes}
\gpsetdashtype{gp dt axes}
\gpsetlinewidth{0.50}
\draw[gp path] (5.294,0.592)--(5.294,4.401);
\gpcolor{color=gp lt color border}
\gpsetlinetype{gp lt border}
\gpsetdashtype{gp dt solid}
\gpsetlinewidth{1.00}
\draw[gp path] (5.294,0.592)--(5.294,0.772);
\draw[gp path] (5.294,4.401)--(5.294,4.221);
\node[gp node center,font={\fontsize{7.0pt}{8.4pt}\selectfont}] at (5.294,0.407) {3072};
\draw[gp path] (1.010,4.401)--(1.010,0.592)--(5.294,0.592)--(5.294,4.401)--cycle;
\node[gp node center,rotate=-270] at (0.194,2.496) {elapsed time in s};
\node[gp node center] at (3.152,0.068) {number of cores};
\node[gp node right] at (1.560,4.519) {tsi};
\gpcolor{rgb color={0.000,0.392,0.000}}
\draw[gp path] (1.670,4.519)--(2.290,4.519);
\draw[gp path] (1.010,3.374)--(1.575,2.944)--(2.140,2.616)--(2.704,2.387)--(3.269,2.051)%
  --(3.834,1.794)--(4.399,1.638)--(4.964,1.583)--(5.294,1.719);
\gpsetpointsize{4.00}
\gp3point{gp mark 7}{}{(1.010,3.374)}
\gp3point{gp mark 7}{}{(1.575,2.944)}
\gp3point{gp mark 7}{}{(2.140,2.616)}
\gp3point{gp mark 7}{}{(2.704,2.387)}
\gp3point{gp mark 7}{}{(3.269,2.051)}
\gp3point{gp mark 7}{}{(3.834,1.794)}
\gp3point{gp mark 7}{}{(4.399,1.638)}
\gp3point{gp mark 7}{}{(4.964,1.583)}
\gp3point{gp mark 7}{}{(5.294,1.719)}
\gp3point{gp mark 7}{}{(1.980,4.519)}
\gpcolor{color=gp lt color border}
\node[gp node right] at (2.950,4.519) {dd};
\gpcolor{rgb color={0.753,0.251,0.000}}
\draw[gp path] (3.060,4.519)--(3.680,4.519);
\draw[gp path] (3.834,3.244)--(4.399,2.637)--(4.964,2.261)--(5.294,2.024);
\gp3point{gp mark 2}{}{(3.834,3.244)}
\gp3point{gp mark 2}{}{(4.399,2.637)}
\gp3point{gp mark 2}{}{(4.964,2.261)}
\gp3point{gp mark 2}{}{(5.294,2.024)}
\gp3point{gp mark 2}{}{(3.370,4.519)}
\gpcolor{color=gp lt color border}
\node[gp node right] at (4.340,4.519) {ideal};
\gpcolor{rgb color={0.000,0.000,0.000}}
\draw[gp path] (4.450,4.519)--(5.070,4.519);
\draw[gp path] (1.010,3.374)--(1.033,3.354)--(1.056,3.334)--(1.078,3.314)--(1.101,3.294)%
  --(1.124,3.273)--(1.147,3.253)--(1.170,3.233)--(1.193,3.213)--(1.215,3.193)--(1.238,3.173)%
  --(1.261,3.153)--(1.284,3.132)--(1.307,3.112)--(1.329,3.092)--(1.352,3.072)--(1.375,3.052)%
  --(1.398,3.032)--(1.421,3.012)--(1.444,2.992)--(1.466,2.971)--(1.489,2.951)--(1.512,2.931)%
  --(1.535,2.911)--(1.558,2.891)--(1.581,2.871)--(1.603,2.851)--(1.626,2.830)--(1.649,2.810)%
  --(1.672,2.790)--(1.695,2.770)--(1.717,2.750)--(1.740,2.730)--(1.763,2.710)--(1.786,2.690)%
  --(1.809,2.669)--(1.832,2.649)--(1.854,2.629)--(1.877,2.609)--(1.900,2.589)--(1.923,2.569)%
  --(1.946,2.549)--(1.968,2.528)--(1.991,2.508)--(2.014,2.488)--(2.037,2.468)--(2.060,2.448)%
  --(2.083,2.428)--(2.105,2.408)--(2.128,2.388)--(2.151,2.367)--(2.174,2.347)--(2.197,2.327)%
  --(2.219,2.307)--(2.242,2.287)--(2.265,2.267)--(2.288,2.247)--(2.311,2.226)--(2.334,2.206)%
  --(2.356,2.186)--(2.379,2.166)--(2.402,2.146)--(2.425,2.126)--(2.448,2.106)--(2.470,2.086)%
  --(2.493,2.065)--(2.516,2.045)--(2.539,2.025)--(2.562,2.005)--(2.585,1.985)--(2.607,1.965)%
  --(2.630,1.945)--(2.653,1.924)--(2.676,1.904)--(2.699,1.884)--(2.722,1.864)--(2.744,1.844)%
  --(2.767,1.824)--(2.790,1.804)--(2.813,1.784)--(2.836,1.763)--(2.858,1.743)--(2.881,1.723)%
  --(2.904,1.703)--(2.927,1.683)--(2.950,1.663)--(2.973,1.643)--(2.995,1.622)--(3.018,1.602)%
  --(3.041,1.582)--(3.064,1.562)--(3.087,1.542)--(3.109,1.522)--(3.132,1.502)--(3.155,1.482)%
  --(3.178,1.461)--(3.201,1.441)--(3.224,1.421)--(3.246,1.401)--(3.269,1.381);
\gpcolor{color=gp lt color border}
\draw[gp path] (1.010,4.401)--(1.010,0.592)--(5.294,0.592)--(5.294,4.401)--cycle;
\gpdefrectangularnode{gp plot 1}{\pgfpoint{1.010cm}{0.592cm}}{\pgfpoint{5.294cm}{4.401cm}}
\end{tikzpicture}
\label{MS_time1}}
	\subfloat[Second computation]{
\begin{tikzpicture}[gnuplot]
\tikzset{every node/.append style={scale=0.60}}
\path (0.000,0.000) rectangle (5.625,4.812);
\gpcolor{color=gp lt color border}
\gpsetlinetype{gp lt border}
\gpsetdashtype{gp dt solid}
\gpsetlinewidth{1.00}
\draw[gp path] (0.825,0.592)--(0.915,0.592);
\draw[gp path] (5.294,0.592)--(5.204,0.592);
\draw[gp path] (0.825,0.723)--(0.915,0.723);
\draw[gp path] (5.294,0.723)--(5.204,0.723);
\draw[gp path] (0.825,0.834)--(0.915,0.834);
\draw[gp path] (5.294,0.834)--(5.204,0.834);
\draw[gp path] (0.825,0.930)--(0.915,0.930);
\draw[gp path] (5.294,0.930)--(5.204,0.930);
\draw[gp path] (0.825,1.015)--(0.915,1.015);
\draw[gp path] (5.294,1.015)--(5.204,1.015);
\gpcolor{color=gp lt color axes}
\gpsetlinetype{gp lt axes}
\gpsetdashtype{gp dt axes}
\gpsetlinewidth{0.50}
\draw[gp path] (0.825,1.090)--(5.294,1.090);
\gpcolor{color=gp lt color border}
\gpsetlinetype{gp lt border}
\gpsetdashtype{gp dt solid}
\gpsetlinewidth{1.00}
\draw[gp path] (0.825,1.090)--(1.005,1.090);
\draw[gp path] (5.294,1.090)--(5.114,1.090);
\node[gp node right] at (0.715,1.090) {$100$};
\draw[gp path] (0.825,1.589)--(0.915,1.589);
\draw[gp path] (5.294,1.589)--(5.204,1.589);
\draw[gp path] (0.825,1.880)--(0.915,1.880);
\draw[gp path] (5.294,1.880)--(5.204,1.880);
\draw[gp path] (0.825,2.087)--(0.915,2.087);
\draw[gp path] (5.294,2.087)--(5.204,2.087);
\draw[gp path] (0.825,2.247)--(0.915,2.247);
\draw[gp path] (5.294,2.247)--(5.204,2.247);
\draw[gp path] (0.825,2.378)--(0.915,2.378);
\draw[gp path] (5.294,2.378)--(5.204,2.378);
\draw[gp path] (0.825,2.489)--(0.915,2.489);
\draw[gp path] (5.294,2.489)--(5.204,2.489);
\draw[gp path] (0.825,2.585)--(0.915,2.585);
\draw[gp path] (5.294,2.585)--(5.204,2.585);
\draw[gp path] (0.825,2.670)--(0.915,2.670);
\draw[gp path] (5.294,2.670)--(5.204,2.670);
\gpcolor{color=gp lt color axes}
\gpsetlinetype{gp lt axes}
\gpsetdashtype{gp dt axes}
\gpsetlinewidth{0.50}
\draw[gp path] (0.825,2.746)--(5.294,2.746);
\gpcolor{color=gp lt color border}
\gpsetlinetype{gp lt border}
\gpsetdashtype{gp dt solid}
\gpsetlinewidth{1.00}
\draw[gp path] (0.825,2.746)--(1.005,2.746);
\draw[gp path] (5.294,2.746)--(5.114,2.746);
\node[gp node right] at (0.715,2.746) {$1000$};
\draw[gp path] (0.825,3.244)--(0.915,3.244);
\draw[gp path] (5.294,3.244)--(5.204,3.244);
\draw[gp path] (0.825,3.535)--(0.915,3.535);
\draw[gp path] (5.294,3.535)--(5.204,3.535);
\draw[gp path] (0.825,3.742)--(0.915,3.742);
\draw[gp path] (5.294,3.742)--(5.204,3.742);
\draw[gp path] (0.825,3.903)--(0.915,3.903);
\draw[gp path] (5.294,3.903)--(5.204,3.903);
\draw[gp path] (0.825,4.034)--(0.915,4.034);
\draw[gp path] (5.294,4.034)--(5.204,4.034);
\draw[gp path] (0.825,4.145)--(0.915,4.145);
\draw[gp path] (5.294,4.145)--(5.204,4.145);
\draw[gp path] (0.825,4.241)--(0.915,4.241);
\draw[gp path] (5.294,4.241)--(5.204,4.241);
\draw[gp path] (0.825,4.325)--(0.915,4.325);
\draw[gp path] (5.294,4.325)--(5.204,4.325);
\gpcolor{color=gp lt color axes}
\gpsetlinetype{gp lt axes}
\gpsetdashtype{gp dt axes}
\gpsetlinewidth{0.50}
\draw[gp path] (0.825,4.401)--(5.294,4.401);
\gpcolor{color=gp lt color border}
\gpsetlinetype{gp lt border}
\gpsetdashtype{gp dt solid}
\gpsetlinewidth{1.00}
\draw[gp path] (0.825,4.401)--(1.005,4.401);
\draw[gp path] (5.294,4.401)--(5.114,4.401);
\node[gp node right] at (0.715,4.401) {$10000$};
\gpcolor{color=gp lt color axes}
\gpsetlinetype{gp lt axes}
\gpsetdashtype{gp dt axes}
\gpsetlinewidth{0.50}
\draw[gp path] (0.825,0.592)--(0.825,4.401);
\gpcolor{color=gp lt color border}
\gpsetlinetype{gp lt border}
\gpsetdashtype{gp dt solid}
\gpsetlinewidth{1.00}
\draw[gp path] (0.825,0.592)--(0.825,0.772);
\draw[gp path] (0.825,4.401)--(0.825,4.221);
\node[gp node center,font={\fontsize{7.0pt}{8.4pt}\selectfont}] at (0.825,0.407) {16};
\gpcolor{color=gp lt color axes}
\gpsetlinetype{gp lt axes}
\gpsetdashtype{gp dt axes}
\gpsetlinewidth{0.50}
\draw[gp path] (1.414,0.592)--(1.414,4.401);
\gpcolor{color=gp lt color border}
\gpsetlinetype{gp lt border}
\gpsetdashtype{gp dt solid}
\gpsetlinewidth{1.00}
\draw[gp path] (1.414,0.592)--(1.414,0.772);
\draw[gp path] (1.414,4.401)--(1.414,4.221);
\node[gp node center,font={\fontsize{7.0pt}{8.4pt}\selectfont}] at (1.414,0.407) {32};
\gpcolor{color=gp lt color axes}
\gpsetlinetype{gp lt axes}
\gpsetdashtype{gp dt axes}
\gpsetlinewidth{0.50}
\draw[gp path] (2.003,0.592)--(2.003,4.401);
\gpcolor{color=gp lt color border}
\gpsetlinetype{gp lt border}
\gpsetdashtype{gp dt solid}
\gpsetlinewidth{1.00}
\draw[gp path] (2.003,0.592)--(2.003,0.772);
\draw[gp path] (2.003,4.401)--(2.003,4.221);
\node[gp node center,font={\fontsize{7.0pt}{8.4pt}\selectfont}] at (2.003,0.407) {64};
\gpcolor{color=gp lt color axes}
\gpsetlinetype{gp lt axes}
\gpsetdashtype{gp dt axes}
\gpsetlinewidth{0.50}
\draw[gp path] (2.593,0.592)--(2.593,4.401);
\gpcolor{color=gp lt color border}
\gpsetlinetype{gp lt border}
\gpsetdashtype{gp dt solid}
\gpsetlinewidth{1.00}
\draw[gp path] (2.593,0.592)--(2.593,0.772);
\draw[gp path] (2.593,4.401)--(2.593,4.221);
\node[gp node center,font={\fontsize{7.0pt}{8.4pt}\selectfont}] at (2.593,0.407) {128};
\gpcolor{color=gp lt color axes}
\gpsetlinetype{gp lt axes}
\gpsetdashtype{gp dt axes}
\gpsetlinewidth{0.50}
\draw[gp path] (3.182,0.592)--(3.182,4.401);
\gpcolor{color=gp lt color border}
\gpsetlinetype{gp lt border}
\gpsetdashtype{gp dt solid}
\gpsetlinewidth{1.00}
\draw[gp path] (3.182,0.592)--(3.182,0.772);
\draw[gp path] (3.182,4.401)--(3.182,4.221);
\node[gp node center,font={\fontsize{7.0pt}{8.4pt}\selectfont}] at (3.182,0.407) {256};
\gpcolor{color=gp lt color axes}
\gpsetlinetype{gp lt axes}
\gpsetdashtype{gp dt axes}
\gpsetlinewidth{0.50}
\draw[gp path] (3.771,0.592)--(3.771,4.401);
\gpcolor{color=gp lt color border}
\gpsetlinetype{gp lt border}
\gpsetdashtype{gp dt solid}
\gpsetlinewidth{1.00}
\draw[gp path] (3.771,0.592)--(3.771,0.772);
\draw[gp path] (3.771,4.401)--(3.771,4.221);
\node[gp node center,font={\fontsize{7.0pt}{8.4pt}\selectfont}] at (3.771,0.407) {512};
\gpcolor{color=gp lt color axes}
\gpsetlinetype{gp lt axes}
\gpsetdashtype{gp dt axes}
\gpsetlinewidth{0.50}
\draw[gp path] (4.360,0.592)--(4.360,4.401);
\gpcolor{color=gp lt color border}
\gpsetlinetype{gp lt border}
\gpsetdashtype{gp dt solid}
\gpsetlinewidth{1.00}
\draw[gp path] (4.360,0.592)--(4.360,0.772);
\draw[gp path] (4.360,4.401)--(4.360,4.221);
\node[gp node center,font={\fontsize{7.0pt}{8.4pt}\selectfont}] at (4.360,0.407) {1024};
\gpcolor{color=gp lt color axes}
\gpsetlinetype{gp lt axes}
\gpsetdashtype{gp dt axes}
\gpsetlinewidth{0.50}
\draw[gp path] (4.949,0.592)--(4.949,4.401);
\gpcolor{color=gp lt color border}
\gpsetlinetype{gp lt border}
\gpsetdashtype{gp dt solid}
\gpsetlinewidth{1.00}
\draw[gp path] (4.949,0.592)--(4.949,0.772);
\draw[gp path] (4.949,4.401)--(4.949,4.221);
\node[gp node center,font={\fontsize{7.0pt}{8.4pt}\selectfont}] at (4.949,0.407) {2048};
\gpcolor{color=gp lt color axes}
\gpsetlinetype{gp lt axes}
\gpsetdashtype{gp dt axes}
\gpsetlinewidth{0.50}
\draw[gp path] (5.294,0.592)--(5.294,4.401);
\gpcolor{color=gp lt color border}
\gpsetlinetype{gp lt border}
\gpsetdashtype{gp dt solid}
\gpsetlinewidth{1.00}
\draw[gp path] (5.294,0.592)--(5.294,0.772);
\draw[gp path] (5.294,4.401)--(5.294,4.221);
\node[gp node center,font={\fontsize{7.0pt}{8.4pt}\selectfont}] at (5.294,0.407) {3072};
\draw[gp path] (0.825,4.401)--(0.825,0.592)--(5.294,0.592)--(5.294,4.401)--cycle;
\node[gp node center] at (3.059,0.068) {number of cores};
\node[gp node right] at (2.255,4.519) {tsi (var 1\%)};
\gpcolor{rgb color={0.000,1.000,0.000}}
\draw[gp path] (2.365,4.519)--(2.985,4.519);
\draw[gp path] (0.825,2.938)--(1.414,2.443)--(2.003,2.035)--(2.593,1.713)--(3.182,1.329)%
  --(3.771,0.934)--(4.360,0.738)--(4.949,0.635)--(5.294,0.681);
\gpsetpointsize{4.00}
\gp3point{gp mark 3}{}{(0.825,2.938)}
\gp3point{gp mark 3}{}{(1.414,2.443)}
\gp3point{gp mark 3}{}{(2.003,2.035)}
\gp3point{gp mark 3}{}{(2.593,1.713)}
\gp3point{gp mark 3}{}{(3.182,1.329)}
\gp3point{gp mark 3}{}{(3.771,0.934)}
\gp3point{gp mark 3}{}{(4.360,0.738)}
\gp3point{gp mark 3}{}{(4.949,0.635)}
\gp3point{gp mark 3}{}{(5.294,0.681)}
\gp3point{gp mark 3}{}{(2.675,4.519)}
\gpcolor{color=gp lt color border}
\node[gp node right] at (4.525,4.519) {dd (var 1\%)};
\gpcolor{rgb color={1.000,0.647,0.000}}
\draw[gp path] (4.635,4.519)--(5.255,4.519);
\draw[gp path] (3.771,3.136)--(4.360,2.508)--(4.949,2.139)--(5.294,1.906);
\gp3point{gp mark 2}{}{(3.771,3.136)}
\gp3point{gp mark 2}{}{(4.360,2.508)}
\gp3point{gp mark 2}{}{(4.949,2.139)}
\gp3point{gp mark 2}{}{(5.294,1.906)}
\gp3point{gp mark 2}{}{(4.945,4.519)}
\gpcolor{rgb color={0.000,0.000,0.000}}
\draw[gp path] (0.825,2.938)--(0.867,2.902)--(0.908,2.867)--(0.950,2.832)--(0.992,2.797)%
  --(1.033,2.762)--(1.075,2.726)--(1.117,2.691)--(1.158,2.656)--(1.200,2.621)--(1.242,2.585)%
  --(1.283,2.550)--(1.325,2.515)--(1.367,2.480)--(1.408,2.444)--(1.450,2.409)--(1.492,2.374)%
  --(1.533,2.339)--(1.575,2.303)--(1.617,2.268)--(1.658,2.233)--(1.700,2.198)--(1.742,2.163)%
  --(1.783,2.127)--(1.825,2.092)--(1.867,2.057)--(1.908,2.022)--(1.950,1.986)--(1.991,1.951)%
  --(2.033,1.916)--(2.075,1.881)--(2.116,1.845)--(2.158,1.810)--(2.200,1.775)--(2.241,1.740)%
  --(2.283,1.704)--(2.325,1.669)--(2.366,1.634)--(2.408,1.599)--(2.450,1.564)--(2.491,1.528)%
  --(2.533,1.493)--(2.575,1.458)--(2.616,1.423)--(2.658,1.387)--(2.700,1.352)--(2.741,1.317)%
  --(2.783,1.282)--(2.825,1.246)--(2.866,1.211)--(2.908,1.176)--(2.950,1.141)--(2.991,1.106)%
  --(3.033,1.070)--(3.075,1.035)--(3.116,1.000)--(3.158,0.965)--(3.200,0.929)--(3.241,0.894)%
  --(3.283,0.859)--(3.325,0.824)--(3.366,0.788)--(3.408,0.753)--(3.450,0.718)--(3.491,0.683)%
  --(3.533,0.647)--(3.575,0.612)--(3.598,0.592);
\gpcolor{color=gp lt color border}
\draw[gp path] (0.825,4.401)--(0.825,0.592)--(5.294,0.592)--(5.294,4.401)--cycle;
\gpdefrectangularnode{gp plot 1}{\pgfpoint{0.825cm}{0.592cm}}{\pgfpoint{5.294cm}{4.401cm}}
\end{tikzpicture}
\label{MS_time2}}
	\subfloat[Elapsed time ratio]{
\begin{tikzpicture}[gnuplot]
\tikzset{every node/.append style={scale=0.60}}
\path (0.000,0.000) rectangle (5.750,4.812);
\gpcolor{color=gp lt color border}
\gpsetlinetype{gp lt border}
\gpsetdashtype{gp dt solid}
\gpsetlinewidth{1.00}
\draw[gp path] (0.680,0.592)--(0.770,0.592);
\draw[gp path] (5.419,0.592)--(5.329,0.592);
\draw[gp path] (0.680,0.760)--(0.770,0.760);
\draw[gp path] (5.419,0.760)--(5.329,0.760);
\draw[gp path] (0.680,0.902)--(0.770,0.902);
\draw[gp path] (5.419,0.902)--(5.329,0.902);
\draw[gp path] (0.680,1.025)--(0.770,1.025);
\draw[gp path] (5.419,1.025)--(5.329,1.025);
\draw[gp path] (0.680,1.133)--(0.770,1.133);
\draw[gp path] (5.419,1.133)--(5.329,1.133);
\gpcolor{color=gp lt color axes}
\gpsetlinetype{gp lt axes}
\gpsetdashtype{gp dt axes}
\gpsetlinewidth{0.50}
\draw[gp path] (0.680,1.230)--(5.419,1.230);
\gpcolor{color=gp lt color border}
\gpsetlinetype{gp lt border}
\gpsetdashtype{gp dt solid}
\gpsetlinewidth{1.00}
\draw[gp path] (0.680,1.230)--(0.860,1.230);
\draw[gp path] (5.419,1.230)--(5.239,1.230);
\node[gp node right] at (0.570,1.230) {$1$};
\draw[gp path] (0.680,1.868)--(0.770,1.868);
\draw[gp path] (5.419,1.868)--(5.329,1.868);
\draw[gp path] (0.680,2.241)--(0.770,2.241);
\draw[gp path] (5.419,2.241)--(5.329,2.241);
\draw[gp path] (0.680,2.506)--(0.770,2.506);
\draw[gp path] (5.419,2.506)--(5.329,2.506);
\draw[gp path] (0.680,2.712)--(0.770,2.712);
\draw[gp path] (5.419,2.712)--(5.329,2.712);
\draw[gp path] (0.680,2.880)--(0.770,2.880);
\draw[gp path] (5.419,2.880)--(5.329,2.880);
\draw[gp path] (0.680,3.021)--(0.770,3.021);
\draw[gp path] (5.419,3.021)--(5.329,3.021);
\draw[gp path] (0.680,3.144)--(0.770,3.144);
\draw[gp path] (5.419,3.144)--(5.329,3.144);
\draw[gp path] (0.680,3.253)--(0.770,3.253);
\draw[gp path] (5.419,3.253)--(5.329,3.253);
\gpcolor{color=gp lt color axes}
\gpsetlinetype{gp lt axes}
\gpsetdashtype{gp dt axes}
\gpsetlinewidth{0.50}
\draw[gp path] (0.680,3.350)--(5.419,3.350);
\gpcolor{color=gp lt color border}
\gpsetlinetype{gp lt border}
\gpsetdashtype{gp dt solid}
\gpsetlinewidth{1.00}
\draw[gp path] (0.680,3.350)--(0.860,3.350);
\draw[gp path] (5.419,3.350)--(5.239,3.350);
\node[gp node right] at (0.570,3.350) {$10$};
\draw[gp path] (0.680,3.988)--(0.770,3.988);
\draw[gp path] (5.419,3.988)--(5.329,3.988);
\draw[gp path] (0.680,4.361)--(0.770,4.361);
\draw[gp path] (5.419,4.361)--(5.329,4.361);
\draw[gp path] (0.680,4.626)--(0.770,4.626);
\draw[gp path] (5.419,4.626)--(5.329,4.626);
\gpcolor{color=gp lt color axes}
\gpsetlinetype{gp lt axes}
\gpsetdashtype{gp dt axes}
\gpsetlinewidth{0.50}
\draw[gp path] (0.680,0.592)--(0.680,4.626);
\gpcolor{color=gp lt color border}
\gpsetlinetype{gp lt border}
\gpsetdashtype{gp dt solid}
\gpsetlinewidth{1.00}
\draw[gp path] (0.680,0.592)--(0.680,0.772);
\draw[gp path] (0.680,4.626)--(0.680,4.446);
\node[gp node center,font={\fontsize{7.0pt}{8.4pt}\selectfont}] at (0.680,0.407) {16};
\gpcolor{color=gp lt color axes}
\gpsetlinetype{gp lt axes}
\gpsetdashtype{gp dt axes}
\gpsetlinewidth{0.50}
\draw[gp path] (1.305,0.592)--(1.305,4.626);
\gpcolor{color=gp lt color border}
\gpsetlinetype{gp lt border}
\gpsetdashtype{gp dt solid}
\gpsetlinewidth{1.00}
\draw[gp path] (1.305,0.592)--(1.305,0.772);
\draw[gp path] (1.305,4.626)--(1.305,4.446);
\node[gp node center,font={\fontsize{7.0pt}{8.4pt}\selectfont}] at (1.305,0.407) {32};
\gpcolor{color=gp lt color axes}
\gpsetlinetype{gp lt axes}
\gpsetdashtype{gp dt axes}
\gpsetlinewidth{0.50}
\draw[gp path] (1.930,0.592)--(1.930,4.626);
\gpcolor{color=gp lt color border}
\gpsetlinetype{gp lt border}
\gpsetdashtype{gp dt solid}
\gpsetlinewidth{1.00}
\draw[gp path] (1.930,0.592)--(1.930,0.772);
\draw[gp path] (1.930,4.626)--(1.930,4.446);
\node[gp node center,font={\fontsize{7.0pt}{8.4pt}\selectfont}] at (1.930,0.407) {64};
\gpcolor{color=gp lt color axes}
\gpsetlinetype{gp lt axes}
\gpsetdashtype{gp dt axes}
\gpsetlinewidth{0.50}
\draw[gp path] (2.554,0.592)--(2.554,4.626);
\gpcolor{color=gp lt color border}
\gpsetlinetype{gp lt border}
\gpsetdashtype{gp dt solid}
\gpsetlinewidth{1.00}
\draw[gp path] (2.554,0.592)--(2.554,0.772);
\draw[gp path] (2.554,4.626)--(2.554,4.446);
\node[gp node center,font={\fontsize{7.0pt}{8.4pt}\selectfont}] at (2.554,0.407) {128};
\gpcolor{color=gp lt color axes}
\gpsetlinetype{gp lt axes}
\gpsetdashtype{gp dt axes}
\gpsetlinewidth{0.50}
\draw[gp path] (3.179,0.592)--(3.179,4.626);
\gpcolor{color=gp lt color border}
\gpsetlinetype{gp lt border}
\gpsetdashtype{gp dt solid}
\gpsetlinewidth{1.00}
\draw[gp path] (3.179,0.592)--(3.179,0.772);
\draw[gp path] (3.179,4.626)--(3.179,4.446);
\node[gp node center,font={\fontsize{7.0pt}{8.4pt}\selectfont}] at (3.179,0.407) {256};
\gpcolor{color=gp lt color axes}
\gpsetlinetype{gp lt axes}
\gpsetdashtype{gp dt axes}
\gpsetlinewidth{0.50}
\draw[gp path] (3.804,0.592)--(3.804,4.626);
\gpcolor{color=gp lt color border}
\gpsetlinetype{gp lt border}
\gpsetdashtype{gp dt solid}
\gpsetlinewidth{1.00}
\draw[gp path] (3.804,0.592)--(3.804,0.772);
\draw[gp path] (3.804,4.626)--(3.804,4.446);
\node[gp node center,font={\fontsize{7.0pt}{8.4pt}\selectfont}] at (3.804,0.407) {512};
\gpcolor{color=gp lt color axes}
\gpsetlinetype{gp lt axes}
\gpsetdashtype{gp dt axes}
\gpsetlinewidth{0.50}
\draw[gp path] (4.429,0.592)--(4.429,4.626);
\gpcolor{color=gp lt color border}
\gpsetlinetype{gp lt border}
\gpsetdashtype{gp dt solid}
\gpsetlinewidth{1.00}
\draw[gp path] (4.429,0.592)--(4.429,0.772);
\draw[gp path] (4.429,4.626)--(4.429,4.446);
\node[gp node center,font={\fontsize{7.0pt}{8.4pt}\selectfont}] at (4.429,0.407) {1024};
\gpcolor{color=gp lt color axes}
\gpsetlinetype{gp lt axes}
\gpsetdashtype{gp dt axes}
\gpsetlinewidth{0.50}
\draw[gp path] (5.054,0.592)--(5.054,4.626);
\gpcolor{color=gp lt color border}
\gpsetlinetype{gp lt border}
\gpsetdashtype{gp dt solid}
\gpsetlinewidth{1.00}
\draw[gp path] (5.054,0.592)--(5.054,0.772);
\draw[gp path] (5.054,4.626)--(5.054,4.446);
\node[gp node center,font={\fontsize{7.0pt}{8.4pt}\selectfont}] at (5.054,0.407) {2048};
\gpcolor{color=gp lt color axes}
\gpsetlinetype{gp lt axes}
\gpsetdashtype{gp dt axes}
\gpsetlinewidth{0.50}
\draw[gp path] (5.419,0.592)--(5.419,4.626);
\gpcolor{color=gp lt color border}
\gpsetlinetype{gp lt border}
\gpsetdashtype{gp dt solid}
\gpsetlinewidth{1.00}
\draw[gp path] (5.419,0.592)--(5.419,0.772);
\draw[gp path] (5.419,4.626)--(5.419,4.446);
\node[gp node center,font={\fontsize{7.0pt}{8.4pt}\selectfont}] at (5.419,0.407) {3072};
\draw[gp path] (0.680,4.626)--(0.680,0.592)--(5.419,0.592)--(5.419,4.626)--cycle;
\node[gp node center,rotate=-270] at (0.194,2.609) {ratio (elapsed time)/(best elapsed time) };
\node[gp node center] at (3.049,0.068) {number of cores};
\gpcolor{rgb color={0.000,0.392,0.000}}
\draw[gp path] (0.680,1.789)--(1.305,1.872)--(1.930,1.974)--(2.554,2.093)--(3.179,2.155)%
  --(3.804,2.331)--(4.429,2.383)--(5.054,2.444)--(5.419,2.560);
\gpsetpointsize{4.00}
\gp3point{gp mark 7}{}{(0.680,1.789)}
\gp3point{gp mark 7}{}{(1.305,1.872)}
\gp3point{gp mark 7}{}{(1.930,1.974)}
\gp3point{gp mark 7}{}{(2.554,2.093)}
\gp3point{gp mark 7}{}{(3.179,2.155)}
\gp3point{gp mark 7}{}{(3.804,2.331)}
\gp3point{gp mark 7}{}{(4.429,2.383)}
\gp3point{gp mark 7}{}{(5.054,2.444)}
\gp3point{gp mark 7}{}{(5.419,2.560)}
\gpcolor{rgb color={0.000,1.000,0.000}}
\draw[gp path] (0.680,1.230)--(1.305,1.230)--(1.930,1.230)--(2.554,1.230)--(3.179,1.230)%
  --(3.804,1.230)--(4.429,1.230)--(5.054,1.230)--(5.419,1.230);
\gp3point{gp mark 3}{}{(0.680,1.230)}
\gp3point{gp mark 3}{}{(1.305,1.230)}
\gp3point{gp mark 3}{}{(1.930,1.230)}
\gp3point{gp mark 3}{}{(2.554,1.230)}
\gp3point{gp mark 3}{}{(3.179,1.230)}
\gp3point{gp mark 3}{}{(3.804,1.230)}
\gp3point{gp mark 3}{}{(4.429,1.230)}
\gp3point{gp mark 3}{}{(5.054,1.230)}
\gp3point{gp mark 3}{}{(5.419,1.230)}
\gpcolor{rgb color={0.753,0.251,0.000}}
\draw[gp path] (3.804,4.188)--(4.429,3.662)--(5.054,3.312)--(5.419,2.950);
\gp3point{gp mark 2}{}{(3.804,4.188)}
\gp3point{gp mark 2}{}{(4.429,3.662)}
\gp3point{gp mark 2}{}{(5.054,3.312)}
\gp3point{gp mark 2}{}{(5.419,2.950)}
\gpcolor{rgb color={1.000,0.647,0.000}}
\draw[gp path] (3.804,4.050)--(4.429,3.497)--(5.054,3.155)--(5.419,2.800);
\gp3point{gp mark 2}{}{(3.804,4.050)}
\gp3point{gp mark 2}{}{(4.429,3.497)}
\gp3point{gp mark 2}{}{(5.054,3.155)}
\gp3point{gp mark 2}{}{(5.419,2.800)}
\gpcolor{color=gp lt color border}
\draw[gp path] (0.680,4.626)--(0.680,0.592)--(5.419,0.592)--(5.419,4.626)--cycle;
\gpdefrectangularnode{gp plot 1}{\pgfpoint{0.680cm}{0.592cm}}{\pgfpoint{5.419cm}{4.626cm}}
\end{tikzpicture}
\label{MS_time3}}
	\caption{Micro structural performances: Elapsed time and performance ration for the first and second computation.\label{MS_curves}.}
\end{figure}

The elapsed times for this first computation are given in the figure \ref{MS_time1}.
It shows that, as in the equivalent test case of section \ref{UOL}, the \tS solver provides better performance compared to the "dd" solver.
The strong scaling efficiency curves (not presented here)  show that the \tS solver  also has a smaller slope than the "dd" solver.

Let us now  consider a broader application framework and imagine that this problem is the subject of a micro-structure optimization study with respect to a displacement field based objective.
In this case, many calculations from a given set of parameters must be performed to obtain sensitivities and/or new configurations.
To mimic such a case, the Young's modulus values  are randomly perturbed by $\pm1\%$ and the computation starts from the first result obtained above.
For the \tS solver, this implies to use the r-set solution of  the previous calculation,  to redo the factorization of the fine-scale problems and to enter the scale loop.
For the "dd" solver, the domains are condensed again, the preconditioners recomputed, and since the Schur complement relies on the same space, the iterative solver starts from the previous boundary problem solution.
Figure \ref{MS_time2} shows the performance for this second computation ("tsi (var 1\%)" and "dd ( var 1\%)" for the \tS and domain decomposition solvers respectively).
In figure \ref{MS_time3}, the relative performance of the two computations is presented as a ratio of the elapsed time of one computation to the best elapsed time among all computations and solvers.
The \tS solver has a reduced elapsed time because it iterated only  37 times, compared to 93 times for the initial calculation.
This can be explained by the fact that the initial solution for this computation is already of better quality than what the unenriched coarse  field can give ($resi=1.5e^{-2}$ compare to $resi=1.6$ for the first computation) even if the Young's moduli have changed. 
Moreover, the local variation of E may only influence the modulus and not the shape of the enrichment functions.
But starting from a better solution also reduces the number of \tS iterations where the full rank solver is used at the  coarse scale.
In the first computation $resi<\epsilon\times 10000$ (rule from section \ref{gsolv}) is met after 25 iterations.
For the second computation, this condition is met only after 4 iterations.
This explains why the ratio of "ts" in figure \ref{MS_time3} grows with the number of processes: as already observed in the analysis of figure\ref{UOL_curve_percent}, when a high number of cores is reached, the reduced set of $nbp_{max}$ processes has a negative impact on the factorization performance.
Thus, this poor performance has an impact only 4 times on the second computation against  25 times for the first. 
For the "dd" solver, restarting from the previous global boundary solution reduces the number of iterations in solving the global problem but not enough to  significantly decrease  elapsed time.
The second "dd" computation cost almost the same as the first one, which gives an advantage to \tS solver in this case.

To obtain an even worse conditioned system,  exactly the same scenario is replayed with a wider range of E values: from 36.5 GPa  to 36500 GPa (maximum ten times larger than previous one).
Results  with 2048 cores and both ranges, are given in the table \ref{MS_tab_bad_cond}.   
\begin{figure}[h]
	\subfloat[Curves]{

\label{MS_converge}}
	\subfloat[table]{\small\begin{tabular}[b]{|c|c|c|c|c|c|c|}
  \hline
  \multirow{3}{*}{E  (GPa)}&\multirow{2}{*}{run}&\multicolumn{3}{c}{\tS} & \multicolumn{2}{|c|}{dd}  \\
\cline{3-7}&&\multicolumn{2}{c|}{nb of iterations}&\multirow{2}{*}{ time (s)} & \multirow{2}{*}{nb of iterations}&\multirow{2}{*}{ time (s)} \\
\cline{3-4}&&all& $resi\geqslant10^{-3}$&&& \\
\hline
  \multirow{2}{*}{$[36.5,3650]$}&first&93&25&198.6&1868&509.9\\
   \cline{2-7}&second&37&4&53.1&1487&429.8\\
  \hline
  \multirow{2}{*}{$[36.5,36500]$}&first&265&96&691.4&2191&621.4\\
  \cline{2-7}&second&134&7&128.8&1487&519.4\\
  \hline
E  (GPa)&run&all& $resi\geqslant10^{-1}$&time (s)&\multicolumn{2}{l}{~} \\
  \cline{1-5}
  \multirow{2}{*}{$[36.5,36500]$}&first&265&27&452.8&\multicolumn{2}{l}{~}\\
  \cline{2-5}&second&134&0&149.4&\multicolumn{2}{l}{~}\\
  \cline{1-5}
	\end{tabular}\label{MS_tab_bad_cond}}
\caption{With 2048 cores, comparison of the first and second micro-structure calculations with different Young's modulus ranges."nb of iterations" is the number of iteration of the loop of the algorithm \ref{TS_algebra} and \ref{domain_decomposition_algo}. "time" is the elapsed time to solve the problem. \label{MS_bad_cond}}
\end{figure}
With the new E-range, the number of \tS iterations for the first computation (the one  starting from  the coarse unenriched solution) increases from 93 to 265.
This increase is related to a worse starting point and  conditioning, which degraded the residual convergence rate, as can be seen in figure \ref{MS_converge}.
In this figure, the slope (after a few iterations) of the first computation with a smaller E-range (thus with better matrix conditioning) is stronger than that of the first computation with a larger E-range (thus with worse matrix conditioning).
This is confirmed by the second calculation where the curves start with a smaller residual error but have the same slope as their  counterpart in the first calculation.
This brings explanations to the observation made in section \ref{UOL} where a finer discretization of the problem at the global scale  gives a slower residual convergence (see the  slopes of the curves in figure \ref{UOL_curve_error}).
In fact, we can now add  that it is the conditioning of the global matrix (related to the level of discretization) that has an impact on the convergence.

Regarding  performance, like the \tS solver, the "dd" solver needs more iterations to converge (the worse conditioning is not fully corrected by domain condensation).
For first calculation, "dd" needs 2191 iterations with the new E-range against  1868 with the original E-range.
The elapsed time is moderately impacted and increase to 621.4s from 509.9s.
The second calculation proceeds in the same way and there is a moderate gain  starting with the first calculation (for both E-range).
For \tS solver, this increasing number of iterations has a greater impact on the first calculation.
The elapsed time goes from 198.6s to 691.4s which gives an advantage to "dd" in this case.
This  increase in time is in fact closely related to the number of iterations where the global problem is computed by the full rank solver.
With \TSI version, the iterative resolution is activated when $resi<\epsilon \times 1000$  (rule of the section \ref{gsolv}) and with $\epsilon=1e^{-7}$ the condition is $resi<1e^{-3}$.
The table \ref{MS_tab_bad_cond} gives the number of iterations that do not meet this condition and therefore only correspond to a full rank resolution.
For the first calculation, 96 iterations are performed in full rank mode with the new E-range compared to 25.
If we change the condition to $resi<0.1$, only 27 iterations are spent in full rank mode and the elapsed time drops to 452.8s which is now better than "dd".
This shows that the \TSI setting can be improved in some cases.
But even with this new setting, some computations not reported here, with a larger E-range, shows that the \tS solver is in trouble compared to "dd" when the matrix conditioning degrades further.
In all cases, the second calculation with the \tS solver, with both settings and both E-ranges, has a better elapsed time than the first calculation.
This confirms the interest of using the \tS solver when a history can be kept from one computation to another.

\subsection{Pull-out test case}\label{PO}
\begin{figure}[h]
	\centering
	\subfloat[geometry]{\def\svgwidth{67mm}
\begingroup%
  \makeatletter%
  \providecommand\color[2][]{%
    \errmessage{(Inkscape) Color is used for the text in Inkscape, but the package 'color.sty' is not loaded}%
    \renewcommand\color[2][]{}%
  }%
  \providecommand\transparent[1]{%
    \errmessage{(Inkscape) Transparency is used (non-zero) for the text in Inkscape, but the package 'transparent.sty' is not loaded}%
    \renewcommand\transparent[1]{}%
  }%
  \providecommand\rotatebox[2]{#2}%
  \ifx\svgwidth\undefined%
    \setlength{\unitlength}{449.82714844bp}%
    \ifx\svgscale\undefined%
      \relax%
    \else%
      \setlength{\unitlength}{\unitlength * \real{\svgscale}}%
    \fi%
  \else%
    \setlength{\unitlength}{\svgwidth}%
  \fi%
  \global\let\svgwidth\undefined%
  \global\let\svgscale\undefined%
  \makeatother%
  \begin{picture}(1,0.86133345)%
    \lineheight{1}%
    \setlength\tabcolsep{0pt}%
    \put(0,0){\includegraphics[width=\unitlength,page=1]{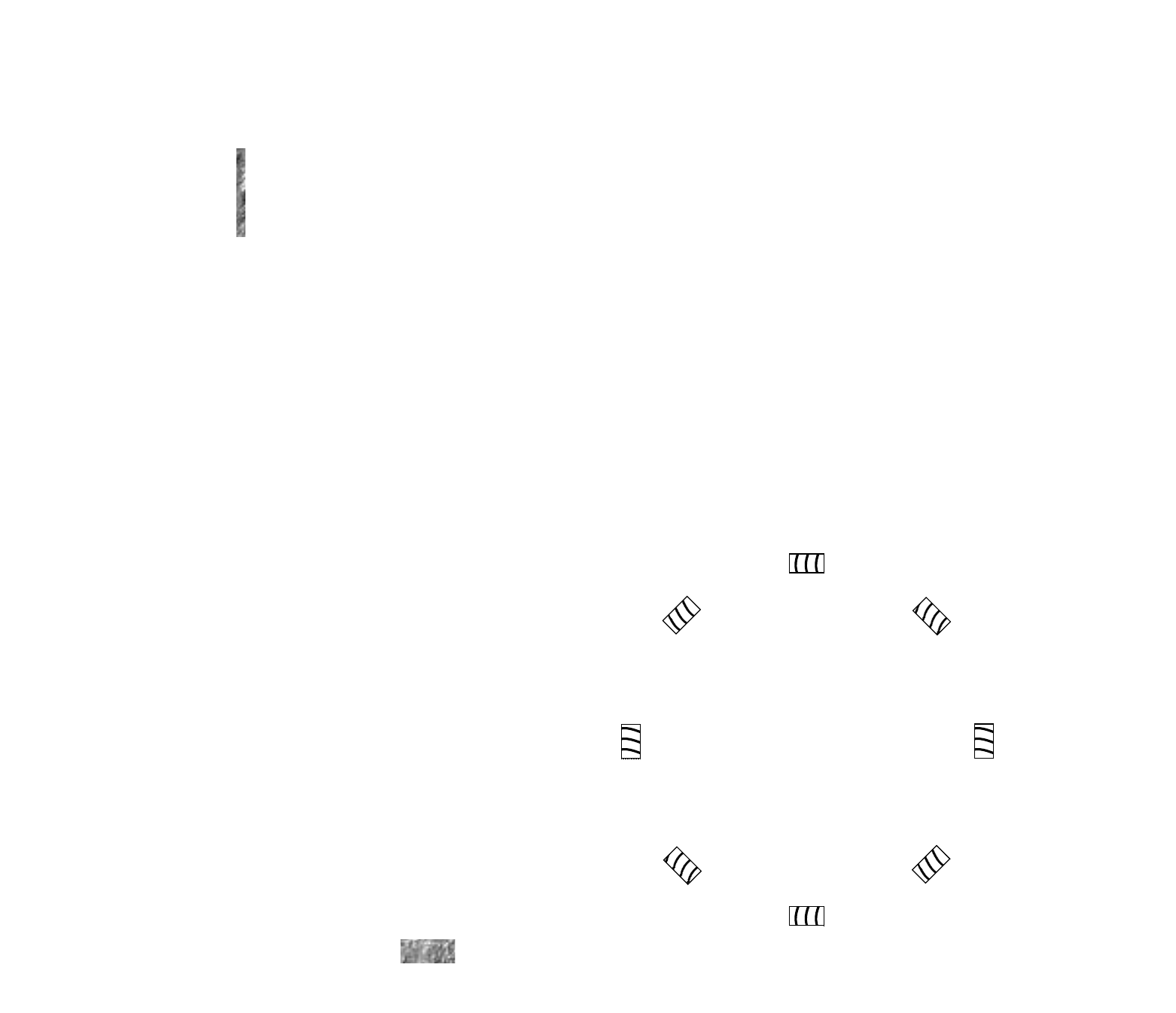}}%
    \put(0.38973457,0.04431764){\color[rgb]{0,0,0}\makebox(0,0)[lt]{\lineheight{0}\smash{\begin{tabular}[t]{l}\tiny{Dirichlet: $U_x=U_z=0$}\end{tabular}}}}%
    \put(0,0){\includegraphics[width=\unitlength,page=2]{Fig19a.pdf}}%
    \put(0.38951107,0.00637104){\color[rgb]{0,0,0}\makebox(0,0)[lt]{\lineheight{0}\smash{\begin{tabular}[t]{l}\tiny{Neuman: $1000.\overrightarrow{e_y}$}\end{tabular}}}}%
    \put(0,0){\includegraphics[width=\unitlength,page=3]{Fig19a.pdf}}%
    \put(0.38973457,0.08426818){\color[rgb]{0,0,0}\makebox(0,0)[lt]{\lineheight{0}\smash{\begin{tabular}[t]{l}\tiny{Dirichlet: $U_y=0$}\end{tabular}}}}%
    \put(0,0){\includegraphics[width=\unitlength,page=4]{Fig19a.pdf}}%
    \put(0.92403746,0.74879901){\color[rgb]{0,0,0}\makebox(0,0)[lt]{\lineheight{0}\smash{\begin{tabular}[t]{l}\tiny{$x$}\end{tabular}}}}%
    \put(0.99999997,0.50719398){\color[rgb]{0,0,0}\makebox(0,0)[lt]{\lineheight{0}\smash{\begin{tabular}[t]{l} \end{tabular}}}}%
    \put(0,0){\includegraphics[width=\unitlength,page=5]{Fig19a.pdf}}%
    \put(0.70208066,0.84661467){\color[rgb]{0,0,0}\makebox(0,0)[lt]{\lineheight{0}\smash{\begin{tabular}[t]{l}\tiny{$y$}\end{tabular}}}}%
    \put(0,0){\includegraphics[width=\unitlength,page=6]{Fig19a.pdf}}%
    \put(0.6994439,0.0075298){\color[rgb]{0,0,0}\makebox(0,0)[lt]{\lineheight{0}\smash{\begin{tabular}[t]{l}\tiny{$z$}\end{tabular}}}}%
    \put(0,0){\includegraphics[width=\unitlength,page=7]{Fig19a.pdf}}%
    \put(0.69121177,0.7394422){\color[rgb]{0,0,0}\makebox(0,0)[lt]{\lineheight{0}\smash{\begin{tabular}[t]{l}\tiny{$0$}\end{tabular}}}}%
    \put(0.70015711,0.2353667){\color[rgb]{0,0,0}\makebox(0,0)[lt]{\lineheight{0}\smash{\begin{tabular}[t]{l}\tiny{$0$}\end{tabular}}}}%
    \put(0.91950671,0.20661995){\color[rgb]{0,0,0}\makebox(0,0)[lt]{\lineheight{0}\smash{\begin{tabular}[t]{l}\tiny{$x$}\end{tabular}}}}%
    \put(0,0){\includegraphics[width=\unitlength,page=8]{Fig19a.pdf}}%
    \put(0.21023509,0.73842732){\color[rgb]{0,0,0}\makebox(0,0)[lt]{\lineheight{0}\smash{\begin{tabular}[t]{l}\tiny{$0$}\end{tabular}}}}%
    \put(0.16948071,0.84276947){\color[rgb]{0,0,0}\makebox(0,0)[lt]{\lineheight{0}\smash{\begin{tabular}[t]{l}\tiny{$y$}\end{tabular}}}}%
    \put(0.00441041,0.75262728){\color[rgb]{0,0,0}\makebox(0,0)[lt]{\lineheight{0}\smash{\begin{tabular}[t]{l}\tiny{$z$}\end{tabular}}}}%
    \put(0.28879544,0.37643839){\color[rgb]{0,0,0}\makebox(0,0)[lt]{\lineheight{0}\smash{\begin{tabular}[t]{l}\tiny{M}\end{tabular}}}}%
    \put(0.28685339,0.31195628){\color[rgb]{0,0,0}\makebox(0,0)[lt]{\lineheight{0}\smash{\begin{tabular}[t]{l}\tiny{P}\end{tabular}}}}%
    \put(0.1128242,0.37710312){\color[rgb]{0,0,0}\makebox(0,0)[lt]{\lineheight{0}\smash{\begin{tabular}[t]{l}\tiny{Q}\end{tabular}}}}%
    \put(0.84131994,0.56411315){\color[rgb]{0,0,0}\makebox(0,0)[lt]{\lineheight{0}\smash{\begin{tabular}[t]{l}\tiny{R}\end{tabular}}}}%
    \put(0.83857454,0.74351386){\color[rgb]{0,0,0}\makebox(0,0)[lt]{\lineheight{0}\smash{\begin{tabular}[t]{l}\tiny{S}\end{tabular}}}}%
    \put(0,0){\includegraphics[width=\unitlength,page=9]{Fig19a.pdf}}%
    \put(0.32703438,0.30079148){\color[rgb]{0,0,0}\makebox(0,0)[lt]{\lineheight{0}\smash{\begin{tabular}[t]{l}\textbf{\footnotesize{D}}\end{tabular}}}}%
    \put(0,0){\includegraphics[width=\unitlength,page=10]{Fig19a.pdf}}%
    \put(0.03056902,0.30079148){\color[rgb]{0,0,0}\makebox(0,0)[lt]{\lineheight{0}\smash{\begin{tabular}[t]{l}\textbf{\footnotesize{D}}\end{tabular}}}}%
    \put(0,0){\includegraphics[width=\unitlength,page=11]{Fig19a.pdf}}%
    \put(0.02483816,0.22899531){\color[rgb]{0,0,0}\makebox(0,0)[lt]{\lineheight{0}\smash{\begin{tabular}[t]{l}\textbf{\footnotesize{D-D}}\end{tabular}}}}%
    \put(0,0){\includegraphics[width=\unitlength,page=12]{Fig19a.pdf}}%
    \put(0.70211637,0.5036235){\color[rgb]{0,0,0}\makebox(0,0)[lt]{\lineheight{0}\smash{\begin{tabular}[t]{l}\textbf{\footnotesize{C}}\end{tabular}}}}%
    \put(0.03988759,0.79549295){\color[rgb]{0,0,0}\makebox(0,0)[lt]{\lineheight{0}\smash{\begin{tabular}[t]{l}\textbf{\footnotesize{C-C}}\end{tabular}}}}%
    \put(0.70211637,0.81557564){\color[rgb]{0,0,0}\makebox(0,0)[lt]{\lineheight{0}\smash{\begin{tabular}[t]{l}\textbf{\footnotesize{C}}\end{tabular}}}}%
    \put(0,0){\includegraphics[width=\unitlength,page=13]{Fig19a.pdf}}%
  \end{picture}%
\endgroup%
}
	\subfloat[general behavior]{\includegraphics[width=67mm]{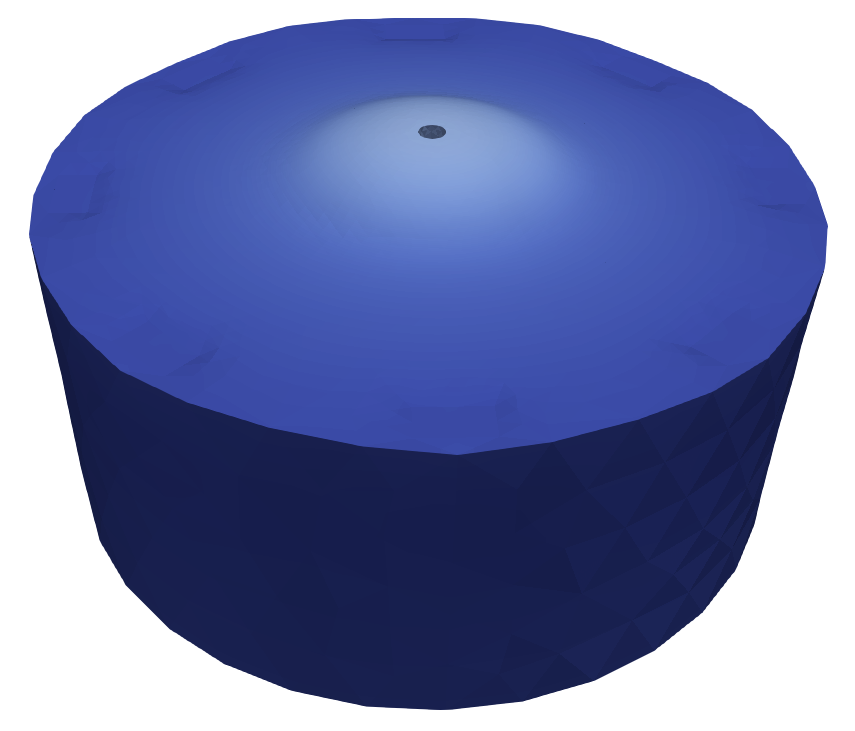}\label{PO_general}}
	\caption{Pull-out test case geometry, boundary conditions and general behavior (obtained by one simulation of this section): \\M(0,-450,-60),P(0,-470,-60),Q(0,-450,36),R(1000,-900,0),S(1000,0,0) unit=mm \label{PO-geo}}
\end{figure}

This test case is inspired by the crack propagation simulation studied in \cite{Eligehausen1989AFM,Ozbolt1999,Gasser2005,Bordas2008:01}.
It represents the pull-out of a steel anchor embedded in an unreinforced  concrete (E=26.4GPa, $\nu$=0.193) cylinder.
The chosen  geometry and boundary conditions used in this test are given in figure \ref{PO-geo}.
The steel anchor, not shown in this figure, is located in the hollowing out in the center of the disk and is pulled in the $\overrightarrow{e_y}$ direction.
Its action  is simply modeled as a force applied to the contact surface between the compression surface of the anchor head  and the concrete.
Other  interactions of the anchor with the concrete are simplified because contact, decohesion or other mechanical bond phenomena are beyond the scope of this study.
Thus, all remaining surfaces of the anchor head are considered to have no interaction with the concrete which is left as free surfaces in these areas.
The surface of the anchor body  is considered to impose a sliding contact interface along the $\overrightarrow{e_y}$ axis.

\begin{figure}[h]
	\subfloat[step 0]{\includegraphics[width=40mm]{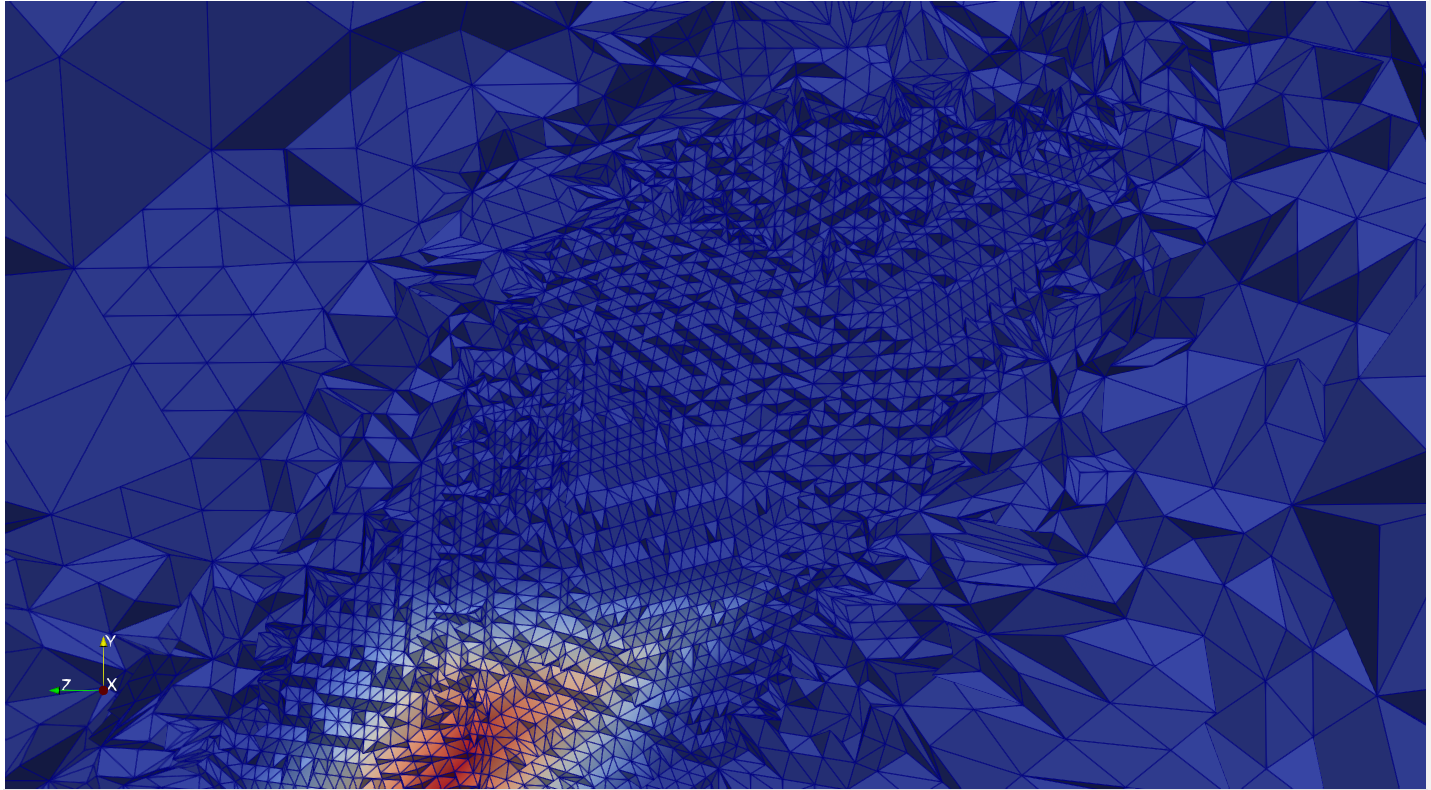}\label{PO_dam_step0}}
	\subfloat[step 1]{\includegraphics[width=40mm]{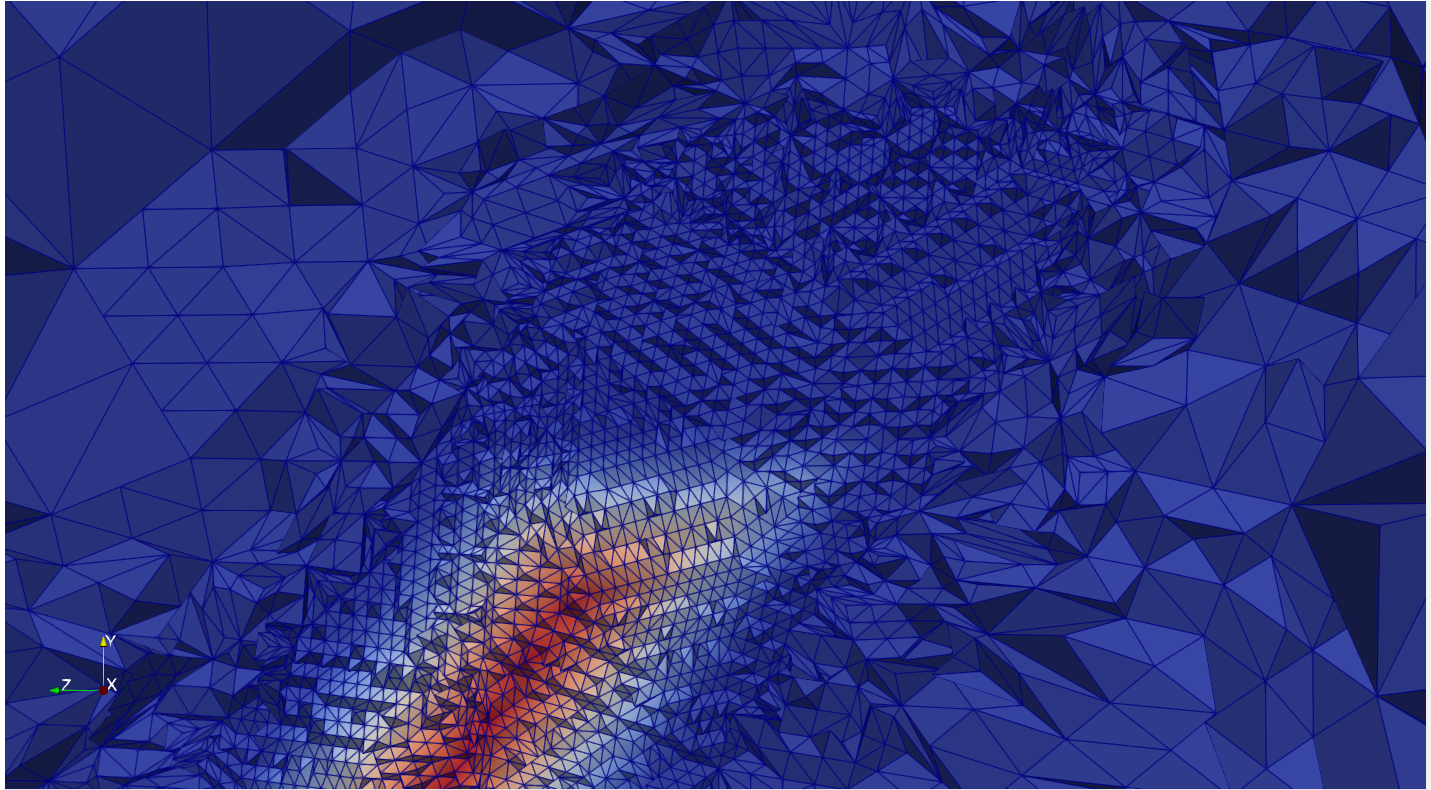}\label{PO_dam_step1}}
	\subfloat[step 2]{\includegraphics[width=40.1mm]{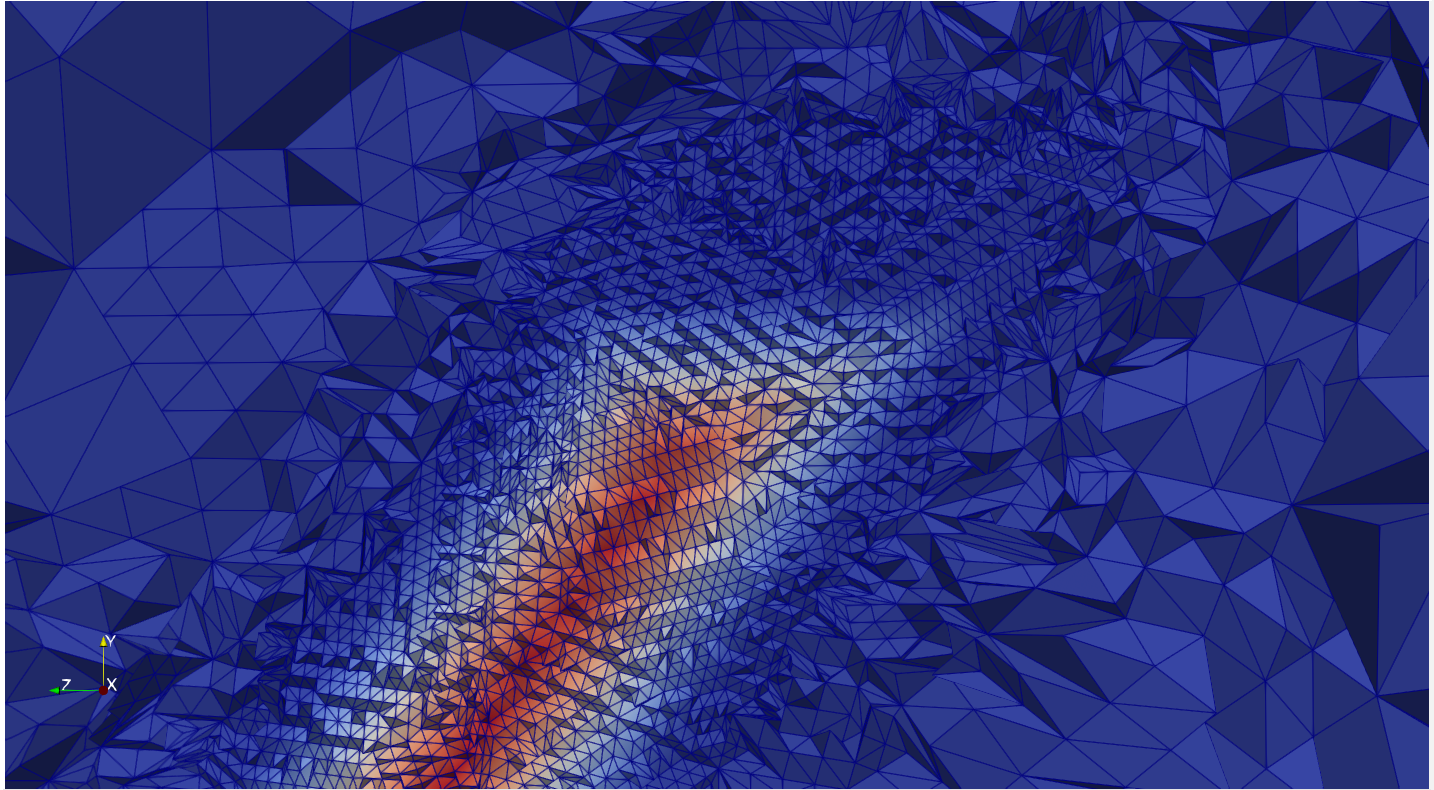}\label{PO_dam_step2}}
	\subfloat[step 3]{\includegraphics[width=45.4mm]{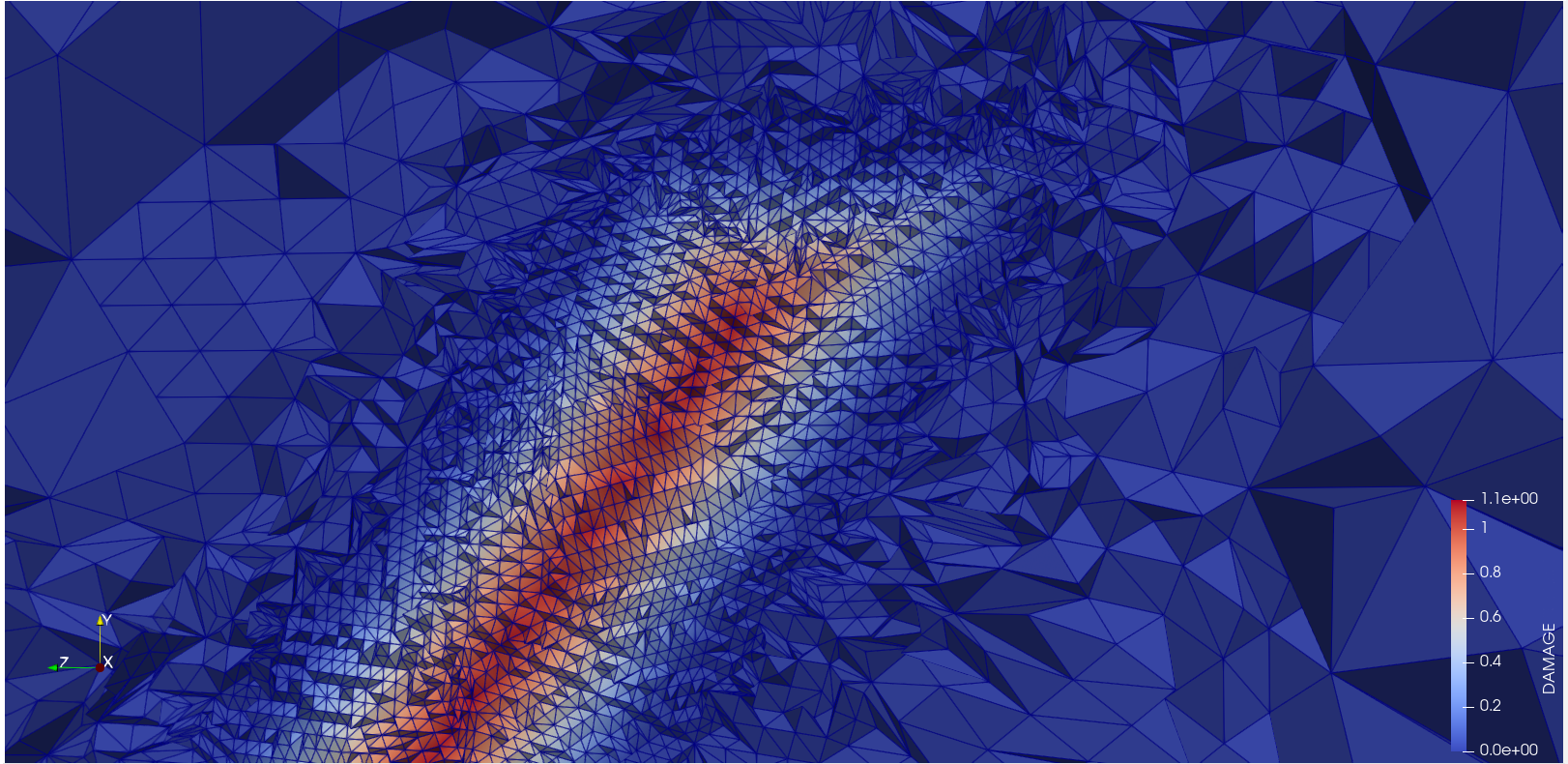}\label{PO_dam_step3}}
	\caption{Evolution of the imposed damage field for the $h=-400$ adaptation. Cross-section of the SP. Zoom of the crack front (with a rather unrealistic shape as mentioned in \ref{PO_cone_anexe}).\label{PO_dam_steps}}	
\end{figure}

The crack is represented by a damage field $d$.
The behavior of the material thus follows the linear elastic damage potential $\varphi(\tens{\epsilon},d)$ which is written,  using the Hooke tensor $\tens{C}$, as follows:
\begin{equation}
	\varphi(\tens{\epsilon},d)=\frac{1}{2}(1- d)\tens{\epsilon}:\tens{C}:\tens{\epsilon}
\end{equation}
When $d=1$, the material has lost all its rigidity and can be considered as a crack.
In the pull-out test, the crack starts under the head of the anchor  and develops at a specific angle in the form of a cone.
The more you pull on the steel anchor, the more the crack continues to grow in a cone.
The classic scenario for simulating such evolutionary phenomena is to use mesh adaptation to follow the movement of the crack front so that the damaged area is always entirely within a finer mesh region.
By always refining  a little ahead of the front, the number of mesh adaptations  can be limited,  keeping the same mesh for a few evolutionary iterations called step below.
To study such a scenario with the \tS method, we consider three arbitrary evolution stages  that correspond to three coarse mesh adaptations.
And for each of them, four fictitious steps advance the crack front evenly according to the cone angle in the adapted fixed mesh.
The damage field is  calculated by considering the conical envelope described in  \ref{PO_cone_anexe} with the parameter $h$ in \eqref{PO_annex_cone_eq}  controlling the crack cone size.
The coarse mesh is adapted fo $h=-400$ (50 874
c-set dofs), $h=-300$ (119 829
c-set dofs ) and $h=-100$ (362 904
c-set dofs) considering that $h$  during the steps will fluctuate between $[-410,-400]$, $[-310,-300]$ and $[-110,-100]$ respectively.
The mesh sizes of the fine problems follow the equation \eqref{ts_jump_level_eq}.
This corresponds approximately to 4, 16 and 58.6 million  r-set dofs  for the $h=-400$, $h=-300$ and $h=-100$ stages respectively.
The \SP is built at each stage to cover only the damaged area during the 4 evolutionary iterations as shown in figure \ref{PO_dam_steps}.
So, in this test case, there are NSP elements and transition patches (the one with hanging node treatment).
The \tS resolution is performed 4 times.
The first one starts from a zero displacement field (in a real simulation, it would be a displacement field obtained by projection of the result of the last step of the last stage).
The next three \tS resolutions will start with the displacement field calculated in the previous iteration.
Here, since  the \SP discretization remains the same, the displacement fields at both scales are on a stable space.
Only some patches need to be re-factored  to  account for material degradation related to damage propagation.
The precision chosen,  as in the other test cases, is $\epsilon=1.e^{-7}$.

\begin{figure}[h!]
	\subfloat[Global scale,h=-400]{\includegraphics[width=55mm]{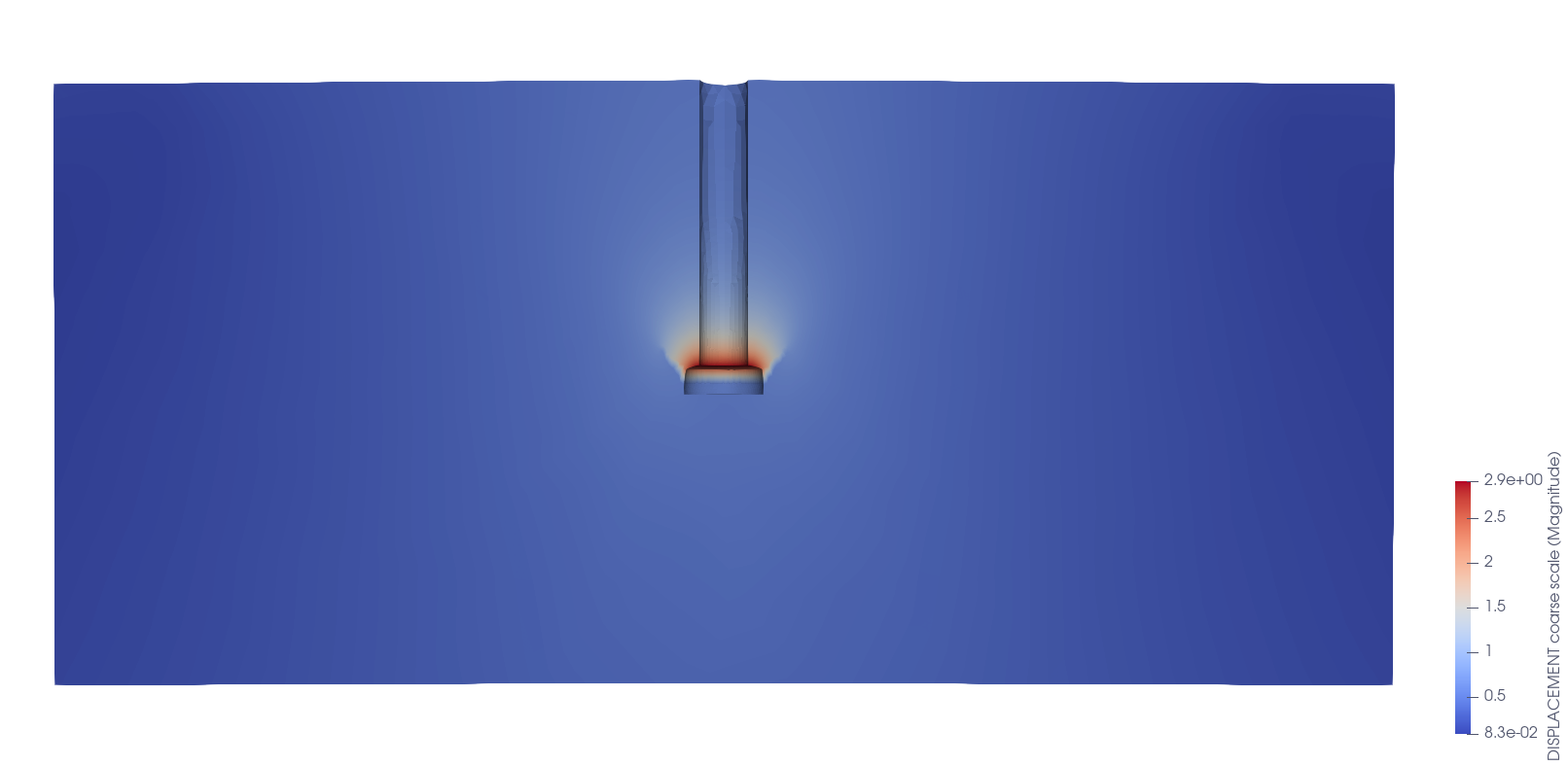}\label{PO_disp_gset_400}}
	\subfloat[Global scale,h=-300]{\includegraphics[width=55mm]{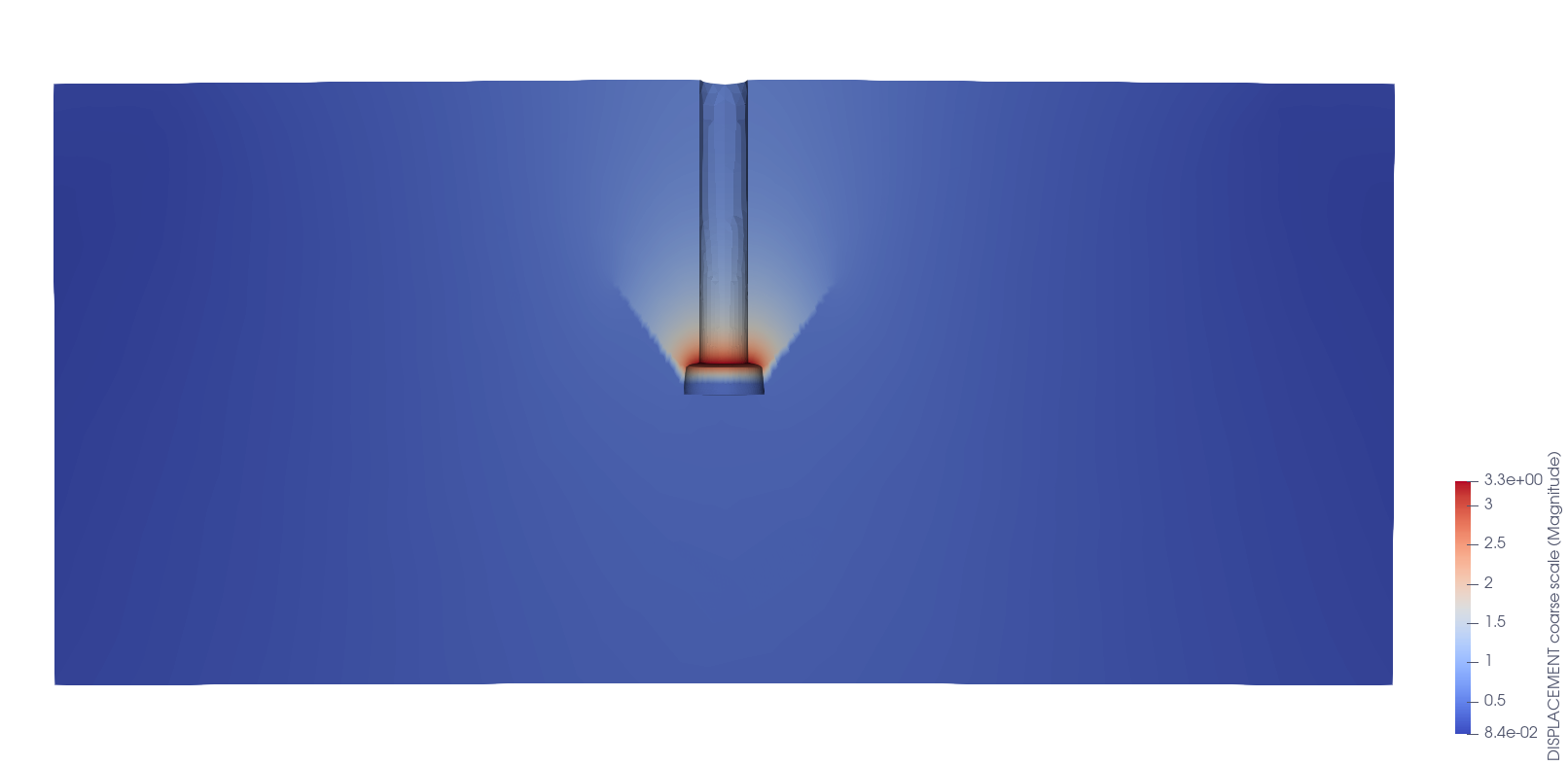}\label{PO_disp_gset_300}}
	\subfloat[Global scale,h=-100]{\includegraphics[width=55mm]{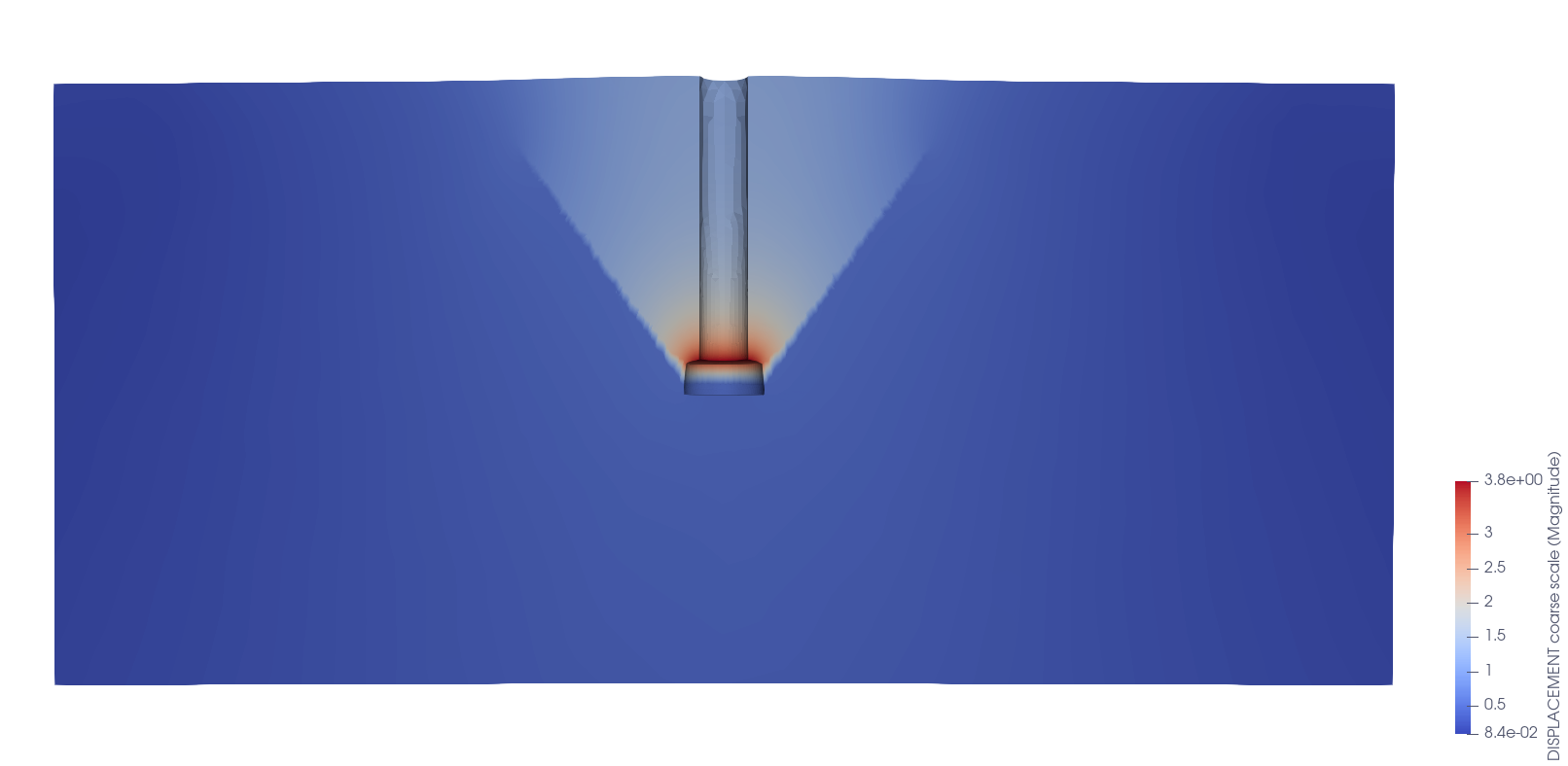}\label{PO_disp_gset_100}}\\
	\subfloat[Local scale,h=-400]{\includegraphics[width=55mm]{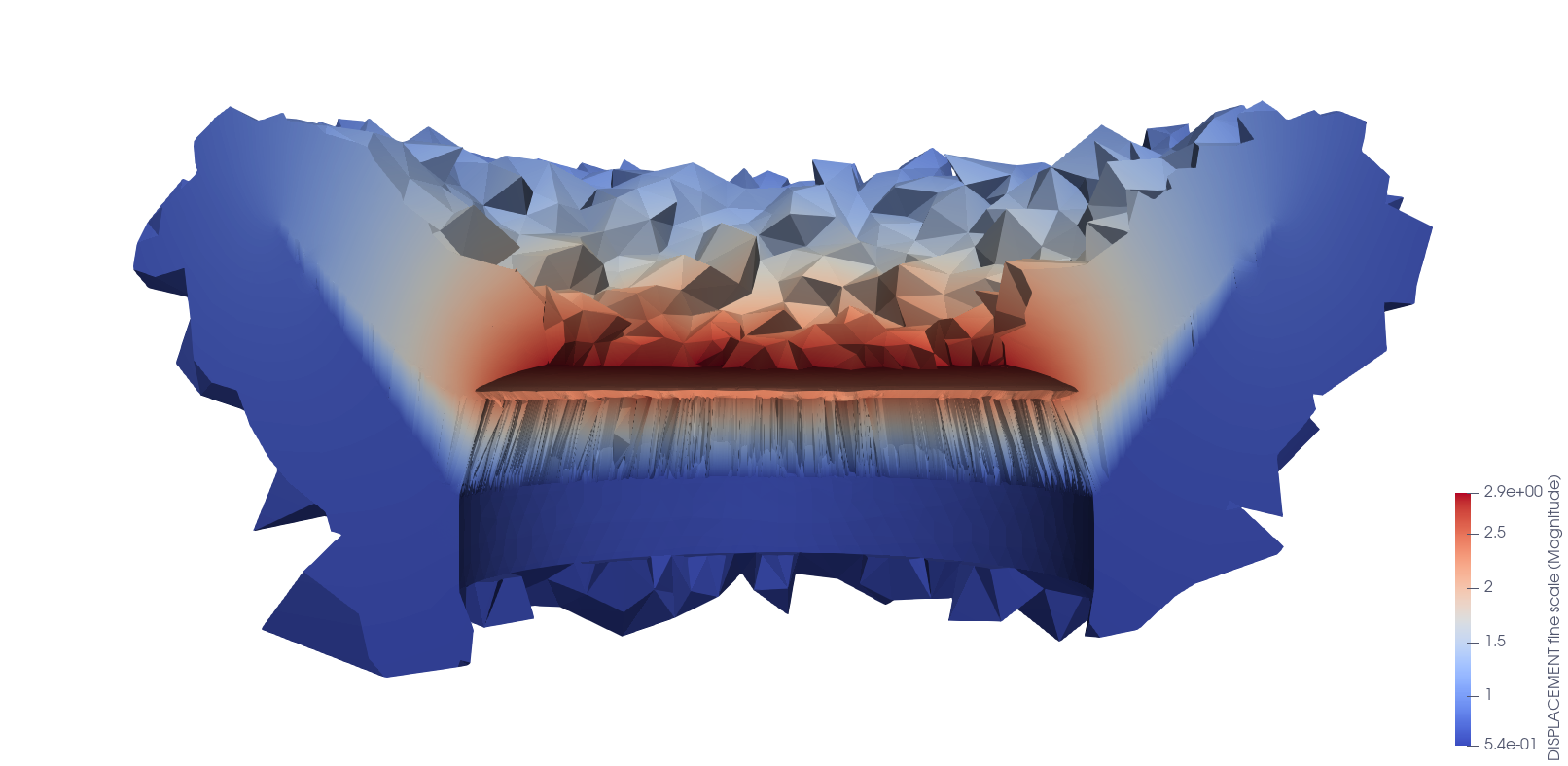}\label{PO_disp_rset_400}}
	\subfloat[Local scale,h=-300]{\includegraphics[width=55mm]{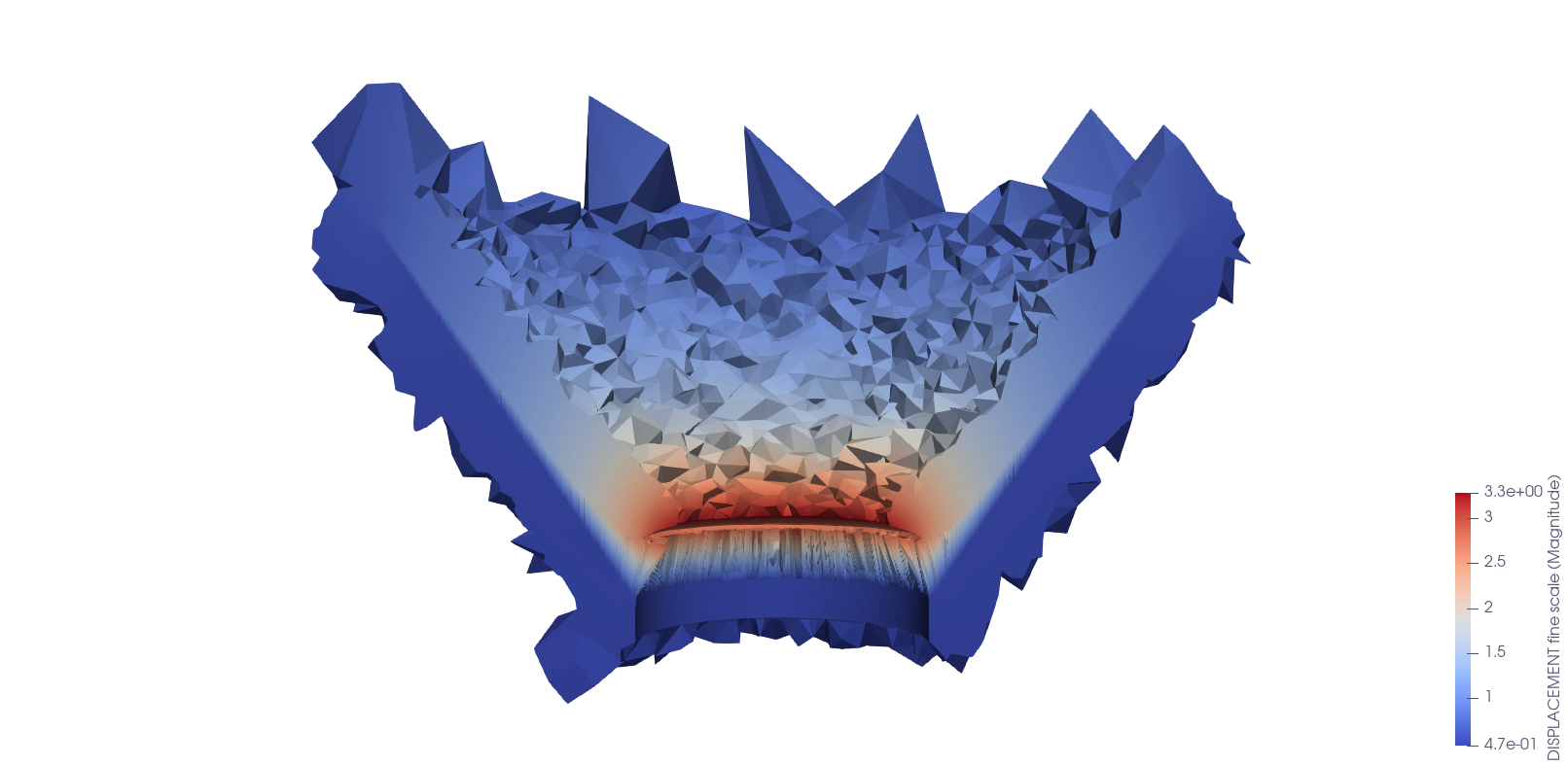}\label{PO_disp_rset_300}}
	\subfloat[Local scale,h=-100]{\includegraphics[width=55mm]{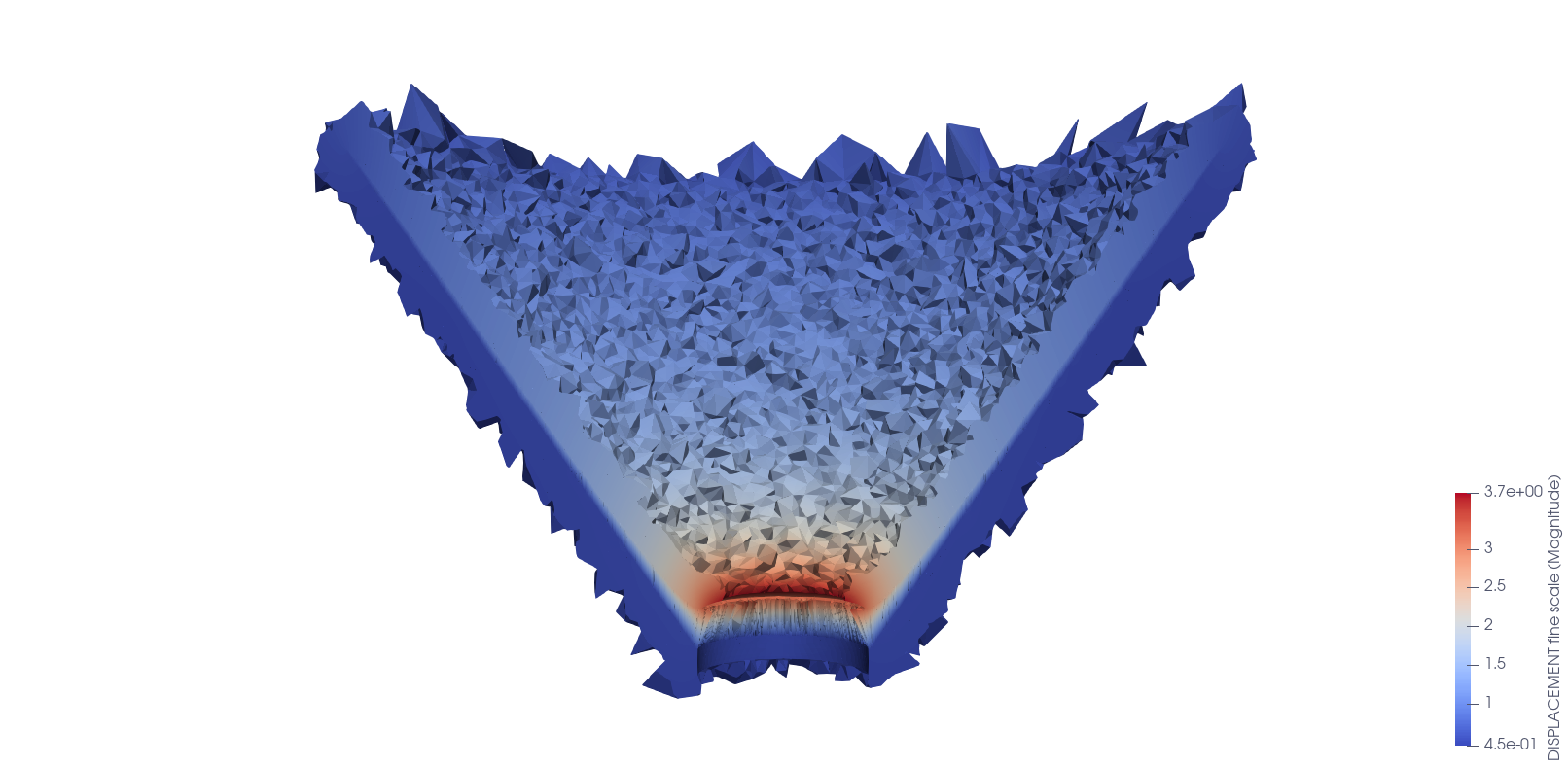}\label{PO_disp_rset_100}}
	\caption{Displacement field for h=-400, h=-300 and h=-100.  Cross-section of the SP discretization and the coarse mesh. Displacement multiplied by 10.\label{PO_disp}}	
\end{figure}

The displacement results of the \tS solver, obtained  at both scales for all stages, are shown in figure \ref{PO_disp}.
The displacements field  are presented on the coarse mesh (using only C-set dofs) and on the SP fine discretization (F-set dofs).
Figure \ref{PO_disp_rset_400}, \ref{PO_disp_rset_300} and \ref{PO_disp_rset_100} show at each adaptation the region covered by the SP.
At the crack, we can see distorted elements confirming the correct displacement jump introduced by the fully damaged material.
This jump is also evident in the clear separation of colors between the cone that is torn off and the rest of the disk that hardly moves.
This is even  clearer in figures \ref{PO_disp_gset_400}, \ref{PO_disp_gset_300} and \ref{PO_disp_gset_100} where the crack appears in context of the complete disk.
The pulled out cone becomes more and more visible as the damage grows.
These results, confirmed by simulation with other solvers, fully validate the proposed \tS implementation (in particular the treatment of the hanging nodes (section \ref{TS_reffield}) and the enrichment of mixed patch  (figure \ref{enriche_construct})).

\begin{figure}[h!]
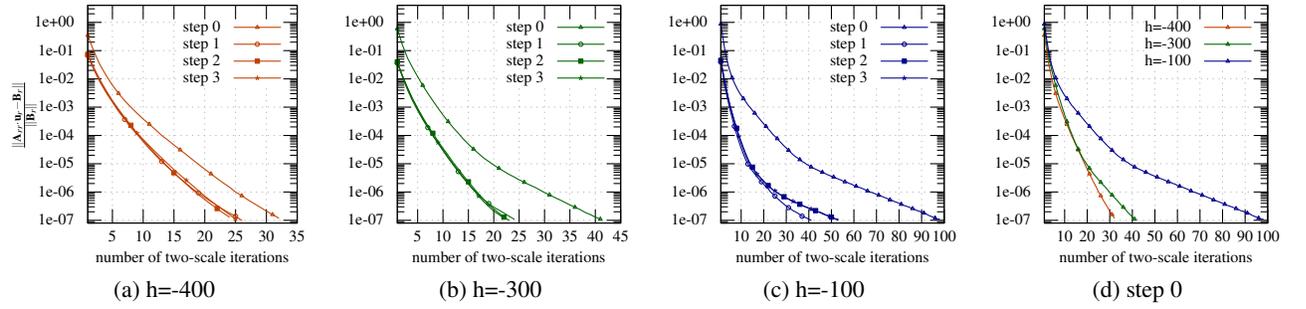

	\subfloat[h=-400]{

\label{PO_curve_resi_s0}}
	\caption{Convergence of residual errors for the 3 adaptation stages and associated steps. Lin/Log scale\label{PO_curve_resi}.}
\end{figure}

The \tS solver (the \TSI version) is again compared  to the "fr",  "blr"  and "dd" solvers.
The cumulative elapsed times of the different steps (per adaptation) give the performance curves shown in the figure \ref{PO_curve_times}.
\begin{figure}[h]
	\subfloat[h=-400]{

\label{PO_curve_times_100_l}}
	\caption{Accumulation over the 4 steps for each adaptation of the elapsed time in second. For the full computation in (\protect\subref*{PO_curve_times_400},\protect\subref*{PO_curve_times_300},\protect\subref*{PO_curve_times_100}) and only for   the \tS loop in (\protect\subref*{PO_curve_times_400_l},\protect\subref*{PO_curve_times_300_l},\protect\subref*{PO_curve_times_100_l}).  Log/log scale.\label{PO_curve_times}.}
\end{figure}
For "fr" and "blr" it is simply 4 times almost the same cost (mainly factoring and solving, the symbolic phase being done only once).
But for "dd" and "tsi", each execution, after the first one, starts from the previous one and gains iterations in their looping part.
For "dd" only the first step requires iterating in the conjugate gradient to solve the global boundary problem.
The other steps just enter the loop of the algorithm \ref{CG_algo} once and stop.
They  cost  less than the first step (almost only the condensation and  preconditionner are updated and do count, but both depends only on the domain size).
For "tsi", the restart process reduces the number of scale iterations as can be seen in figure \ref{PO_curve_resi} by starting  with a solution with smaller residual error.  
As the crack grows, the matrix conditioning decreases, and we can again observe  that the  rate of convergence of the residual error decreases (figure \ref{PO_curve_resi_s0}).

As in the other test cases, the \tS solver provides consistent performance for all numbers of cores tested.
Only the "dd" solver achieves equal or better performance when used with a high number of processes.
This decrease in scalability is due, when using many cores, on the one hand to the low number of patches per process and  the high number of distributed patches and on the other hand to the low scalability of the coarse-scale resolution (even with the use of the  \TSI version).
This can be seen in figures \ref{PO_curve_times_400_l}, \ref{PO_curve_times_300_l} and \ref{PO_curve_times_100_l} where the times spent in the \tS loop to solve the fine-scale problems ("loop fine-scale solv") and the global-scale problem ("loop global-scale solv") are shown.
We can observe the drop in scalability when reaching $nbp_{max}$ processes (see section \ref{gsolv})  with the coarse scale resolution.
To prove that this last point, if addressed, can significantly improve the scalability of the \tS solver, a new \tS version ("tsdd" hereafter) has been   implemented with  the domain decomposition solver of \ref{DDRES} used at the coarse scale level ($\dagger$ of algorithm \ref{TS_algebra}).
This change brings the scalability of the "dd" solver to the coarse scale.
And by always using the last solution of the global boundary problem (at coarse scale) as the starting point for the "dd" conjugate gradient resolution (algorithm \ref{domain_decomposition_algo}), the evolutionary aspect of the simulation is now taken into account in coarse scale resolution.
In addition, during the \tS loop, restarting from the previous solution of the global boundary problem also reduces the  number of iterations of the  conjugate gradient. 
The iteration gains appears in the figure \ref{PO_curve_itertsd}.
Note the oscillation that appears for step 0 around the 35th iteration just as the residual curve in the figure \ref{PO_curve_resi_100} stabilizes on another slope. 
This must be understood in a future work on the subject. 
This new "tsdd" version gives, as expected, better results (see figure \ref{PO_curve_times_100}) since the use of the "dd" solver at coarse scale is more efficient than a direct/iterative resolution as can be seen in \ref{PO_curve_times_100_l}.
The slope at both scales is correct, giving a nice overall slope for the "tsdd" solver.
This proves, that any solver that scales well to a large number of cores and  benefits from iterative context resolution will be a good candidate  to solve large coarse-scale problem.
It is natural to think of the \tS distributed solver for this, as it can provide such capabilities, as has been proven in this work.
This is left as a future prospect. 

Since the patches do not have the same "weight", the local problem scheduling proposed in section \ref{scheduling} can be fully validated by this test case.
In table \ref{tabsortingversion} tree variants of the algorithms   \ref{sequencing_algo:general}, \ref{sequencing_algo:compute_G}, \ref{sequencing_algo:pick_first} and \ref{sequencing_algo:pick}, called VO, V1  and V2,  are tested for h=-300 and h=-100 with 64 and 512 processes.
\begin{figure}[h!]
	\begin{minipage}{0.34\textwidth}
\begin{tikzpicture}[gnuplot]
\tikzset{every node/.append style={scale=0.70}}
\path (0.000,0.000) rectangle (5.625,4.375);
\gpcolor{color=gp lt color axes}
\gpsetlinetype{gp lt axes}
\gpsetdashtype{gp dt axes}
\gpsetlinewidth{0.50}
\draw[gp path] (0.925,0.691)--(5.237,0.691);
\gpcolor{color=gp lt color border}
\gpsetlinetype{gp lt border}
\gpsetdashtype{gp dt solid}
\gpsetlinewidth{1.00}
\draw[gp path] (0.925,0.691)--(1.105,0.691);
\draw[gp path] (5.237,0.691)--(5.057,0.691);
\node[gp node right] at (0.796,0.691) {$0$};
\gpcolor{color=gp lt color axes}
\gpsetlinetype{gp lt axes}
\gpsetdashtype{gp dt axes}
\gpsetlinewidth{0.50}
\draw[gp path] (0.925,1.186)--(5.237,1.186);
\gpcolor{color=gp lt color border}
\gpsetlinetype{gp lt border}
\gpsetdashtype{gp dt solid}
\gpsetlinewidth{1.00}
\draw[gp path] (0.925,1.186)--(1.105,1.186);
\draw[gp path] (5.237,1.186)--(5.057,1.186);
\node[gp node right] at (0.796,1.186) {$50$};
\gpcolor{color=gp lt color axes}
\gpsetlinetype{gp lt axes}
\gpsetdashtype{gp dt axes}
\gpsetlinewidth{0.50}
\draw[gp path] (0.925,1.682)--(5.237,1.682);
\gpcolor{color=gp lt color border}
\gpsetlinetype{gp lt border}
\gpsetdashtype{gp dt solid}
\gpsetlinewidth{1.00}
\draw[gp path] (0.925,1.682)--(1.105,1.682);
\draw[gp path] (5.237,1.682)--(5.057,1.682);
\node[gp node right] at (0.796,1.682) {$100$};
\gpcolor{color=gp lt color axes}
\gpsetlinetype{gp lt axes}
\gpsetdashtype{gp dt axes}
\gpsetlinewidth{0.50}
\draw[gp path] (0.925,2.177)--(5.237,2.177);
\gpcolor{color=gp lt color border}
\gpsetlinetype{gp lt border}
\gpsetdashtype{gp dt solid}
\gpsetlinewidth{1.00}
\draw[gp path] (0.925,2.177)--(1.105,2.177);
\draw[gp path] (5.237,2.177)--(5.057,2.177);
\node[gp node right] at (0.796,2.177) {$150$};
\gpcolor{color=gp lt color axes}
\gpsetlinetype{gp lt axes}
\gpsetdashtype{gp dt axes}
\gpsetlinewidth{0.50}
\draw[gp path] (0.925,2.672)--(5.237,2.672);
\gpcolor{color=gp lt color border}
\gpsetlinetype{gp lt border}
\gpsetdashtype{gp dt solid}
\gpsetlinewidth{1.00}
\draw[gp path] (0.925,2.672)--(1.105,2.672);
\draw[gp path] (5.237,2.672)--(5.057,2.672);
\node[gp node right] at (0.796,2.672) {$200$};
\gpcolor{color=gp lt color axes}
\gpsetlinetype{gp lt axes}
\gpsetdashtype{gp dt axes}
\gpsetlinewidth{0.50}
\draw[gp path] (0.925,3.167)--(3.380,3.167);
\draw[gp path] (5.108,3.167)--(5.237,3.167);
\gpcolor{color=gp lt color border}
\gpsetlinetype{gp lt border}
\gpsetdashtype{gp dt solid}
\gpsetlinewidth{1.00}
\draw[gp path] (0.925,3.167)--(1.105,3.167);
\draw[gp path] (5.237,3.167)--(5.057,3.167);
\node[gp node right] at (0.796,3.167) {$250$};
\gpcolor{color=gp lt color axes}
\gpsetlinetype{gp lt axes}
\gpsetdashtype{gp dt axes}
\gpsetlinewidth{0.50}
\draw[gp path] (0.925,3.663)--(3.380,3.663);
\draw[gp path] (5.108,3.663)--(5.237,3.663);
\gpcolor{color=gp lt color border}
\gpsetlinetype{gp lt border}
\gpsetdashtype{gp dt solid}
\gpsetlinewidth{1.00}
\draw[gp path] (0.925,3.663)--(1.105,3.663);
\draw[gp path] (5.237,3.663)--(5.057,3.663);
\node[gp node right] at (0.796,3.663) {$300$};
\gpcolor{color=gp lt color axes}
\gpsetlinetype{gp lt axes}
\gpsetdashtype{gp dt axes}
\gpsetlinewidth{0.50}
\draw[gp path] (0.925,4.158)--(5.237,4.158);
\gpcolor{color=gp lt color border}
\gpsetlinetype{gp lt border}
\gpsetdashtype{gp dt solid}
\gpsetlinewidth{1.00}
\draw[gp path] (0.925,4.158)--(1.105,4.158);
\draw[gp path] (5.237,4.158)--(5.057,4.158);
\node[gp node right] at (0.796,4.158) {$350$};
\gpcolor{color=gp lt color axes}
\gpsetlinetype{gp lt axes}
\gpsetdashtype{gp dt axes}
\gpsetlinewidth{0.50}
\draw[gp path] (0.925,0.691)--(0.925,4.158);
\gpcolor{color=gp lt color border}
\gpsetlinetype{gp lt border}
\gpsetdashtype{gp dt solid}
\gpsetlinewidth{1.00}
\draw[gp path] (0.925,0.691)--(0.925,0.871);
\draw[gp path] (0.925,4.158)--(0.925,3.978);
\node[gp node center] at (0.925,0.475) {$0$};
\gpcolor{color=gp lt color axes}
\gpsetlinetype{gp lt axes}
\gpsetdashtype{gp dt axes}
\gpsetlinewidth{0.50}
\draw[gp path] (1.356,0.691)--(1.356,4.158);
\gpcolor{color=gp lt color border}
\gpsetlinetype{gp lt border}
\gpsetdashtype{gp dt solid}
\gpsetlinewidth{1.00}
\draw[gp path] (1.356,0.691)--(1.356,0.871);
\draw[gp path] (1.356,4.158)--(1.356,3.978);
\node[gp node center] at (1.356,0.475) {$10$};
\gpcolor{color=gp lt color axes}
\gpsetlinetype{gp lt axes}
\gpsetdashtype{gp dt axes}
\gpsetlinewidth{0.50}
\draw[gp path] (1.787,0.691)--(1.787,4.158);
\gpcolor{color=gp lt color border}
\gpsetlinetype{gp lt border}
\gpsetdashtype{gp dt solid}
\gpsetlinewidth{1.00}
\draw[gp path] (1.787,0.691)--(1.787,0.871);
\draw[gp path] (1.787,4.158)--(1.787,3.978);
\node[gp node center] at (1.787,0.475) {$20$};
\gpcolor{color=gp lt color axes}
\gpsetlinetype{gp lt axes}
\gpsetdashtype{gp dt axes}
\gpsetlinewidth{0.50}
\draw[gp path] (2.219,0.691)--(2.219,4.158);
\gpcolor{color=gp lt color border}
\gpsetlinetype{gp lt border}
\gpsetdashtype{gp dt solid}
\gpsetlinewidth{1.00}
\draw[gp path] (2.219,0.691)--(2.219,0.871);
\draw[gp path] (2.219,4.158)--(2.219,3.978);
\node[gp node center] at (2.219,0.475) {$30$};
\gpcolor{color=gp lt color axes}
\gpsetlinetype{gp lt axes}
\gpsetdashtype{gp dt axes}
\gpsetlinewidth{0.50}
\draw[gp path] (2.650,0.691)--(2.650,4.158);
\gpcolor{color=gp lt color border}
\gpsetlinetype{gp lt border}
\gpsetdashtype{gp dt solid}
\gpsetlinewidth{1.00}
\draw[gp path] (2.650,0.691)--(2.650,0.871);
\draw[gp path] (2.650,4.158)--(2.650,3.978);
\node[gp node center] at (2.650,0.475) {$40$};
\gpcolor{color=gp lt color axes}
\gpsetlinetype{gp lt axes}
\gpsetdashtype{gp dt axes}
\gpsetlinewidth{0.50}
\draw[gp path] (3.081,0.691)--(3.081,4.158);
\gpcolor{color=gp lt color border}
\gpsetlinetype{gp lt border}
\gpsetdashtype{gp dt solid}
\gpsetlinewidth{1.00}
\draw[gp path] (3.081,0.691)--(3.081,0.871);
\draw[gp path] (3.081,4.158)--(3.081,3.978);
\node[gp node center] at (3.081,0.475) {$50$};
\gpcolor{color=gp lt color axes}
\gpsetlinetype{gp lt axes}
\gpsetdashtype{gp dt axes}
\gpsetlinewidth{0.50}
\draw[gp path] (3.512,0.691)--(3.512,3.078);
\draw[gp path] (3.512,3.978)--(3.512,4.158);
\gpcolor{color=gp lt color border}
\gpsetlinetype{gp lt border}
\gpsetdashtype{gp dt solid}
\gpsetlinewidth{1.00}
\draw[gp path] (3.512,0.691)--(3.512,0.871);
\draw[gp path] (3.512,4.158)--(3.512,3.978);
\node[gp node center] at (3.512,0.475) {$60$};
\gpcolor{color=gp lt color axes}
\gpsetlinetype{gp lt axes}
\gpsetdashtype{gp dt axes}
\gpsetlinewidth{0.50}
\draw[gp path] (3.943,0.691)--(3.943,3.078);
\draw[gp path] (3.943,3.978)--(3.943,4.158);
\gpcolor{color=gp lt color border}
\gpsetlinetype{gp lt border}
\gpsetdashtype{gp dt solid}
\gpsetlinewidth{1.00}
\draw[gp path] (3.943,0.691)--(3.943,0.871);
\draw[gp path] (3.943,4.158)--(3.943,3.978);
\node[gp node center] at (3.943,0.475) {$70$};
\gpcolor{color=gp lt color axes}
\gpsetlinetype{gp lt axes}
\gpsetdashtype{gp dt axes}
\gpsetlinewidth{0.50}
\draw[gp path] (4.375,0.691)--(4.375,3.078);
\draw[gp path] (4.375,3.978)--(4.375,4.158);
\gpcolor{color=gp lt color border}
\gpsetlinetype{gp lt border}
\gpsetdashtype{gp dt solid}
\gpsetlinewidth{1.00}
\draw[gp path] (4.375,0.691)--(4.375,0.871);
\draw[gp path] (4.375,4.158)--(4.375,3.978);
\node[gp node center] at (4.375,0.475) {$80$};
\gpcolor{color=gp lt color axes}
\gpsetlinetype{gp lt axes}
\gpsetdashtype{gp dt axes}
\gpsetlinewidth{0.50}
\draw[gp path] (4.806,0.691)--(4.806,3.078);
\draw[gp path] (4.806,3.978)--(4.806,4.158);
\gpcolor{color=gp lt color border}
\gpsetlinetype{gp lt border}
\gpsetdashtype{gp dt solid}
\gpsetlinewidth{1.00}
\draw[gp path] (4.806,0.691)--(4.806,0.871);
\draw[gp path] (4.806,4.158)--(4.806,3.978);
\node[gp node center] at (4.806,0.475) {$90$};
\gpcolor{color=gp lt color axes}
\gpsetlinetype{gp lt axes}
\gpsetdashtype{gp dt axes}
\gpsetlinewidth{0.50}
\draw[gp path] (5.237,0.691)--(5.237,4.158);
\gpcolor{color=gp lt color border}
\gpsetlinetype{gp lt border}
\gpsetdashtype{gp dt solid}
\gpsetlinewidth{1.00}
\draw[gp path] (5.237,0.691)--(5.237,0.871);
\draw[gp path] (5.237,4.158)--(5.237,3.978);
\node[gp node center] at (5.237,0.475) {$100$};
\draw[gp path] (0.925,4.158)--(0.925,0.691)--(5.237,0.691)--(5.237,4.158)--cycle;
\node[gp node center,rotate=-270] at (0.204,2.424) {number of conjugate gradient iterations};
\node[gp node center] at (3.081,0.151) {number of \tS iterations};
\node[gp node right] at (4.154,3.865) {step 0};
\gpcolor{rgb color={0.753,0.251,0.000}}
\draw[gp path] (4.283,3.865)--(4.979,3.865);
\draw[gp path] (0.925,3.801)--(0.968,3.940)--(1.011,3.970)--(1.054,3.613)--(1.097,3.564)%
  --(1.141,3.455)--(1.184,3.177)--(1.227,3.128)--(1.270,3.108)--(1.313,3.088)--(1.356,3.058)%
  --(1.399,3.029)--(1.442,3.009)--(1.486,2.999)--(1.529,2.989)--(1.572,2.969)--(1.615,2.950)%
  --(1.658,2.930)--(1.701,2.900)--(1.744,2.880)--(1.787,2.860)--(1.831,2.831)--(1.874,2.821)%
  --(1.917,2.801)--(1.960,2.791)--(2.003,2.771)--(2.046,2.751)--(2.089,2.741)--(2.132,2.732)%
  --(2.175,2.722)--(2.219,2.712)--(2.262,2.702)--(2.305,2.682)--(2.348,2.672)--(2.391,2.652)%
  --(2.434,2.642)--(2.477,1.652)--(2.520,2.593)--(2.564,1.592)--(2.607,2.524)--(2.650,1.533)%
  --(2.693,2.434)--(2.736,1.503)--(2.779,2.098)--(2.822,1.493)--(2.865,2.058)--(2.909,1.474)%
  --(2.952,1.493)--(2.995,2.038)--(3.038,1.444)--(3.081,1.454)--(3.124,2.028)--(3.167,1.394)%
  --(3.210,1.355)--(3.253,1.424)--(3.297,1.999)--(3.340,1.325)--(3.383,1.305)--(3.426,1.315)%
  --(3.469,1.325)--(3.512,1.969)--(3.555,1.256)--(3.598,1.246)--(3.642,1.236)--(3.685,1.236)%
  --(3.728,1.236)--(3.771,1.256)--(3.814,0.978)--(3.857,1.335)--(3.900,0.978)--(3.943,1.929)%
  --(3.987,0.949)--(4.030,1.166)--(4.073,0.929)--(4.116,1.166)--(4.159,0.909)--(4.202,1.147)%
  --(4.245,0.879)--(4.288,0.939)--(4.331,1.137)--(4.375,0.859)--(4.418,0.919)--(4.461,1.137)%
  --(4.504,0.840)--(4.547,0.869)--(4.590,0.909)--(4.633,1.117)--(4.676,0.810)--(4.720,0.840)%
  --(4.763,0.859)--(4.806,0.869)--(4.849,0.770)--(4.892,1.127)--(4.935,0.780)--(4.978,0.780)%
  --(5.021,0.790)--(5.065,0.731)--(5.108,0.840)--(5.151,0.731);
\gpsetpointsize{2.00}
\gp3point{gp mark 8}{}{(0.925,3.801)}
\gp3point{gp mark 8}{}{(1.141,3.455)}
\gp3point{gp mark 8}{}{(1.356,3.058)}
\gp3point{gp mark 8}{}{(1.572,2.969)}
\gp3point{gp mark 8}{}{(1.787,2.860)}
\gp3point{gp mark 8}{}{(2.003,2.771)}
\gp3point{gp mark 8}{}{(2.219,2.712)}
\gp3point{gp mark 8}{}{(2.434,2.642)}
\gp3point{gp mark 8}{}{(2.650,1.533)}
\gp3point{gp mark 8}{}{(2.865,2.058)}
\gp3point{gp mark 8}{}{(3.081,1.454)}
\gp3point{gp mark 8}{}{(3.297,1.999)}
\gp3point{gp mark 8}{}{(3.512,1.969)}
\gp3point{gp mark 8}{}{(3.728,1.236)}
\gp3point{gp mark 8}{}{(3.943,1.929)}
\gp3point{gp mark 8}{}{(4.159,0.909)}
\gp3point{gp mark 8}{}{(4.375,0.859)}
\gp3point{gp mark 8}{}{(4.590,0.909)}
\gp3point{gp mark 8}{}{(4.806,0.869)}
\gp3point{gp mark 8}{}{(5.021,0.790)}
\gp3point{gp mark 8}{}{(4.631,3.865)}
\gpcolor{color=gp lt color border}
\node[gp node right] at (4.154,3.640) {step 1};
\gpcolor{rgb color={0.000,0.392,0.000}}
\draw[gp path] (4.283,3.640)--(4.979,3.640);
\draw[gp path] (0.925,2.672)--(0.968,2.355)--(1.011,2.216)--(1.054,2.078)--(1.097,1.999)%
  --(1.141,1.959)--(1.184,1.691)--(1.227,1.632)--(1.270,1.592)--(1.313,1.553)--(1.356,1.523)%
  --(1.399,1.384)--(1.442,1.365)--(1.486,1.355)--(1.529,1.335)--(1.572,1.315)--(1.615,1.295)%
  --(1.658,0.988)--(1.701,1.285)--(1.744,0.958)--(1.787,1.196)--(1.831,0.919)--(1.874,1.176)%
  --(1.917,0.879)--(1.960,0.929)--(2.003,1.137)--(2.046,0.840)--(2.089,0.859)--(2.132,0.889)%
  --(2.175,1.087)--(2.219,0.800)--(2.262,0.810)--(2.305,0.840)--(2.348,0.849)--(2.391,0.760)%
  --(2.434,1.067)--(2.477,0.770)--(2.520,0.770)--(2.564,0.731)--(2.607,0.790)--(2.650,0.721)%
  --(2.693,0.879)--(2.736,1.087)--(2.779,0.790)--(2.822,0.800)--(2.865,0.849)--(2.909,0.859)%
  --(2.952,0.750)--(2.995,1.087)--(3.038,0.770)--(3.081,0.770)--(3.124,0.731)--(3.167,0.800)%
  --(3.210,0.721);
\gp3point{gp mark 6}{}{(0.925,2.672)}
\gp3point{gp mark 6}{}{(1.184,1.691)}
\gp3point{gp mark 6}{}{(1.442,1.365)}
\gp3point{gp mark 6}{}{(1.701,1.285)}
\gp3point{gp mark 6}{}{(1.960,0.929)}
\gp3point{gp mark 6}{}{(2.219,0.800)}
\gp3point{gp mark 6}{}{(2.477,0.770)}
\gp3point{gp mark 6}{}{(2.736,1.087)}
\gp3point{gp mark 6}{}{(2.995,1.087)}
\gp3point{gp mark 6}{}{(4.631,3.640)}
\gpcolor{color=gp lt color border}
\node[gp node right] at (4.154,3.415) {step 2};
\gpcolor{rgb color={0.000,0.000,0.545}}
\draw[gp path] (4.283,3.415)--(4.979,3.415);
\draw[gp path] (0.925,2.702)--(0.968,3.118)--(1.011,2.256)--(1.054,2.177)--(1.097,2.078)%
  --(1.141,2.008)--(1.184,1.969)--(1.227,1.691)--(1.270,1.662)--(1.313,1.592)--(1.356,1.563)%
  --(1.399,1.533)--(1.442,1.503)--(1.486,1.464)--(1.529,1.365)--(1.572,1.345)--(1.615,1.325)%
  --(1.658,1.315)--(1.701,1.295)--(1.744,1.275)--(1.787,1.256)--(1.831,1.226)--(1.874,1.196)%
  --(1.917,1.186)--(1.960,1.176)--(2.003,0.978)--(2.046,1.176)--(2.089,0.958)--(2.132,1.166)%
  --(2.175,0.929)--(2.219,1.157)--(2.262,0.899)--(2.305,1.147)--(2.348,0.869)--(2.391,0.939)%
  --(2.434,1.137)--(2.477,0.849)--(2.520,0.879)--(2.564,1.117)--(2.607,0.830)--(2.650,0.849)%
  --(2.693,0.869)--(2.736,1.127)--(2.779,0.800)--(2.822,0.830)--(2.865,0.840)--(2.909,0.849)%
  --(2.952,0.750)--(2.995,1.117)--(3.038,0.780)--(3.081,0.780)--(3.124,0.721)--(3.167,0.820)%
  --(3.210,0.721);
\gp3point{gp mark 5}{}{(0.925,2.702)}
\gp3point{gp mark 5}{}{(1.227,1.691)}
\gp3point{gp mark 5}{}{(1.529,1.365)}
\gp3point{gp mark 5}{}{(1.831,1.226)}
\gp3point{gp mark 5}{}{(2.132,1.166)}
\gp3point{gp mark 5}{}{(2.434,1.137)}
\gp3point{gp mark 5}{}{(2.736,1.127)}
\gp3point{gp mark 5}{}{(3.038,0.780)}
\gp3point{gp mark 5}{}{(4.631,3.415)}
\gpcolor{color=gp lt color border}
\node[gp node right] at (4.154,3.190) {step 3};
\gpcolor{rgb color={0.545,0.000,0.000}}
\draw[gp path] (4.283,3.190)--(4.979,3.190);
\draw[gp path] (0.925,2.672)--(0.968,2.484)--(1.011,2.345)--(1.054,2.246)--(1.097,2.197)%
  --(1.141,2.048)--(1.184,1.999)--(1.227,1.989)--(1.270,1.672)--(1.313,1.652)--(1.356,1.583)%
  --(1.399,1.553)--(1.442,1.533)--(1.486,1.523)--(1.529,1.493)--(1.572,1.483)--(1.615,1.365)%
  --(1.658,1.325)--(1.701,1.315)--(1.744,1.305)--(1.787,1.295)--(1.831,1.275)--(1.874,1.256)%
  --(1.917,0.978)--(1.960,1.315)--(2.003,0.968)--(2.046,1.315)--(2.089,0.949)--(2.132,1.186)%
  --(2.175,0.929)--(2.219,1.166)--(2.262,0.899)--(2.305,1.166)--(2.348,0.869)--(2.391,0.929)%
  --(2.434,1.157)--(2.477,0.859)--(2.520,0.879)--(2.564,1.147)--(2.607,0.840)--(2.650,0.849);
\gp3point{gp mark 3}{}{(0.925,2.672)}
\gp3point{gp mark 3}{}{(1.270,1.672)}
\gp3point{gp mark 3}{}{(1.615,1.365)}
\gp3point{gp mark 3}{}{(1.960,1.315)}
\gp3point{gp mark 3}{}{(2.305,1.166)}
\gp3point{gp mark 3}{}{(2.650,0.849)}
\gp3point{gp mark 3}{}{(4.631,3.190)}
\gpcolor{color=gp lt color border}
\draw[gp path] (0.925,4.158)--(0.925,0.691)--(5.237,0.691)--(5.237,4.158)--cycle;
\gpdefrectangularnode{gp plot 1}{\pgfpoint{0.925cm}{0.691cm}}{\pgfpoint{5.237cm}{4.158cm}}
\end{tikzpicture}
	\caption{Number of iterations for solving the "dd" conjugate gradient  at the coarse scale ("tsdd" version) for h=-100 with 1024 cores.\label{PO_curve_itertsd}.}
	\end{minipage}
	\begin{minipage}{0.64\textwidth}
	\footnotesize
	\begin{tabular}{c|l|l|c|c|c|}
		\cline{4-6}
		\multicolumn{3}{l|}{} &           \rotatebox{90}{ h=-300 }            &      \rotatebox{90}{ h=-100 }        &       \rotatebox{90}{ h=-100 }        \\ \cline{2-6}
		\multicolumn{1}{l|}{}                   & \multicolumn{2}{l|}{Number of processes} & 64 & 64 & 512 \\ \cline{2-6} 
		\multicolumn{1}{l|}{\multirow{4}{*}{}}  & \multicolumn{2}{l|}{Total number of patches} &   26609      &       94388       &      94388         \\ \cline{2-6} 
		\multicolumn{1}{l|}{}                   & \multicolumn{2}{l|}{Min number of patches per process}      &     456         &       1457          &     225        \\ \cline{2-6} 
		\multicolumn{1}{l|}{}                   & \multicolumn{2}{l|}{Max number of patches per process}      &        610      &   1802   &   316   \\ \cline{2-6} 
		\multicolumn{1}{l|}{}                   & \multicolumn{2}{l|}{Average number of patches per process}   &    535.5       &   1702.1         &  272.6          \\ \cline{2-6} 
		\multicolumn{1}{l|}{}                   & \multicolumn{2}{l|}{Average number of distributed patches per process}   &    226.6       &   442.1          &  160.4       \\ \hline
		\multicolumn{1}{|c|}{\multirow{5}{*}{V0}} & \multicolumn{2}{l|}{Number of sequence}  &        303.3      &       586.3           &     221.3     \\ \cline{2-6} 
		\multicolumn{1}{|c|}{}                  & \multirow{2}{*}{Average observed cost per process} &elapsed time (s)&   158.4        &    908.7     &    184.0          \\ \cline{3-6} 
		\multicolumn{1}{|c|}{}                  & &vs V2 (\%)     &   8.1     &    3.7 &    7.3     \\ \cline{2-6} 
		\multicolumn{1}{|c|}{}                  & \multirow{2}{*}{Observed cost standard deviation}    & in s &    12.6          &   57.0  &       9.0          \\ \cline{3-6} 
		\multicolumn{1}{|c|}{}                  &  & vs V2 (\%)     &   24.9     &         38.3      &    30.5     \\ \hline
		\multicolumn{1}{|c|}{\multirow{5}{*}{V1}} & \multicolumn{2}{l|}{Number of sequence  }      &  301.3      &       591.7           &    216.7      \\ \cline{2-6} 
		\multicolumn{1}{|c|}{}                  & \multirow{2}{*}{Average observed cost per process} &elapsed time (s)&   167.7     &     986.1 & 186.0            \\ \cline{3-6} 
		\multicolumn{1}{|c|}{}                  &   & vs V2 (\%)      &    14.4    &   12.5            &      8.4   \\ \cline{2-6} 
		\multicolumn{1}{|c|}{}                  & \multirow{2}{*}{Observed cost standard deviation}    & in s &    17.0          &         81.2       &      11.0           \\ \cline{3-6} 
		\multicolumn{1}{|c|}{}                  &   & vs V2 (\%)      &    68.4    &          97.0     &    59.3     \\ \hline
		\multicolumn{1}{|c|}{\multirow{3}{*}{V2}} & \multicolumn{2}{l|}{Number of sequence}        & 303.3      &       580.0           &     213.7     \\ \cline{2-6} 
		\multicolumn{1}{|c|}{}                  & Average observed cost per process & elapsed time (s) &      146.5                 &        876.4  &   171.6  \\ \cline{2-6} 
		\multicolumn{1}{|c|}{}                  & Observed cost standard deviation    & in s &      10.1        &               41.2 &    6.9   \\  \hline
	\end{tabular}
	\captionof{table}{Comparison of distributed patches sequencing algorithms. V2 corresponds to the algorithm \ref{sequencing_algo:general}, \ref{sequencing_algo:compute_G}, \ref{sequencing_algo:pick_first} and \ref{sequencing_algo:pick}. V0 corresponds to the algorithm \ref{sequencing_algo:general}, \ref{sequencing_algo:pick_first} and \ref{sequencing_algo:pick} without taking into account  the weight (The $\mathcal{D}$ and $\mathcal{L}$ are not sorted, the $mxwg$ and $mxwl$ are not computed nor used, the first available patch in these groups is selected). V1 corresponds to the algorithm \ref{sequencing_algo:general}, \ref{sequencing_algo:pick_first} and \ref{sequencing_algo:pick} with a patch selection, neither using nor computing, $mxwg$ and $mxwl$ (the first available patch in the group is selected). The results given  are the average values of 3 runs with the same parameters. The observed cost corresponds to the creation, factoring and  resolution of all patches for all \tS iterations and steps.}
	\label{tabsortingversion}
	\end{minipage}
\end{figure} 

The V0 version corresponds to the algorithms \ref{sequencing_algo:general}, \ref{sequencing_algo:pick_first} and \ref{sequencing_algo:pick} where the weight is not used: $\mathcal{D}$ and $\mathcal{L}$ are not sorted, $mxwg$ and $mxwl$ are not calculated or used, the first available patch in the group is selected (i.e. the one with a "random" weight that depends only on how the patches were  constructed).
The V1 version uses the descending weight ordered groups $\mathcal{D}$ and $\mathcal{L}$ but does not use or compute  $mxwg$ and $mxwl$.
The first available patch in the group is selected (i.e. the one with the highest weight).
The V2 version corresponds to the proposed algorithms \ref{sequencing_algo:general}, \ref{sequencing_algo:pick_first} and \ref{sequencing_algo:pick}.
To compare these versions, a metric called "observed cost per process" counts the elapsed time (averaged over 3 runs) per process related to all tasks connected with distributed patch sequencing (i.e. all creations, factorizations, and  resolutions of all patches for all \tS iterations and steps).
The  mean and standard deviation of this metric, for all processes,  are given in absolute value and relative to the V2 performance.
The V0 and V1 versions are   3.7\% to 8.1\%  and  8.4\% to 14.4\% slower than the V2 version, respectively, in terms of average performance.
Surprisingly, the "random" V0 version  is not the worst.
Certainly because during construction, the groups were built in a favorable order by chance. 
The standard deviation of this measure reflects the load imbalance introduced by the sequencing algorithm and the coarse mesh distribution.
Between V0,V1 and V2, only the sequencing varies and the coarse mesh distribution remain the same.
Thus, the variation in standard deviation from V2 somewhat reflects the load imbalance added by V0 and V1, which are  24.9\% to 38.3\%  and  59.3\% to 97.0\% more dispersed than the V2 version, respectively.
In any case, the V2 version is the one that also provides the smallest (or equal) number of sequences compared to the V0 and V1 versions.
These numbers of sequences are as expected higher but close to the average number of distributed patches per process.
All these observations validate the choice of V2 as the best version to use in this paper among V0,V1 and V2.

Also note  that this test case validates the enrichment function and the processing of hanging nodes in patches with a mixture of refined and unrefined elements.
Some test not presented here confirm that these choices give the best \tS convergence compared to some truncated patch solutions.

\section{Conclusions}\label{conclusions}
This work introduces some novelties in the original \GL method:
\begin{itemize}
\item a new criterion to control the iterations of the \tS resolution  loop;
\item a efficient  algebraic computation of the matrices of linear systems at both scales;
\item the use for the enrichment of a linearly interpolated function (with a shift);
\item the treatement of blending element problem with the use of additional enriched nodes corresponding to mixed discretization patches;
\item a distributed implementation of the \tS method;
\item a new algorithm for scheduling the computation of fine-scale problems; 
\item a proof given in sequential that there is an optimal scale jump (confirmed in parallel);
\item the use of specific resolution (\TSDn, \TSIn, iterative domain decomposition solver) at global level to improve scalability of this scale.
\end{itemize}

These new features improve the capability of the \tS method in several ways.
The new stopping criterion is comparable to the error (in energy norm) of the boundary condition imposed on the fine-scale problem.
It thus allows to obtain results with a given precision.
And when  used by the \tS solver in evolutionary phenomena (where the resolution starts from a close previous solution) the number of \tS iterations can be considerably reduced.

With respect to  performance, the proposed distributed implementation is promising and performs efficiently compared to some other parallel solvers for a wide range of used cores.
In particular, it allows, with few cores, to handle larger problems  than with the other solvers tested. 
Parallelism at both scales is clearly to be maintained,  but some elements  still need to be improved (at global scale in particular). 
The scalability of local scale is correct thanks to, among other things, the new scheduling of the local scale problems.

Now, in terms of future works, parallel performance and the use of the new criterion are two topics that can be explored further. The proposed new residual error criterion has already given valuable insights into the convergence of the \tS solver.
In particular, this allowed to highlight the influence of the matrix conditioning of the system at the global scale  on the convergence.
A full analysis of this observation should be done to get a theoretical confirmation.
Regarding enrichment,  a closer comparison with SGFEM and other enrichment techniques could also be conducted in a future work in a quantitative way, thanks to this new convergence criterion. 

Otherwise, for parallel performance, when the global problem becomes large, both   solutions, \TSI and iterative domain decomposition solver, exploited the fact that between \tS iterations, at the global scale, the linear system  evolves relatively slowly.
Therefore, the iterative resolution of this type of problem can gain in performance by using the previous computed solution.
This lead to the conclusion that any solver that scales well over a wide range of cores numbers and that takes advantage of an iterative context in its computation will be a good candidate for solving large coarse scale problem.
And of course, the \tS solver itself can be used for this.
But many questions need to be answered before embarking on such a process.
Would we act independently, with one \tS solver solving the large problem at the coarse scale  of another \tS solver at each of its iterations?
Or should we just do some "fine to coarse" transfer, solve a smaller coarse problem, and do some "coarse to fine" transfers?
And in this last case, what "some" should be: 2 more?
This last idea, which is somehow related to the multigrid methodology ( in \cite{Saad03} grid transfer according  to V or W or other scheme), is an interesting future work.

Still on the subject of parallel performance, one can also consider a hybrid multithreading/MPI implementation of the proposed \tS distributed  implementation to somehow select the right range of cores for the right scale. 
Specifically, in the algorithm \ref{TS_algebra_patch}, as mentioned in section \ref{scheduling}, when all distributed patches are calculated, the remaining local patches are computed independently (processes are not longer synchronized).
For this, in a process, a  sequential computation was naturally  used but as in \cite{Kim2011} a multithreaded computation can also be a solution.
The main impact would be to provide a parameter to adjust the number of distributed patches and the memory consumption per process.
In a multithreaded context, the constraint added by the dynamic scheduling proposed in \cite{Kim2011} (the largest local patch limits the number of threads) will certainly play a role in  choosing the MPI/thread partition.
But this flexibility will certainly also allow to choose  the  right number of processes at the global scale  depending on the solver used at this scale.
As a side effect, many other tasks can also benefit from this hybrid implementation: many loops on the \SP macro elements can be nicely paralyzed across many threads in the algorithm \ref{TS_algebra_init},  \ref{TS_algebra_micro_update}, \ref{TS_algebra_macro_update}, \ref{TS_algebra_residual} and  \ref{TS_algebra_B_norm} due to the per-macro-element  storage and computation of data.
But for large global-scale problem, the question of which solver to choose for that scale will certainly remain open.

Finally, a last perspective on local scale  performance concerns the use of a direct solver to solve patches (distributed or not).
It is possible to consider  that the enrichment functions  vary slowly during the \tS iteration.
Thus, if conditioning is not an issue, using an iterative solver to solve the patches may be more efficient because using the previous solutions at almost all \tS iteration would certainly reduce the number of solver iterations.
And in this context, storing data by macro element could relatively naturally result in a versatile block Jacobi preconditioner (each block can be used for many patches).
This  can also have an impact on the memory footprint (it is no longer necessary to store the factorization of each patch). 
 
\section*{acknowledgements}
The authors would like to thank Nicolas Chevaugeon  and Gregory Legrain with whom the conversations about parallelism and the \GL method were fruitful.

\renewcommand{\thealgorithm}{\Alph{section}.\arabic{algorithm}} 
\setcounter{table}{0}
\setcounter{figure}{0}

\bibliographystyle{abbrv}
\bibliography{bib}
\appendix
\section{notation convention}\label{convention_anexe} 
The following conventions are used in the algorithms in this document:
\begin{itemize}
	\item $\Leftarrow$ implies assembling a vector or matrix of a space(s) or set(s) embedded in the receiving space(s) or set(s)
    \item $\left( A,B\right) \gets\left( C,D\right) $ is equivalent to $A \gets C$ and $B\gets D$ 
    \item $\left( A,B\right) \Leftarrow \left( C,D\right) $ is equivalent to $A \Leftarrow C$ and $B\Leftarrow D$ 
     \item $A \Leftarrow\left( C,D\right) $ is equivalent to $A \Leftarrow C$ and $A\Leftarrow D$ 
    \item $A\gets \left\lbrace B,C\right\rbrace $ indicate that the matrices $B$ and $C$ are merged by column to form a single matrix stored in $A$
    \item $\left( A,B,C\right) \gets$ PROCEDURE... indicates that the called PROCEDURE returns data which are stored in $A$,$B$ and $C$.
       See return declaration in PROCEDURE.
	\item $\triangleright$ represents a task involving communication (point to point or collective)
	\item If $\mathcal{C}$ is a set of entities: 
	\begin{itemize}
		\item $\mathcal{C} \setminus e$ is  removing $e$ from $\mathcal{C}$
		\item $\mathcal{C} \cup e$ is  adding $e$ to the end of $\mathcal{C}$
	\end{itemize}
    \item $weitgh_k$ is the weight of a patch identified by the letter $k$: number of micro-scale elements embedded  in the patch $k$
    \item To simplify the algorithms,  $A^{-1}$ expresses the inverse of the matrix $A$, but it is in many cases the factorization that is actually obtained and used.
       Thus $X=A^{-1}\cdot B$ is calculate as the solution of $L\cdot D\cdot L^{t}\cdot X=B$ where $L$ and $D$ are the lower triangular and diagonal matrices obtained by factoring $A$.
       \item The subscripted letters, corresponding to  the sets of values defined in the table  \ref{tab_dof_set},  give the  dofs  on which the matrix or the vector is defined. Thus,  $\vm{A}_{gg}$ is a matrix of dimension $card(g)\times card(g)$ and  $\vm{U}_{g}$ is a vector of size $card(g)$.
       \item The superscript  $patches$ indicates the set of matrices or vectors of all fine-scale problems:  $X^{patches}=\left\{X^p : p\in patches\right\}$ with $patches=\left\{1,2,...,card(I_e^g)\right\}$ (see section \ref{TSMethodOrig} for definition of $I_e^g$)
\end{itemize}
\section{values set}\label{anex_sets}
\begin{table}[h]
	\small
	\begin{tabular}{|c|c|c|c|}
		\hline
		Set&Values type&Relation&Scale\\
		\hline
		\hline
		$G$&global-scale values&$G=D\cup g,D\cap g=\oslash$&\multirow{9}{*}{\rotatebox{0}{coarse-level}}\\
		\cline{1-3}
		$D$&fixed values (Dirichlet) &$D\subset G$&\\
		\cline{1-3}
		$C$&classical values&$C=D\cup c$&\\
		\cline{1-3}
		$g$&free values &$g=c \cup e$&\\
		\cline{1-3}	
		$c$&classical free values&$c=C\setminus D$, $c=h\cup k$&\\
		\cline{1-3}
		$e$&enriched free values&$e=G\setminus C$&\\
		\cline{1-3}
		$m$&free values  related to SP elements &$m=k\cup e$&\\\cline{1-3}
		$k$&classical free values  related to SP elements &$k\in c,k\cap h=\oslash$&\\
		\cline{1-3}
		$h$&classical free values not related to SP elements &$h=c\setminus k$&\\
		\hline
		\hline
		$Q$&fine-scale values&$Q=L\cup DG\cup d\cup q$&\multirow{5}{*}{ patch-level ($p^{th}$)}\\
		\cline{1-3}
		$L$&fixed values (by linear relation) &$L\subset Q,L\cap DP=\oslash,L\cap d=\oslash$&\\
		\cline{1-3}
		$DP$&fixed values ($D$ restricted to patch $p$)&$DP=D\cap  Q$&\\
		\cline{1-3}
		$d$&fixed values (Dirichlet from fine solution $f$)&$d= f\cap Q$&\\
		\cline{1-3}
		$q$&free values&$q=Q\setminus(d\cup DG \cup L)$&\\
		\hline
		\hline
		$R$&reference values &$R= M\cup L,M\cap L=\oslash$&\multirow{9}{*}{mix fine/coarse level }\\
		\cline{1-3}
		$L$&fixed values (by linear relation) &$L\subset R,L\cap DR=\oslash$&\\
		\cline{1-3}
		$M$&reference values without L &$M=r\cup DR=F\cup H,r\cap DR=\oslash$&\\
		\cline{1-3}
		$DR$&fixed values (inherited from D) &$DR\subset M, D\subseteq DR$&\\
		\cline{1-3}
		$F$&reference values related to SP elements&$F=f\cup (DR\setminus (DR\cap H))$&\\
		\cline{1-3}
		$H$&reference values not related to SP elements&$H=h\cup (DR\setminus (DR\cap F))$&\\
		\cline{1-3}
		$r$&reference free values &$r=h\cup f$&\\
		\cline{1-3}
		$f$&free values related to SP elements&$f\subset r$&\\
		\cline{1-3}
		$h$&free values not related to SP elements &$h=r\setminus f,h\text{ same as }h\text{ at coarse-level}$&\\
		\hline
		\hline
	\end{tabular}
	\caption{Value sets description ( lowercase letter for free values, upper case letter for free and	 fixed values)}
	\label{tab_dof_set}
\end{table}

A large number of sets are put in place to clarify the presentation of the different discretizations.
The value of the discretized field at a node is called the value.
This value can be fixed or free.
When it is free, it can be classical or enriched.
There are three reasons why a value may be fixed: on the physical Dirichlet boundary conditions, on the patch boundary  when solving  fine scale problems or when setting the hanging node value of a micro node.
Regarding the set notations, lower case will be used for free values while upper case will be  used for fixed values or for a set of values gathering both fixed and free values.
All value types are described in table \ref{tab_dof_set}.
\section{Mesh adaptation}\label{rsplit_anexe}
\begin{figure}[h]
	\subfloat[][Starting mesh]{
		\includegraphics[width=0.2\textwidth]{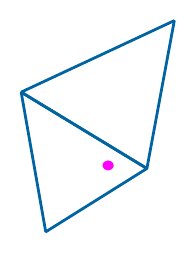}
		\label{split_0}
	}
	\subfloat[][L1]{
		\includegraphics[width=0.2\textwidth]{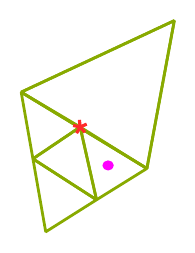}
		\label{split_1}
	}
	\subfloat[][L2: wrong  ]{
		\includegraphics[width=0.2\textwidth]{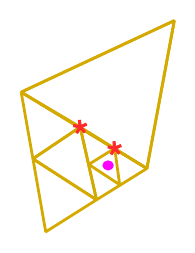}
		\label{split_2w}
	}
	\subfloat[][L2: corect]{
		\includegraphics[width=0.2\textwidth]{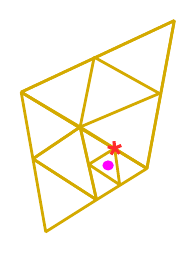}
		\label{split_2}
	}
    \subfloat[][Adapted mesh]{
    	\includegraphics[width=0.2\textwidth]{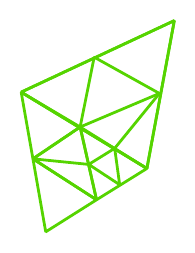}
    	\label{split_2f}
    }
	\caption{mesh refinement principles (presented in 2D but the same  principle applies in 3D): the dot represents the area of interest and the stars the hanging nodes. Only one hanging node is accepted per edge}
	\label{split}
\end{figure}

In this work, the mesh adaptation uses a simple splitting strategy.
For each level of refinement,  we take the impacted elements (i.e. those that meet a certain criterion provided by a user-defined function) from the previous level and divide each of them in an encapsulated manner.
New nodes are added to the middle of the edges of impacted elements.
The affected elements are then replaced by a tetrahedronization based on the old and new nodes, so that the new elements are encapsulated in the destroyed old elements.
The figure \ref{split} illustrates this process, in 2D for simplicity.
In the figure \ref{split_0} 2 triangles represent the starting level of the adaptation with a dot representing the area of interest.
The first level of refinement, L1 (figure \ref{split_1}),  divides only the triangle covering the area of interest.
A hanging node (start in figure \ref{split_1}) appears in the middle of the edge common to the triangles in the  figure \ref{split_0}.
Then, when moving to the next level of refinement, L2, again only the  element covering the area of interest is split.
If we simply split this element (figure \ref{split_2w}) 2 hanging nodes appear on the same edge.
This  is considered too abrupt in terms of mesh transition and is therefore avoided by forcing the splitting of the  second  initial unmodified triangle (figure \ref{split_2}).
This "one hanging node per edge" rule imposes a smoother mesh transition.
When the adaptation is complete (i.e.  for example when the target size of the element  is reached in the area of interest), the hanging nodes are removed by splitting the  element with the hanging edge/face using only that node  and the existing vertices (figure \ref{split_2f}).
Note that in 3D,  an additional node can be added to the center of gravity of the tetrahedron which is modified by removing the hanging node.

The algorithm associated with this process is implemented in parallel.
In short, for one level of refinement,  it reviews all impacted elements.
Next, it examines whether splitting these elements would violate the "one hanging node per edge" rule.
If this is the case, the connected elements (i.e. those with a hanging edge/face) are added to the list of elements to divide.
This verification operation is repeated until  no more elements are added.
There may be some communication  during this task.
Then all selected element are split.
The next level of refinement  may occur or, if it is complete  (i.e. there are no more elements selected for splitting), the final deletion of hanging nodes is performed.
More details can be found in \cite{Salzman2019}. 

Note that a constrained Delaunay tetrahedralization can be an interesting alternative for this work as long as the coarse mesh remains a skeleton of the new discretization:  all faces/edges of the tetrahedra of the coarse mesh  remain defined (possibly divided  several times) in the new mesh.

\section{$\vm{T}_{MG}$ construction}\label{Trg_anexe}

The $\vm{T}_{MG}$ matrix is a linear interpolation operator with the following structure, considering that the M-set is ordered with the H-set first and the F-set second, and that the G-set is ordered with the D-set first, the c-set second and the e-set third with h-set first in the c-set:
\begin{equation}
	\vm{T}_{MG}=\left(  \begin{array}{cccc}
		\vm{T}_{HD}&\vm{T}_{Hh}&\vm{0}&\vm{0}\\
		\vm{T}_{FD}&\vm{0}&\vm{T}_{Fk}&\vm{T}_{Fe}
	\end{array}\right)
	\label{oper_rg}
\end{equation} 
A term of $\vm{T}_{MG}$ at index $(a,b)$  is, according to the equation  \eqref{relation}, either  $N^i(\vm{x}^j)$ or $N^i(\vm{x}^j)\vm{F}^p(\vm{x}^j)\cdot\vm{e_k}$ where $\vm{x}^j$ is the node associated to the dof $a$, $\vm{e_k}$ is the unite-vector associated to the component that dof $a$ represent on the node $\vm{x}^j$,  and  $\vm{x}^i$ is the node associated to dof $b$.
When  $a\in H$ and $b\in m$, as the node $\vm{x}^j$ is not in the support of the node $\vm{x}^i$ all $N^i(\vm{x}^j)$ terms are zero and the H-set$\times$m-set block is null.
And it is the same when $a\in F$ and $b\in h$, the block F-set$\times$h-set is null.
When $b\in D$ only the nodes $\vm{x}^j$ on the derived Dirichlet boundary conditions give nonzero terms, i.e. when  $a\in DR$.
All other terms are zero in $\vm{T}_{HD}$ and $\vm{T}_{FD}$.
For $a\in H$ and $b\in h$ the macro and micro meshes are the same and the only non-zero terms are when $a=b$.
The h-set$\times$h-set block is then an identity matrix and all the other terms of $\vm{T}_{Hh}$ are null.
The remaining blocks $\vm{T}_{Fk}$ and $\vm{T}_{Fe}$ correspond respectively to $N^i(\vm{x}^j)$ for a classical dof ($b\in k$) and to $N^i(\vm{x}^j)\vm{F}^p(\vm{x}^j)\cdot\vm{e_k}$ for an enriched dof ($b\in e$).

\section{Global scale problem construction}\label{Agg_anexe}
The equality \eqref{energy_g} gives the following system:
\begin{equation}
	 \vm{T}_{MG}^t\cdot  \vm{A}_{MM} \cdot \vm{T}_{MG}\cdot \vm{U}_g=\vm{T}_{MG}^t\cdot \vm{B}_{M}
	\label{Agg_begin}
\end{equation} 
It can be  rewritten as follows using \eqref{oper_rg} and the divided matrices on the sets H and F :
\begin{equation}
	\left(  \begin{array}{cc}
		\vm{T}_{HD}^t&\vm{T}_{FD}^t\\
		\vm{T}_{Hh}^t&\vm{0}\\
		\vm{0}&\vm{T}_{Fk}^t\\
		\vm{0}&\vm{T}_{Fe}^t
	\end{array}\right) \cdot\left(  \begin{array}{cc}
		\vm{A}_{HH}&\vm{A}_{HF}\\
		\vm{A}_{FH}&\vm{A}_{FF}
	\end{array}\right) \cdot \left(  \begin{array}{cccc}
		\vm{T}_{HD}&\vm{T}_{Hh}&\vm{0}&\vm{0}\\
		\vm{T}_{FD}&\vm{0}&\vm{T}_{Fk}&\vm{T}_{Fe}
	\end{array}\right)  \cdot\left( \begin{array}{c}
		\vm{u}_D\\
		\vm{u}_h\\
		\vm{u}_k\\
		\vm{u}_e
	\end{array}\right)=\left(  \begin{array}{cc}
		\vm{T}_{HD}^t&\vm{T}_{FD}^t\\
         \vm{T}_{Hh}^t&\vm{0}\\
		\vm{0}&\vm{T}_{Fk}^t\\
		\vm{0}&\vm{T}_{Fe}^t
	\end{array}\right) \cdot \left( \begin{array}{c}
		\vm{B}_H\\
		\vm{B}_F
	\end{array}\right) 
	\label{g_sys_detail}
\end{equation} 

This gives the system:
\begin{equation}
	\left(  \begin{array}{cc}
		\vm{A}_{DD}&\vm{A}_{Dg}\\
		\vm{A}_{gD}&\vm{A}_{gg}
	\end{array}\right)\cdot \left( \begin{array}{c}
	\vm{U}_D\\
	\vm{U}_g
    \end{array}\right)=\left( \begin{array}{c}
    \vm{B}_D\\
    \vm{BN}_g
    \end{array}\right)
	\label{G_sys}
\end{equation}

where, considering here a symmetric system (but there is no restriction for a non-symmetric system) :
\begin{equation}
	\vm{A}_{gg}=\left(  \begin{array}{ccc}
		\vm{T}_{Hh}^t\cdot \vm{A}_{HH}  \cdot \vm{T}_{Hh}&\vm{T}_{Hh}^t\cdot \vm{A}_{HF}\cdot \vm{T}_{Fk}&\vm{T}_{Hh}^t\cdot \vm{A}_{HF}\cdot \vm{T}_{Fe}\\
		\vm{T}_{Fk}^t \cdot \vm{A}_{HF}^t\cdot \vm{T}_{Hh}&\vm{T}_{Fk}^t \cdot \vm{P}_{Fk}&\vm{P}_{Fk}^t \cdot \vm{T}_{Fe} \\
		\vm{T}_{Fe}^t \cdot \vm{A}_{HF}^t\cdot \vm{T}_{Hh}&\vm{T}_{Fe}^t \cdot \vm{P}_{Fk}& \vm{T}_{Fe}^t\cdot \vm{A}_{FF}  \cdot \vm{T}_{Fe}
	\end{array}\right)~\text{with}~\vm{P}_{Fk}=\vm{A}_{FF}\cdot \vm{T}_{Fk}
	\label{Agg_1}
\end{equation} 
\begin{equation}
	\vm{BN}_{g}=\left(  \begin{array}{c}
		\vm{T}_{Hh}^t\cdot \vm{B}_{H}\\
		\vm{T}_{Fk}^t \cdot \vm{B}_{F} \\
		\vm{T}_{Fe}^t \cdot \vm{B}_{F}
	\end{array}\right)
	\label{Bg0}
\end{equation} 
\begin{equation}
	\vm{A}_{gD}=\left(  \begin{array}{c}
		\vm{T}_{Hh}^t\cdot \left( \vm{A}_{HH} \cdot \vm{T}_{HD} +\vm{A}_{HF} \cdot \vm{T}_{FD}\right)\\
		\vm{T}_{Fk}^t\cdot  \vm{A}_{HF}^t \cdot \vm{T}_{HD} +\vm{P}_{Fk}^t \cdot \vm{T}_{FD}\\
		\vm{T}_{Fe}^t\cdot \left( \vm{A}_{HF}^t \cdot \vm{T}_{HD} +\vm{A}_{FF} \cdot \vm{T}_{FD}\right)
	\end{array}\right)
	\label{AgD}
\end{equation}
and $\vm{A}_{DD}$, $\vm{A}_{Dg}$ and $\vm{B}_{D}$ not being  shown here because not used later.
The $\vm{A}_{HF}$ term represents the coupling term of the NSP part with the SP boundary (nodes surrounding SP, around yellow part in the figure \ref{ts_sp}, named SPF in this work).
For this boundary, we know that the nodes are not enriched, so the sub-bloc corresponding to these dofs in  $\vm{T}_{Fe}$ is zero and $\vm{A}_{HF}\cdot \vm{T}_{Fe}=0$.
Moreover, for this boundary, all the fine-scale elements are identical to the global-scale elements, therefore the sub-bloc corresponding to these dofs in  $\vm{T}_{Fk}$ is an identity matrix $\vm{I}_{Fk}$ and $\vm{A}_{HF}\cdot \vm{T}_{Fk}=\vm{A}_{HF}\cdot \vm{I}_{Fk}$.
Otherwise, the  $\vm{T}_{Hh}$ operator, because of its structure (given in \ref{oper_rg})  reduces a block with H-set rows or columns to a block with h-set rows or columns. Thus $\vm{A}_{gg}$ can be written as :

\begin{equation}
	\vm{A}_{gg}=\left(  \begin{array}{ccc}
		\vm{A}_{hh}&\vm{A}_{hF}\cdot \vm{I}_{Fk}&\vm{0}\\
		\vm{I}_{Fk}^t\cdot \vm{A}_{hF}^t&\vm{T}_{Fk}^t \cdot \vm{P}_{Fk}&\vm{P}_{Fk}^t \cdot \vm{T}_{Fe} \\
		\vm{0}&\vm{T}_{Fe}^t \cdot \vm{P}_{Fk}& \vm{T}_{Fe}^t\cdot \vm{A}_{FF}  \cdot \vm{T}_{Fe}
	\end{array}\right)
	\label{Aggs}
\end{equation}

And the  $\vm{BN}_{g}$ and $\vm{A}_{gD}$  matrices simplify as follows:
\begin{equation}
	\vm{BN}_{g}=\left(  \begin{array}{c}
		\vm{B}_{h}\\
		\vm{T}_{Fk}^t \cdot \vm{B}_{F} \\
		\vm{T}_{Fe}^t \cdot \vm{B}_{F}
	\end{array}\right)
	\label{BN1}
\end{equation} 
\begin{equation}
	\vm{A}_{gD}=\left(  \begin{array}{c}
		\vm{A}_{hH}\cdot \vm{T}_{HD}  +\vm{A}_{hF} \cdot \vm{T}_{FD}\\
		\vm{I}_{Fk}^t\cdot \vm{A}_{HF}^t \cdot \vm{T}_{HD}+\vm{P}_{Fk}^t \cdot \vm{T}_{FD} \\
		\vm{T}_{Fe}^t\cdot  \vm{A}_{FF} \cdot \vm{T}_{FD}
	\end{array}\right)
	\label{AgD1}
\end{equation} 
By eliminating the Dirichlet boundary condition of the system \eqref{G_sys}, the final system to solve is \eqref{g_sys} with
\begin{equation}
	\vm{B}_{g}=\vm{BN}_{g}-\vm{A}_{gD}\cdot \vm{U}_D
	\label{Bg1}
\end{equation} 
and with full terms:
\begin{equation}
	\vm{B}_{g}=\left(  \begin{array}{c}
		\vm{B}_{h} -\vm{A}_{hH}\cdot \vm{T}_{HD}  \cdot \vm{U}_{D} -\vm{A}_{hF} \cdot \vm{T}_{FD}\cdot \vm{U}_{D}\\
		\vm{T}_{Fk}^t \cdot \vm{B}_{F}- \vm{I}_{Fk}^t\cdot \vm{A}_{HF}^t \cdot \vm{T}_{HD}\cdot \vm{U}_{D} -\vm{P}_{Fk}^t \cdot \vm{T}_{FD}\cdot \vm{U}_{D} \\
		\vm{T}_{Fe}^t \cdot \left( \vm{B}_{F}- \vm{A}_{FF} \cdot \vm{T}_{FD}\cdot \vm{U}_{D}\right)
	\end{array}\right)
	\label{Bg}
\end{equation} 
In this work, the  Dirichlet boundary conditions are zero imposed displacements in all tested cases.
Thus, to simplify the presentation,  $\vm{U}_{D}$ is fixed at $\vm{0}$ and \ref{Bg} becomes:
\begin{equation}
	\vm{B}_{g}=\left(  \begin{array}{c}
		\vm{B}_{h}\\
		\vm{T}_{Fk}^t \cdot \vm{B}_{F} \\
		\vm{T}_{Fe}^t \cdot \vm{B}_{F}
	\end{array}\right)
	\label{Bgs}
\end{equation} 

\section{\TS procedures}\label{anexe_algo}
\subsection{INIT
}The INIT procedure is given by algorithm \ref{TS_algebra_init}.
\begin{algorithm}
	\footnotesize
	\begin{algorithmic}
		\Procedure{INIT}{}
		\State create empty $\vm{A}_{gg}^{ini}$ and $\vm{B}_{g}^{ini}$
		\For{$e_{macro}\in SP$}
		\State Initialize $\vm{u}_{R}^{e_{macro}}$ fine-scale dofs covered by $e_{macro}$ 
		\State Eliminate linear relation if any from  $\vm{u}_{R}^{e_{macro}}$ to create  $\vm{u}_{F}^{e_{macro}}$
		\State create empty $\vm{A}_{FF}^{e_{macro}}$  matrix and $\vm{B}_{F}^{e_{macro}}$ vector
		\For{$e_{micro}\in e_{macro}$}
		\State compute $\vm{A}_{FF}^{e_{micro}}$ matrix and $\vm{B}_{F}^{e_{micro}}$ vector
		\State $\left( \vm{A}_{FF}^{e_{macro}},\vm{B}_{F}^{e_{macro}}\right)  \Leftarrow \left( \vm{A}_{FF}^{e_{micro}},\vm{B}_{F}^{e_{micro}}\right) $ 
		\EndFor
		\State Create $\vm{T}_{Fk}^{e_{macro}}$ operator based on classical form function of $e_{macro}$ at $\vm{u}_{F}^{e_{macro}}$ dof location
		\State $\vm{P}_{Fk}^{e_{macro}}\gets \vm{A}_{FF}^{e_{macro}}\cdot \vm{T}_{Fk}^{e_{macro}}$
		\State $\vm{A}_{kk}^{e_{macro}}\gets \vm{T}_{Fk}^{t~e_{macro}}\cdot \vm{P}_{Fk}^{e_{macro}}$
		\State $\vm{B}_{k}^{e_{macro}}\gets \vm{T}_{Fk}^{t~e_{macro}}\cdot \vm{B}_{F}^{e_{macro}}$
		\State $\left(\vm{A}_{gg}^{ini},\vm{B}_{g}^{ini}\right) \Leftarrow \left(  \vm{A}_{kk}^{e_{macro}},\vm{B}_{k}^{e_{macro}}\right)$
		\State Create a zero value $\vm{T}_{Fe}^{e_{macro}}$ operator, component product of shifted enrichment function with classical form function of $e_{macro}$ at $\vm{u}_{F}^{e_{macro}}$ dof location
		\State $\vm{T}_{Fm}^{e_{macro}}\gets \left\{\vm{T}_{Fk}^{e_{macro}},\vm{T}_{Fe}^{e_{macro}}\right\}$
		\EndFor
		\For{$e_{macro}\in NSP$}
		\State compute $\vm{A}_{hh}^{e_{macro}}$,$\vm{A}_{hF}^{e_{macro}}$ matrix and $\vm{B}_{h}^{e_{macro}}$ vector
		\State $\left( \vm{A}_{gg}^{ini},\vm{B}_{g}^{tmp}\right) \Leftarrow\left( \left( \vm{A}_{hh}^{e_{macro}},\vm{A}_{hF}^{e_{macro}} \cdot \vm{I}_{Fk}^{e_{macro}}\right) ,\vm{B}_{h}^{e_{macro}} \right) $
		\EndFor
		\State $\left\|\vm{B}_r\right\|\gets$\Call{COMPUTE\_B\_NORM}{$\vm{B}_{g}^{tmp}$, $\vm{B}_F$}\Comment{~} 
		\State $\vm{B}_{g}^{ini} \gets \vm{B}_{g}^{ini}+\vm{B}_{g}^{tmp}$
		\For{$p\in patches$}
		\State  Create $\vm{u}_{Q}^{p}$ fine-scale dofs covered by patch $p$\Comment{~}
		\State Eliminate fixed values from $\vm{u}_{Q}^{p}$ to create $\vm{u}_{q}^{p}$ \Comment{~} 
		\State Create empty $\vm{A}_{qq}^{p}$ matrix and $\vm{BI}_{q}^{p}$ vector
		\For{$e_{macro}\in p$}
		\State $\left( \vm{A}_{qq}^{p}~\text{and}~\vm{A}_{qd}^{p}, \vm{BI}_{q}^{p} \right) \Leftarrow \left( \vm{A}_{FF}^{e_{macro}},\vm{B}_{F}^{e_{macro}}\right) $
		\EndFor
		\State   $\vm{D}_{qd}^{p}\gets -\vm{A}_{qd}^{p}$
		\EndFor
		\State \Return
		$\vm{A}_{gg}^{ini},\vm{B}_{g}^{ini},\vm{A}_{FF},\vm{B}_F,\left\|\vm{B}_r\right\|,\vm{A}_{qq}^{patches},\vm{BI}_q^{patches},\vm{D}_{qd}^{patches}$
		\EndProcedure
	\end{algorithmic}
	\caption{\TS algorithm: initialization part. Procedure COMPUTE\_B\_NORM is given by algorithm \ref{TS_algebra_B_norm}}
	\label{TS_algebra_init}
\end{algorithm} 
In this procedure, the constant sub-blocks of \ref{Aggs} and \ref{Bgs} are stored in memory as $\vm{A}_{gg}^{ini}$ and $\vm{B}_{g}^{ini}$ which will be reused in the \tS loop to compute the final matrices $\vm{A}_{gg}$ and $\vm{B}_g$. 
In INIT, the first loop on SP macro-elements computes and attaches to a  macro-element $e_{macro}$ the  matrices $\vm{A}_{FF}^{e_{macro}}$ and $\vm{B}_{F}^{e_{macro}}$,  computed by integration and assembly of \eqref{equilibR} over $\omega^{e_{macro}}$ with elimination of \eqref{hanging_eq} but keeping Dirichlet values.
In this loop, the vectors $\vm{u}_F^{e_{macro}}$,  dedicated to the storage of the $e_{macro}$ part of the $\vm{u}_F$ vector, are also allocated and attached to $e_{macro}$.
Still in this loop, $\vm{T}_{Fk}^{e_{macro}}$  and $\vm{P}_{Fk}^{e_{macro}}$ are computed and attached to the macro-elements because they will also be  reused in the \tS loop (in fact it is $\vm{T}_{Fm}^{e_{macro}}$ operator with a  $\vm{T}_{Fe}^{e_{macro}}$ null at this stage that is saved in memory).
Finally, the $\vm{T}_{Fk}^{t~e_{macro}} \cdot \vm{P}_{Fk}^{e_{macro}}$ and $\vm{T}_{Fk}^{t~e_{macro}}\cdot \vm{B}_{F}^{e_{macro}}$ terms are calculated and assembled into $\vm{A}_{gg}^{ini}$ and $\vm{B}_{g}^{ini}$.
The second loop on the NSP elements (if any) computes the terms  $\vm{A}_{hh}$, $\vm{A}_{hF}$ and $\vm{B}_h$ and assembles them into $\vm{A}_{gg}^{ini}$ and $\vm{B}_{g}^{tmp}$.
This last temporary vector contains the boundary condition for the NSP elements and is used as an argument to the COMPUTE\_B\_NORM procedure ( algorithm \ref{TS_algebra_B_norm} below) called by INIT to obtain $\left \|\vm{B}_{r} \right \|$, the constant denominator of $resi$.
Once this is done, $\vm{B}_{g}^{tmp}$ is added to $\vm{B}_{g}^{ini}$ to save it.
The INIT procedure, with a loop over the patches, is also responsible for defining and partially creating the fine-scale linear systems \eqref{p_sys}.
Construction of  one of these systems  for $p\in I_e^g$, is simply obtained by assembling   $\vm{A}_{FF}^{e_{macro}}$ and $\vm{B}_{F}^{e_{macro}}$ with $e_{macro}\in J^p$ to form the following system:
\begin{equation}
	\left(  \begin{array}{ccc}
		\vm{A}_{DR\,DR}^p&\vm{A}_{DRd}^p&\vm{A}_{DRq}^p\\
		\vm{A}_{dDR}^p&\vm{A}_{dd}^p&\vm{A}_{dq}^p\\
		\vm{A}_{qDR}^p&\vm{A}_{qd}^p&\vm{A}_{qq}^p
	\end{array}\right)\cdot \left( \begin{array}{c}
		\vm{u}_{DR}^p\\
		\vm{u}_d^p\\
		\vm{u}_q^p
	\end{array}\right)=\left( \begin{array}{c}
		\vm{BN}_{DR}^p\\
		\vm{BN}_d^p\\
		\vm{BN}_q^p
	\end{array}\right)
	\label{QMF_sys}
\end{equation}
Eliminating  from \eqref{QMF_sys} the Dirichlet boundary condition on $\partial\omega_p$  we obtain:
\begin{equation}
	\vm{A}_{qq}^p
	\cdot \vm{u}_q^p =
		\vm{BN}_q^p-\vm{A}_{qDR}^p\cdot\vm{u}_{DR}^p-\vm{A}_{qd}^p\cdot\vm{u}_d^p
	\label{QMF_sys_red}
\end{equation}
Right-hand side of \eqref{QMF_sys_red} (which corresponds to  the vector $\vm{B}_q^p$) is divided in two contributions: $\vm{BI}_q^p=\vm{BN}_q^p-\vm{A}_{qDR}^p\cdot\vm{u}_{DR}^p$ which is constant during the \tS loop and $\vm{D}_{qd}^p\cdot\vm{u}_d^p$ (with $\vm{D}_{qd}^p=-\vm{A}_{qd}^p$) which varies during the \tS loop as does $\vm{u}_d^p$ (subset of $\vm{u}_F$).
Thus, the INIT procedure creates  for all patches $\vm{A}_{qq}^p$, $\vm{BI}_q^p$ and $\vm{D}_{qd}^p$.

\subsection{UPDATE\_MICRO\_DOFS}
The  UPDATE\_MICRO\_DOFS procedure (algorithm \ref{TS_algebra_micro_update}) updates  the values of $\vm{u}_F$  from a $\vm{U}_G$ vector given as argument according to the equation \eqref{u_ru_g}. 
This is done by a simple loop over the SP elements where for  $e_{macro}\in SP$ the following is computed :
 $\vm{u}_F^{e_{macro}}=\vm{T}_{Fm}^{e_{macro}}\cdot \vm{U}_m^{e_{macro}}$ where $\vm{U}_m^{e_{macro}}$ is the restriction of $\vm{U}_G$ to $e_{macro}$ dofs. 
The $\vm{T}_{Fm}^{e_{macro}}$ is supposed to be correctly computed when we enter this procedure (i.e. its $\vm{T}_{Fe}^{e_{macro}}$ part either null or updated by the  UPDATE\_MACRO\_PRB procedure).
 
\subsection{MICRO-SCALE\_RESOLUTION}
The  MICRO-SCALE\_RESOLUTION procedure (algorithm \ref{TS_algebra_patch})  loops over the patches to compute, for each  patch $p$, the  term $\vm{D}_{qd}^p\cdot\vm{u}_d^p$  (given in \eqref{QMF_sys_red}) with  $\vm{u}_d^p$ being the  current $\vm{u}_F$ restricted  to $\partial\omega_p\setminus\partial\Omega^D$.  
This term added to $\vm{BI}_q^p$ gives  the right-hand side $\vm{B}_q^p$ of the system \eqref{p_sys}.
Then, as mentioned earlier,  this system of equations  is solved with a direct solver ($\ddagger$ in algorithm \ref{TS_algebra_patch})  that first  processes a factorization of $\vm{A}_{qq}^p$ before doing a backward/forward resolution.
When the factorization is finished, as $\vm{A}_{qq}^p$ is constant during the \tS loop, only the factors are kept in memory, $\vm{A}_{qq}^p$ being freed thus gaining memory.
The following resolutions use the factors directly to solve the systems, which saves factoring time. 
This procedure, when completed, provides all $\vm{u}_q^p$  solutions. 
\subsection{UPDATE\_MACRO\_PRB}
The UPDATE\_MACRO\_PRB procedure (algorithm \ref{TS_algebra_macro_update} ) is dedicated to the update of the problem at the global-scale (local to global transfer).
It starts by copying $\vm{A}_{gg}^{ini}$ and $\vm{B}_{g}^{ini}$ into $\vm{A}_{gg}$ and $\vm{B}_{g}$ respectively.
Then, a loop on the SP elements finishes the computation of $\vm{A}_{gg}$ and $\vm{B}_{g}$  as follows.
The $\vm{T}_{Fe}^{e_{macro}}$, based on all $\vm{u}_q^p$ covering $e_{macro}$, is constructed by following the enrichment function defined in \eqref{shift_enrich} and the kinematic equation \eqref{kinematicnew}.
This operator is used to update the $\vm{T}_{Fm}^{e_{macro}}$ operator (for the algorithm \ref{TS_algebra_micro_update}) and compute the terms $\vm{T}_{Fe}^{t~e_{macro}} \cdot \vm{B}_{F}^{e_{macro}}$, $\vm{T}_{Fe}^{t~e_{macro}} \cdot \vm{P}_{Fk}^{e_{macro}}$ and $\vm{T}_{Fe}^{t~e_{macro}}\cdot \vm{A}_{FF}^{e_{macro}}\cdot \vm{T}_{Fe}^{e_{macro}}$  which are assembled into $\vm{A}_{gg}$ and $\vm{B}_{g}$.
\vspace{2mm}

\noindent\begin{minipage}{\textwidth}
	\captionsetup[algorithm]{format=bfcn} 
	\begin{minipage}{0.46\textwidth}
		\algva
		\captionof{algorithm}{\TS algorithm : update macro problem. Here $\vm{1}_F^{e_{macro}}$ is a unit vector corresponding to the F-set dofs (restricted to the dofs covered by $e_{macro}$ ). The  $e^P$ corresponds to the index of the enriched dof for the patch $p$ and $k^{p}$ to the index of the classical dof related to $e^P$ }
		\label{TS_algebra_macro_update}
		\vspace{-0.2em}\rule{\textwidth}{0.4pt}\vspace{-0.2em}
		\footnotesize 
		\begin{algorithmic}
			\Procedure{UPDATE\_MACRO\_PRB}{$\vm{A}_{gg}^{ini}$, $\vm{B}_{g}^{ini}$, $\vm{u}_q^{patches}$}
			\State $\left( \vm{A}_{gg},\vm{B}_{g} \right) \gets \left( \vm{A}_{gg}^{ini},\vm{B}_{g}^{ini}\right) $
			\For{$e_{macro}\in SP$}
			\For{$p\in$ patches including $e_{macro}$}
			\State $\vm{X}_F^{p,e_{macro}}\gets$ restriction of $\vm{u}_{q}^{p}$ to dofs covered by $e_{macro}$
			\State $shift^p \gets \vm{X}_{F}^{p,e_{macro}}[e^{p}]$
			\State $\vm{T}_{Fe^{p}}^{e_{macro}} \gets \vm{T}_{Fk^{p}}^{e_{macro}}*\left( \vm{X}_F^{p,e_{macro}} -\vm{1}_F^{e_{macro}}\times shift^{p} \right)$
			\EndFor
			\State $\vm{T}_{Fm}^{e_{macro}}\gets \left\{\vm{T}_{Fk}^{e_{macro}},\vm{T}_{Fe}^{e_{macro}}\right\}$
			\State $\vm{A}_{ee}^{e_{macro}}\gets \vm{T}_{Fe}^{t~e_{macro}}\cdot \vm{A}_{FF}^{e_{macro}}\cdot \vm{T}_{Fe}^{e_{macro}}$
			\State $\vm{A}_{ek}^{e_{macro}}\gets \vm{T}_{Fe}^{t~e_{macro}}\cdot \vm{P}_{Fk}^{e_{macro}}$
			\State $ \vm{A}_{gg} \Leftarrow\left( \vm{A}_{ee}^{e_{macro}},\vm{A}_{ek}^{e_{macro}} \right) $
			\State $\vm{B}_{g} \Leftarrow \vm{T}_{Fe}^{t~e_{macro}}\cdot \vm{B}_{F}^{e_{macro}} $
			\EndFor
			\State \Return $\vm{A}_{gg},\vm{B}_{g}$
			\EndProcedure
		\end{algorithmic}
		\vspace{-0.8em}\rule{\textwidth}{0.4pt}
	\end{minipage}\hfill\noindent\begin{minipage}{0.5\textwidth}
		\footnotesize
		\algva
		\captionof{algorithm}{\raggedright \TS algorithm: update micro-scale field.}
		\label{TS_algebra_micro_update}
		\vspace{-0.2em}\rule{\textwidth}{0.4pt}\vspace{-0.2em}
		\begin{algorithmic}
			\Procedure{UPDATE\_MICRO\_DOFS}{$\vm{U}_g$}
			\For{$e_{macro}\in SP$}
			\State $U_m^{e_{macro}}\gets \text{restriction~of~}U_{g}~\text{to}~e_{macro}~\text{dofs}$
			\State $\vm{u}_F^{e_{macro}} \gets \vm{T}_{Fm}^{e_{macro}}\cdot U_m^{e_{macro}}$
			\EndFor
			\State \Return $\vm{u}_F$
			\EndProcedure
		\end{algorithmic}
		\vspace{-0.8em}\rule{\textwidth}{0.4pt}
		\algva
		\captionof{algorithm}{\raggedright \TS algorithm: fine-scale problem computation.}
		\label{TS_algebra_patch}
		\vspace{-0.2em}\rule{\textwidth}{0.4pt}\vspace{-0.2em}
		\begin{algorithmic}
			\Procedure{MICRO-SCALE\_RESOLUTION}{$\vm{u}_F$, $\vm{A}_{qq}^{patches}$, $\vm{BI}_q^{patches}$, $\vm{D}_{qd}^{patches}$}
			\For{$p\in patches$}
			\For{$e_{macro}\in p$}
			\State $\vm{u}_{d}^{p} \Leftarrow \vm{u}_{F}^{e_{macro}}$ 
			\EndFor
			\State $\vm{B}_{q}^{p}\gets \vm{BI}_{q}^{p}+\vm{D}_{qd}^{p}.\vm{u}_{d}^{p}$
			\State $\vm{u}_{q}^{p}\gets \vm{A}_{qq}^{-1~p}\cdot \vm{B}_{q}^{p}$  {\large$\ddagger$}\Comment{~}
			\EndFor
			\State \Return $\vm{u}_q^{patches}$
			\EndProcedure
		\end{algorithmic}
		\vspace{-0.8em}\rule{\textwidth}{0.4pt}
	\end{minipage}
	\vspace{1em}
\end{minipage}

\subsection{COMPUTE\_RESIDUAL}
The COMPUTE\_RESIDUAL procedure (algorithm \ref{TS_algebra_residual}) computes only the term $\left \| \vm{A}_{fr}\cdot \vm{u}_{r}-\vm{B}_{f} \right \|$   of   \eqref{residual_error_simp}  considering that it is  mainly a scalar product of the vector $\vm{A}_{fr}\cdot \vm{u}_{r}-\vm{B}_{f}$  by itself.
In this procedure, since the information is stored on the F-set,  the computations are performed on this set. 
But since the organization of the data is by macro-element, the residual vector is not assembled into a complete F-set vector.
During the first loop on SP, its contribution per macro-element ($\vm{A}_{FF}^{e_{macro}}\cdot \vm{u}_{F}^{e_{macro}}-\vm{B}_{F}^{e_{macro}}$) is stored in an accumulation buffer called $\vm{VR}_{F}^{e_{macro}}$.
It is an accumulation buffer because the residual contributions of $e_{macro}$ boundaries from other adjacent macro-elements (remote or local to the process)  are added to it.
This accumulation allows to compute a part of the dot product of the residual F-set vector as a local dot product (appearing in the second loop on the SP element).
Their sum is the numerator of $resi$ to the power of two.
But since this local computation is redundant at the boundaries of the macro-elements, a diagonal scaling matrix, applied to the local dot product, is computed once for the entire \tS loop, so that the boundary terms are counted only once.
This scaling matrix is also used, with a zero scaling factor, to remove the Dirichlet DR-set contributions in the local dot product. 
This design maintains a semi-independent computation of these accumulation buffers and local dot product which can be treated in parallel (MPI and even with threads as proposed in conclusion).
Note that, with the enrichment strategy adopted in section \ref{TS_enrich_strategy}, as the h-set part of the residual vector that has been eliminated,  the rows associated with the SPF nodes (see \ref{Agg_anexe}) can also be considered null.
They are eliminated by using a zero scale factor in the local scalar product for the SPF node dofs.
The  $resi$ value is obtained by dividing the computed numerator by $\left \|\vm{B}_{r} \right \|$ term given as argument to the procedure. 
\subsection{COMPUTE\_B\_NORM}
The  COMPUTE\_B\_NORM procedure (algorithm \ref{TS_algebra_B_norm} ) provides the term $\left \|\vm{B}_{r} \right \|$  from the vectors $\vm{B}_g^{tmp}$ and $\vm{B}_F$ containing the Neumann boundary condition and volume loading related to the NSP  and SP elements respectively.
It uses an algorithm  very similar  to that of \ref{TS_algebra_residual} in that the  scalar product of  vector $\vm{B}_r$  by itself is transformed into  a sum of local dot products.
The scalar product of $\vm{B}_g^{tmp}$ initializes the sum of the local scalar products of  $\vm{B}_r$.
The same diagonal scaling matrix as in \ref{TS_algebra_residual}, is applied to the local dot product to correct for the redundant contribution to the macro-element boundaries.

\vspace{2mm}
\noindent\begin{minipage}{\textwidth}
	\captionsetup[algorithm]{format=bfcn} 
	\begin{minipage}{0.45\textwidth}
		\algva
		\captionof{algorithm}{\TS algorithm: compute residual value.
			Return $ \frac{\left \| \vm{A}_{rr}\cdot \vm{u}_{r}-\vm{B}_{r} \right \|}{\left \|\vm{B}_{r} \right \|} $}
		\label{TS_algebra_residual}
		\vspace{-0.2em}\rule{\textwidth}{0.4pt}\vspace{-0.2em}
		\footnotesize 
		\begin{algorithmic}
			\Procedure{COMPUTE\_RESIDUAL}{$\vm{u}_F$, $\vm{A}_{FF}$, $\vm{B}_F$, $\left\|\vm{B}_r\right\|$}
			\State $resi \gets 0$
			\For{$e_{macro}\in SP$}
			\State $\vm{V}_{F}^{e_{macro}}\gets \vm{A}_{FF}^{e_{macro}}\cdot \vm{u}_{F}^{e_{macro}}-\vm{B}_{F}^{e_{macro}}$
			\State $\vm{VR}_{F}^{e_{macro}}\gets \vm{VR}_{F}^{e_{macro}}+\vm{V}_{F}^{e_{macro}}$
			\For{$e_{adj}\in e_{macro}$ Adjacency }
			\State $\vm{VR}_{F}^{e_{adj}}\gets \vm{VR}_{F}^{e_{adj}}+\vm{V}_{F}^{e_{macro}~\cap ~e_{adj}}$ \Comment{~}
			\EndFor
			\EndFor
			\For{$e_{macro}\in SP$}
			\State $\vm{V}_{F}^{e_{macro}}\gets \text{scaled}~\vm{VR}_{F}^{e_{macro}}$
			\State $resi\gets resi+ \vm{V}_{F}^{t~e_{macro}}\cdot \vm{VR}_{F}^{e_{macro}}$
			\EndFor
			\State reduce $resi$ on all processes\Comment{~}
			\State $resi\gets \sqrt{resi}/\left\|\vm{B}_r\right\|$
			\State \Return $resi$
			\EndProcedure
		\end{algorithmic}
		\vspace{-0.8em}\rule{\textwidth}{0.4pt}
	\end{minipage}\hfill
	\begin{minipage}{0.45\textwidth}
		\algva
		\captionof{algorithm}{\raggedright \TS algorithm: compute $\left\|\vm{B}_r\right\|$.}
		\label{TS_algebra_B_norm}
		\vspace{-0.2em}\rule{\textwidth}{0.4pt}\vspace{-0.2em}
		\footnotesize	
		\begin{algorithmic}
			\Procedure{COMPUTE\_B\_NORM}{$\vm{B}_{g}^{tmp}$, $\vm{B}_F$}
			\State $\left\|\vm{B}_r\right\| \gets \vm{B}_{g}^{t~tmp}\cdot \vm{B}_{g}^{tmp}$\Comment{~}
			\For{$e_{macro}\in SP$}
			\State $\vm{VR}_{F}^{e_{macro}}\gets \vm{VR}_{F}^{e_{macro}}+\vm{B}_{F}^{e_{macro}}$
			\For{$e_{adj}\in e_{macro}$ Adjacency }
			\State $\vm{VR}_{F}^{e_{adj}}\gets \vm{VR}_{F}^{e_{adj}}+\vm{B}_{F}^{e_{macro}~\cap ~e_{adj}}$ \Comment{~}
			\EndFor
			\EndFor
			\For{$e_{macro}\in SP$}
			\State $\vm{V}_{F}^{e_{macro}}\gets \text{scaled}~\vm{VR}_{F}^{e_{macro}}$
			\State $\left\|\vm{B}_r\right\| \gets \left\|\vm{B}_r\right\| + \vm{V}_{F}^{t~e_{macro}}\cdot \vm{VR}_{F}^{e_{macro}}$
			\EndFor
			\State reduce $\left\|\vm{B}_r\right\| $ on all processes\Comment{~}
			\State $\left\|\vm{B}_r\right\| \gets \sqrt{\left\|\vm{B}_r\right\| }$
			\State \Return $\left\|\vm{B}_r\right\|$
			\EndProcedure
		\end{algorithmic}
		\vspace{-0.8em}\rule{\textwidth}{0.4pt}
	\end{minipage}
	\vspace{1em}
\end{minipage}

\section{Patch cost estimate}\label{TAP_anexe}
In the example of the cube, we have 4 types of patches depending on whether the enriched dof is on the corners, edges, faces or volumes of the cube.
Following the octree refinement strategy, for a coarse level $L_c$ and a final level $L$  ($L_c\leqslant L$), the numbers of these patches and their number of dofs  are given by the following:
\begin{subequations}
	\label{TAP_anexe_nbp}
	\begin{gather}
		nbpatch_{corner}\left( L_c \right)=8\\
		nbpatch_{edge}\left( L_c\right) =12 \left( {{2}^{L_c}}-1\right)\\
		nbpatch_{face}\left( L_c\right) =6 {{\left( {{2}^{L_c}}-1\right) }^{2}}\\
		nbpatch_{vol}\left( L_c\right) ={{\left( {{2}^{L_c}}-1\right) }^{3}}
	\end{gather}
\end{subequations} 
\begin{subequations}
	\label{TAP_anexe_nbd}
	\begin{gather}
	nbdofs_{corner}\left( L_c,L\right) = 3 {{\left( {{2}^{L-L_c}}+1\right) }^{3}} \\
	nbdofs_{edge}\left( L_c,L\right) =3 \left( 2 {{\left( {{2}^{L-L_c}}+1\right) }^{3}}-{{\left( {{2}^{L-L_c}}+1\right) }^{2}}\right)\\
    nbdofs_{face}\left( L_c,L\right) =3 \left( 4 {{\left( {{2}^{L-L_c}}+1\right) }^{3}}-4 {{\left( {{2}^{L-L_c}}+1\right) }^{2}}+\left( {{2}^{L-L_c}}+1\right) \right)\\
    nbdofs_{vol}\left( L_c,L\right) =3 \left( 8 {{\left( {{2}^{L-L_c}}+1\right) }^{3}}-12 {{\left( {{2}^{L-L_c}}+1\right) }^{2}}+6 \left( {{2}^{L-L_c}}+1\right) -1\right)
\end{gather}
\end{subequations} 
Solving a pacth with $N_d$ dofs and  a sparse ratio $SR_p$, for $N_L$ iterations using \eqref{TAP:costr} can be expressed as 
\begin{equation}
	\label{TAP_anexe_c1}
	cost_{1\_patch}(N_d,N_l,SR_p)=count_{fact}(N_d,SR_p)+NL\times count_{bf}(N_d,SR)
\end{equation}
since only one factorization is performed and a backward/forward resolution  is done at each \tS iteration.

The final cost for all patches using (\ref{TAP_anexe_nbp},~\ref{TAP_anexe_nbd},~\ref{TAP_anexe_c1}) will be the following:
\begin{equation}
\label{TAP_anexe_cp}
\begin{split}
cost_{patch}(L_c,L,N_l,SR_p)=nbpatch_{corner}\left( L_c \right)\times cost_{1\_patch}(nbdofs_{corner}\left( L_c,L\right),N_l,SR_p)+\\ nbpatch_{edge}\left( L_c \right)\times cost_{1\_patch}(nbdofs_{edge}\left( L_c,L\right),N_l,SR_p)+\\ nbpatch_{face}\left( L_c \right)\times cost_{1\_patch}(nbdofs_{face}\left( L_c,L\right),N_l,SR_p)+\\ nbpatch_{vol}\left( L_c \right)\times cost_{1\_patch}(nbdofs_{vol}\left( L_c,L\right),N_l,SR_p)
\end{split}
\end{equation}
\section{Static scheduling}\label{static_scheduling}
The scheduling algorithm, formalized in the algorithms \ref{sequencing_algo:general},\ref{sequencing_algo:pick_first} and \ref{sequencing_algo:pick}, consists in  arbitrarily selecting a vertex for $\mathcal{S}$ by ascending order of process.
\begin{algorithm}
	\footnotesize
	\begin{algorithmic}
		\State Create $\mathcal{D}$ and  $\mathcal{L}$  sets of distributed and local patches sorted in decreasing weight order
		\State Create $\mathcal{D}_o$ and  $\mathcal{L}_o$ empty sets that will store, in sequence ordering, distributed and local patches 
		\State Create color vector $\vm{col}$ of size $nbpid-pid$  
		\State $nbs\gets card\left( \mathcal{D}\right)$
		\State $nbs\gets \max_{i\in \mathcal{P}} {nbs}_i$  \Comment{~}
		\State $sid \gets  1 $
		\State $nbe \gets nbs $
		\Repeat
		\State $nbs\gets nbe$
		\While{$sid\leqslant nbs$}
		\State  $mxwg \gets$\Call{COMPUTE\_WG}{$\mathcal{D}$,$pid$}\Comment{~}
		\If {$pid=0$}
		\State  $\left( mxwl,\mathcal{D}, \mathcal{L},\mathcal{D}_o,\mathcal{L}_o, \vm{col}\right)  \gets$\Call{PICK\_FIRST}{$\mathcal{D}$, $\mathcal{L}$,$\mathcal{D}_o$,$\mathcal{L}_o$, $\vm{col}$, $pid$, $mxwg$}	
		\EndIf
		\If {$pid>0$}
		\State Receive from $pid-1$ $\vm{col}$ and $mxwl$ \Comment{~}
		\If {$\vm{col}[0]<0$}
		\State $\left( mxwl,\mathcal{D}, \mathcal{L},\mathcal{D}_o,\mathcal{L}_o, \vm{col}\right)  \gets$\Call{PICK}{$\mathcal{D}$, $\mathcal{L}$,$\mathcal{D}_o$,$\mathcal{L}_o$, $\vm{col}$, $pid$, $mxwg$,$mxwl$}
		\Else
		\State $patch_s\gets$ pick in $\mathcal{D}$ the patch with identifier$=\vm{col}[0]$
		\State  $\mathcal{D}_o\gets \mathcal{D}_o \cup  patch_s$
		\State $\mathcal{D}\gets \mathcal{D} \setminus  patch_s$
		\EndIf
		\EndIf
		\If {$pid<nbpid-1$}
		\State Send to $pid+1$ last $nbpid-1-pid$ components of $\vm{col}$ and $mxwl$ \Comment{~}
		\EndIf
		\State Store $\vm{col}$ as a the new last column of $M$
		\State $sid\gets sid+1$
		\EndWhile
		\State $nbe\gets nbe+card\left( \mathcal{D}\right)$
		\State $nbe\gets \max_{i\in \mathcal{P}} {nbe}_i$  \Comment{~}
		\Until{$nbs < nbe$}
		\State  $\mathcal{L}_o\gets \mathcal{L}_o \cup  \mathcal{L}$
		\State Create communicators based on colored sequence of $M$  \Comment{~}
	\end{algorithmic}
	\caption{Sequencing algorithm: Let $\mathcal{P}$ be the set of process identifiers (starting at 0) and $nbpid= card\left( \mathcal{P}\right) $.
		This algorithm is executed on each process $pid\in \mathcal{P}$ containing  patches  either  local  (entirely defined in $pid$) or distributed (on $pid$ and other processes ).
		In this algorithm, $nbs$ is the maximum number of sequences for all processes in  $\mathcal{P}$, $nbe$ is the maximum number of enlarged sequences for all processes in  $\mathcal{P}$, $mxwg$ is the maximum weight of all distributed patches of all process with an identifier greater than $pid$ (given by the algorithm \ref{sequencing_algo:compute_G}), $mxwl$ is the maximum weight of all selected patches of all processes with an identifier less than $pid$, $M$ is a matrix storing per column each $\mathcal{S}$ as a set of colors for all sequences, and  PICK\_FIRST, PICK, COMPUTE\_WG are procedures described respectively in the algorithms  \ref{sequencing_algo:pick_first}, \ref{sequencing_algo:pick},\ref{sequencing_algo:compute_G}.
		See \ref{convention_anexe} for other conventions.}
	\label{sequencing_algo:general}
\end{algorithm}

\clearpage
\noindent\begin{minipage}{\textwidth}
	\captionsetup[algorithm]{format=bfcn} 
	\begin{minipage}{0.45\textwidth}
		\algva
		\captionof{algorithm}{Pick first algorithm: On $pid=0$ the choice of the candidate is independent of the choice made in the other processes. First try in $\mathcal{D}$ then in $\mathcal{L}$. If both are empty, this process does  not participate in the current sequence. See the algorithm \ref{sequencing_algo:general} for the notation. }
		\label{sequencing_algo:pick_first}
		\vspace{-0.2em}\rule{\textwidth}{0.4pt}\vspace{-0.2em}
		\footnotesize
		\begin{algorithmic}
			\Function{PICK\_FIRST}{$\mathcal{D}$, $\mathcal{L}$,$\mathcal{D}_o$,$\mathcal{L}_o$, $\vm{col}$, $pid$, $mxwg$}
			\State $\vm{col}\gets -\vm{1}$
			\If {$\mathcal{D}\neq\varnothing$}
			\State $patch_s\gets$ pick out first patch in $\mathcal{D}$
			\State $\mathcal{A}_s \gets \{n,n\in \mathcal{P} , n\geqslant pid,  n~\text{participate to}~ patch_s\}$
			\State $\forall l \in \mathcal{A}_s ~\vm{col}[l-pid]\gets {patch_s}$ Identifier
			\State $\mathcal{D}_o\gets \mathcal{D}_o \cup patch_s$
			\State $\mathcal{D}\gets \mathcal{D} \setminus patch_s$
			\ElsIf {$\mathcal{L}\neq\varnothing$}
			\State $patch_s\gets$ the first patch $s$ in $\mathcal{L}$ with $weight_{s}<mxwg$ 
			\If {$patch_s=\varnothing$} 				\State$patch_s\gets$ pick out last patch in $\mathcal{L}$
			\EndIf
			\State $\vm{col}[0]\gets {patch_s}$ Identifier
			\State $\mathcal{L}_o\gets \mathcal{L}_o \cup patch_s$
			\State $\mathcal{L}\gets \mathcal{L} \setminus patch_s$
			
			\Else
			\State $patch_s\gets \varnothing$
			\State $\vm{col}[0]\gets -1$
			\EndIf
			\State $mxwl\gets weight_{patch_s}$
			\State \Return 	$mxwl$,$\mathcal{D}$, $\mathcal{L}$,$\mathcal{D}_o$,$\mathcal{L}_o$, $\vm{col}$
			\EndFunction
		\end{algorithmic}
		\vspace{-0.8em}\rule{\textwidth}{0.4pt}
	\end{minipage}\hfill
	\begin{minipage}{0.45\textwidth}
		\algva
		\captionof{algorithm}{Pick algorithm: First, try in $\mathcal{D}$. The choice of the candidate depends on the choices made by all the processes having an identifier $< pid$.  If no candidate is available, try in $\mathcal{L}$ respecting a certain condition on the weights. If no candidate is found, this process do  not participate in the current sequence. See algorithm \ref{sequencing_algo:general} for  the notation.}
		\label{sequencing_algo:pick}
		\vspace{-0.2em}\rule{\textwidth}{0.4pt}\vspace{-0.2em}
		\footnotesize
		\begin{algorithmic}
			\Procedure{PICK}{$\mathcal{D}$, $\mathcal{L}$,$\mathcal{D}_o$,$\mathcal{L}_o$, $\vm{col}$, $pid$, $mxwg$,$mxwl$}
			\State $patch_s\gets \varnothing$
			\If {$\mathcal{D}\neq\varnothing$}
			\State $patch_s\gets$ the first patch $s$, if any, in $\mathcal{D}$ not already selected in $\vm{col}$
			\If {$patch_s \neq\varnothing$}
			\State $\mathcal{A}_s \gets \{n,n\in \mathcal{P} , n\geqslant pid,  n~\text{participate to}~ patch_s\}$
			\State $\forall l \in \mathcal{A}_s ~\vm{col}[l-pid]\gets {patch_s}$ Identifier
			\State $\mathcal{D}_o\gets \mathcal{D}_o \cup patch_s$
			\State $\mathcal{D}\gets \mathcal{D} \setminus patch_s$
			\EndIf
			\EndIf
			\If {$patch_s=\varnothing~\text{and}~ \mathcal{L}\neq\varnothing$}
			\State $crit\gets \max(mxwg,mxwl)$
			\State $patch_s\gets$  the first patch s in $\mathcal{L}$ with $weight_{s}<crit$ 
			\If {$patch_s=\varnothing$} \State$patch_s\gets$ pick out last patch in $\mathcal{L}$ \EndIf
			\State $\vm{col}[0]\gets {patch_s}$ identifier
			\State $\mathcal{L}_o\gets \mathcal{L}_o \cup patch_s$
			\State $\mathcal{L}\gets \mathcal{L} \setminus patch_s$
			\Else
			\State $patch_s\gets \varnothing$
			\State $\vm{col}[0]\gets -1$
			\EndIf
			\State $mxwl\gets \max (mxwl,weight_{patch_s})$
			\State \Return 	$mxwl$,$\mathcal{D}$, $\mathcal{L}$,$\mathcal{D}_o$,$\mathcal{L}_o$, $\vm{col}$
			\EndProcedure
		\end{algorithmic}
		\vspace{-0.8em}\rule{\textwidth}{0.4pt}
	\end{minipage}
	\vspace{1em}
\end{minipage}

A color is imposed on the distributed vertex with the highest weight in the lowest process (see figure \ref{seq:4}\subref*{seq:4pseq1}-\subref*{seq:4pseq8} for illustration).
In the algorithms \ref{sequencing_algo:pick_first} and \ref{sequencing_algo:pick} this color is the identifier of the selected patch.
All its adjacent vertices are frozen to avoid choosing one of them for a new color (it corresponds to fill $\vm{col}$ by a color in algorithm \ref{sequencing_algo:pick_first} and \ref{sequencing_algo:pick}).
Then, if a vertex remains uncolored and  unfrozen, it can be assigned a new color by choosing again the one with  the highest weight in the lowest process.
It will freeze its own set of vertices.
The selection continues until all vertices of the graph have a color or are frozen.
The set of colored vertices in $\mathcal{S}$  can be used directly to create a colored  MPI communicator\footnote{with MPI\_Comm\_split}.
During selection, the lowest process may not have a  distributed vertex available.
In this case, a local patch with an appropriate weight is chosen if it is present, otherwise this process will not participate in this sequence.

At each sequence construction, two weights, $mxwl$ and $mxwg$ (in the algorithm \ref{sequencing_algo:general}), are exchanged  between the processes respectively, from the lowest identifier to the highest identifier and in reverse direction.
The  $mxwg$ weight is the maximum weight of all distributed patches stored in  all processes whose  process identifier is greater than the current one.
It is obtained with the algorithm \ref{sequencing_algo:compute_G}.
\begin{algorithm}
	\footnotesize
	\begin{algorithmic}
		\Procedure{COMPUTE\_WG}{$\mathcal{D}$,$pid$}
		\State $mxwg\gets 0$
		\If {$pid<nbpid-1$}
		\State $mxwg\gets mxwg_{pid+1}$   \Comment{~}
		\EndIf
		\State  $mxwg\gets \max(mxwg,\max_{k\in \mathcal{D}} (weight_k))$
		\If {$pid>0$}\State $mxwg_{pid-1}\gets mxwg$ \Comment{~}
		\EndIf
		\State \Return $mxwg$
		\EndProcedure
	\end{algorithmic}
	\caption{Compute the maximum weight of all distributed patches of all process whose identifier is greater than $pid$. See algorithm \ref{sequencing_algo:general} for the notation.}
	\label{sequencing_algo:compute_G}
\end{algorithm} 
It is used directly in the algorithm \ref{sequencing_algo:pick_first} when a local patch has to be select (i.e. when all the distributed patches of process 0 have been consumed): the selected local patch must have a weight lower than $mxwg$ in order not to  introduce an imbalance with respect to the weight of the remaining distributed patches  on all the other processes that will  potentially be selected when the remaining sequence is selected.
In the algorithm \ref{sequencing_algo:pick}, $mxwg$ is used in conjunction with $mxwl$ which corresponds to the maximum weight of all selected patches in all processes whose identifier is less than the current one.
In this algorithm, again when it comes to selecting a local patch, the one whose  weight is lower than  $max(mxwg,mxwl)$ is selected.
Thus, the selected local patch will have a cost lower than the worst cost between the already selected patches and those that will be potentially selected.
\section{Computational condition}\label{cluster}
This work was performed by using HPC resources of the Centrale Nantes Supercomputing Center ICI-CNSC on the cluster Liger (France).
It is an INTEL-based computer composed of 252  nodes, with 24 cores (2 x 12 cores Xeon E5-2680v3 at 2.5GHz) per node.
Those nodes use a Gpfs network drive for I/O and 128GB of  memory.
The software used on this platform are eXlibris 2018 (https://git.gem.ec-nantes.fr/explore), the GCC 9.2 compiler, OpenMPI 4.03, Mumps 5.4.1,  Intel MKL Scalapack, ParMetis 4.0.3,  Intel MKL lapack, Intel MKL blas.
All Intel software comes from the parallel studio 2020 suites.
The integers are compiled  with a 4 bytes precision.
\section{Micro structure plan definition}\label{MS_plan_anexe}
The equations of the 64 planes are given below:\\
$\footnotesize \begin{array}{cc}
0.8609265\times x+0.1956651\times y+0.4695963\times z-1.748384=0&
0.7195847\times x+0.5102510\times y+0.4710009\times z-2.943090=0\\
0.5688037\times x+0.06691808\times y+0.8197465\times z-1.924178=0&
0.1071920\times x+0.6029550\times y+0.7905410\times z-1.271770=0\\
0.7416569\times x+0.2076639\times y+0.6378250\times z-2.580966=0&
0.4863807\times x+0.4607817\times y+0.7423705\times z-1.215735=0\\
0.2568314\times x+0.9470659\times y+0.1926236\times z-1.424000=0&
0.7151591\times x+0.3916348\times y+0.5789384\times z-1.921521=0\\
0.1054744\times x+0.9844276\times y+0.1406325\times z-1.712618=0&
0.6326948\times x+0.2613305\times y+0.7289744\times z-1.181465=0\\
0.2716518\times x+0.9055059\times y+0.3259821\times z-1.135105=0&
0.6398444\times x+0.3745431\times y+0.6710563\times z-2.180073=0\\
0.5836437\times x+0.5275241\times y+0.6173154\times z-2.052645=0&
0.5684090\times x+0.1764028\times y+0.8036127\times z-1.657721=0\\
0.8197048\times x+0.1024631\times y+0.5635471\times z-2.360125=0&
0.1274946\times x+0.7139698\times y+0.6884709\times z-2.370968=0\\
0.7305503\times x+0.1398926\times y+0.6683758\times z-0.9220794=0&
0.4882044\times x+0.5858453\times y+0.6468708\times z-0.5461269=0\\
0.7241275\times x+0.5738369\times y+0.3825579\times z-2.412523=0&
0.3979635\times x+0.3684848\times y+0.8401452\times z-1.930610=0\\
0.8650248\times x+0.1441708\times y+0.4805693\times z-1.795212=0&
0.8572357\times x+0.4653565\times y+0.2204320\times z-1.902939=0\\
0.8669178\times x+0.1284323\times y+0.4816210\times z-2.635038=0&
0.1543033\times x+0.7715167\times y+0.6172134\times z-1.961483=0\\
0.4793697\times x+0.7989495\times y+0.3631589\times z-1.759905=0&
0.1659080\times x+0.2156804\times y+0.9622663\times z-0.6748800=0\\
0.2468854\times x+0.7715167\times y+0.5863527\times z-0.6423204=0&
0.5649452\times x+0.6013933\times y+0.5649452\times z-1.464101=0\\
0.4231527\times x+0.8711968\times y+0.2489134\times z-1.657172=0&
0.6388813\times x+0.5525460\times y+0.5352789\times z-1.462725=0\\
0.6333450\times x+0.4446890\times y+0.6333450\times z-2.730715=0&
0.2433962\times x+0.2920754\times y+0.9249055\times z-2.773892=0\\
0.3298492\times x+0.7985822\times y+0.5034540\times z-2.098559=0&
0.1632993\times x+0.4082483\times y+0.8981462\times z-2.001109=0\\
0.5423839\times x+0.6693248\times y+0.5077637\times z-0.9228155=0&
0.6018227\times x+0.5249942\times y+0.6018227\times z-1.695434=0\\
0.9093977\times x+0.3247849\times y+0.2598279\times z-2.058806=0&
0.8230470\times x+0.5534282\times y+0.1277142\times z-0.7484869=0\\
0.4433384\times x+0.1313595\times y+0.8866768\times z-0.8443746=0&
0.6734445\times x+0.6884099\times y+0.2693778\times z-1.324568=0\\
0.2972254\times x+0.9145396\times y+0.2743619\times z-2.464607=0&
0.7218661\times x+0.3925938\times y+0.5698943\times z-2.207026=0\\
0.5392394\times x+0.6564654\times y+0.5275168\times z-2.088212=0&
0.3743731\times x+0.06606583\times y+0.9249217\times z-1.333261=0\\
0.1427762\times x+0.1665723\times y+0.9756376\times z-0.5186730=0&
0.4939317\times x+0.2798946\times y+0.8232196\times z-1.744109=0\\
0.7648147\times x+0.07114555\times y+0.6403100\times z-1.913092=0&
0.8026276\times x+0.2390806\times y+0.5464699\times z-2.138411=0\\
0.7220829\times x+0.6804243\times y+0.1249759\times z-1.493591=0&
0.5806682\times x+0.5806682\times y+0.5706566\times z-1.237695=0\\
0.8945864\times x+0.4388537\times y+0.08439495\times z-1.095704=0&
0.3212124\times x+0.4534764\times y+0.8313734\times z-1.009433=0\\
0.9747546\times x+0.05415304\times y+0.2166121\times z-1.855889=0&
0.5572679\times x+0.6868651\times y+0.4665499\times z-2.181626=0\\
0.4943023\times x+0.8687737\times y+0.02995771\times z-0.7418343=0&
0.4652615\times x+0.6203487\times y+0.6314263\times z-1.781062=0\\
0.3647265\times x+0.4103173\times y+0.8358315\times z-0.9762100=0&
0.5807795\times x+0.5915347\times y+0.5592691\times z-1.467443=0\\
0.9473874\times x+0.2368468\times y+0.2153153\times z-1.995499=0&
0.9486833\times x+0.3162278\times y-1.061239=0\\
0.8602915\times y+0.5098024\times z-0.3024251=0&
0.9363292\times x+0.3511234\times z-0.9363292=0\\
0.9805807\times x+0.1961161\times z-1.668649=0&
0.1240347\times y+0.9922779\times z-0.9292094=0
\end{array}$

\section{Pull-out damage computation}\label{PO_cone_anexe}
The damage field  is computed using level-set technologies where the closest distance to the surfaces  given by \eqref{PO_annex_cone_eq} is interpreted, after scaling, as a damage between $[0,1]$.
The conical envelope where the damage is calculated (non-zero) is defined by the surface  given by the following series of equations:
\begin{equation}
	 \begin{array}{c}
		(x^2+z^2).(\cos\theta)^2 - (y-o_1)^2 {(\sin\theta)}^2 =0\\
		(x^2+z^2).(\cos\theta)^2 - (y-o_2)^2 {(\sin\theta)}^2 =0\\
		y+469.00=0\\
		y-h=0
	\end{array}
	\label{PO_annex_cone_eq}
\end{equation}
with \begin{itemize}
	\item $\theta=35$° is the angle at the apex between the axes $\overrightarrow{e_y}$ and the generating line of the lateral surface.
	\item $o_1=-545.08$ and $o_2=-531.31$ are $y$ coordinates of the two apexes.
	\item $h$ is the parameter that controls the  stage of disk damage.
\end{itemize}
The envelope of the conic is then any point $P$ such that:
\begin{equation}
P(x,y,z)\in \text{conic envelope if}~\left\lbrace \begin{array}{c}
		(x^2+z^2).(\cos\theta)^2 < (y-o_1)^2 {(\sin\theta)}^2\\
		(x^2+z^2).(\cos\theta)^2 > (y-o_2)^2 {(\sin\theta)}^2\\
		-469.<y\\
		y<h
	\end{array}\right.
	\label{PO_annex_cone_env}
\end{equation}

Note that using planes to stop the conical envelope is a simple and convenient way to obtain a bounded region of easily adjustable size.
But for sure, it gives, a completely unrealistic shape of the damage front.

\section{Domain decomposition resolution algorithm}\label{DDRES}
The  domain decomposition method used, given by the algorithm \ref{domain_decomposition_algo}, is a distributed non-overlapping Schur complement method.
Each process holds a  domain that is  condensed on the  process boundary (using the incomplete LDL$^t$ MUMPS factorization).
The global boundary problem (i.e. the global Schur complement problem) is solved with a distributed preconditionned (block Jacobi) conjugate gradient given by the algorithm \ref{CG_algo}.
This distributed version of the conjugate gradient  works only with local contributions and  communicates mainly when computing the scalar product.
The block Jacobi preconditioner  simply uses the factorization of the global boundary matrix diagonal block owned by the current process.
For a given process, this diagonal block is arbitrary chosen as its local boundary dofs,  excluding all dofs present in any process whose identifier is less than that of the current process.
Communication only occurs  at the time  of the construction of the diagonal block,   when all the contributions of the other processes are added to the process that owns the block.
The factorization of this block is then local to each process.
This arbitrary choice, easy to implement, induces that some processes (at least the last one) may not have a diagonal block to treat, which induces some imbalance in the preconditionning task.
Note that this choice ensures that the diagonal block can always be factorized as long as the global matrix is not ill-conditioned.

\noindent\begin{minipage}{\textwidth}
	\captionsetup[algorithm]{format=bfcn} 
	\begin{minipage}{0.45\textwidth}
		\algva
\captionof{algorithm}{Block Jacobi preconditionned iterative parallel domain decomposition algorithm:  
	The CONJGRAD procedure is given in the algorithm \ref{CG_algo}
}
\label{domain_decomposition_algo}
\vspace{-0.2em}\rule{\textwidth}{0.4pt}\vspace{-0.2em}
	\footnotesize
	\begin{algorithmic} 
		\State $\vm{S}_{bb} \gets \vm{A}_{bb}-\vm{A}_{bI}\cdot \vm{A}_{II}^{-1}\cdot \vm{A}_{Ib}$ 
		\State $\vm{B}_{b} \gets \vm{B}_{b}-\vm{A}_{bI}\cdot \vm{A}_{II}^{-1}\cdot \vm{B}_{I}$ 
		\State $ \vm{M}_{jj} \Leftarrow \bigcup\limits_{\forall p \in \mathcal{P}~\text{with}~t_p=b_p\cap j\neq \oslash} \vm{S}_{t_pt_p} $   \Comment{~}
		\State Factorize $\vm{M}_{jj}$ to form the block Jacobi preconditionner
		\State initialize a null vector $\vm{X}_s^0$ (or use a previously computed one)
		\State $\left( \vm{X}_s,crit\right) \gets$\Call{CONJGRAD}{$\vm{X}_s^0$,$\vm{S}_{bb}$,$\vm{B}_b$,$\vm{M}_{jj}^{-1}$,$\epsilon$} \Comment{~}
		\State $\vm{X}_b\gets \vm{X}_{s \cap b}$
		\State $\vm{X}_I\gets \vm{A}_{II}^{-1} \cdot \left( \vm{B}_I-\vm{A}_{Ib}\cdot \vm{X}_b \right) $		
	\end{algorithmic}
\vspace{-0.8em}\rule{\textwidth}{0.4pt}
\end{minipage}\hfill
\begin{minipage}{0.45\textwidth}
\algva
\captionof{algorithm}{\raggedright Parallel conjugate gradient algorithm.}\label{CG_algo}
\vspace{-0.2em}\rule{\textwidth}{0.4pt}\vspace{-0.2em}
	\footnotesize
	\begin{algorithmic}
		\Procedure{CONJGRAD}{$\vm{X}_s^0$,$\vm{S}_{bb}$,$\vm{B}_b$,$\vm{M}_{jj}^{-1}$,$\epsilon$}
		\State $\vm{X}_s \gets \vm{X}_s^0$
		\State $\vm{R}_s \Leftarrow \vm{B}_b$
		\State $stop \gets \vm{R}_s\cdot \vm{R}_s\times \epsilon^2$\Comment{~}
		\State $\vm{SP}_s \Leftarrow \vm{S}_{bb}\cdot \vm{X}_{s\cap b}$
		\State $\vm{R}_s \gets \vm{R}_s-\vm{SP}_s$
		\State $\vm{Z}_s \Leftarrow \vm{M}_{jj}^{-1}\cdot \vm{R}_{s\cap j}$
		\State $res_o \gets \vm{R}_s\cdot \vm{Z}_s$ \Comment{~}
		\State $\vm{P}_s \gets \vm{Z}_s$
		\State $alt \gets 0$
		\Repeat
		\State $res_n \gets res_o$
		\State $\vm{SP}_s \Leftarrow \vm{S}_{bb}\cdot \vm{P}_{s\cap b}$
		\State $\alpha \gets \frac{res_n}{\vm{P}_s\cdot \vm{SP}_s}$ \Comment{~}
		\State $\vm{X}_s \gets \vm{X}_s + \alpha\times \vm{P}_s$
		\State $\vm{R}_s \gets \vm{R}_s - \alpha\times \vm{SP}_s$
		\State $crit^2 \gets \vm{R}_s\cdot \vm{R}_s$\Comment{~}
		\If{$crit^2 \geqslant stop$}
		\State $\vm{Z}_s \Leftarrow \vm{M}_{jj}^{-1}\cdot \vm{R}_{s\cap j}$
		\State $res_o \gets \vm{R}_s\cdot \vm{Z}_s$ \Comment{~}
		\State $\beta \gets \frac{res_o}{res_n}$
		\State $\vm{P}_s \gets \vm{Z}_s + \beta\times \vm{P}_s$
		\State $iter \gets iter+1$
		\Else
		\State $alt \gets 1$
		\EndIf
		\Until{$alt=1$ or $iter>iter_{max}$}
		\State $crit \gets \epsilon.\sqrt{\frac{crit^2}{stop}}$
		\State \Return $\vm{X}_s$,crit
		\EndProcedure
	\end{algorithmic}
\vspace{-0.8em}\rule{\textwidth}{0.4pt}
\end{minipage}
\vspace{1em}
\end{minipage}

In the algorithm \ref{domain_decomposition_algo} and \ref{CG_algo} let $\mathcal{P}$ be the set of process identifiers (starting at 0) and $nbpid= card\left( \mathcal{P}\right)$.
These algorithms are executed on each process $pid\in \mathcal{P}$ with:
\begin{itemize}
	\item $K$  the set of dofs of the domain held by the process $pid$, $K= I\cup b$ , $I\cap b=\oslash$
	\item $I$ the set of  dofs eliminated by condensation in the $pid$ process
	\item $b$ the set of boundary dofs of the domain $pid$ (i.e. the Schur complement dofs)
	\item $s$ the set of all  domain boundary dofs, $s=\bigcup\limits_{p\in \mathcal{P}} b_{p}$
	\item $j$ the set of boundary dofs owned by the $pid$ process: $j=b\setminus \left( \bigcup\limits_{p=0,pid-1} b_{p}\right) $
\end{itemize}  

The following applies for these sets:
\begin{itemize}
	\item $\forall p \in \mathcal{P}$ and $\forall q \in \mathcal{P}$,$p \neq q$ then $j_p\cap j_q = \oslash$ 
	\item $\forall p \in \mathcal{P}$ and $\forall q \in \mathcal{P}$,$p \neq q$ then 
	$b_p\cap b_q \neq \oslash$ if $p$ and $q$ are connected by at least one mesh node
	\item for $pid=0$ $j=b$
	\item for $pid=nbid-1$ $j=\oslash$
	\item $j$ can also be $\oslash$ for any process $p\in \mathcal{P}\setminus (nbpid-1)$  if the domain $p$ is surrounded by domains of lower process identifier
\end{itemize}

\section{Mumps block low-rank resolution}\label{BLRRES}
The block low-rank (BLR) feature  has been introduced  in Mumps based on  \cite{Amestoy2015}.
We choose to activated it with the automatic choice of the BLR option by the library (ICNTL(35)=1 see Mumps 5.4.1 Users' guide ).
The BLR factorization variant is the default (ICNTL(36)=0, UFSC variant).
The dropping parameter used during BLR compression (CNTL(7)), is chosen to be identical to residual threshold $\epsilon$ of the algorithms \ref{TS_algebra} and \ref{domain_decomposition_algo}.
This choice is guided  by the fact that the dropping factor and residual error are expected to be strongly related, as shown in \cite{Higham2022}.
To enforce the condition on the residual error, the BLR resolution is followed by a parallel preconditioned  conjugate gradient resolution (algorithm \ref{CG_algo}) using both the factorization (as a preconditioner) and the solution (as a starting point) of the low-rank resolution (algorithm \ref{blr_algo}).
\begin{algorithm}
	\footnotesize
	\begin{algorithmic} 
		\State $\left( \tilde{\vm{A}}_{rr}^{-1},\tilde{\vm{X}}_r\right) \gets \text{Mumps BLR resolution of } \vm{A}_{rr}\cdot \vm{X}_{r}=\vm{B}_r$ \Comment{~}
		\State $\left( \vm{X}_r,crit\right)  \gets$\Call{CONJGRAD}{$\tilde{\vm{X}}_r$,$\vm{A}_{rr}$,$\vm{B}_r$,$\tilde{\vm{A}}_{rr}^{-1}$,$\epsilon$} \Comment{~}
	\end{algorithmic}
	\caption{Block low-rank algorithm:  
		the CONJGRAD procedure is given in the algorithm \ref{CG_algo}. Here $\tilde{.}$ represent the approximate factorization and the solution of the low-rank resolution. 
	}
	\label{blr_algo}
\end{algorithm} 
This second resolution is supposed to iterate very little to  just force the residual error to be less than $\epsilon$.

\end{document}